\documentclass[leqno]{amsbook}
\usepackage{amsmath, amsthm, amssymb, latexsym}
\newcounter{theorem}
\newcounter{proposition}
\newcounter{lemma}
\newcounter{corollary}
\newcounter{conjecture}
\newtheorem{definition}{Definition}[section]
\newtheorem{theorem}[definition]{Theorem}
\newtheorem{proposition}[definition]{Proposition}
\newtheorem{lemma}[definition]{Lemma}

\newtheorem{lemma-appendix}{Lemma}[chapter]

\theoremstyle{remark}
\newtheorem{example}[definition]{Example}
\newtheorem{remark}[definition]{Remark}
\newtheorem*{claim}{Claim}

\newtheorem*{acknowledgements}{Acknowledgements}
 
\numberwithin{section}{chapter}
\numberwithin{equation}{section}
\newcommand{\op}[1]{\operatorname{#1}}
\newcommand{\acou}[2]{\ensuremath{\langle #1 , #2 \rangle}}

\newcommand{\brak}[1]{\ensuremath{\langle\! #1\!\rangle}}

\newcommand{\Tr}{\ensuremath{\op{Tr}}}
\newcommand{\tr}{\op{tr}}
\newcommand{\Tra}{\ensuremath{\op{Trace}}}

\newcommand{\Hol}{\op{Hol}}

\newcommand{\C}{\ensuremath{\mathbb{C}}} 
\newcommand{\bH}{\ensuremath{\mathbb{H}}} 
\newcommand{\N}{\ensuremath{\mathbb{N}}} 
\newcommand{\R}{\ensuremath{\mathbb{R}}} 
\newcommand{\Z}{\ensuremath{\mathbb{Z}}}

\newcommand{\Rd}{\ensuremath{\R^{d+1}}}
\newcommand{\Rdd}{\ensuremath{\R^{d+2}}}
\newcommand{\Rdo}{\R^{d+1}\!\setminus\! 0}

\newcommand{\URd}{U\times\R^{d+1}}
\newcommand{\URdd}{U\times\R^{d+2}}
\newcommand{\URdo}{U\times(\R^{d+1}\!\setminus\! 0)}

\newcommand{\Ca}[1]{\ensuremath{\mathcal{#1}}}
\newcommand{\cA}{\Ca{A}}
\newcommand{\cB}{\Ca{B}}
\newcommand{\cD}{\ensuremath{\mathcal{D}}}
\newcommand{\cE}{\Ca{E}}

\newcommand{\cG}{\ensuremath{\mathcal{G}}}
\newcommand{\cH}{\ensuremath{\mathcal{H}}}

\newcommand{\cK}{\ensuremath{\mathcal{K}}}
\newcommand{\cL}{\ensuremath{\mathcal{L}}}

\newcommand{\cN}{\ensuremath{\mathcal{N}}}
\newcommand{\cS}{\ensuremath{\mathcal{S}}}

\newcommand{\cU}{\ensuremath{\mathcal{U}}}

\newcommand{\fg}{\ensuremath{\mathfrak{g}}}
\newcommand{\fh}{\ensuremath{\mathfrak{h}}}

\newcommand{\psivdo}{$\Psi_{H}$DO}
\newcommand{\psivdos}{$\Psi_{H}$DO's}

\newcommand{\pvdo}{\ensuremath{\Psi_{H}}} 
\newcommand{\pvhdo}{\ensuremath{\Psi_{H,\op{v}}}}

\newcommand{\psido}{$\Psi$DO} 
\newcommand{\psidos}{$\Psi$DO's} 
\newcommand{\psinf}{\ensuremath{\Psi^{-\infty}}} 
\newcommand{\vo}{\op{v}}

\newcommand{\Svb}{S_{\scriptscriptstyle{| |}}}
\newcommand{\SvbU}[1]{S_{\scriptscriptstyle{| |}}^{#1}(U\times\Rd)}
\newcommand{\ah}{\text{ah}}
\newcommand{\reg}{{\text{reg}}}
\newcommand{\ord}{{\op{ord}}}

\newcommand{\xiy}{{\xi\rightarrow y}}

\newcommand{\yxi}{{y\rightarrow\xi}}

\newcommand{\supp}{\op{supp}}
\newcommand{\loc}{\op{loc}}

\newcommand{\rk}{\op{rk}}
\newcommand{\im}{\op{im}}
 
\newcommand{\dom}{\op{dom}}

\newcommand{\End}{\ensuremath{\op{End}}}

\newcommand{\hotimes}{\hat\otimes}

\renewcommand{\Box}{\square}
\newcommand{\subsubset}{\subset\!\subset}

\newcommand{\Sp}{\op{Sp}}

\newcommand{\dbarb}{\overline{\partial}_{b}}
\newcommand{\dbarbpq}{\overline{\partial}_{b;p,q}}

\newcommand{\Boxb}{\Box_{b}}
\newcommand{\Boxbpq}{\Box_{b;p,q}}
 
\begin{document}
\frontmatter
\title{Heisenberg calculus and spectral theory of hypoelliptic operators on Heisenberg manifolds.} 

\author{Rapha\"el Ponge}

\address{Department of Mathematics, Ohio State University, Columbus, USA.}
\email{ponge@math.ohio-state.edu}
 \keywords{Heisenberg calculus, complex powers, heat equation, spectral asymptotics, analysis on CR and contact manifolds, hypoelliptic operators}
 \subjclass[2000]{Primary 58J40, 58J50; Secondary 58J35, 32V10, 35H10, 53D10}
\thanks{Research partially supported by NSF grant DMS 0409005}

\begin{abstract}
This memoir deals with the hypoelliptic calculus on Heisenberg manifolds, including CR and contact 
manifolds. In this context the main differential operators at stake include the H\"ormander's sum of squares, the Kohn Laplacian, the horizontal sublaplacian, 
the CR conformal operators of Gover-Graham and the contact Laplacian. These operators cannot be elliptic and the relevant pseudodifferential calculus to study them is provided 
by the Heisenberg calculus of Beals-Greiner and Taylor. 

The Heisenberg manifolds generalize CR and contact manifolds and their name stems from the fact that the relevant notion of tangent space in this 
setting is rather that of a  bundle of graded two-step nilpotent Lie groups. Therefore, the idea behind the Heisenberg calculus, which goes back to Stein, 
is to construct a pseudodifferential calculus modelled on homogeneous left-invariant convolution operators on nilpotent groups.  

The aim of this monograph is threefold. First, we give an intrinsic 
approach to the Heisenberg calculus by finding an intrinsic notion of principal symbol in this setting, in connection with the construction of the 
tangent groupoid in~\cite{Po:Pacific1}. This framework allows us to prove that the pointwise invertibility of a principal symbol, which can be restated in terms of the 
so-called Rockland condition, actually implies its global invertibility. 

Second, we study complex powers of hypoelliptic operators on Heisenberg manifolds in terms of the Heisenberg calculus. In particular, we show that complex 
powers of such operators give rise to holomorphic families in the Heisenberg calculus. To this end, due to the lack of microlocality of the 
Heisenberg calculus, we cannot make use of the standard approach of Seeley, so we rely on an alternative approach based on the pseudodifferential representation of 
the heat kernel in~\cite{BGS:HECRM}. This has some interesting consequences related to hypoellipticity and allows us to construct a scale of weighted
Sobolev spaces providing us with sharp estimates for the operators in the Heisenberg calculus. 

Third, we make use of the Heisenberg calculus and of the results of this monograph to derive spectral asymptotics for hypoelliptic operators on 
Heisenberg manifolds. The advantage of using the Heisenberg calculus is illustrated by reformulating in a geometric fashion these asymptotics for the main 
geometric operators on CR and contact manifolds, namely, the Kohn Laplacian and the horizontal sublaplacian in the CR 
setting and the horizontal sublaplacian and the contact Laplacian in the contact setting. 
\end{abstract}

\maketitle 

\tableofcontents

\mainmatter 

% \include{Intro}
% 
% \include{Chap2}
%        
% \include{Chap3}
%        
% \include{Chap4}
%        
% \include{Chap5}
%        
% \include{Chap6}
% 
% \include{Appendix}

% \part{Introduction}

\chapter{Introduction} 
This memoir deals with the hypoelliptic calculus on Heisenberg manifolds, including CR and contact 
manifolds. In this context the main differential operators at stake include the H\"ormander's sum of squares, the Kohn Laplacian, the horizontal 
sublaplacian, the CR conformal  operators of Gover-Graham and the contact Laplacian. 
These operators cannot be elliptic and the relevant pseudodifferential calculus to study them is provided 
by the Heisenberg calculus of of Beals-Greiner~\cite{BG:CHM} and Taylor~\cite{Ta:NCMA}. 

The Heisenberg manifolds generalize CR and contact manifolds and their name stems from the fact that the relevant notion of tangent space in this 
setting is rather that of a  bundle of graded two-step nilpotent Lie groups. Therefore, the idea behind the Heisenberg calculus, which goes back to Stein, 
is to construct a pseudodifferential calculus modelled on homogeneous left-invariant convolution operators on nilpotent groups.  

Our aim in this monograph is threefold. First, we give an intrinsic 
approach to the Heisenberg calculus by defining an intrinsic notion of principal symbol in this setting, in connection with the construction of the 
tangent groupoid in~\cite{Po:Pacific1}. This framework allows us to prove that the pointwise invertibility of a principal symbol, which can be restated in terms of the 
so-called Rockland condition, actually implies its global invertibility. 

These results have been already used in~\cite{Po:AofM1} to produce new invariants for CR and contact manifolds, 
extending previous results of Hirachi~\cite{Hi:LSSKGISPD} and Boutet de Monvel~\cite{BdM:LTTP}. Moreover, since our approach to the principal symbol 
connects nicely with 
the construction of the tangent groupoid of a Heisenberg manifold in~\cite{Po:Pacific1}, this presumably allows us to make use of global $K$-theoretic arguments in 
the Heisenberg setting, as those involved in the proof of the (full) 
Atiyah-Singer index theorem~(\cite{AS:IEO1}, \cite{AS:IEO3}). Therefore, this part of the memoir can also be seen as a step towards a reformulation 
of the Index Theorem for hypoelliptic operators on Heisenberg manifolds.

Second, we study complex powers of hypoelliptic operators on Heisenberg manifolds in terms of the Heisenberg calculus. In particular, we show that complex 
powers of such operators give rise to holomorphic families in the Heisenberg calculus. To this end, due to the lack of microlocality of the 
Heisenberg calculus, we cannot make use of the standard approach of Seeley, so we rely on an alternative approach based on the pseudodifferential representation of 
the heat kernel in~\cite{BGS:HECRM}. This has some interesting consequences related to hypoellipticity and allows us to construct a scale of weighted
Sobolev spaces providing us with sharp estimates for the operators in the Heisenberg calculus. 

These results are important ingredients in~\cite{Po:GAFA1} to construct an analogue for the Heisenberg calculus of the noncommutative residue 
trace of  Wodzicki~(\cite{Wo:LISA}, \cite{Wo:NCRF}) and Guillemin~\cite{Gu:NPWF} and to study the zeta and eta functions of hypoelliptic operators. 
In turn this has several geometric consequences. In particular, this allows us to make use of the framework of Connes' noncommutative geometry, 
including the local index formula of~\cite{CM:LIFNCG}. 

Third, we make use of the Heisenberg calculus and of the results of this monograph to derive spectral asymptotics for hypoelliptic operators on 
Heisenberg manifolds. The advantage of using the Heisenberg calculus is illustrated by reformulating in a geometric fashion these asymptotics for the main 
geometric operators on CR and contact manifolds, namely, the Kohn Laplacian, the horizontal sublaplacian and the Gover-Graham operators in the CR 
setting and the horizontal sublaplacian and the contact Laplacian in the contact setting.

On the other hand, although the setting of this monograph is the hypoelliptic calculus on Heisenberg calculus, it is believed that the results herein can be extended to 
more general settings such as the hypoelliptic calculus on Carnot-Carath\'eodory manifolds which are equiregular in the sense of~\cite{Gr:CCSSW}.

Following is a more detailed description of the contents of this memoir.

\section{Heisenberg manifolds and their main differential operators}
A Heisenberg manifold $(M,H)$ consists of a manifold $M$ together with a distinguished hyperplane bundle $H\subset TM$. This definition 
 covers many examples: Heisenberg group and its quotients by cocompact lattices, (codimension 1) foliations, CR and contact manifolds and the 
 confolations of Elyahsberg and Thurston. 
 
 In this setting 
 the relevant tangent structure for a Heisenberg manifold $(M,H)$ is rather that of a bundle $GM$ of two-step nilpotent Lie 
 groups %whose Lie group structure is encoded by an intrinsic Levi form $\cL:H\times H\rightarrow TM/H$  
 (see~\cite{BG:CHM}, \cite{Be:TSSRG}, \cite{EMM:HAITH}, \cite{EMM:RLSPD}, \cite{FS:EDdbarbCAHG}, \cite{Gr:CCSSW}, \cite{Po:Pacific1}, 
 \cite{Ro:INA}).  
 
 The main examples of differential operators on Heisenberg manifolds are the following. 
 
(a) H\"ormander's sum of squares on a Heisenberg manifold $(M,H)$ of the form,
% We will also consider differential 
% operators of the form, 
\begin{equation}
    \Delta=\nabla_{X_{1}}^{*}\nabla_{X_{1}}+\ldots+\nabla_{X_{m}}^{*}\nabla_{X_{m}},
     \label{eq:Intro.sum-of-squares}
\end{equation}
where the (real) vector fields $X_{1},\ldots,X_{m}$ span $H$ and $\nabla$ is a connection on a vector bundle $\cE$ over $M$ and 
the adjoint is taken with respect to a smooth positive measure on $M$ 
and a Hermitian metric on $\cE$.\smallskip

(b) Kohn Laplacian $\Box_{b;p,q}$ acting on $(p,q)$-forms on a CR manifold $M^{2n+1}$ endowed with a CR compatible Hermitian metric (not necessarily a 
Levi metric).\smallskip

(c) Horizontal sublaplacian $\Delta_{b;k}$ acting on horizontal differential forms of degree $k$ on a Heisenberg manifold $(M,H)$. When $M^{2n+1}$ 
is a CR manifold the horizontal sublaplacian preserves the bidegree and so we can consider its restriction $\Delta_{b;p,q}$ to foms of bidegree 
$(p,q)$.\smallskip

(d) Gover-Graham operators $\boxdot_{\theta}^{(k)}$, $k=1,\ldots,n+1, n+2, n+4,\ldots$ on a strictly pseudoconvex CR manifold 
$M^{2n+1}$ endowed with a CR compatible contact form $\theta$ (so that $\theta$ defines a pseudohermitian structure on $M$). These operators have been 
constructed by Gover-Graham~\cite{GG:CRIPSL} as the CR analogues 
of the conformal GJMS operators of~\cite{GMJS:CIPLIE}. In particular, they are differential operators which tranforms conformally under a conformal change of contact 
form and for $k=1$ we recover the conformal sublaplacian of Jerison-Lee~\cite{JL:YPCRM}.\smallskip
% of the (scalar) horizontal sublaplacian $\Delta_{b;0}$  These operators 
% were recently constructed by  in such way that  has same principal 
% part as $\Delta_{b;0}^{k}$ and it transforms conformally  (this generalizes the conformal sublaplacian of 

(e) Contact Laplacian on a contat manifold $M^{2n+1}$ associated to the contact complex of Rumin~\cite{Ru:FDVC}. This complex acts between sections 
of a graded subbunbdle $(\oplus_{k\neq n}\Lambda^{k}_{R})\oplus \Lambda_{R,1}^{n}\oplus \Lambda_{R,2}^{n}$ of horizontal forms.  The contact 
Laplacian is a differential operator of order 2 in degree $k\neq n$ and of order 4 in degree $n$.\smallskip

In the examples (a)--(c) the operators are instances of sublaplacians. More precisely, a sublaplacian is a second order differential 
$\Delta:C^{\infty}(M,\cE)\rightarrow C^{\infty}(M,\cE)$ which near any point $a\in M$ is of the form, 
\begin{equation}
     \Delta=-\sum_{j=1}^{d}X_{j}^{2}-i\mu(x)X_{0}+\sum_{j=1}^{d}a_{j}(x)X_{j}+b(x),
      \label{eq:Intro.sublaplacian}
 \end{equation}
 where  $X_{0},X_{1},\ldots,X_{d}$ is a local frame of $TM$ such that $X_{1},\ldots,X_{d}$ span $H$ and the coefficients $\mu(x)$ and 
 $a_{1}(x),\ldots,a_{d}(x), b(x)$ are local sections of $\End \cE$.  
%  and $\mu(x)$ and $a_{j}(x)$ 
%  the term $\mu(x)$ are local sections of $\End \cE$ 
% and the notation $\op{O}_{H}(1)$ stands for a differential operator of Heisenberg order~$\leq 1$.  
 
 \section{Intrinsic approach to the Heisenberg calculus}
  Although the differential operators above may be hypoelliptic under some conditions, they are definitely not elliptic.  Therefore, we cannot rely 
  on the standard pseudodifferential 
  calculus to study these operators. 
 
The substitute to the standard pseudodifferential calculus is provided by the Heisenberg calculus, independently introduced by
 Beals-Greiner~\cite{BG:CHM} and Taylor~\cite{Ta:NCMA} (see also \cite{BdM:HODCRPDO},  \cite{CGGP:POGD}, \cite{Dy:POHG}, \cite{Dy:APOHSC}, 
 \cite{EMM:HAITH}, \cite{FS:EDdbarbCAHG}, \cite{RS:HDONG}).  
The idea in the Heisenberg calculus, which goes back to Elias Stein, is the following. Since the relevant notion of tangent structure for a Heisenberg 
manifold $(M,H)$ is that of a bundle $GM$ of 2-step nilpotent graded Lie groups, it stands for reason to construct a pseudodifferential calculus which at every point 
$x \in M$ is well modelled by the calculus of convolution operators on the nilpotent tangent group $G_{x}M$. 

The result is a class of pseudodifferential operators, the \psivdos, which are locally \psidos\ of type $(\frac{1}{2},\frac{1}{2})$, but unlike the latter possess a full 
symbolic calculus and makes sense on a general Heisenberg manifold. In particular, a \psivdo\ admits a parametrix in the Heisenberg calculus if, and only if, its principal 
symbol is invertible, and then the \psivdo\ is hypoelliptic with a gain of derivatives controlled by its order (see Section~\ref{sec:PsiHDO} for a detailed 
review of the Heisenberg calculus). 

\subsection{Intrinsic notion of principal symbol} 
  In~\cite{BG:CHM} and~\cite{Ta:NCMA} the principal symbol of a \psivdo\ is defined in local coordinates only, so the definition 
\emph{a priori} depends on the choice of these coordinates. In the special case of a contact manifold, an intrinsic definition have been given 
in~\cite{EMM:HAITH} and~\cite{EMM:RLSPD} as a section over a bundle of jets of vector fields representing the tangent group bundle of the contact manifold. 
This approach is similar to that of Melrose~\cite{Me:APSIT} in the setting of the $b$-calculus for manifolds with boundary.   

In this paper we give an intrinsic definition of the principal symbol, valid for an arbitrary Heisenberg manifold, using the results of~\cite{Po:Pacific1}. 

Let $(M^{d+1},H)$ be a Heisenberg manifold. As shown in~\cite{Po:Pacific1} the tangent Lie group bundle of $(M,H)$ can be described as the 
bundle $(TM/H)\oplus H$ together with the grading and group law such that, for sections $X_{0}$, $Y_{0}$ of $TM/H$ and sections $X'$, $Y'$ of $H$, we have 
\begin{gather}
    t.(X_{0}+X')=t^{2}X_{0}+tX', \quad t \in \R,\\
    (X_{0}+X').(Y_{0}+Y')=X_{0}+Y_{0}+\frac{1}{2}\cL(X',Y')+X'+Y',
%     \label{eq:¥}
\end{gather}
where $\cL:H\times H \longrightarrow TM/H$ is the intrinsic Levi form such that 
\begin{equation}
    \cL(X',Y')=[X',Y'] \qquad \text{in $TM/H$}.
     \label{eq:Intro.Levi-form}
\end{equation}

In suitable coordinates, called Heisenberg coordinates, this description of $GM$  is equivalent to previous descriptions of $GM$ in terms of the 
Lie group of a nilpotent Lie algebra of jets of vector fields. A consequence of this equivalence is a tangent approximation result for Heisenberg 
diffeormorphisms stating that in Heisenberg coordinates such a diffeormorphism is well approximated by the induced isomorphisms between the tangent 
groups (see~\cite[Prop.~2.21]{Po:Pacific1}). 

Let $\cE$ be a vector bundle over $M$. For $m \in \C$ let $\pvdo^{m}(M,\cE)$ denote the class of 
\psivdos\ of order $m$ acting on the sections of $\cE$. Furthermore, let $\fg^{*}M$ be the linear dual of the Lie 
algebra bundle $\fg M$ of $GM$,  with canonical projection $\pi:\fg^{*}M\rightarrow M$, and let us define
$S_{m}(\fg^{*} M,\cE)$ as the space of sections $p_{m}(x,\xi)$ in $C^{\infty}(\fg^{*}M\setminus 0,\End \pi_{*}\cE)$ such that 
$p_{m}(x,\lambda.\xi)=\lambda^{m}p_{m}(x,\xi)$ for any $\lambda>0$.

The key to our definition of the principal symbol is the aforementioned approximation result for Heisenberg diffeomorphisms of~\cite{Po:Pacific1}. 
More precisely,  it allows us to carry out in Heisenberg coordinates a proof of the invariance by Heisenberg diffeomorphisms of 
the Heisenberg calculus allong similar lines as that in~\cite{BG:CHM} (see Appendix~\ref{chap.Appendix-Invariance}). The upshot is that it yields  
a change of variable formula for the principal symbol in Heisenberg coordinates showing that the latter can be intrinsically defined as a section 
over $\fg^{*}M\setminus 0$. Therefore, we obtain: 

\begin{proposition}\label{prop:Intro.principal-symbol}
    For any  $P\in \pvdo^{m}(M,\cE)$, $m \in \C$ there exists a unique symbol 
    $\sigma_{m}(P)\in S_{m}(\fg^{*}M, \cE)$ such that, for any $a \in M$, the symbol $\sigma_{m}(P)(a,.)$ agrees in trivializing Heisenberg coordinates centered at 
    $a$ with the principal symbol of $P$ at $x=0$.  %in the sense of~\cite{BG:CHM}. 
\end{proposition}

The symbol $\sigma_{m}(P)(x,\xi)$ is called the principal symbol of $P$. In local coordinates it can be explicitly related to the principal symbol in 
the sense of~\cite{BG:CHM} (see Eqs.~(\ref{eq:PsiHDO.local-global-principal-symbol})--(\ref{eq:PsiHDO.isomorphism-symbols})). 
In general the two definitions don't agree,  but they do when $P$ is a differential operator or the 
bundle $H$ is integrable. In any case we have a linear isomorphism 
$\sigma_{m}:  \pvdo^{m}(M,\cE)/\pvdo^{m-1}(M,\cE)\stackrel{\sim}{\rightarrow}  S_{m}(\fg^{*}M,\cE)$.
% %     \label{eq:}
% \end{equation}

% In order to distinguish it from the other definition of the principal symbol 
% in local coordinates, we will sometimes refer to it as the \emph{global} principal symbol, while the other principal symbol will be called the 
% \emph{local} principal symbol. 

As a consequence of this intrinsic definition of principal symbol we can define the model operator of a \psivdo\ $P\in \pvdo^{m}(M,\cE)$ at a point 
$a\in M$ as the left-invariant \psivdo\ $P^{a}$  on $G_{a}M$ with symbol $\sigma_{m}(a,.)$, that is, the left-convolution operator with 
the inverse Fourier transform of $\sigma_{m}(a,.)$ (see Definition~\ref{def:PsiHDO.model-operator}). 

These notions principal symbol and of model operators show that the 
Heisenberg calculus is well approximated by the calculus of left-invariant pseudodifferential operators on the tangent groups $G_{a}M$, $a\in M$. 

First, as it follows from the results of~\cite{BG:CHM}, for any $a \in M$ the convolution product on the group $G_{a}M$ defines a bilinear product, 
\begin{equation}
    *^{a}:S_{m_{1}}(\fg^{*}_{a}M)\times S_{m_{2}}(\fg^{*}_{a}M) \longrightarrow S_{m_{1}+m_{2}}(\fg^{*}_{a}M).
\end{equation}
This product depends smoothly on $a$ in such way  to give rise to a bilinear product, 
\begin{gather}
    *:S_{m_{1}}(\fg^{*}M,\cE)\times S_{m_{2}}(\fg^{*}M,\cE) \longrightarrow S_{m_{1}+m_{2}}(\fg^{*}M,\cE),
    \label{eq:Intro.product-symbols}\\
      p_{m_{1}}*p_{m_{2}}(a,\xi)=[p_{m_{1}}(a,.)*^{a}p_{m_{2}}(a,.)](\xi) \qquad \forall p_{m_{j}}\in S_{m_{j}}(\fg^{*}M,\cE).
\end{gather}
It then can be shown that the above product correspond to the product of \psivdos\ at the level of principal symbols and that the model operator of a 
product of two \psivdos\ is the product of the corresponding model operators (see Proposition~\ref{prop:PsiHDO.composition2}). 

On the other hand, we also can carry out in Heisenberg coordinates versions of the proofs that the transpose and adjoints of \psivdos\ are again 
\psivdos~(see Appendix~\ref{chap:Appendix-transpose}). As a consequence we can identify their principal symbols 
and see that the model operators at a point of the transpose and the adjoint of a \psivdo\ are respectively the transpose and 
the adjoint of its model operator (see Propositions~\ref{prop:PsiHDO.transpose-global} and~\ref{prop:PsiHDO.adjoint-manifold}).  

\subsection{Rockland condition, parametrices and hypoellipticity}
It follows from the results of~\cite{BG:CHM} that for a \psivdo\ $P\in \pvdo^{m}(M,\cE)$  the existence of a parametrix in 
$\pvdo^{-m}(M,\cE)$ is equivalent to the invertibility of its principal symbol $\sigma_{m}(P)$. Moreover, when $\Re m\geq 0$ this implies that $P$ is 
hypoelliptic with gain of $\frac{1}{2}\Re m$ derivatives.  

% It is shown in~\cite{BG:CHM} that, in local coordinates, the invertibility of the local principal symbol is equivalent to the existence of a \psivdo\ 
% parametrix. Thanks to Proposition~\ref{prop:Intro.composition} and  the relationship in local coordinates between the local and global principal symbols, we can 
% reformulate this result in a global fashion. More precisely, we have:  %using the global principal symbol: 
% 
% \begin{proposition}\label{thm:Intro.hypoellipticity}
%    Let $P:C^{\infty}_{c}(M,\cE)\rightarrow C^{\infty}(M,\cE)$ be a \psivdo\ of order $m$. Then the following are equivalent:\smallskip 
% 
%    1) The principal symbol $\sigma_{m}(P)$ of $P$ is invertible with respect to the convolution product for homogeneous 
%     symbols;\smallskip 
%     
%     2) The operator $P$ admits a parametrix $Q$ in $\pvdo^{-m}(M,\cE)$.\smallskip 
%  
%  \noindent Furthermore, if 1) and 2) hold then $P$ is hypoelliptic with gain of $\frac{k}{2}$-derivatives.
% \end{proposition}

In general it may be difficult to determine the invertibility of the principal symbol of a \psivdo, because the product~(\ref{eq:Intro.product-symbols}) 
for symbols is not anymore the pointwise product of symbols. Nevertheless, this problem can be understood in terms of a representation theoretic criterion, the 
so-called Rockland condition. 

If $P$ is a homogeneous left invariant \psido\ on a nilpotent graded group $G$ then to any unitary representation $\pi$ we can associate an  
(unbounded) operator $\pi_{P}$ on the representation space $\cH_{\pi}$ such that the  domain of its closure contains the space $C^{\infty}(\pi)$ of 
smooth vectors of $\pi$. The Rockland condition then requires that for 
any non-trivial irreducible unitary representation $\pi$ of $G$ the 
closure $\overline{\pi_{P}}$ is injective on $C^{\infty}(\pi)$. 

It is a remarkable result that the Rockland condition for $P$ is equivalent to its hypoellipticity 
(see~\cite{Ro:HHGRTC}, \cite{HN:HGNRN3}, \cite{HN:COHHIGGLG}, \cite{CGGP:POGD}). 
Moreover, it can be shown that $P$ is hypoelliptic iff it admits a 
left \psido\ inverse (see~\cite{Fo:SEFNLG}, \cite{Ge:LTHG}, \cite{CGGP:POGD}). 
It then follows that $P$ admits a two-sided \psido\ inverse if, and only if, $P$ and $P^{t}$ 
satisfies the Rockland condition. 

In the setting of the Heisenberg calculus we say that a \psivdo\ 
$P\in \pvdo^{m}(M,\cE)$ satisfies the Rockland condition at a point $a$ when the model operator $P^{a}$ 
satisfies the Rockland condition on $G_{a}M$. It then follows that the principal symbol $\sigma_{m}(P)$ is invertible at $x=a$, i.e., 
$\sigma_{m}(P)(a,.)$ admits an inverse in $S_{-m}(\fg^{*}_{a}M,\cE_{a})$ with respect to the product $*^{a}$, if, and only if, $P$ and $P^{t}$ satisfies the 
Rockland condition at $a$. 

If $P\in \pvdo^{m}(M,\cE)$ is such that $P$ and $P^{t}$ satisfy the Rockland condition at every point then, as mentioned above, for each point $a\in M$ we get an 
inverse $q^{a}\in S_{-m}(\fg^{*}_{a}M,\cE_{a})$ for $\sigma_{m}(P)(a,.)$. However, in order to obtain an inverse for $\sigma_{m}(P)$ 
in $S_{m}(\fg^{*}M,\cE)$ we still have to check that the family $(q^{a})_{a \in M}$ varies smoothly with $a$. 

By using an idea of Christ~\cite{Ch:ISASIO} it has been shown in~\cite{CGGP:POGD} that given a smooth family of homogeneous left invariant \psidos\ 
$(P_{u})_{u \in U}$ on a fixed nilpotent homogeneous group $G$ such that $P^{u}$ and $(P^{u})^{t}$ satisfy the Rockland condition for every $u$ then the family 
of inverses depend smoothly on $u$. 

In this memoir we show that in the Heisenberg setting this result is also true when the group varies from point to point. Namely, we prove: 
% we yet have to check 
% that the inverse depends smoothly enough on $a$ to yield an element of $S_{-m}(\fg^{*}M,\cE)$. This is a difficult issue whose solution seems to be 
% only known in the case a family of operators on fixed nilpotent homogeneous group (see~\cite{CGGP:POGD}). We prove that in the Heisenberg setting the 
% result is true in general, that is, we obtain: 

\begin{theorem}\label{thm:Intro.Rockland-invertibility}
  Let $P:C^{\infty}(M,\cE)\rightarrow C^{\infty}(M,\cE)$ be a \psivdo\ of order $m$. Then the following are equivalent:\smallskip

  (i) $P$ and $P^{t}$ satisfy the Rockland condition at every point of $M$;\smallskip 
    
  (ii) The principal symbol of $P$ is invertible.\smallskip 
  
  \noindent Moreover, when $m=0$ both (i) and (ii) are equivalent to:\smallskip
  
  (iii) For any $a \in M$ the model operator $P^{a}$ is invertible on $L^{2}(G_{a}M,\cE_{a})$.
\end{theorem}

In substance this theorem states that in the Heisenberg the pointwise invertibility of a principal symbol is equivalent to its invertibility. The 
proof elaborates on the ideas of~\cite{Ch:ISASIO} and~\cite{CGGP:POGD} and is divided into two steps. 

In the first step we prove the theorem in the case $m=0$. In this case, the equivalence of (i) and (iii) follows from a result of G{\l}owacki~\cite{Gl:RCNDO2} 
and it is immediate that (ii) implies (iii). Therefore, we only have to prove that (iii) implies (ii). The arguments are based on the ideas 
of~\cite{Ch:ISASIO} as in~\cite{CGGP:POGD}, but instead of relying on the result of Christ~\cite{Ch:ORISIO} 
on the $L^{p}$ boundedness of zero'th order 
convolutions operators on nilpotent graded groups, we rely on its earlier version due to Knapp-Stein~\cite{KS:IOSSG}, which can be more conveniently 
generalized to the setting of families of groups. 
% 
% generalizations to 
% family of groups more convenient.  
% 
% The latter require more 
% off-diagonal regularity but this not issue in our setting and, more importantly, its proof makes it more convenient to study 
% 
% holds in the setting of a family 
% of groups where the law group varies from point to point, but not the grading. 

The second step is the reduction to the case $m=0$. This is similar to what is done in~\cite{CGGP:POGD}, but instead of 
making use of the commutative approximation of the identity on a fixed nilpotent group of~\cite{Gl:SSGMCAINGHG}, 
we make use of integer powers of a sublaplacian with an invertible principal symbol (such an operator always exists thanks to the results 
of~\cite{BG:CHM}).

On the other hand, Theorem~\ref{thm:Intro.Rockland-invertibility} 
has several interesting consequences. First, if $P$ satisfies the Rockland condition at every  point and we have $\Re m\geq 0$ 
then $P$ is hypoelliptic with  gain of $\frac{1}{2}\Re m$ derivatives (Proposition~\ref{prop:Heisenberg.Rockland-hypoellipticity}). 

Second, even though the representation theory of $G_{a}M$ may vary as $a$ ranges over points of $M$ the Rockland condition is an open condition. More 
precisely, we prove: 

\begin{proposition}\label{prop:Intro.Rockland.open}
    Let $P:C^{\infty}_{c}(M,\cE)\rightarrow C^{\infty}(M,\cE)$ be a \psivdo\ of integer order $m$ with principal symbol $p_{m}(x,\xi)$
    and let $a\in M$.\smallskip
    
    1) If $P$ satisfies the Rockland condition at $a$ then there exists an open neighborhood $V$ of $a$ such that $P$ satisfies the Rockland 
    condition at every point of $V$.\smallskip
    
    2) If $p_{m}(a,\xi)$ is invertible in $S_{m}(\fg^{*}_{a}M,\cE_{a})$ then there exists an open neighborhood $V$ of $a$ such that 
    $p_{m|_{V}}$ is invertible on $S_{m}(\fg^{*}V,\cE)$.
\end{proposition}
% 
% 
% $P$ satisfies the Rockland condition at a point $a$ then there exists an open neighborhood of $a$ such that $P$ satisfies the Rockland condition at 
% every point of this neighborhood (Proposition~\ref{prop:Chap3.Rockland.open}). 
% Likewise, if the principal symbol of $P$ is invertible at a point $a$ then there exists a neighborhood $V$ of $a$ such 
% that $\sigma_{m}(P)_{|_{V}}$ is invertible.

Finally, if $(p_{\nu\in B})_{\nu \in B}$ is a smooth family with values in $S_{m}(\fg^{*}M,\cE)$ parametrized by a manifold $B$ 
such that $p_{\nu}$ admits an inverse 
$p^{(-1)}_{\nu}$ in $S_{-m}(\fg^{*}M,\cE)$ for any $\nu\in B$, then the family $(p^{(-1)}_{\nu})_{\nu \in B}$ too depends smoothly on $B$ 
(see Proposition~\ref{prop:Chap3.Rockland.family}).
% In particular, when $P$ is selfadfoint and the Levi form has constant rank the Rockland condition for $P$ is equivalent to the invertibility of the 
% principal symbol of $P$. 
% 
% Finally, if $2n$ denotes the rank of the Levi form at $a$ then we have $G_{a}M\simeq \bH^{2n+1}\times \R^{d-2n}$. Therefore, the irreducible 
% unitary representations of $G_{a}M$ are two kinds, one dimensional representations on $\C$ and infinite dimensional representations on 
% $L^{2}(\R^{n})$. In the former case the Rockland condition corresponds to the invertibility of the principal symbol along $H^{*}$, while in the latter 
% case it is enough to look at the representations coming from that of $\bH^{2n+1}$ with Planck constants $\pm 1$ (see Section~\ref{sec:hypoellipticity}). 

\subsection{Invertibility criteria for sublaplacians}
As alluded to above the sublaplacians cover the important examples that are the H\"ormander's sum of 
squares, the Kohn Laplacian or the horizontal sublaplacian. If $\Delta:C^{\infty}(M,\cE)\rightarrow C^{\infty}(M,\cE)$ is a sublaplacian 
then as shown in~\cite{BG:CHM} the Rockland condition can be formulated in terms of the Levi form~(\ref{eq:Intro.Levi-form}) 
as follows. 

For $a \in M$ let $2n$ be the rank of the Levi form $\cL_{a}$ and, using the same notation as in~(\ref{eq:Intro.sublaplacian}), consider the singular set 
  \begin{gather} 
  \Lambda_{a}=     (-\infty, -\frac12 \sum_{j=1}^{d}|\lambda_{j}|]\cup [\frac12 \sum_{j=1}^{d}|\lambda_{j}|
    |L(a)|,\infty) \qquad \text{if $2n<d$}, 
    \label{eq:Intro.singular-set1}\\
     \Lambda_{a}=\{\pm\frac12 \sum_{j=1}^{d} (1+2\alpha_{j})|\lambda_{j}|; \alpha_{j}\in \N^{d}\}\qquad \text{if $2n=d$},
     \label{eq:Intro.singular-set2}
  \end{gather}
where $\lambda_{1},\ldots,\lambda_{d}$ are the eigenvalues of $\cL(a)$ with respect to the frame $X_{0},\ldots,X_{d}$ in~(\ref{eq:Intro.sublaplacian}). 
Then the Rockland conditions at $a$ for $\Delta$ and $\Delta^{t}$ are both equivalent to the single condition, 
  \begin{equation}
     \Sp \mu(a) \cap \Lambda_{a}=\emptyset.
      \label{eq:Intro.sublaplacian.condition}
  \end{equation}

  In fact, when the condition~(\ref{eq:Intro.sublaplacian.condition}) holds at every point, we evan can derive an explicit formula for the inverse of 
  the principal symbol of $\Delta$.  This is carried out in~\cite{BG:CHM} in the scalar case only, but we really need to deal with the system case in order to study 
  sublaplacians acting on forms. For instance, the Kohn Laplacian locally is scalar modulo lower order terms if, and only if, the Levi form diagonalizes in a smooth 
  eigenframe, which needs not exist in general. 

In Section~\ref{sec:sublaplacian}, after having recalled the arguments of~\cite{BG:CHM} in the scalar case, we explain how to extend them for systems of 
sublaplacians. In particular, this allows us to complete the treatment of the Kohn Laplacian in~\cite{BG:CHM} (see below). 

\subsection{Invertibility criteria for the main differential operators on Heisenberg manifolds}
In Section~\ref{sec:Examples} we work out the previous invertibility criteria for the principal symbols of the 
main examples of operators on Heisenberg manifolds. In particular, 
we recover in a unified fashion several known hypoellipticity results. 

\subsubsection*{(a) H\"ormander's sum of squares} For a sum of squares as in~(\ref{eq:Intro.sum-of-squares}) 
the condition~(\ref{eq:Intro.sublaplacian}) is equivalent to have $\rk \cL_{a}\neq 0$, so that the 
invertibility of the principal symbol of $\Delta$ is equivalent to the condition, 
\begin{equation}
    H+[H,H]=TM. 
    \label{eq:Intro.bracket-condition}
\end{equation}
This is exactly the bracket condition of H\"ormander~\cite{Ho:HSODE} for a codimension 1 distribution $H\subset TM$. 
% When $X_{1},\ldots,X_{m}$ span $H$ we get a sublaplacian and 
% Proposition~\ref{prop:Intro.sublaplacian.Rockland-bundle}  allows us to 
% recover, in this special case, the celebrated result of H\"ormander~\cite{Ho:HSODE} about the hypoellipticity of sum of squares under the bracket 
% condition. 

\subsubsection*{(b) Kohn Laplacian} In the case of the Kohn Laplacian  acting on $(p,q)$-forms on a CR manifold $M^{2n+1}$ 
the condition~(\ref{eq:Intro.sublaplacian})  reduces to Kohn's 
$Y(q)$-condition. For instance when $M$ is $\kappa$-strictly pseudoconvex this reduces to have 
$q\neq \kappa$ and $q\neq n-\kappa$. 
% 
% The Kohn Laplacian is the Laplacian associated to the tangential Cauchy-Riemann complex, or $\dbarb$-complex, on 
% a CR manifold~(\cite{KR:EHFBCM}, \cite{Ko:BCM}). It was shown by Kohn~\cite{Ko:BCM} that under the condition $Y(q)$ the Kohn Laplacian acting on 
% $(p,q)$-forms is hypoelliptic with gain of one derivative. 
% 
% It was proved by Beals-Greiner~\cite{BG:CHM} that for the Kohn Laplacian acting on $(p,q)$-forms the 
% condition~(\ref{eq:Intro.sublaplacian.condition})  reduces to the condition 
% $Y(q)$, so we may apply Proposition~\ref{prop:Intro.sublaplacian.Rockland-bundle} to recover Kohn's result. This allows us to complete the treatment  of the Kohn 
% in~\cite{BG:CHM}, because the initial argument there is not quite complete (see Remark~\ref{rem:Examples.Boxb}). 

\subsubsection*{(c) Horizontal sublaplacian} 
For the horizontal sublaplacian $\Delta_{b;k}$ acting on horizontal forms of degree $k$ on a Heisenberg manifold the relevant condition to look at is 
a condition that we call condition $X(k)$: given a point $a\in M$ and letting $2n$ be 
the rank of  the Levi form $\cL$ at $a$,  we say that the condition $X(k)$ is satisfied at $a$ when we have 
\begin{equation}
     k\not \in\{n,n+1,\ldots,d-n\}.
%     \label{eq:¥}
\end{equation}
More precisely, we show that the condition~(\ref{eq:Intro.sublaplacian}) reduces to the condition $X(k)$ and so $\Delta_{b;k}$  has an invertible 
principal symbol if, and only if, the condition $X(k)$ holds at every point (see Proposition~\ref{prop:Examples.horizontal-sublaplacian}). 

For $k=0$ we get $\rk \cL_{a}\neq 0$ which is equivalent to the condition~(\ref{eq:Intro.bracket-condition}) 
(in fact $\Delta_{b;0}$ is a sum of squares modulo lower 
order terms). When $M^{2n+1}$ is a contact manifold the 
condition $X(k)$ exactly means that we must have $k\neq n$, so that we recover the hypoellipticity results of~\cite{Ta:DGSSPCM} and \cite{Ru:FDVC}, 
but in the non-contact case our invertibility criterion for the horizontal sublaplacian seems to be new.

When $M^{2n+1}$ is a CR manifold and we consider the horizontal sublaplacian $\Delta_{b;p,q}$ acting on $(p,q)$-forms we can refine the $X(k)$ 
condition into the $X(p,q)$ condition (see Proposition~\ref{prop:Examples.Deltabpq}). 
For instance when $M^{2n+1}$ is $\kappa$-strictly pseudoconvex it means that we must have 
$(p,q)\neq (\kappa,n-\kappa)$ and 
$(p,q)\neq (n-\kappa,\kappa)$. 

\subsubsection*{(d) Gover-Graham operators} On a strictly pseudoconvex CR manifold $M^{2n+1}$ 
the Gover-Graham operators $\boxdot_{\theta}^{(k)}$, $k=1,2,\ldots,n+1, n+2, n+4,\ldots$, are products of sublaplacians modulo 
lower order terms. Except for $k=n+1$ all the sublaplacians that are involved have invertible principal symbols, so
except for the value $k=n+1$ the principal symbol of $\boxdot_{\theta}^{(k)}$ is invertible (see Proposition~\ref{prop:Examples.Gover-Graham}). 

\subsubsection*{(e) Contact Laplacian} It has been shown by Rumin~\cite{Ru:FDVC} that in every degree the contact Laplacian satisfies the Rockland condition 
at every point, so it follows from Theorem~\ref{thm:Intro.Rockland-invertibility} that in every degree its principal symbol is invertible.

\section{Holomorphic families of \psivdos} 
In order to deal with complex powers of hypoelliptic operators we define holomorphic families of \psivdos\ and check their main properties in 
Chapter~\ref{chap.HolPHDO}. 
 
In a local Heisenberg chart $U\subset \Rd$ the definition of a holomorphic family of \psivdos\ parametrized by an open $\Omega\subset \C$ 
is similar to that of the definition of a holomorphic family 
of \psidos\ in~\cite[7.14]{Wo:LISA} and \cite[p.~189]{Gu:GLD} (see also~\cite{KV:GDEO}). In particular, we allow the order of the family of \psivdos\ to vary 
analytically. 

Most of the properties of \psivdos\ extend \emph{mutatis mutandis} to the setting of holomorphic families of \psivdos. 
In particular, the product of two holomorphic families of \psivdos\ is again a holomorphic family of \psivdos\ (Proposition~\ref{prop:HolPHDO.composition}). 
% 
% \begin{proposition}
%     For $j=1,2$ let $(P_{j,z})_{z \in \Omega}\subset \pvdo^{*}(M,\cE)$ be a holomorphic family of \psivdos\ and supposed that $(P_{1,z})_{z \in 
%     \Omega}$ or  $(P_{2,z})_{z \in \Omega}$ is uniformly properly supported with respect to $z$. Then the family of products $(P_{1,z}P_{2,z})_{z \in 
%     \Omega}$ is a holomorphic family of \psivdos. 
% \end{proposition}

There is, however, a difficulty when trying to extend the definition to general Heisenberg manifolds. More precisely, 
the proof of the invariance of the Heisenberg calculus by Heisenberg diffeomorphisms relies on a 
characterization of the distribution kernels of \psivdos\ by means of a suitable class of distributions 
$\cK^{*}(\URd)= \sqcup_{m\in \C}\cK^{m}(\URd)\subset \cD'(\URd)$. Each distribution $K \in \cK^{m}(\URd)$ 
admits an asymptotic expansion, in the sense of distributions, 
\begin{equation}
    K\sim \sum_{j\geq 0} K_{m+j}, \qquad K_{l}\in \cK_{l}(\URd),
\end{equation}
where  $\cK_{l}(\URd)$ consists of distributions that are smooth for $y\neq 0$ and homogeneous of degree $l$ if $l \not \in\N$ 
and are homogeneous of degree $l$ up to logarithmic terms otherwise (see~\cite{BG:CHM} and Chapter~\ref{chap:Heisenberg-calculus}). 
In particular, the definition of $\cK_{l}(\URd)$ depends upon whether $l$ is an 
integer or in not, which causes trouble for defining holomorphic families with values in $\cK^{*}(\URd)$ when the order crosses integers. 

This issue is resolved by means of a new description of the class $\cK^{*}(\URd)$ in terms of what we call \emph{almost homogeneous} 
distributions. The latter are homogenous modulo smooth terms and under the Fourier transform they correspond to the almost homogeneous 
symbols considered in~\cite{BG:CHM}. 

Since the definition of an almost homogeneous dsitribution of degree $l$ does not depend on whether $l$ is an integer or not, there is no trouble anymore 
to define holomorphic families of almost homogeneous kernels. Therefore, we can make use of the characterization of $\cK^{*}(\URd)$ in terms of almost 
homogenous distributions
to define holomorphic families with values in $\cK^{*}(\URd)$ (see~Definition~\ref{def:HolPHDO.kernels-families}).
% for the precise definition). 

We show that the distribution of kernel holomorphic families of \psivdos\ can be characterized in terms of holomorphic families with values in 
$\cK^{*}(\URd)$. This allows us to extend the arguments in the proof of the invariance by Heisenberg diffeomorphisms of the Heisenberg calculus 
to prove that holomorphic families of \psivdos\ too are invariant under Heisenberg diffeomorphisms (Proposition~\ref{prop:HolPHDO.invariance}). 
As a consequence we can define holomorphic families of \psivdos\ on an arbitrary Heisenberg manifold independently of the choice of a covering by 
Heisenberg charts. 

Let us also mention that the almost homogeneous approach to the Heisenberg calculus can also be used to constructing a class of \psivdos\ with 
parameter containing the resolvents of hypoelliptic \psivdos\ (see~\cite{Po:CPDE1}).  

\section{Heat equation and complex powers of hypoelliptic operators} 
% of hypoelliptic operators}
One of the main goals of this memoir is to obtain complex powers of hypoelliptic operators on  Heisenberg manifolds as holomorphic families of 
\psivdos\ along with some applications to hypoellipticity.  

It has been shown by Mohammed~\cite{Mo:ESOHCM2} that the complex powers of invertible positive hypoelliptic operators with
multicharacteristics are \psidos\ in the class constructed in~\cite{BdMGH:POPDCM}, but in the Heisenberg setting we would like to obtain them as 
holomorphic families of \psivdos. To this end we cannot follow the standard approach of Seeley~\cite{Se:CPEO} due to the lack of microlocality 
of the Heisenberg calculus. Instead we make use of the pseudodifferential representation of the heat kernel of~\cite{BGS:HECRM}, which is especially 
suitable for dealing with positive differential operators (we will deal with the general case in~\cite{Po:CPDE1} using another approach). 

Let us also mention that a similar approach to complex powers has been used independently by Mathai-Melrose-Singer~\cite{MMS:FAI} 
and Melrose~\cite{Me:SPLLB} in the context of projective pseudodifferential operators on Azamaya bundles. 

From now on we let $(M^{d+1},H)$ be a compact Heisenberg manifold equipped with a smooth density~$>0$ and let $\cE$ be a Hermitian vector bundle 
over $M$. 

\subsection{Pseudodifferential representation of the heat kernel}
Consider a selfadjoint differential operator $P:C^{\infty}(M,\cE)\rightarrow C^{\infty}(M,\cE)$ which is bounded from below and has an invertible 
principal symbol, so that the heat kernel $k_{t}(x,y)$ of $P$ is smooth for $t>0$.  

Recall that the heat semigroup $e^{-tP}$ allows us to invert the heat operator $P+\partial_{t}$. Conversely, 
constructing a suitable pseudodifferential calculus nesting parametrices for $P+\partial_{t}$ allows us to derive the small time heat kernel asymptotics for $P$. 

In the elliptic setting this approach was carried out by Greiner~\cite{Gr:AEHE} and 
the relevant pseudodifferential calculus is the Volterra calculus (see \cite{Gr:AEHE}, \cite{Pi:COPDTV}). 
The latter consists only in a modification of the classical pseudodifferential calculus 
in order to take into account the parabolicity and the Volterra property with respect to the time variable of the heat equation. In particular,  
Greiner's approach holds in fairly greater generality 
and has many applications~(see, e.g., \cite{BGS:HECRM}, \cite{BS:HEPNP1}, \cite{BS:HEPNP2}, \cite{Gr:AEHE}, \cite{Kr:PhD}, \cite{Kr:VFPDO}, \cite{KS:IPSPDE}, 
\cite{Me:APSIT}, \cite{Pi:COPDTV}, \cite{Po:CMP1}, \cite{Po:PAMS2}). 

The Greiner's approach has been extended to the Heisenberg calculus in~\cite{BGS:HECRM}, with the purpose of deriving the small time heat 
kernel asymptotics for the Kohn Laplacian on CR manifolds. In particular, a class of Volterra \psivdos\ is obtained which contains parametrices for the heat
operator $P+\partial_{t}$. As a consequence, once the principal symbol of $P+\partial_{t}$  is invertible in this calculus,
the inverse of $P+\partial_{t}$ is a Volterra \psivdo\ which, in turn, yields a pseudodifferential representation of the heat 
kernel of $P$.  More precisely, we have: 

\begin{theorem}[\cite{BGS:HECRM}]\label{thm:Intro1.heat}
Let $P:C^{\infty}(M,\cE)\rightarrow C^{\infty}(M,\cE)$ be a selfadjoint differential operator of even Heisenberg order $v$ which is bounded from below 
and such that the principal symbol of $P+\partial_{t}$ is invertible in the Volterra-Heiseneberg calculus. Then:\smallskip

1)  The inverse $(P+\partial_{t})^{-1}$ is a Volterra \psivdo;\smallskip 
   
   2) The heat kernel $k_{t}(x,y)$ of $P$ has an asymptotics in  $C^{\infty}(M,(\End \cE)\otimes|\Lambda|(M))$ of the form
    \begin{equation}
     k_{t}(x,x) \sim_{t\rightarrow 0^{+}} t^{-\frac{d+2}{v}} \sum t^{\frac{2j}{v}} a_{j}(P)(x),
    \label{eq:Intro1.Rockland-Heat.heat-kernel-asymptotics}
    \end{equation}
where the density $a_{j}(P)(x)$ is locally computable in terms of the symbol $q_{-v-2j}(x,\xi,\tau)$ of degree $-v-2j$ of any Volterra-\psivdo\ parametrix 
for $P+\partial_{t}$.
  \end{theorem}

This framework is recalled in Section~\ref{sec.Volterra-PsiHDO-calculus} and we can extend to this setting the intrinsic approach of 
Chapter~\ref{chap:Heisenberg-calculus}. 

\subsection{Heat equation and sublaplacians} 
Let $\Delta:C^{\infty}(M,\cE)\rightarrow C^{\infty}(M,\cE)$ be a selfadjoint sublaplacian which is bounded from below. 
In~\cite{BGS:HECRM} the authors construct explicitly an inverse in the Volterra-Heisenberg calculus for the principal symbol of 
$\Delta+\partial_{t}$ when $\cE$ is the trivial line bundle and when at every point $a \in M$,  with the 
notation of~(\ref{eq:Intro.singular-set1})--(\ref{eq:Intro.singular-set2}), we have
\begin{equation}
    |\mu(a)| < \frac12 \sum_{j=1}^{d}|\lambda_{j}|. 
%     \label{eq:¥}
\end{equation}
Since $\Delta$ is selfadjoint, and so $\mu(a)$ is real, 
the above condition is the same as~(\ref{eq:Intro.sublaplacian.condition}) 
when $\rk \cL_{a}<d$, but when $\rk \cL_{a}=d$ this is a stronger condition. 
% Moreover, as in~\cite{BG:CHM} 
% this is done by explicitly constructing the inverse of the principal symbol of $\Delta+\partial_{t}$. 

In fact, the explicit formulas of~\cite{BGS:HECRM} can be extended to the case where $\cE$ is an arbitrary vector bundle and where $\Delta$ satisfies the 
weaker condition~(\ref{eq:Intro.sublaplacian}) at every point (see Proposition~\ref{thm:Powers1.heat-sublaplacians}). 
As a consequence Theorem~\ref{thm:Intro1.heat} holds for the Kohn Laplacian even when the Levi form is not diagonalizable.

Moreover, as we actually can invert the principal symbol of $\Delta+\partial_{t}$ in a refined class of symbols, for any integer $k=2,3,\ldots$  
we can invert the principal 
symbol of $\Delta^{k}+\partial_{t}$ in the Volterra-Heisenberg calculus (see~Proposition~\ref{prop:Heat1.Delta^k}). 
% In particular, Theorem~\ref{thm:Intro1.heat} is valid for the  of the horizontal sublaplacian on a strictly pseudoconvex CR manifold as defined 
% in~\cite{GG:CRIPSL}.  
% 
% On the other hand, the pseudodifferential of the heat kernel in~\cite{BGS:CHM} is powerful enough to yield for almost free the small time heat kernel 
% asymptotics for hypoelliptic operators. It has been shown in~\cite{BGS:CHM} could be applied to (selfadjoint) sublaplacians. We also show here that, 
% partly as a consequence of our results on complex powers. 
% 
% using a different approach using a new kind of pseudodifferential representation 
% of the resolvent. 
% 
% the main aims of this memoir is to realize complex powers of hypoelliptic operators as holomorphic families of \psivdos. As alluded to above we cannot carry out the 
% standard approach of Seeley~\cite{Se:CPEO} in the Heisenberg setting. Therefore, we have to rely on another approach based on the pseudodifferential 
% representation of the heat kernel of~\cite{BGS:HECRM}. 

\subsection{Complex powers}
% we rely on a new approach based on the the pseudodifferential representation 
% in~\cite{BGS:HECRM} of the heat kernel of a hypoelliptic differential operator.  
% 
% In this paper, and its sequel~\cite{Po:CPDE1}, we aim to get Unfortunately, 
% the non-microlocality of the composition of principal symbols does not allow us to extend to the Heisenberg setting 
% the standard approach of Seeley~\cite{Se:CPEO} (see also~\cite{Sh:POST}, \cite{Gr:FCPDBP}) to complex powers of elliptic operators. Instead, we rely on 
% new approaches. In the present paper, we deal with positive hypoelliptic differential operators using an approach based 
% on the pseudodifferential representation of the heat kernel of a hypoelliptic operator in~\cite{BGS:HECRM}. In~\cite{Po:CPDE1} we will deal with general 
% hypoelliptic \psivdo\ using another approach. 
% 
% On the other hand, the lack of microlocality of the Heisenberg calculus does not allow us to carry out in the Heisenberg setting the standard approach 
% of Seeley~\cite{Se:CPEO} to complex powers of elliptic operators. Instead, we rely on two new approaches. The first one, used in this paper, 
% is quite suitable for positive differential operators and similar approaches have also been used  by 
% Mathai-Melrose-Singer~\cite{MMS:FAI} and Melrose~\cite{Me:SPLLB} in the context of projective pseudodifferential operators on Azamaya bundles. 
% The second approach is  carried out in~\cite{Po:CPDE1} and allows us to deal with general \psivdos.
Let $P:C^{\infty}(M,\cE)\rightarrow C^{\infty}(M,\cE)$ be a positive selfadjoint differential operator of even Heisenberg order $v$ and assume that the principal 
symbol of $P$ is invertible, i.e., $P$ satisfies the Rockland condition at every point. Thanks to the spectral theorem we can define the complex powers 
$P^{s}$, $s\in \C$,  of $P$ as unbounded operators on $L^{2}(M,\cE)$ 
which are bounded for $\Re s\leq 0$. 

Moreover, for $\Re s<0$ the Mellin formula holds, 
\begin{equation}
    P^{s}=\Gamma(s)^{-1}\int_{0}^{\infty}t^{s}(1-\Pi_{0}(P))e^{-tP}\frac{dt}{t},
    \label{eq:Intro1.Mellin-formula}
\end{equation}
where $\Pi_{0}(P)$ denotes the orthogonal projection onto the kernel of $P$. 
Combining this  formula with the pseudodifferential representation of the heat kernel of $P$ in terms of the Volterra-Heisenberg calculus allows us to prove: 

\begin{theorem}\label{thm:Intro1.complex-powers}
Assume that the principal symbol of $P+\partial_{t}$ is an invertible Volterra-Heisenberg symbol. Then the 
complex powers  $P^{s}$, $s\in \C$, of $P$ form a holomorphic 1-parameter group of \psivdos\ such that $\ord P^{s}=ms$ $\forall s\in \C$.  
\end{theorem}
 % In~\cite{BGS:HECRM} it was shown that for a selfadjoint sublaplacian $\Delta:C^{\infty}(M,\cE)\rightarrow C^{\infty}(M,\cE)$
% the principal symbol of $\Delta+\partial_{t}$ is invertible under a condition closely related to the Rockland condition and the invertibility of the 
% principal symbol of $\Delta$ (see~\cite[5.23]{BGS:HECRM} and Section~\ref{sec.powers1}). Therefore, in the case of a sublaplacian we get: 
% 
% \begin{theorem}\label{thm:Intro1.complex-powers-sublaplacian}
%     Let $\Delta:C^{\infty}(M,\cE)\rightarrow C^{\infty}(M,\cE)$ be a positive sublaplacian satisfying the condition of~\cite[5.23]{BGS:HECRM}. 
%     Then  the family $(\Delta^{s})_{s\in \C}$ of the complex powers of $\Delta$ is a holomorphic 1-parameter group 
%     of \psivdos\ such that $\ord \Delta^{s}=2s$ for any $s\in \C$.  
% \end{theorem}

In particular, this theorem holds for the following sublaplacians:\smallskip 

(a) A sum of squares of the form~(\ref{eq:Intro.sum-of-squares}), 
% $\nabla_{X_{1}}^{*}\nabla_{X_{1}}+\ldots+\nabla_{X_{m}}^{*}\nabla_{X_{m}}$, 
% where the vector fields $X_{1},\ldots,X_{m}$ linearly span $H$ and $\nabla$ is a connection on $\cE$, 
provided that the bracket condition~(\ref{eq:Intro.bracket-condition}) holds;\smallskip 

(b) The Kohn Laplacian on a CR manifold acting on $(p,q)$-forms under condition $Y(q)$;\smallskip 

(c) The horizontal Laplacian on a Heisenberg manifold acting on horizontal forms of degree $k$  under condition $X(k)$;\smallskip 

(d) The horizontal Laplacian on a CR manifold acting  on $(p,q)$-forms under condition $X(p,q)$.\smallskip 

In fact, partly by making use of Theorem~\ref{thm:Intro1.complex-powers},
we will show that when the  condition~(\ref{eq:Intro.bracket-condition}) holds the principal symbol of $P+\partial_{t}$ is 
automatically invertible in the Volterra-Heisenberg symbol (see below). Therefore, we obtain: 

\begin{theorem}\label{thm:Intro1.complex-powers-Rockland}
    If the bracket condition~(\ref{eq:Intro.bracket-condition})  holds then the 
complex powers  $P^{s}$, $s\in \C$, of $P$ form a holomorphic 1-parameter group of \psivdos\ such that $\ord P^{s}=ms$ $\forall s\in \C$.  
\end{theorem}

In particular, Theorem~\ref{thm:Intro1.complex-powers-Rockland} is valid for the contact Laplacian on a contact manifold. In this context 
this allows us to fill a technical gap in~\cite{JK:OKTGSU} concerning the proof of the fact that the complex powers of the contact Laplacian give rise to 
\psivdos\  which is an important step in the proof there of the Baum-Connes conjecture for $SU(n,1)$ (see~\cite{Po:Crelle1}). 

\subsection{Rockland condition and heat equation} 
Theorem~\ref{thm:Intro1.complex-powers} 
has several interesting applications related to hypoellipticity.
% and in particular we can make use of it in connection with the Rockland condition when the bracket 
% condition $H+[H,H]=TM$  holds.  

First, Theorem~\ref{thm:Intro.Rockland-invertibility} can be extended to \psivdos\ with non-integer orders as follows. 

\begin{theorem}\label{thm:Intro.Heat-Rockland.invertibility-non-integer-order}
 Assume that the bracket condition~(\ref{eq:Intro.bracket-condition}) holds. Then for any $P\in \pvdo^{m}(M,\cE)$, $m\in \C$, 
%  and let $P:C^{\infty}_{c}(M,\cE)\rightarrow C^{\infty}(M,\cE)$ be a \psivdo\ of order $m \in 
%  \C$. Then 
 the following are equivalent:\smallskip
  
  (i) The principal symbol of $P$ is invertible;\smallskip
  
  (ii) $P$ and $P^{t}$ satisfy the Rockland condition at every point $a\in M$.\smallskip

  (iii) $P$ and $P^{*}$ satisfy the Rockland condition at every point $a\in M$.
\end{theorem}

As a consequence of this theorem we can prove that when the condition~(\ref{eq:Intro.bracket-condition})  
holds any $P\in \pvdo^{m}(M,\cE)$ with $\Re m \geq 0$ satisfying the 
Rockland condition at every point is hypoelliptic with gain of $\frac{1}{2}\Re m$ derivative(s) (see Proposition~\ref{prop:Heat-Rockland.Rockland-hypoellipticity}). 

Next, let $P:C^{\infty}(M,\cE)\rightarrow C^{\infty}(M,\cE)$ be a selfadjoint differential operator of even Heisenberg order $v$. We shall say that 
the principal symbol of $P$ is positive when it can be put into the form $\overline{q_{\frac{v}{2}}}*q_{\frac{v}{2}}$ for some symbol 
$q_{\frac{v}{2}}$ homogeneous of degree $\frac{v}{2}$. Then, by making use of Theorem~\ref{thm:Intro1.complex-powers} 
and by extending to the Volterra-Heisenberg setting the 
arguments of the proof of Theorem~\ref{thm:Intro.Rockland-invertibility} we prove: 

\begin{theorem}\label{thm:Intro.Rockland-heat}
 Assume that the  bracket condition~(\ref{eq:Intro.bracket-condition}) holds and that $P$ satisfies the Rockland condition at every point.\smallskip
 
 1) $P$ is bounded from below if, and only if, it has a positive principal symbol.\smallskip
 
 2) If $P$ has a positive principal symbol, then the principal symbol $P+\partial_{t}$ is invertible in the Volterra-Heisenberg calculus. 
% $S_{\vo,-v}(\fg^{*}M\times \R_{(v)},\cE)$.  
%  Hence Theorems~\ref{thm:volterra.inverse} and~\ref{thm:Powers1.heat-kernel-asymptotics} hold for~$P$.  
\end{theorem}

This proves that, when the condition~(\ref{eq:Intro.bracket-condition})  holds, in Theorem~\ref{thm:Intro1.heat} we can replace the invertibility 
condition on the principal symbol of $P+\partial_{t}$ by the validity of the Rockland condition for $P$ at every point. 
Consequently, we see that the results of~\cite{BGS:HECRM} actually hold for a wide class of operators. In particular, 
Theorem~\ref{thm:Intro1.heat} is valid for the contact Laplacian on a contact manifold.

\subsection{Weighted Sobolev spaces}
As another application of Theorem~\ref{thm:Intro1.complex-powers}, under the bracket condition~(\ref{eq:Intro.bracket-condition}) 
we can construct a scale of Weighted Sobolev spaces $W_{H}^{s}(M,\cE)$, $s \in \C$, providing us with sharp regularity estimates for \psivdos. 

Let $\Delta_{\nabla,X}:C^{\infty}(M,\cE)\rightarrow C^{\infty}(M,\cE)$ 
be a sum of squares as in~(\ref{eq:Intro.sum-of-squares}).
% Consider a sum of squares, 
% \begin{equation}
%     \Delta_{\nabla,X}= \nabla_{X_{1}}^{*} \nabla_{X_{1}} +\ldots+ \nabla_{X_{m}}^{*} \nabla_{X_{m}}, 
% %     \label{eq:¥}
% \end{equation}
% where the vector fields $X_{1},\ldots,X_{m}$ span $H$ and $\nabla$ is a connection on $\cE$. 
Since  the bracket condition~(\ref{eq:Intro.bracket-condition}) holds,  
Theorem~\ref{thm:Intro1.complex-powers} tells us that the complex powers $(1+\Delta_{X})^{s}$, $s \in \C$,  give rise to an analytic 1-parameter group 
of invertible \psivdos. 

For $s \in \R$ the weighted Sobolev space $W_{H}^{s}(M,\cE)$ is defined as the space of distributional sections $u \in 
\cD'(M,\cE)$ such that $(1+\Delta_{\nabla,X})^{\frac{s}{2}}u$ is in $L^{2}(M,\cE)$ together with the Hilbertian norm, 
\begin{equation}
    \|u\|_{W_{H}^{s}}=\|(1+\Delta_{\nabla, X})^{\frac{s}{2}}u\|_{L^{2}}, \qquad u \in W_{H}^{s}(M,\cE). 
\end{equation}
It can be shown that, up to the choice of an equivalent Hilbertian norm, this definition does not depend on the choices of the vector fields 
$X_{1},\ldots,X_{m}$ and of the connection $\nabla$ and that when $s$ is a positive integer it agrees with the previous definition of 
the Weighted Sobolev spaces of Folland-Stein~\cite{FS:EDdbarbCAHG} (see Section~\ref{sec.Sobolev}). 

Moreover, the spaces $W_{H}^{s}(M,\cE)$ can be nicely compared to the standard Sobolev spaces $L^{2}_{s}(M,\cE)$. More precisely, we show that we have 
the following continuous embeddings, 
   \begin{equation}
   \begin{array}{rcl}
       L^{2}_{s}(M)  & \hookrightarrow W_{H}^{s}(M) \hookrightarrow L^{2}_{s/2}(M) & \qquad \text{if $s\geq 0$},\\
       L^{2}_{s/2}(M)  & \hookrightarrow W_{H}^{s}(M) \hookrightarrow L^2_{s}(M) & \qquad \text{if $s< 0$}.
   \end{array}
         \label{eq:Intro1.embeddings}
    \end{equation}

On the other hand, these Sobolev spaces are suitable for studying \psivdos, for we have: 

\begin{proposition}\label{prop:Intro.Sobolev.regularity-PsiHDOs1}
    Let $P:C^{\infty}(M,\cE)\rightarrow C^{\infty}(M,\cE)$ be a \psivdo\ of order $m$ and set $k=\Re m$. 
    Then, for any $s\in \R$, the operator $P$ extends to a continuous linear mapping from 
    $W_{H}^{s+k}(M,\cE)$ to $W_{H}^{s}(M,\cE)$. 
\end{proposition}

As a consequence we get sharp regularity results for \psivdos\ satisfying the Rockland condition: 

\begin{proposition}\label{prop:Intro.Sobolev.regularity-PsiHDOs2}
  Let $P:C^{\infty}(M,\cE)\rightarrow C^{\infty}(M,\cE)$ be a \psivdo\ of order $m$ such that $P$ satisfies the Rockland condition at every point and 
  set $k=\Re m$. 
  Then for any $u \in \cD'(M,\cE)$ we have
    \begin{equation}
        Pu \in W_{H}^{s}(M,\cE) \Longrightarrow u \in W_{H}^{s+k}(M,\cE).
         \label{eq:intro.Sobolev.hypoellipricity-WHs}
    \end{equation}
   In fact, for any $s'\in \R$ we have the estimate,
    \begin{equation}
        \|u\|_{W^{s+k}_{H}} \leq C_{ss'}(\|Pu\|_{W^{s}_{H}}+\|u\|_{W^{s'}_{H}}), \qquad u \in W^{s+k}_{H}(M,\cE).
        \label{eq:Intro.Sobolev.hypoelliticity-WHs-estimates}
    \end{equation}
\end{proposition}

When $P$ is a differential operator of Heisenberg order $v$ the 
properties~(\ref{eq:intro.Sobolev.hypoellipricity-WHs}) and~(\ref{eq:Intro.Sobolev.hypoelliticity-WHs-estimates}) 
correspond to the maximal hypoellipticity of~\cite{HN:HMOPCV}. 

% Therefore, if we extend the definition of maximal hypoellipticty in terms of~(\ref{eq:intro.Sobolev.hypoellipricity-WHs}) 
% and~(\ref{eq:Intro.Sobolev.hypoelliticity-WHs-estimates})   then 
% Proposition~\ref{prop:Intro.Sobolev.regularity-PsiHDOs2}  
% says that any \psivdo\ satisfying the Rockland condition at every point is maximal hypoelliptic.

In addition, the weighted Sobolev spaces $W_{H}^{s}(M,\cE)$ can be localized, so that it makes sense to say that a distributional section is $W_{H}^{s}$  
near a point, and we can prove a localized version of Proposition~\ref{prop:Intro.Sobolev.regularity-PsiHDOs2} 
(see Proposition~\ref{prop:Sobolev.hypoellipricity-WHs-localized}). 

Finally, we also give a version of Proposition~\ref{prop:Intro.Sobolev.regularity-PsiHDOs1} 
for holomorphic families of \psivdos\ and in particular for complex powers of positive 
differential operator satisfying the Rockland condition (Propositions~\ref{prop:Sobolev.regularity-PsiHDOs-families} 
and~\ref{prop:Sobolev.regularity-complex-powers} for the detailed statements).

\section{Spectral asymptotics for hypoelliptic operators} 
Another main goal of this monograph is to make use of the Heisenberg calculus to derive spectral asymptotics for hypoelliptic operators on 
Heisenberg manifolds and in particular to get explicit geometric expressions for the leading terms of these asymptotics for the main geometric differential operators on 
CR and contact manifolds. 
% that is, the Kohn Laplacian, the horizontal sublaplacian and its  and the contact Laplacian. 

% already known,  (see below) 
% 
% 
% interesting application of the heat kernel asymptotics in~\cite{BGS:HECRM} and in~(\ref{eq:Intro1.Rockland-Heat.heat-kernel-asymptotics}) 
% is to allows us to derive spectral asymptotics for hypoelliptic operators as follows. 

\subsection{Heat equation and spectral asymptotics}
Consider a selfadjoint differential operator  $P:C^{\infty}(M,\cE)\rightarrow C^{\infty}(M,\cE)$ of even Heisenberg order $v$ which is bounded from 
below and such that the principal symbol of $P+\partial_{t}$ is  invertible in the Volterra-Heisenberg calculus. Then 
% by Theorem~\ref{thm:Intro.Rockland-invertibility} 
% when the bracket 
% condition $H+[H,H]=TM$ holds it is enough to assume that $P$ satisfies the Rockland condition at every point. 
 the heat kernel asymptotics~(\ref{eq:Intro1.Rockland-Heat.heat-kernel-asymptotics}) holds at the level of densities, so that 
 as $t\rightarrow 0^{+}$ we have 
 \begin{equation}
         \Tr e^{-tP} \sim t^{-\frac{d+2}{m}} \sum t^{\frac{2j}{m}} A_{j}(P), \qquad  
   A_{j}(P)=\int_{M}\tr_{\cE} a_{j}(P)(x).
          \label{eq:Intro1.heat-trace-asymptotics}
     \end{equation}
%  where the density $a_{j}(P)(x)$ is the coefficient of $t^{\frac{j-d+2}{m}}$ in the 
%  asymptotics~(\ref{eq:Rockland-Heat.heat-kernel-asymptotics}). 

Next, let $\lambda_{0}(P)\leq \lambda_{1}(P)\leq \ldots$ denote the eigenvalues of $P$ counted with multiplicity and  
let $N(P;\lambda)$ denote its counting function, that is, 
\begin{equation}
     N(P;\lambda)=\#\{k\in \N;\ \lambda_{k}(P)\leq \lambda \}, \qquad \lambda\geq 0.
\end{equation}
In addition, define
\begin{equation}
    \nu_{0}(P)=\Gamma(1+\frac{d+2}{m})^{-1}A_{0}(P).
    \label{eq:Intro1.counting-function-asymptotics}
\end{equation}
Then we obtain: %the following Weyl asymptotics:  

\begin{proposition}\label{thm:Intro1.spectral-asymptotics}
1) We have $\nu_{0}(P)>0$.\smallskip 
     
2) As $\lambda\rightarrow \infty$ we have $N(P;\lambda) \sim \nu_{0}(P)\lambda^{\frac{d+2}{m}}$.\smallskip 

3) As $k\rightarrow \infty$ we have $\lambda_{k}(P)\sim \left(\frac{k}{\nu_{0}(P)}\right)^{\frac{m}{d+2}}$. 
 \end{proposition}

Once it is proved that $\nu_{0}(P)$ is~$>0$ we can make use of Karamata's Tauberian theorem to deduce from~(\ref{eq:Intro1.heat-trace-asymptotics}) 
the asymptotics for $N(P;\lambda)$ and $\lambda_{k}(P)$. Thus the bulk the proof is to establish the positivity of $\nu_{0}(P)$, which is carried out
via spectral theoretic considerations.

 By relying on other pseudodifferential calculi several authors have also obtained Weyl asymptotics
in the more general setting of hypoelliptic operators with multicharacteristics~(see \cite{II:PDPEAABSFSP}, \cite{Me:HOCVC2WE}, \cite{MS:ECH2}, 
\cite{Mo:ESOHCM1}, \cite{Mo:ESOHCM2}). Nevertheless, as 
far as the Heisenberg setting is concerned, the approach using the Volterra-Heisenberg calculus has two main advantages. 
 
First, the pseudodifferential analysis is significantly simpler. In particular, the Volterra-Heisenberg calculus yields for free 
 the heat kernel asymptotics once the principal symbol of the heat operator is shown to be invertible, for which it is enough to use the Rockland condition 
 when the condition~(\ref{eq:Intro.bracket-condition}) holds. 
 
 Second, since the Volterra-Heisenberg calculus fully takes into account the underlying Heisenberg geometry of the manifold and is invariant 
 by change of Heisenberg coordinates, we can get explicit geometric expressions for the coefficient $\nu_{0}(P)$ in the case of the main geometric differential 
 operators on CR and contact manifolds (see below for the precise formulas). 
%  we can very effectively deal with operators admitting normal forms (see Proposition~\ref{prop:Spectral.normal-form-nu0P} on this point). 
%  In particular, as explained below, we can reformulate in geometric terms the Weyl asymptotics~(\ref{eq:Intro1.counting-function-asymptotics}) for the main 
%  differential operators on CR and contact manifolds. 

\section{Weyl asymptotics and CR geometry} 
Let $M^{2n+1}$ be a compact $\kappa$-strictly pseudoconvex CR manifold and let $\theta$ be a contact form whose associated Levi form has 
signature $(n-\kappa,\kappa,0)$, so that $\theta$ defines a pseudohermitian structure on $M$. We endow $M$ with a Levi metric compatible with $\theta$.  
Then the volume of $M$ with respect to this Levi metric is independent of the choice of the Levi form and is equal to 
\begin{equation}
        \op{vol}_{\theta}M=\frac{(-1)^{\kappa}}{n!} \int_{M} \theta \wedge d\theta^{n}.
\end{equation}
We call $\op{vol}_{\theta}M$ the pseudohermitian volume of $(M,\theta)$ and we relate it to the Weyl asymptotics~(\ref{eq:Intro1.counting-function-asymptotics}) 
for the Kohn Laplacian and the horizontal sublaplacian as follows. 

For $\mu\in (-n,n)$ we let 
\begin{equation}
    \nu(\mu)= (2\pi)^{-(n+1)} \int_{-\infty}^{\infty}e^{-\mu\xi_{0}}(\frac{\xi_{0}}{\sinh \xi_{0}})^{n}d\xi_{0}.
\end{equation}
Then for the Kohn Laplacian  we prove: 
\begin{theorem}
 Let $\Boxbpq$ be the Kohn Laplacian
  acting on $(p,q)$ forms with  $q\neq \kappa$ and $q\neq n-\kappa$. Then  
    as $\lambda \rightarrow \infty$ we have 
\begin{equation}
        N(\Boxbpq ;\lambda) \sim \alpha_{n\kappa pq}(\op{vol}_{\theta}M)\lambda^{n+1},
%     \label{eq:}
\end{equation}
where $ \alpha_{n\kappa pq}$ is equal to  
\begin{equation}
   \binom{n}{p} \sum_{\max(0,q-\kappa)\leq  k\leq \min(q,n-\kappa)} \frac{1}{2} \binom{n-\kappa}{k}\binom{\kappa}{q-k} 
        \nu(n-2(\kappa-q+2k)).
         \label{eq:Intro.Heat1.alphapq}  
%     \label{eq:}
\end{equation}
In particular $\alpha_{n\kappa pq}$ is a universal constant depending only on $n$, $\kappa$, $p$ and $q$.
\end{theorem}
In the strictly pseudoconvex case, i.e., when $\kappa=0$, this theorem follows from the computation of $A_{0}(\Boxbpq)$ 
in~\cite{BGS:HECRM}, but for the case $\kappa \geq 1$ 
this seems to be a new result. 

Next, in the CR setting the horizontal sublaplacian preserves the bidegree and, in the same way as with the Kohn Laplacian, we prove: 

\begin{theorem}
 Let $\Delta_{b;p,q}:C^{\infty}(M,\Lambda^{p,q})\rightarrow C^{\infty}(M,\Lambda^{p,q})$ be the  horizontal sublaplacian acting on $(p,q)$-forms with  
  $(p,q)\neq (\kappa,n-\kappa)$ and $(p,q)\neq (n-\kappa,\kappa)$. Then as $\lambda  \rightarrow \infty$ we have 
\begin{equation}
             N(\Delta_{b;p,q};\lambda) \sim \beta_{n\kappa pq}(\op{vol}_{\theta}M) \lambda^{n+1},
    \label{eq:Intro1.Weyl-Deltab-CR}
\end{equation}
 where $\beta_{n\kappa pq}$ is equal to
\begin{equation}
        \!  \!  \!  \!   \sum_{\substack{\max(0,q-\kappa)\leq  k\leq \min(q,n-\kappa)\\ \max(0,p-\kappa)\leq l\leq \min(p,n-\kappa)}}  \!  \!  \!  \!
        2^{n}\binom{n-\kappa}{l}\binom{\kappa}{p-l} \binom{n-\kappa}{k}\binom{\kappa}{q-k} 
        \nu(2(q-p)+4(l-k)).
         \label{eq:Intro.Heat1.betapq}
%     \label{eq:}
\end{equation}
In particular $\beta_{n\kappa pq}$  is a universal constant depending only on $n$, $\kappa$, $p$ and $q$.
\end{theorem}

 Finally, suppose that $M$ is strictly pseudoconvex, i.e., $\kappa=0$, and for $k=1,\ldots,n+1,n+2, n+4,\ldots$ let 
 $\boxdot_{\theta}^{(k)}$ be the Gover-Graham operator of order $k$. Then we have:
 
\begin{theorem}
 Assume $k\neq n+1$. Then there exists a  universal constant $\nu_{n}^{(k)}>0$ depending only on $n$ and $k$  such that as $\lambda 
     \rightarrow \infty$ we have 
     \begin{equation}
          N(\boxdot_{\theta}^{(k)};\lambda) \sim  \nu_{n}^{(k)}(\op{vol}_{\theta}M) \lambda^{\frac{n+1}{k}}. 
     \end{equation}
\end{theorem}
% we also deal with the Gover-Graham operators. We obtain: 
% as the value of $\nu_{0}(P)$ in~(\ref{eq:Intro1.counting-function-asymptotics}) depends only on the principal symbol of $P$, in the 
% % strictly pseudoconvex case we also can deal with the  of the horizontal Laplacian.
% % 
%  \begin{theorem}
%     Assume $M$ is strictly pseudoconvex and 
%     for $k=1,2,\ldots,n+1$ let $\boxdot_{\theta}^{(k)}:C^{\infty}(M)\rightarrow C^{\infty}(M)$ be a $k$'th conformal power of the sublaplacian. 
%     Then as $\lambda  \rightarrow \infty$ we have
%     \begin{equation}
%         N(\boxdot_{\theta}^{(k)};\lambda) \sim \nu(0) (\op{vol}_{\theta}M) \lambda^{\frac{n+1}{k}}.
% % %        \label{eq:}
%     \end{equation}
% \end{theorem}

 \section{Weyl asymptotics and contact geometry}
 Let $(M^{2n+1},H)$ be a compact orientable contact manifold. Let $\theta$ be a contact form and let $J$ be a calibrated almost complex structure on 
 $H$ so that  $d\theta(X,JX)=-d\theta(JX,X)>0$ for any section $X$ of $H\setminus 0$. We then endow $M$ with the Riemannian metric 
 $g_{\theta, J}=d\theta(.,J.)+\theta^{2}$. The volume of $M$ with respect to $g_{\theta,J}$ depends only on $\theta$ and is equal to:
 \begin{equation}
     \op{vol}_{\theta}M=\frac{1}{n!}\int_{M}d\theta^{n}\wedge \theta.
 \end{equation}
 We call $\op{vol}_{\theta}M$ the contact volume of $M$. 
  
 We can relate the Weyl asymptotics for the horizontal sublaplacian to the contact volume to get:
 
 \begin{theorem}
    Let $\Delta_{b;k}:C^{\infty}(M,\Lambda^{k}_{\C}H^{*})\rightarrow C^{\infty}(M,\Lambda^{k+1}_{\C}H^{*})$ be the horizontal sublaplacian on $M$ in 
    degree $k$ with $k \neq n$. Then as $\lambda  \rightarrow \infty$ we have
    \begin{gather}
         N(\Delta_{b;k};\lambda) \sim \gamma_{nk} (\op{vol}_{\theta}M) \lambda^{n+1},\qquad 
          \label{eq:Intro.Spectral.Weyl-Deltab-contact}
 \gamma_{nk}=\sum_{p+q=k}2^{n}\binom{n}{p} \binom{n}{q}\nu(p-q).         
    \end{gather}
In particular $\gamma_{nk}$ is universal constant depending on $n$ and $k$ only.
\end{theorem}
 
Note that  when $M$ is a strictly pseudoconvex CR manifold the asymptotics~(\ref{eq:Intro.Spectral.Weyl-Deltab-contact}) 
is compatible with~(\ref{eq:Intro1.Weyl-Deltab-CR}) because the contact volume differs from the pseudohermitian volume by a factor of $2^{-n}$.

Finally, we can also deal with the contact Laplacian as follows. 
\begin{theorem}
    1)  Let  $\Delta_{R;k}:C^{\infty}(M,\Lambda^{k})\rightarrow C^{\infty}(M,\Lambda^{k})$ be the contact Laplacian in degree $k$ with 
    $k \neq n$. Then  
    there exists a  universal constant $\nu_{nk}>0$ depending only on $n$ and $k$ 
    such that as $\lambda \rightarrow \infty$ we have
     \begin{equation}
           N(\Delta_{R;k};\lambda)\sim \nu_{nk} (\op{vol}_{\theta}M)\lambda^{n+1}. 
                \label{eq:Intro.Spectral.Weyl-contact-Laplacian1}
      \end{equation}
    
    2) For $j=1,2$ consider the contact Laplacian $\Delta_{R;n}:C^{\infty}(M,\Lambda^{n}_{j})\rightarrow C^{\infty}(M, \Lambda^{n}_{j})$.  
    Then there exists a  universal constant $\nu_{n}^{(j)}>0$ depending only on $n$ and $j$  such that as $\lambda \rightarrow \infty$ we have 
    \begin{equation}
      N(\Delta_{R;nj};\lambda)\sim \nu_{n}^{(j)} (\op{vol}_{\theta}M)\lambda^{\frac{n+1}{2}}. 
         \label{eq:Intro.Spectral.Weyl-contact-Laplacian2}
    \end{equation}
\end{theorem}

\section{Organization of the memoir} The rest of the memoir is organized as follows. In Chapter~\ref{chap:Heisenberg-manifolds-operators}
 we start by recalling the main  definitions and examples concerning Heisenberg manifolds and their tangent Lie group bundles. Then we review the constructions 
 of the main differential operators on Heisenberg manifolds: sum of squares, Kohn 
Laplacian, horizontal sublaplacian, Gover-Graham operators and the contact Laplacian.

In Chapter~\ref{chap:Heisenberg-calculus} 
after a detailed review of the main known facts about the Heisenberg calculus, we give an intrinsic definition of the principal symbol and model 
operators of a \psivdo\ and check their main properties. Then we prove Theorem~\ref{thm:Intro.Rockland-invertibility} and its consequences. 
We conclude the chapter by a closer look at the main differential operators on Heisenberg manifolds.  

In Chapter~\ref{chap.HolPHDO} we define holomorphic families of \psivdo\ and study their main properties. In particular, we make use of an almost homogeneous 
approach to the Heisenberg calculus. 
% 
% 
% n Section~\ref{sec:PsiHDO} we give a detailed overview of the Heisenberg 
% calculus of~\cite{BG:CHM} and~\cite{Ta:NCMA}, following closely the expositions of~\cite{BG:CHM} and~\cite{Po:BSM1}. 

Chapter~\ref{chap.complex-powers} is devoted to complex powers of positive hypoelliptic differential operators in connection with the heat 
equation. After having recalled the pseudodifferential represention of the heat kernel of such an operator in terms of the Volterra-Heisenberg calculus 
of~\cite{BGS:HECRM}, we use it to establish Theorem~\ref{thm:Intro1.complex-powers}.  
Then we make use of Theorem~\ref{thm:Intro1.complex-powers} to extend Theorem~\ref{thm:Intro.Rockland-invertibility} to \psivdos\ with 
non-integer orders  and to prove Theorem~\ref{thm:Intro.Rockland-heat}. Eventually, we construct the weighted Sobolev spaces $W_{H}^{s}(M,\cE)$, $s\in 
\R$,  and check their main properties. In particular, we prove that they yield sharp regularity results for \psivdos. 

% In Section~\ref{sec.powers1}, after 
% having recalled the pseudodifferential representation of the heat kernel of~\cite{BGS:HECRM}, we make use of the latter to 
% construct complex powers of 
% positive hypoelliptic operators as holomorphic families of \psivdos.  
% In Section~\ref{sec.Sobolev} we define the weighted Sobolev spaces $W_{H}(M,\cE)$ and show that they give 
% sharp regularity results for \psivdos. In Section~\ref{sec:Rockland-heat} we show that if $P$ is a selfadjoint differential operator bounded from below, then 
% the fact that $P$ satisfies the Rockland condition at every point is 
% enough to insure us the invertibility of the principal symbol of $P+\partial_{t}$. This show that the results of~\cite{BGS:HECRM} and Section~\ref{sec.powers1} 
% are valid for a wide class of hypoelliptic operators. 

In Chapter~\ref{chap:Spectral}, 
we deal with spectral asymptotics for hypoelliptic operators on Heisenberg manifolds. First, we derive 
general spectral asymptotics for such operators on a general Heisenberg manifold. We then express these asymptotics in a geometric fashion. We first 
proceed with the Kohn Laplacian and the horizontal sublaplacian on a CR manifold. Then we deal with
the horizontal sublaplacian and the contact Laplacian on a contact manifold. 

Finally, two appendices are included. In Appendix~\ref{chap.Appendix-Invariance} 
we give a version in Heisenberg coordinates of the proof of the invariance by Heisenberg 
diffeomorphisms of the Heisenberg calculus, which is used in the intrinsic definition of the principal symbol of a \psivdo. 
In Appendix~\ref{chap:Appendix-transpose} we 
similarly give a version in Heisenberg coordinates of the proof that the transpose of a \psivdo\ is again a \psivdo. Both proofs will also be useful for 
generalizing the aforementionned results to the setting of holomorphic families of \psivdos. 

\begin{acknowledgements} 
I am grateful  to Alain Connes, Charles Epstein, Colin Guillermou, Bernard Helffer, Henri Moscovici, Michel Rumin and Elias Stein for helpful and 
 stimulating discussions. I would like also to thank for their hospitality the mathematics departments of Princeton 
 University, Harvard University and University of California at Berkeley where the memoir was finally completed. 
 
 On the other hand, some of the results of this memoir were announced in~\cite{Po:CRAS1} and presented as part of 
 the author's PhD thesis at University of Paris-Sud (Orsay, France) made under the supervision of Professor Alain Connes.  
 \end{acknowledgements}

%%%%%%%%%%%%%%%%%%%%%%%%%%%%%%%%%
%%%%%%%%%%%%%  Chap2.tex %%%%%%%%%%%%%%%
%%%%%%%%%%%%%%%%%%%%%%%%%%%%%%%%%%%

\chapter{Heisenberg manifolds and their main differential operators} 
\label{chap:Heisenberg-manifolds-operators}
In this chapter we recall the main definitions and properties of Heisenberg manifolds 
and we review the construction of the main examples of differential operators on such manifolds. 
% which will motivate and illustrate the main reson which we will focus throughout the memoir. 

\section{Heisenberg manifolds}
\label{sec:Heisenberg-manifolds}
In this section we gather the main facts about Heisenberg manifolds and their tangent Lie group bundles. 
% The exposition here follows closely that of~\cite{Po:Pacific1}.

\begin{definition}
   1) A Heisenberg manifold is a smooth manifold $M$ equipped with a distinguished hyperplane bundle $H \subset TM$. \smallskip 
   
   2) A Heisenberg diffeomorphism $\phi$ from a Heisenberg manifold $(M,H)$ onto another Heisenberg manifold 
   $(M,H')$ is a diffeomorphism $\phi:M\rightarrow M'$ such that $\phi^{*}H = H'$. 
\end{definition}

\begin{definition}
   Let $(M^{d+1},H)$ be a Heisenberg manifold. Then:\smallskip 
   
   1) A (local) $H$-frame for $TM$ is  a (local)  frame $X_{0}, X_{1}, \ldots, X_{d}$ of $TM$ so that $X_{1}, \ldots, 
   X_{d}$ span~$H$.\smallskip  
   
   2) A local Heisenberg chart is a  local chart with a local $H$-frame of $TM$ over its domain.
\end{definition}

The main examples of Heisenberg manifolds are the following.\smallskip 

\emph{a) Heisenberg group}. The $(2n+1)$-dimensional Heisenberg group
$\bH^{2n+1}$ is $\R^{2n+1}=\R \times \R^{2n}$ equipped with the 
group law, 
\begin{equation}
    x.y=(x_{0}+y_{0}+\sum_{1\leq j\leq n}(x_{n+j}y_{j}-x_{j}y_{n+j}),x_{1}+y_{1},\ldots,x_{2n}+y_{2n}).  
\end{equation}
A left-invariant basis for its Lie algebra $\fh^{2n+1}$ is 
provided by the vector-fields, 
\begin{equation}
    X_{0}=\frac{\partial}{\partial x_{0}}, \quad X_{j}=\frac{\partial}{\partial x_{j}}+x_{n+j}\frac{\partial}{\partial 
    x_{0}}, \quad X_{n+j}=\frac{\partial}{\partial x_{n+j}}-x_{j}\frac{\partial}{\partial 
    x_{0}}, %\quad j=1,\ldots,n,
     \label{eq:Examples.Heisenberg-left-invariant-basis}
\end{equation}
with $j=1,\ldots,n$. For $j,k=1,\ldots,n$ and $k\neq j$ we have the relations,
\begin{equation}
    [X_{j},X_{n+k}]=-2\delta_{jk}X_{0}, \qquad [X_{0},X_{j}]=[X_{j},X_{k}]=[X_{n+j},X_{n+k}]=0.
     \label{eq:Examples.Heisenberg-relations}
\end{equation}
In particular, the subbundle spanned by the vector fields 
$X_{1},\ldots,X_{2n}$ defines a left-invariant Heisenberg structure on 
$\bH^{2n+1}$.\smallskip

 \emph{(b) Codimension 1 foliations.} These are the Heisenberg manifolds $(M,H)$ such that $H$ is integrable in Fr\"obenius' sense, i.e.,
 $C^{\infty}(M,H)$ is closed under the Lie bracket of vector fields.\smallskip

 \emph{(c) Contact manifolds}. A contact manifold is a Heisenberg manifold $(M^{2n+1},H)$ such that near any point of $M$ there exists a contact form 
 anihilating $H$, i.e., a 1-form $\theta$ such that $d\theta_{|_{H}}$ is non-degenerate. When $M$ is orientable it is equivalent to require the 
 existence of a globally defined contact form on $M$  anihilating $H$. More specific examples of contact manifolds include the Heisenberg group 
 $\bH^{2n+1}$, boundaries of strictly pseudoconvex domains $D\subset \C^{2n+1}$, like the sphere $S^{2n+1}$, or even the cosphere bundle 
 $S^{*}M$ of a Riemannian manifold $M^{n+1}$.\smallskip

 \emph{d) Confoliations}. The confoliations of Elyashberg and Thurston in~\cite{ET:C} interpolate between contact manifolds and foliations. They can be seen 
 as oriented Heisenberg manifolds $(M^{2n+1},H)$ together with  a non-vanishing $1$-form $\theta$ on $M$ anihilating $H$ and such that 
$(d\theta)^{n}\wedge \theta\geq 0$.\smallskip

 \emph{e) CR manifolds.} If $D\subset \C^{n+1}$ a bounded domain with boundary $\partial D$ then the maximal complex structure, or CR structure, 
of $T(\partial D)$ is given 
by $T_{1,0}=T(\partial D)\cap T_{1,0}\C^{n+1}$, where $T_{1,0}$ denotes the holomorphic tangent bundle of $\C^{n+1}$. 
More generally, a CR structure on an orientable manifold $M^{2n+1}$ is given by a complex rank $n$ integrable subbundle $T_{1,0}\subset T_{\C}M$ such 
that $T_{1,0}\cap \overline{T_{1,0}}=\{0\}$. 
% Equivalently, the subbundle $H=\Re (T_{1,0}\oplus T_{0,1})$ has the structure of a complex vector bundle 
% of dimension $n$. 
Besides on boundaries of complex domains, and more generally 
such structures naturally appear on real hypersurfaces in $\C^{n+1}$, quotients of the Heisenberg group $\bH^{2n+1}$ by cocompact lattices, 
boundaries of complex hyperbolic spaces, and circle bundles over complex manifolds.  

A real hypersurface $M=\{r=0\}\subset \C^{n+1}$ is  strictly pseudoconvex when the Hessian $\partial \overline{\partial} r$ is positive definite. In general, 
to a CR manifold $M$ we can associate a Levi form $L_{\theta}(Z,W)=-id\theta(Z,\overline{W})$ on the CR tangent bundle $T_{1,0}$ by picking a 
non-vanishing real 1-form $\theta$ anihilating $T_{1,0}\oplus T_{0,1}$.  We then say that $M$ is strictly pseudoconvex 
(resp.~$\kappa$-strictly pseudoconvex) when we can choose $\theta$ so that $L_{\theta}$ is positive definite (resp.~is nondegenerate with $\kappa$ negative 
eigenvalues) at every point. In particular, when this happens  $\theta$ is non-degenerate on $H=\Re (T_{1,0}\oplus 
T_{0,1})$ and so $(M,H)$ is a contact manifold. 

\subsection{Tangent Lie group bundle of a Heisenberg manifold}
A simple description of the tangent Lie group bundle of a Heisenberg manifold $(M^{d+1},H)$ is given as follows.

\begin{lemma}[\cite{Po:Pacific1}]
The Lie bracket of vector fields induces on $H$ a 2-form with values in $TM/H$, 
\begin{equation}
    \cL: H\times H \longrightarrow TM/H,
     \label{eq:Bundle.Levi-form1}
\end{equation}
so that for any sections $X$ and $Y$ of $H$ near a point $a\in M$ we have
\begin{equation}
    \cL_{a}(X(a),Y(a)) = [X,Y](a) \quad \bmod H_{a}.
     \label{eq:Bundle.Levi-form2}
\end{equation}
\end{lemma}

\begin{definition}
 The $2$-form  $\cL$ is called the Levi form of $(M,H)$.
\end{definition}

The Levi form $\cL$ allows us to define a bundle $\fg M$ of graded Lie algebras  by endowing the vector bundle $(TM/H)\oplus H$ 
with the smooth fields of Lie brackets and gradings such that, for sections $X_{0}$, $Y_{0}$ of $TM/H$ and $X'$, $Y'$ 
of $H$ and for $t\in \R$, we have 
% for $a\in M$ and $X_{0}$, $Y_{0}$ in $T_{a}M/H_{a}$ and $X'$, $Y'$ in $H_{a}$, 
\begin{equation}
    [X_{0}+X',Y_{0}+Y']_{a}=\cL_{a}(X',Y'),  \qquad t.(X_{0}+X')=t^{2}X_{0}+tX'. %\quad t \in \R,
    \label{eq:Heisenberg.intrinsic-Lie-algebra-structure}
\end{equation} 

% \begin{definition}
%     The bundle $\fg M$ is the tangent Lie algebra bundle of $M$.
% \end{definition}

As we can easily check $\fg M$ is a bundle of $2$-step nilpotent Lie algebras which contains the normal bundle $TM/H$ in its center.
Therefore, its associated 
graded Lie group bundle $GM$ can be described as follows. As a bundle $GM$ is $(TM/H)\oplus H$ and the exponential 
map is merely the identity. In particular, the grading of $GM$ is as in~(\ref{eq:Heisenberg.intrinsic-Lie-algebra-structure}). 
Moreover, since  $\fg M$ is 
2-step nilpotent the Campbell-Hausdorff formula shows that, for  sections $X$, $Y$ of $\fg M$, we have
\begin{equation}
    (\exp X)(\exp Y)= \exp(X+Y+\frac{1}{2}[X,Y]).
\end{equation}
From this we deduce that the product on $GM$ is such that  
\begin{equation}
    (X_{0}+X').(Y_{0}+X')=X_{0}+Y_{0}+\frac{1}{2}\cL(X',Y')+X'+Y',    
    \label{eq:Bundle.Lie-group-law}
\end{equation}
for  sections $X_{0}$, $Y_{0}$ of $TM/H$  and sections $X'$, $Y'$ of $H$.

\begin{definition}
    The bundles  $\fg M$ and $GM$ are respectively called the tangent Lie group bundle and the tangent Lie group of $M$. 
\end{definition}

In fact, the fibers of $GM$ are classified by the Levi form $\cL$ as follows.

\begin{proposition}[\cite{Po:Pacific1}]\label{prop:Bundle.intrinsic.fiber-structure}
  1) Let $a\in M$. Then $\cL_{a}$ has rank $2n$ if, and only if, as a 
  graded Lie group $G_{a}M$ is isomorphic to $\bH^{2n+1}\times \R^{d-2n}$.\smallskip 
  
  2) The Levi form $\cL$ has constant rank $2n$ if, and only if, $GM$ is  a fiber bundle with typical fiber 
  $\bH^{2n+1}\times \R^{d-2n}$.
\end{proposition}
 
Now, let $\phi:(M,H)\rightarrow (M',H')$ be a Heisenberg diffeomorphism from $(M, H)$ onto another Heisenberg manifold 
$(M',H')$. Since $\phi_{*}H=H'$ we see that $\phi'$ induces a smooth vector bundle isomorphism 
$\overline{\phi}:TM/H\rightarrow TM'/H'$. 

\begin{definition}
We let  $\phi_{H}':(TM/H)\oplus 
  H \rightarrow (TM'/H')\oplus H'$ denote the vector bundle isomorphism such that
    \begin{equation}
    \phi'_{H}(a)(X_{0}+X')=\overline{\phi}'(a)X_{0}+\phi'(a)X',
     \label{eq:Bundle.Intrinsic.Phi'H}
\end{equation}
for any $a\in M$ and any $X_{0}\in T_{a}/H_{a}$ and $X'\in H_{a}$.
\end{definition}

\begin{proposition}[\cite{Po:Pacific1}]\label{prop:Bundle.Intrinsic.Isomorphism}
The vector bundle isomorphism  $\phi'_{H}$ is an isomorphism of graded Lie group bundles from $GM$ onto $GM'$. In particular, the Lie group bundle isomorphism 
class of $GM$ depends only on the Heisenberg diffeomorphism class of $(M,H)$.  
\end{proposition}

\subsection{Heisenberg coordinates and nilpotent approximation of vector fields}
It is interesting to relate the intrinsic description of $GM$ above with the more extrinsic description of~\cite{BG:CHM} (see also~\cite{Be:TSSRG}, 
\cite{EMM:HAITH}, \cite{EMM:RLSPD}, \cite{FS:EDdbarbCAHG},  \cite{Gr:CCSSW}, \cite{Ro:INA}) in terms of the Lie group 
associated to a nilpotent Lie algebra of model vector fields. 

First, let $a\in M$ and let us describe $\fg_{a}M$ as the graded Lie algebra of left-invariant vector fields on $G_{a}M$  
by identifying any $X \in \fg_{a}M$ with the left-invariant vector fields $L_{X}$ on $G_{a}M$ given by 
\begin{equation}
    L_{X}f(x)= \frac{d}{dt}f[x.(t\exp X)]_{|_{t=0}}= \frac{d}{dt}f[x.(tX)]_{|_{t=0}}, \qquad f \in C^{\infty}(G_{a}M).
\end{equation}
This allows us to associate to any vector fields $X$ near $a$ a unique left-invariant vector fields $X^{a}$ on $G_{a}M$ 
such that 
\begin{equation}
    X^{a}= \left\{ 
    \begin{array}{ll}
        L_{X_{0}(a)} & \text{if $X(a)\not \in H_{a}$},  \\
        L_{X(a)} & \text{otherwise,} 
    \end{array}\right.
     \label{eq:Bundle.intrinsic.model-vector-fields}
\end{equation}
where $X_{0}(a)$ denotes the class of $X(a)$ modulo $H_{a}$. 

\begin{definition}
    The left-invariant vector field $X^{a}$ is called the model vector field of $X$ at $a$.
\end{definition}

Let us look at the above construction in terms of a $H$-frame $X_{0},\ldots,X_{d}$ near 
$a$, i.e.,~of a local trivialization of the vector bundle $(TM/H)\oplus H$. For $j,k=1,\ldots,d$ we let 
\begin{equation}
    \cL(X_{j},X_{k})=[X_{j},X_{k}]X_{0}=L_{jk}X_{0} \quad \bmod H.
\end{equation}
With respect to the coordinate system $(x_{0},\ldots,x_{d})\rightarrow x_{0}X_{0}(a)+\ldots+x_{d}X_{d}(a)$ we can 
write the product law of $G_{a}M$ as 
\begin{equation}
    x.y=(x_{0}+\frac{1}{2}\sum_{j,k=1}^{d}L_{jk}x_{j}x_{k},x_{1},\ldots,x_{d}).
     \label{eq:Heisenberg.productGmM-coordinates}
\end{equation}
Then the vector fields $X_{j}^{a}$, $j=1,\ldots,d$, in~(\ref{eq:Bundle.intrinsic.model-vector-fields}) 
are just the left-invariant vector fields corresponding to the vector $e_{j}$ of the canonical basis 
of $\Rd$, that is, we have
\begin{equation}
    X_{0}^{a}=\frac{\partial}{\partial x_{0}} \quad \text{and}  \quad X_{j}^{a}=\frac{\partial}{\partial x_{j}} 
    -\frac{1}{2}\sum_{k=1}^{d}L_{jk}x_{k}\frac{\partial}{\partial x_{0}}, \quad 1\leq j\leq d.
     \label{eq:Heisenberg.Xjm.coordinates}
\end{equation}
In particular, for $j,k=1,\ldots,d$ we have the relations, 
\begin{equation}
    [X_{j}^{a},X_{k}^{a}]=L_{jk}(a)X_{0}^{a}, \qquad [X_{j}^{a},X_{0}^{a}]=0.
     \label{eq:Heisenberg.constant-structures.Gm}
\end{equation}

Now, let $\kappa:\dom \kappa \rightarrow U$ be a Heisenberg chart near $a=\kappa^{-1}(u)$ and let 
$X_{0},\ldots,X_{d}$ be the associated $H$-frame of $TU$.  
Then there is a unique affine coordinate change $x \rightarrow \psi_{u}(x)$ such that 
$\psi_{u}(u)=0$ and $\psi_{u*}X_{j}(0)=\frac{\partial}{\partial x_{j}}$ for 
$j=0,1,\ldots,d$. Indeed, if for $j=1,\ldots,d$ we set $X_{j}(x)=\sum_{k=0}^{d}B_{jk}(x)\frac{\partial}{\partial x_{k}}$ then 
we have
\begin{equation}
    \psi_{u}(x)=A(u)(x-u), \qquad A(u)=(B(u)^{t})^{-1}.
\end{equation}

\begin{definition}\label{def:Heisenberg.extrinsic.u-coordinates}
1) The coordinates provided by $\psi_{u}$ are called the privileged coordinates at $u$ 
with respect to the $H$-frame $X_{0},\ldots,X_{d}$. 

2) The map $\psi_{u}$ is called the privileged-coordinate map with respect to the $H$-frame $X_{0},\ldots,X_{d}$.
\end{definition}
\begin{remark}
    The privileged coordinates at $u$ are called $u$-coordinates in~\cite{BG:CHM}, but they correspond to the privileged coordinates 
    of~\cite{Be:TSSRG} and \cite{Gr:CCSSW} in the special case of a Heisenberg manifold. 
\end{remark}

Next, on $\Rd$ we consider the dilations, 
\begin{equation}
    \delta_{t}(x)=t.x=(t^{2}x_{0},tx_{1}, \ldots, tx_{d}), \qquad t \in \R,
    \label{eq:Heisenberg.dilations}
\end{equation}
with respect to which $\frac{\partial}{\partial_{x_{0}}}$ is 
homogeneous of degree $-2$ and $\frac{\partial}{\partial_{x_{1}}},\ldots,\frac{\partial}{\partial_{x_{d}}}$ 
is homogeneous of degree~$-1$. 

Since in the privileged coordinates at $u$ we have $X_{j}(0)=\frac{\partial}{\partial x_{j}}$, we can write
\begin{equation}
    X_{j}= \frac{\partial}{\partial_{x_{j}}}+ \sum_{k=0}^{d} a_{jk}(x) \frac{\partial}{\partial_{x_{k}}}, 
    \qquad j=0,1,\ldots d,
\end{equation}
where the $a_{jk}$'s are smooth functions such that $a_{jk}(0)=0$. Thus, we can let
\begin{gather}
    X_{0}^{(u)}= \lim_{t\rightarrow 0} t^{2}\delta_{t}^{*}X_{0}= \frac{\partial}{\partial_{x_{0}}},
    \label{eq:Heisenberg.X0u}\\
     X_{j}^{(u)}= \lim_{t\rightarrow 0} t^{-1}\delta_{t}^{*}X_{j}= 
     \frac{\partial}{\partial_{x_{j}}}+\sum_{k=1}^{d}b_{jk}x_{k} \frac{\partial}{\partial_{x_{0}}}, \quad 
     j=1,\ldots,d, \label{eq:Heisenberg.Xju}     
\end{gather}
where for $j,k=1,\ldots,d$ we have set $b_{jk}= \partial_{x_{k}}a_{j0}(0)$. 

Observe that $X_{0}^{(u)}$ is homogeneous of degree $-2$ and $X_{1}^{(u)},\ldots,X_{d}^{(u)}$ are homogeneous of degree $-1$.
 Moreover, for $j,k=1,\ldots,d$ we have 
\begin{equation}
    [X_{j}^{(u)},X_{0}^{(u)}]=0 \quad \text{and} \quad [X_{j}^{(u)},X_{0}^{(u)}]=(b_{kj}-b_{jk})X_{0}^{(u)}.
     \label{eq:Heisenberg.constant-structures.Gu1}
\end{equation}
Thus, the linear space spanned by $X_{0}^{(u)},X_{1}^{(u)}, \ldots, X_{d}^{(u)}$ is a graded 2-step nilpotent 
Lie algebra $\fg^{(u)}$. In particular, $\fg^{(u)}$ is the Lie algebra of left-invariant vector fields over the graded Lie group $G^{(u)}$ 
consisting of $\Rd$ equipped with the grading~(\ref{eq:Heisenberg.dilations}) and the group law,
\begin{equation}
    x.y=(x_{0}+\sum_{j,k=1}^{d}b_{kj}x_{j}x_{k},x_{1},\ldots,x_{d}).
\end{equation}

Now, if near $a$ we let $\cL(X_{j},X_{k})=[X_{j},X_{k}]=L_{jk}(x)X_{0}\bmod H$ then we have 
\begin{multline}
    [X_{j}^{(u)},X_{k}^{(u)}]=\lim_{t\rightarrow 0}[t\delta_{t}^{*}X_{j},t\delta_{t}^{*}X_{k}] = 
    \lim_{t\rightarrow 0} t^{2}\delta_{t}^{*}(L_{jk}(\circ \kappa^{-1}(x))X_{0})\\ =L_{jk}(a)X_{0}^{(u)}.
     \label{eq:Heisenberg.constant-structures.Gu2}
\end{multline}
Comparing this with~(\ref{eq:Heisenberg.constant-structures.Gm}) and~(\ref{eq:Heisenberg.constant-structures.Gu1}) 
then shows that $\fg^{(u)}$ has the same the constant structures as those of 
$\fg_{a}M$, hence is isomorphic to $\fg_{a}M$. Consequently, the Lie groups $G^{(u)}$ and $G_{a}M$ are isomorphic. 
In fact, as it follows from~\cite{BG:CHM} and~\cite{Po:Pacific1} an explicit isomorphism is given by 
\begin{equation}
     \phi_{u}(x_{0},\ldots,x_{d})= (x_{0}-\frac{1}{4}\sum_{j,k=1}^{d}(b_{jk}+b_{kj})x_{j}x_{k},x_{1},\ldots,x_{d}).
     \label{eq:Bundle.Extrinsic.Phiu}
\end{equation}

\begin{definition}\label{def:Bundle.extrinsic.normal-coordinates}
Let $\varepsilon_{u}=\phi_{u}\circ \psi_{u}$. Then:\smallskip 

1) The new coordinates provided by $\varepsilon_{u}$  are called Heisenberg 
coordinates at $u$ with respect to the $H$-frame $X_{0},\ldots,X_{d}$.\smallskip  

2) The map $\varepsilon_{u}$ is called the $u$-Heisenberg coordinate map.
\end{definition}

\begin{remark}
       The Heisenberg coordinates at $u$ have been also considered in~\cite{BG:CHM} as a technical tool 
       for inverting the principal symbol of a hypoelliptic sublaplacian.
\end{remark}

  Next, as it follows from~\cite[Lem.~1.17]{Po:Pacific1} we also have
\begin{gather}
    \phi_{*}X_{0}^{(u)}=\frac{\partial}{\partial x_{0}}=X_{0}^{a}, 
     \label{eq:Heisenberg.Xu-Xm1} \\
     \phi_{*}X_{j}^{(u)}=\frac{\partial}{\partial x_{j}}-\frac{1}{2}\sum_{k=1}^{d}L_{jk}x_{k}\frac{\partial}{\partial x_{0}}=X_{j}^{a}, \quad 
    j=1,\ldots,d. 
    \label{eq:Heisenberg.Xu-Xm2}
\end{gather}
Since $\phi_{u}$ commutes with the dilations~(\ref{eq:Heisenberg.dilations}) 
using~(\ref{eq:Heisenberg.X0u})--(\ref{eq:Heisenberg.Xju}) we get
\begin{equation}
    \lim_{t\rightarrow 0} t^{2}\delta_{t}^{*}\phi_{u*}X_{0}^{(u)}=X^{a}_{0} \quad \text{and} \quad 
    \lim_{t\rightarrow 0} t\delta_{t}^{*}\phi_{u*}X_{j}^{(u)}=X^{a}_{j}, \quad j=1,\ldots,d.
\end{equation}

In fact, as shown in~\cite{Po:Pacific1} 
for any vector fields $X$ near $a$, as $t\rightarrow 0$ and in Heisenberg coordinates at $a$, we have
\begin{equation}
   \delta_{t}^{*}X=  \left\{ 
   \begin{array}{ll}
       t^{-2}X^{a} +\op{O}(t^{-1})& \text{if $X(a)\in H_{a}$},\\
       t^{-1}X^{a} +\op{O}(1) & \text{otherwise}. 
   \end{array}\right. 
   \label{eq:Bundle.Extrinsic.approximation-normal}
\end{equation}
Therefore, we obtain:

\begin{proposition}[\cite{Po:Pacific1}]\label{prop:Bundle.equivalent-descriptions}
    In the Heisenberg coordinates centered at $a=\kappa^{-1}(u)$ the tangent Lie group $G_{a}M$ coincides with $G^{(u)}$ and for any vector fields 
    $X$ the model vector fields $X^{a}$ approximates $X$ near $a$ in the sense of~(\ref{eq:Bundle.Extrinsic.approximation-normal}).
\end{proposition}

One consequence of the equivalence between the two approaches to $GM$ is a tangent approximation for Heisenberg diffeomorphisms as follows. 

Let  $\phi:(M,H)\rightarrow (M',H')$ be a Heisenberg diffeomorphism
from $(M,H)$ to another Heisenberg manifold $(M',H')$. We also endow $\Rd$ with the pseudo-norm,
\begin{equation}
    \|x\|= (x_{0}^{2}+(x_{1}^{2}+\ldots+x_{d}^{2})^{2})^{1/4}, \qquad x\in \Rd,
\end{equation}
so that for any $x \in \Rd$ and any $t \in \R$ we have 
\begin{equation}
    \|t.x\|=|t|\, \|x\| . 
     \label{eq:Bundle.homogeneity-pseudonorm}
\end{equation}

\begin{proposition}[{\cite[Prop.~2.21]{Po:Pacific1}}]\label{prop:Heisenberg.diffeo}
   Let $a\in M$ and set $a'=\phi(a)$. Then, in Heisenberg coordinates at $a$ and at $a'$   
   the diffeomorphism $\phi(x)$ has a behavior near 
   $x=0$ of the form 
   \begin{equation}
       \phi(x)= \phi_{H}'(0)x+(\op{O}(\|x\|^{3}), \op{O}(\|x\|^{2}),\ldots,\op{O}(\|x\|^{2})).
        \label{eq:Bundle.Approximation-diffeo}
   \end{equation}
   In particular, there is no term of the form $x_{j}x_{k}$, $1\leq j,k\leq d$, in the Taylor expansion of $\phi_{0}(x)$ 
   at $x=0$.
\end{proposition}
   \begin{remark}
  An asymptotics similar to~(\ref{eq:Bundle.Approximation-diffeo}) is given  
  in~\cite[Prop.~5.20]{Be:TSSRG} in privileged coordinates at $u$ and 
   $u'=\kappa_{1}(a')$, 
  but the leading term there is only a Lie algebra isomorphism from $\fg^{(u)}$ onto $\fg^{(u')}$. This is only 
  in Heisenberg coordinates that we recover the Lie group isomorphism $\phi'_{H}(a)$ as the leading term of the asymptotics.
\end{remark}
   
\begin{remark}
    An interesting application  of Proposition~\ref{prop:Heisenberg.diffeo} in~\cite{Po:Pacific1} is the construction of the tangent groupoid 
    $\cG_{H}M$ of $(M,H)$ as the differentiable groupoid encoding the smooth deformation of  $M\times M$ to $GM$. 
   This groupoid is the analogue in the Heisenberg setting of Connes' tangent groupoid (see~\cite[II.5]{Co:NCG}, \cite{HS:MKOEFFTK}) and it shows 
   that $GM$ is tangent to $M$ in a differentiable  fashion (compare~\cite{Be:TSSRG}, \cite{Gr:CCSSW}). 
\end{remark}

\section{Main differential operators on Heisenberg manifolds} 
\label{sec:Operators}
In this section we recall the definitions of the most common operators on a Heisenberg manifold. With the exception of the contact Laplacian, 
all these operators are sublaplacians or are product of such operators up to lower order terms. 

A sublaplacian on a Heisenberg manifold $(M^{d+1},H)$ acting on the sections of a vector bundle $\cE$ over $M$ is a 
differential operator $\Delta:C^{\infty}(M,\cE)\rightarrow C^{\infty}(M,\cE)$ such that, near any $a\in M$, there exists a 
$H$-frame $X_{0},X_{1},\ldots,X_{d}$ of $TM$ so that $\Delta$ takes the form
 \begin{equation}
    \Delta=-\sum_{j=1}^{d} X_{j}^{2} + \sum_{j=1}^{d}a_{j}(x)X_{j}+c(x),
    \label{eq:Heisenberg.sublaplacian}
\end{equation}
for some local sections $a_{1}(x),\ldots,a_{d}(x)$  and $c(x)$ of $\End \cE$.
% where $\mu(x)$ is a smooth section of $\End \cE$ and $\op{O}_{H}(1)$ is a differential operator of Heisenberg order 1, i.e., it 
% 
% $\lambda (x),b_{1}(x),\ldots,b_{d}(x)$ and $c(x)$ are smooth  complex-valued functions.
% \end{definition}
% 
% We can also define sublaplacians acting on the sections of a 
% \begin{definition}\label{def:Operators.sublaplacian-bundle}
% A differential operator $\Delta:C^{\infty}(M,\cE)\rightarrow C^{\infty}(M,\cE)$ is a sublaplacian if, near any $y\in M$, there exists a Heisenberg chart 
%  with $H$-frame $X_{0},X_{1},\ldots,X_{d}$ and a local trivialization of $\cE$ such that in the local coordinates provided by this chart and this 
%  trivialization $\Delta$ takes the form
%  \begin{equation}
%     \Delta=-\sum_{j=1}^{d} X_{j}^{2} + i\lambda (x) X_{0}+ \sum_{j=1}^{d} b{j}(x)X_{j}+ c(x), \quad \lambda(x)=(\lambda_{1}(x),\ldots,\lambda_{r}(x)).
%     \label{eq:Heisenberg.sublaplacian.bundle}
% \end{equation}
% where $\lambda_{1}(x),\ldots,\lambda_{d}(x)$ are smooth complex-valued functions and $b_{1}(x),\ldots,b_{d}(x)$ and $c(x)$ are smooth functions with 
% values in $M_{r}(\C)$.
% \end{definition}

\subsection{H\"ormander's sum of squares.}
 Let $X_{1},\ldots,X_{m}$ be (real) vector fields on a manifold $M^{d+1}$ and consider the sum of squares,
\begin{equation}
    \Delta=-(X_{1}^{2}+\ldots+X_{m}^{2}).
\end{equation}
By a celebrated theorem of H\"ormander~\cite{Ho:HSODE} the operator $\Delta$ is hypoelliptic provided that the following bracket condition is 
satisfied: the vector fields $X_{0},\ldots,X_{m}$ together with their successive Lie brackets $[X_{j_{1}}, [X_{j_{2}},\ldots ,
X_{j_{1}}]\ldots]]$ span the tangent bundle $TM$ at every point. 

When $X_{1},\ldots,X_{m}$ span a hyperplane bundle $H$ the operator $\Delta$ is a sublaplacian with \emph{real} coefficients and the bracket condition 
reduces to $H+[H,H]=TM$ or, equivalently, to the nonvanishing of the Levi form of $(M,H)$. 

In fact, given a vector bundle $\cE$, the theorem of H\"ormander holds more generally for sublaplacians $\Delta:C^{\infty}(M,\cE)\rightarrow 
C^{\infty}(M,\cE)$ of the form
\begin{equation}
    \Delta=-(\nabla_{X_{1}}^{2}+\ldots+\nabla_{X_{m}}^{2})+L,
     \label{eq:Operators.generalized-sum-of-squares}
\end{equation}
where $\nabla$ is a connection on $\cE$ and $L$ is a first order differential operator with \emph{real} coefficients. In particular, if $M$ is endowed with a 
smooth positive density and $\cE$ with a Hermitian metric, this includes the selfadjoint sum of squares,
\begin{equation}
    \Delta=\nabla_{X_{1}}^{*}\nabla_{X_{1}}+\ldots+\nabla_{X_{m}}^{*}\nabla_{X_{m}}.
\end{equation}

\subsection{Kohn Laplacian}  
% In~ Kohn-Rossi showed that the Dolbeault complex on a bounded complex domain  induces on 
% its boundary a horizontal complex of differential forms. This was later extended by Kohn~ to the general setting of a CR 
% manifold $M^{2n+1}$ as follows.
Let $M^{2n+1}$ be an orientable CR manifold with CR tangent bundle $T_{1,0}\subset T_{\C}M$,  
let $\theta$  be a non-vanishing real 1-form annihilating the hyperplane bundle $H=\Re(T_{1,0}\otimes T_{0,1})$ and let $L_{\theta}$ be its 
associated Levi form. 
% Then the $\dbarb$-complex of 
% Kohn-Rossi~can be defined as follows. 

Let $\cN$ be a supplement of $H$ in $TM$. This is an orientable line bundle which gives rise to the splitting, 
\begin{equation}
    T_{\C}M=T_{1,0}\oplus T_{0,1}\oplus (\cN\otimes \C).
    \label{eq:CR-decomposition}
\end{equation}
For $p,q=0,\ldots,n$ let $\Lambda^{p,q}=(\Lambda^{1,0})^{p}\wedge (\Lambda^{0,1})^{q}$ be the bundle of $(p,q)$-forms, where 
$\Lambda^{1,0}$  and $\Lambda^{0,1}$) denote the annihilators 
in $T^{*}_{\C}M$ of $T_{0,1}\oplus (\cN\otimes \C)$ and $T_{1,0}\oplus (\cN\otimes \C)$ respectively. Then we have the splitting, 
\begin{equation}
    \Lambda^{*}T_{\C}^{*}M=(\bigoplus_{p,q=0}^{n}\Lambda^{p,q})\oplus (\theta\wedge  \Lambda^{*}T_{\C}^{*}M).
     \label{eq:CR-Lambda-pq-decomposition}
\end{equation}
Notice that this decomposition does not depend on the choice of $\theta$, but it does depend on that of $\cN$. 

The complex $\dbarb:C^{\infty}(M, \Lambda^{p,*})\rightarrow C^{\infty}(M,\Lambda^{p,*+1})$ of 
Kohn-Rossi~(\cite{KR:EHFBCM},~\cite{Ko:BCM}) is defined as follows. For any $\eta \in C^{\infty}(M, \Lambda^{p,q})$ we can uniquely decompose $d\eta$ as 
\begin{equation}
    d\eta =\dbarbpq \eta + \partial_{b;p,q}\eta + \theta \wedge \cL_{X_{0}}\eta,
%     \label{eq:¥}
\end{equation}
where $\dbarbpq \eta $  and $\partial_{b;p,q}\eta$ are sections of $\Lambda^{p,q+1}$ and $\Lambda^{p+1,q}$ respectively and $X_{0}$ is the section of $\cN$ 
such that $\theta(X_{0})=1$.  Thanks to the integrability of $T_{1,0}$ we have $\overline{\partial}_{b;p,q+1}\circ \dbarbpq=0$, so we really get a chain 
complex. This complex depends only on the CR structure of $M$ and on the choice of $\cN$, but the latter dependence is only up to the intertwinning 
by vector bundle isomorphisms (see, e.g.,~\cite[Lem.~4.1]{Po:AofM1}).

Next, assume that $T_{\C}M$ is endowed with a Hermitian metric compatible with the CR structure in the sense that it commutes with complex 
conjugation and the splitting~(\ref{eq:CR-decomposition}) becomes orthogonal. Let $\dbarbpq^{*}$ be 
 the formal adjoint of $\dbarbpq$. Then the Kohn Laplacian 
$\Boxbpq :C^{\infty}(M,\Lambda^{p,q})\rightarrow C^{\infty}(M,\Lambda^{p,q})$ is 
\begin{equation}
   \Boxbpq=\dbarbpq^{*}\dbarbpq +  \overline{\partial}_{b;p,q-1}\overline{\partial}_{b;p,q-1}^{*}.
\end{equation}

% set $T_{0,1}=\overline{T_{1,0}}$. Then the subbundle   
% $H=\Re (T_{1,0}\oplus T_{0,1})\subset TM$ admits an integrable complex structure and the splitting $H\otimes \C=T_{1,0}\oplus T_{0,1}$ 
% gives rise to a decomposition  $\Lambda H^{*}\otimes \C = \oplus_{0\leq p,q\leq n } \Lambda^{p,q}$ into forms of bidegree $(p,q)$.
% 
% Assume that $T_{\C}M$ is endowed with a Hermitian metric which is compatible with the CR structure in the sense that it  and complex conjugation is 
% an (antilinear) isometry. This Hermitian metric gives rise to a 
% Hermitian metric on $\Lambda^{*}T^{*}_{\C}M$ with respect to which the decomposition  $\Lambda H^{*}\otimes \C = \oplus_{0\leq p,q\leq n } 
% \Lambda^{p,q}$ becomes orthogonal. Let $\Pi_{p,q}:\Lambda^{*}T_{\C}^{*}M \rightarrow \Lambda^{p,q}$ be the orthogonal projection onto $\Lambda^{p,q}$. 
% Then the Kohn-Rossi operator $\bar \partial_{b}:  C^{\infty}(M,\Lambda^{p,q}) \rightarrow C^{\infty}(M,\Lambda^{p,q+1})$ is given by 
% \begin{equation}
%     \bar\partial_{b} \eta = \Pi_{p,q+1}(d\eta),\qquad  \eta \in C^{\infty}(M,\Lambda^{p,q}). 
% \end{equation}
% Since the integrability of $T_{1,0}$ implies that $\bar\partial_{b}^{2}=0$, this yields chain complexes 
% $\bar\partial_{b}:C^{\infty}(M,\Lambda^{p,*}) \rightarrow C^{\infty}(M,\Lambda^{p,*+1})$. 
% 
% 
% and . Then we 
% 
% 
% Endowing $M$ with a smooth density $\rho>0$ we let   Then the Kohn Laplacian is   

The Kohn Laplacian is a sublaplacian (see, e.g.,~\cite[Sect.~13]{FS:EDdbarbCAHG}, \cite[Sect.~20]{BG:CHM}), so is not elliptic. Nevertheless,  
Kohn~\cite{Ko:BCM} proved that 
under a geometric condition on the Levi form~$L_{\theta}$, the so-called condition $Y(q)$, 
the operator $\Boxbpq$ is hypoelliptic with gain of one derivative, 
i.e.,~for any compact $K\subset M$ we have 
\begin{equation}
    \|u \|_{s+1}  \leq C_{Ks}( \|\Boxbpq u \|_{s}  + \|u\|_{0}) \qquad \forall u \in C_{K}^{\infty}(M,\Lambda^{p,q}), 
     \label{eq:Operators.hypoellipticity}
\end{equation} 
where $\|.\|_{s}$ denotes the norm of the Sobolev space $L^{2}_{s}(M,\Lambda^{p,q})$. 

The condition $Y(q)$ at point $x\in M$ means that if we let $(r(x)-\kappa(x),\kappa(x),n-r(x))$ be the signature of $L_{\theta}$ at $x$, so that $r(x)$ is 
the rank of $L_{\theta}$ and $\kappa(x)$ the number of its negative eigenvalues, then we must have
\begin{equation}
    q\not \in \{\kappa(x),\ldots,\kappa(x)+n-r(x)\}\cup \{r(x)-\kappa(x),\ldots,n-\kappa(x)\}.
     \label{eq:Operators.Y(q)-condition}
\end{equation}
For instance, when $M$ is $\kappa$-strictly pseudoconvex, the $Y(q)$-condition exactly means that we must have $q\neq \kappa$ and $q\neq n-\kappa$. 

In general  this condition is equivalent to the existence of a parametrix within the Heisenberg calculus (see~\cite{BG:CHM} for the case of a 
smoothly diagonalizable Levi form and Section~\ref{sec:sublaplacian} for the general case;  
 see also~\cite{BdM:HODCRPDO}, \cite{FS:EDdbarbCAHG}), from which we recover the hypoellipticity of $\Boxbpq$.  

 Finally, the condition $Y(q)$ is only a sufficient condition for the hypoellipticity of the Kohn Laplacian, for  
the latter may be hypoelliptic even when the condition $Y(q)$ fails (see, e.g.,~\cite{Ko:SEPDCRM}, \cite{Ko:OMSHETCROBL},  \cite{Ni:GRDbWPCCRM}).

\subsection{Horizontal sublaplacian}
Let $(M^{d+1},H)$ be a Heisenberg manifold endowed with a Riemannian metric. Identifying $H^{*}$ with the subbundle of $T^{*}M$ annihilating 
the orthogonal supplement $H^{\perp}$, we define the horizontal sublaplacian as the differential operator, 
$\Delta_{b;k}:C^{\infty}(M,\Lambda^{k}_{\C}H^{*})\rightarrow C^{\infty}(M,\Lambda^{k+1}_{\C}H^{*})$ such that
\begin{equation}
    \Delta_{b;k}=d_{b;k}^{*}d_{b;k}+d_{b;k-1}d_{b;k-1}^{*}, \qquad d_{b;k}\alpha=\pi_{b;k+1}(d\alpha), 
    \label{eq:Operators.sublaplacian}
\end{equation}
where $\pi_{b;k+1}$ denotes the orthogonal projection onto $\Lambda^{k+1}_{\C}H^{*}$. 

This operator was first introduced by Tanaka~\cite{Ta:DGSSPCM} in the CR setting, but versions of this operator acting on functions 
were independently defined by Greenleaf~\cite{Gr:FESPM} and Lee~\cite{Le:FMPHI}. Moreover, it can be shown that $d_{b}^{2}=0$ if, and only if, 
the subbundle $H$ is integrable, so in general $\Delta_{b}$ is not the Laplacian of a chain complex. 

 On functions $\Delta_{b;0}$ is a sum of squares modulo a lower order term, hence is hypoelliptic by H\"ormander's theorem. 
 On horizontal forms of higher degree, that is, on sections of $\Lambda^{k}_{\C}H^{*}$ with $k\geq 1$, it is  shown  
in~\cite{Ta:DGSSPCM} and~\cite{Ru:FDVC}, in the contact case, and in Section~\ref{sec:Examples}, in the general case, that $\Delta_{b;k}$ is hypoelliptic 
when some condition, called condition $X(k)$, holds everywhere. More precisely, the condition $X(k)$ is satisfied at a point $x \in M$ when we have
\begin{equation}
    k\not \in\{\frac{1}{2}r(x), \frac{1}{2}r(x)+1,\ldots,d-\frac{1}{2}r(x)\},
%     \label{eq:Operators.X(k)-condition}
\end{equation}
where $r(x)$ denotes the rank of the Levi form $\cL$ at $x$. For instance, if $M^{2n+1}$ is a contact manifold or a nondegenerate CR manifold then the Levi form 
is everywhere nondegenerate, so $r(x)=2n$ and the $X(k)$-condition becomes $k\neq n$. 
% When the condition $\Delta_{b}$ admits a parametrix in the Heisenberg calculus, hence is hypoelliptic 
% (see~\cite{Po:BSM1}). 

Assume now that $M$ is an orientable CR manifold of dimension $2n+1$ with Heisenberg structure $H=\Re(T_{1,0}\oplus T_{0,1})$ and 
let $\theta$ be a global  nonvanishing  section of 
$TM/H$ with associated Levi form $L_{\theta}$. Assume in addition that $T_{\C}M$ is endowed with a Hermitian metric compatible with its CR structure. 

Under these assumptions we have  $d_{b}=\dbarb +\partial_{b}$, where $\partial_{b}$ denotes the conjugate of $\dbarb$, that is, the operator 
such that $\partial_{b}\omega =\overline{\dbarb \overline \omega}$ 
for any $\omega\in C^{\infty}(M,\Lambda_{\C}^{*}H^{*})$. Moreover, one can check that 
$\dbarb \partial_{b}^{*}+\partial_{b}^{*} \dbarb= \dbarb^{*} \partial_{b}+\partial_{b} 
 \dbarb^{*}=0$, from which we get
 \begin{equation}
     \Delta_{b}=\Box_{b}+ \overline{\Box}_{b},
       \label{eq:Operators.Tanaka-Kohn}
 \end{equation}
 where $\overline{\Box}_{b}$ is the conjugate of $\Box_{b}$. In particular, this shows that the horizontal sublaplacian preserves the bidegree, i.e., it 
 acts on $(p,q)$-forms. 
 
  Next, as shown in Section~\ref{sec:Examples}, the operator $\Delta_{b;p,q}$ acting on $(p,q)$-forms is hypoelliptic 
  when at any point $x \in M$ the condition~$X(p,q)$ is satisfied. The latter requires to have  
  \begin{equation}
      \{(p,q), (q,p)\}\cap \{(\kappa(x)+j,r(x)-\kappa(x)+k);\ \max(j,k)\leq n-r(x)\} =\emptyset, %\cup  \{(r(x)-\kappa(x)+j,\kappa(x)+k);\ \max(j,k)\leq n-r(x)\},
%        \label{eq:¥}
  \end{equation}
  where $r(x)$ denotes the rank of the Levi form $L_{\theta}$ at $x$ and $\kappa(x)$ its number of negative eigenvalues. For instance, 
  when $M$ is $\kappa$-strictly pseudoconvex the condition $X(p,q)$ means that we must have $(p,q)\neq (\kappa,n-\kappa)$ 
  and $(p,q)\neq (n-\kappa,\kappa)$.
 
%  The conjugate Kohn Laplacian acting on $(p,q)$-forms is hypoelliptic under the conjugate version of the $Y(q)$-condition, that is, the 
%  $Y(p)$-condition. We then can show that the horizontal Laplacian acting on $(p,q)$-forms is hypoelliptic when at least one of the conditions $Y(p)$ or $Y(q)$ is 
%  satisfied. Thus, when $M$ is $\kappa$-strictly pseudoconvex the horizontal Laplacian acting on $(p,q)$-forms is hypoelliptic whenever we have 
%  $(p,q)\neq (\kappa,n-\kappa)$ and $(p,q)\neq (n-\kappa,\kappa)$. 
 
\subsection{Gover-Graham operators}
Let $M^{2n+1}$ be a strictly pseudoconvex CR manifold. Let $T_{1,0}\subset T_{\C}M$ be the CR tangent bundle of $M$  and set 
$T_{0,1}=\overline{T_{1,0}}$ and $H=\Re(T_{1,0}\oplus T_{0,1})$. Let $\theta$ be a pseudohermitian contact form, i.e., $\theta$ is a 
1-form  annihilating $H$ such that the associated Levi form $L_{\theta}$ is positive definite on $T_{1,0}$. Thus $\theta$ defines a pseudohermitian 
structure on $M$ in the sense of Webster~\cite{We:PHSRH}. 
% so that $\theta$ is a contact 
% form with respect to which the complex structure $J$ of $H$ is calibrated, i.e.,~we have $d\theta(X,JX)>0$ for any non-zero section of $H$. Such a 
% contact form is called a pseudohermitian contac

We extend the Levi form $L_{\theta}$ into a Hermitian metric $h_{\theta}$ on $T_{\C}M$ such that $T_{1,0}$ and $T_{0,1}$ are orthogonal subspaces, 
complex conjugation is an (antilinear) isometry and $h_{\theta|_{H^{\perp}}}=\theta^{2}$. Then as shown by Tanaka~\cite{Ta:DGSSPCM} 
and Webster~\cite{We:PHSRH} there is a unique unitary connection on $T_{\C}M$ preserving the pseudohermitian structure. 
Note that the  contact form $\theta$ is unique up to a conformal change $f \rightarrow e^{2f} \theta$, $f\in C^{\infty}(M,\R)$. 

In order to study 
the analogue in CR geometry of the Yamabe problem  Jerison-Lee~\cite{JL:YPCRM} (see also~\cite{JL:ESIHGCRYP}, \cite{JL:ICRNCCRYP}) 
introduced a conformal version of the horizontal sublaplacian acting on functions as the operator $ \boxdot_{\theta} :C^{\infty}(M) \rightarrow 
C^{\infty}(M)$ such that
\begin{equation}
\boxdot_{\theta}=\Delta_{b;0}+ \frac{n}{n+2} R_{n},
\end{equation}
where $R_{n}$ denotes the scalar curvature of the Tanaka-Webster connection. This is a conformal operator in the sense that we have
\begin{equation}
     \boxdot_{e^{2f}\theta} = e^{-(n+2)f} 
     \boxdot_{\theta} e^{nf}\qquad \forall f\in C^{\infty}(M, \R). 
\end{equation}

The construction of Jerison-Lee 
has been  generalized by Gover-Graham~\cite{GG:CRIPSL}, who produced CR analogues of the conformal operators of~\cite{GMJS:CIPLIE}. 
For $k=1,\ldots,n+1$ and for $k=n+2,n+4,\ldots$ their constructions yield  a selfadjoint differential operator 
 \(
    \boxdot_{\theta}^{(k)}: C^{\infty}(M)\rightarrow C^{\infty}(M)
\)
such that 
\begin{equation}
     \boxdot_{e^{2f}\theta}^{(k)}= e^{-(n+1+k)f}\boxdot_{\theta}^{(k)} 
    e^{(n+1-k)f}, \qquad \forall f\in C^{\infty}(M, \R).
\end{equation}

We make the convention that for $k=1,2,\ldots,n+1$ the operator $ \boxdot_{\theta}^{(k)}$ corresponds to the operator $P_{w,w}$ of~\cite{GG:CRIPSL} 
with $w=\frac{k-1-n}{2}$ under the canonical trivializations of the 
density bundles $\cE(w,w)=|\Lambda^{n,n}|^{w}$ coming from the trivialization of $\Lambda^{n,n}$ provided by $d\theta^{n}$. For $k=n+2,n+4,\ldots$ the 
operator $ \boxdot_{\theta}^{(k)}$ similarly corresponds to the operator $\mathcal{P}_{w,w}$ of~\cite{GG:CRIPSL} 
with $w=\frac{k-1-n}{2}$. 

The operator $P_{w,w}$ is obtained by pushing down to $M$ the GJMS operator of order $k$ on the associated Fefferman bundle, while is 
$\mathcal{P}_{w,w}$ is contructed by making use of a CR geometric version of the tractor calculus of~\cite{BEG:TSCPRS}. In particular, for $k=1$ the 
operator $\boxdot_{\theta}^{(1)}$ agrees with the conformal sublaplacian of Jerison-Lee. In general, if we let $X_{0}$ denote the Reeb vector field of $\theta$, 
so that $\imath_{X_{0}} \theta=1$ and $\imath_{X_{0}} d\theta=0$, then $\boxdot_{\theta}^{(k)}$  has same principal part  (in the Heisenberg sense) as 
\begin{equation}
    (\Delta_{b;0}+i(k-1)X_{0})(\Delta_{b;0}+i(k-3)X_{0})\cdots (\Delta_{b;0}-i(k-1)X_{0}).
%     \label{eq:}
\end{equation}
 In particular, unless for the value $k=n+1$ the operator 
$\boxdot_{\theta}^{(k)}$ is hypoelliptic (see Proposition~\ref{sec:Examples}). 
% This is essentially the operator $P_{w,w'}$ of Theorem~1.1 of~\cite{GG:CRIPSL} with $w=w'=\frac{k-n-1}{2}$, noticing that when $w'=w$ the density bundles 
% $\cE(w,w)$, on which the Gover-Graham operators act, are in fact trivializable. 
% For $k=1$ we recover the operator $\boxdot_{\theta}$ of Jerison-Lee. 

Finally, let $Q_{\theta}$ denote the CR $Q$ curvature as defined by Fefferman-Hirachi~\cite{FH:AMCQCCCRG}. 
This the CR analogue of Branson's $Q$ curvature in 
conformal Riemannian geometry and as with the GJMS operators we have 
\begin{equation}
    e^{2(n+1)f}Q_{e^{2f}\theta}= Q_{\theta} + \boxdot_{\theta}^{n+1} \qquad \forall f \in C^{\infty}(M,\R).
%     \label{}
\end{equation}

\subsection{Contact complex and contact Laplacian}
Let  $(M^{2n+1},H)$ be an orientable contact manifold, let $\theta$ be a contact form on $H$ and let $J\in C^{\infty}(M,\End H)$, $J^{2}=-1$, 
be an almost complex structure on $H$ which is calibrated, i.e., $d\theta(X,JX)=-d\theta(JX,X)>0$ for any section $X$ of $H$. 
Then we can endow $M$ with the Riemannian metric, 
\begin{equation}
    g_{\theta,J}=d\theta(.,J.)+\theta^{2}.
\end{equation}

The contact complex of Rumin~\cite{Ru:FDVC} can be seen as an attempt to get on $M$ 
a complex of horizontal differential forms by forcing up the equalities $d_{b}^{2}=0$ and $(d_{b}^{*})^{2}=0$ as follows.

Let $X_{0}$ be the Reeb field associated to $\theta$, that is, so that $\imath_{X_{0}}\theta=1$ and $\imath_{X_{0}}d\theta=0$. 
Then we have 
\begin{equation}
    d_{b}^{2}=-\cL_{X_{0}}\varepsilon(d\theta)=-\varepsilon(d\theta)\cL_{X_{0}},
    \label{eq:DO.contact.square-db}
\end{equation}
where $\varepsilon(d\theta)$ denotes the exterior multiplication by $d\theta$. 

There are two natural ways of modifying the space $\Lambda^{*}_{\C}H^{*}$ of horizontal forms to get a complex. The first one is to force the equality 
$d_{b}^{2}=0$ by
restricting the operator $d_{b}$ to $\Lambda^{*}_{2}:=\ker \varepsilon(d\theta) \cap \Lambda^{*}_{\C}H^{*}$ since this 
bundle is closed under $d_{b}$ and annihilates $d_{b}^{2}$. 

The second way is to similarly force the equality $(d_{b}^{*})^{2}=0$ by restricting $d_{b}^{*}$ to  
$\Lambda^{*}_{1}:=\ker \iota(d\theta)\cap \Lambda^{*}_{\C}H^{*}=(\im \varepsilon(d\theta))^{\perp}\cap \Lambda^{*}_{\C}H^{*}$, 
where $\iota(d\theta)$ denotes the interior product with $d\theta$. This amounts to replace $d_{b}$ by the operator $d_{b}'=\pi_{1}\circ d_{b}$, 
where $\pi_{1}$ is the orthogonal projection onto $\Lambda^{*}_{1}$.

In fact, as $d\theta_{|_{H}}$ is nondegenerate the operator $\varepsilon(d\theta):\Lambda^{k}_{\C}H^{*}\rightarrow \Lambda^{k+2}_{\C}H^{*}$  is 
injective for $k\leq n-1$ and surjective for $k\geq n-1$. This implies that $\Lambda_{2}^{k}=0$ for $k\leq n-1$ and $\Lambda_{1}^{k}=0$ for $k\geq n+1$. 
Therefore, we only have two halves of complexes, but we can get a whole complex by 
connecting the two halves to each other as follows. 

Consider the differential operator $D_{R;n}:C^{\infty}(M,\Lambda_{\C}^{n}H^{*}) \rightarrow C^{\infty}(M,\Lambda_{\C}^{n}H^{*})$ such that
\begin{equation}
    D_{R;n}=\cL_{X_{0}}+d_{b;n-1}\varepsilon(d\theta)^{-1}d_{b;n},
\end{equation}
where $\varepsilon(d\theta)^{-1}$ is the inverse of $\varepsilon(d\theta):\Lambda^{n-1}_{\C}H^{*}\rightarrow \Lambda^{n+1}_{\C}H^{*}$. Notice that 
$D_{R,n}$ has order 2. Then we have: 

\begin{lemma}[\cite{Ru:FDVC}]
The operator $D_{R;n}$ maps to sections of $\Lambda_{2}^{n}$ and satisfies $d_{b;n}D_{R;n}=0$ and $D_{R;n}d_{b;n}'=0$.   
\end{lemma}
\begin{proof}
 Let $\alpha \in C^{\infty}(M,\Lambda^{n}_{\C}H^{*})$ and set $\beta=\varepsilon(d\theta)^{-1}d_{b;n}\alpha$, so that $d_{b;n}\alpha=d\theta \wedge \beta$.  
 Then we have 
 \begin{multline}
   \theta \wedge D_{R;n}\alpha = \theta\wedge \cL_{X_{0}}\alpha +\theta\wedge d_{b;n-1}\beta=d\alpha-d_{b;n}\alpha+\theta\wedge d\beta \\ 
     =d\alpha - d\theta \wedge \beta+\theta \wedge d\beta =d (\alpha-\theta\wedge \beta).  
 \end{multline}
 Hence $0=d(\theta\wedge D_{R;n}\alpha)=d\theta \wedge D_{R;n}\alpha -\theta\wedge d_{b;n}D_{R;n}\alpha$, which gives $d_{b;n}D_{R;n}\alpha=0$ and $d\theta 
 \wedge D_{R;n}\alpha=0$. In particular $D_{R;n}\alpha$ is a section of $\ker \varepsilon(d\theta)\cap \Lambda^{n}_{\C}H^{*}=\Lambda^{n}_{2}$. 
 
 Next, let $\alpha \in C^{\infty}(M,\Lambda^{n-1}_{\C}H^{*})$. Since $d_{b;n}'\alpha$ is the orthogonal projection of $d_{b;n}\alpha$ onto 
 $\ker \iota(d\theta)=\ker \varepsilon(d\theta)^{\perp}$ we can write $d_{b;n}'\alpha=d_{b;n}\alpha-d\theta\wedge \gamma$ for some  section $\gamma$ 
 of  
 $\Lambda^{n-2}_{\C}H^{*}$. Then using~(\ref{eq:DO.contact.square-db}) 
 we get $d_{b;n+1}d_{b;n}'\alpha=d_{b}^{2}\alpha-d_{b}(\theta \wedge \gamma)= 
 -d\theta \wedge \cL_{X_{0}}\alpha - d\theta \wedge d_{b;n-2}\gamma$. Thus, 
 \begin{equation}
     d_{b;n-1}\varepsilon(d\theta)^{-1}d_{b;n+1}d_{b;n}'\alpha=-d_{b;n}\cL_{X_{0}}\alpha-d_{b}^{2}\alpha=-d_{b;n}\cL_{X_{0}}\alpha+d\theta \wedge \cL_{X_{0}}\gamma.
      \label{eq:DO.contact.Dbdb'}
 \end{equation}
 
 On the other hand, we have 
 \begin{equation}
     \cL_{X_{0}}d_{b;n}'\alpha=\cL_{X_{0}}d_{b;n}\alpha-\cL_{X_{0}}(d\theta\wedge \gamma)= 
      \cL_{X_{0}}d_{b;n}\alpha-d\theta \wedge \cL_{X_{0}}\gamma. 
%      \label{}
 \end{equation}
Furthemore, we see that $\cL_{X_{0}}d_{b;n}\alpha$ is equal to 
\begin{equation}
    \imath_{X_{0}}d[d\alpha- \theta \wedge \cL_{X_{0}}\alpha]= 
     -\imath_{X_{0}}(d\theta \wedge \cL_{X_{0}}\alpha -\theta \wedge 
         d_{b;n}\cL_{X_{0}}\alpha)=d_{b;n}\cL_{X_{0}}\alpha.
%     \label{}
\end{equation}
 Combining this with~(\ref{eq:DO.contact.Dbdb'}) then gives 
 \begin{equation}
     D_{R;n}d_{b;n}'\alpha=\cL_{X_{0}}d_{b;n}'\alpha+ 
      d_{b;n-1}\varepsilon(d\theta)^{-1}d_{b;n+1}d_{b;n}'\alpha= 0. 
%      \label{}
 \end{equation}
 The lemma is thus proved. 
\end{proof}

All this shows that we have the complex, 
\begin{equation}
    C^{\infty}(M)\stackrel{d_{R;0}}{\rightarrow}
    \ldots %\stackrel{d_{R;n-1}}{\rightarrow} 
    C^{\infty}(M,\Lambda^{n})\stackrel{d_{R;n}}{\rightarrow} C^{\infty}(M,\Lambda^{n}) 
%     \stackrel{d_{R;n}}{\rightarrow}
    \ldots \stackrel{d_{R;2n-1}}{\rightarrow} C^{\infty}(M,\Lambda^{2n}),
  \label{eq:Operators.contact-complex}
\end{equation}
where $d_{R;k}$ agrees with $\pi_{1}\circ d_{b;k}$ for $k=0,\ldots,n-1$ and  with $d_{b;k}$ otherwise. This complex is called the contact 
complex. 

The above definition depends on the choice of the contact form $\theta$ and of the calibrated almost complex structure $J$. As shown by 
Rumin~\cite{Ru:FDVC} there is an alternative description of the contact complex depending only on $H=\ker \theta$. For $k=0,1\ldots,2n+1$ consider the vector 
bundles, 
\begin{gather}
    \Lambda^{k}_{H,1}=\Lambda_{\C}^{k}T^{*}M/ [(\im \varepsilon(\theta)+\im \varepsilon(d\theta))\cap \Lambda_{\C}^{k}T^{*}M],\\
    \Lambda^{k}_{H,2}=\ker \varepsilon(\theta) \cap \ker \varepsilon(d\theta)\cap \Lambda_{\C}^{k}T^{*}M.
\end{gather}
Notice that these bundles depend only on $H$ and that $ \Lambda^{k}_{H,1}=\{0\}$ for $k\geq n+1$ and $\Lambda^{k}_{H,2}=\{0\}$ for $k\leq n$. 
Furthermore, there are natural vector bundle isomorphisms 
$\phi_{k}:\Lambda_{1}^{k}\rightarrow   \Lambda^{k}_{H,1}$ and 
$\psi_{k}:\Lambda^{k}_{2}: \rightarrow \Lambda^{k+1}_{H,2}$. The former is obtained by restricting to $\Lambda^{k}_{1}$ the canonical projection from 
$\Lambda^{k}_{\C}T^{*}M$ onto $\Lambda^{k}_{H,1}$ and the latter is given by the restriction to $\Lambda_{2}^{k}$ of 
$\varepsilon(\theta)$, whose inverse is the restriction of $\imath_{X_{0}}$ to $\Lambda^{k+1}_{H,2}$. 

Next, the de Rham differential maps induces differential operators,
\begin{gather}
  d_{k}:C^{\infty}(M,\Lambda^{k}_{H,1})\rightarrow 
C^{\infty}(M,\Lambda^{k+1}_{H,1}),\quad k=0,1,\ldots,n-1,\\
d_{k}:C^{\infty}(M,\Lambda^{k}_{H,2})\rightarrow 
C^{\infty}(M,\Lambda^{k+1}_{H,2}), \quad k=n,\ldots,2n,
\end{gather}
which give rise to two halves of complex. As shown by Rumin these two halves can be connected by means of a 
 differential operator $D_{n}:
C^{\infty}(M,\Lambda^{n}_{H,1})\rightarrow C^{\infty}(M,\Lambda^{n+1}_{H,2})$ 
in such way that we get a full complex, 
\begin{equation}
    C^{\infty}(M)\stackrel{d_{0}}{\rightarrow} %C^{\infty}(M,\Lambda_{H,1}^{1}) \stackrel{d_{1}}{\rightarrow}
    \ldots C^{\infty}(M,\Lambda^{n}_{H,1})\stackrel{D_{n}}{\rightarrow}  C^{\infty}(M,\Lambda^{n+1}_{H,2}) 
%     \stackrel{d_{n+1}}{\rightarrow} 
    \ldots \stackrel{d_{2n}}{\rightarrow} C^{\infty}(M,\Lambda^{2n+1}_{H,2}).
     \label{eq:Operators.contact-complex2}
\end{equation}
In addition, the isomorphisms $\phi_{k}$ and $\psi_{k}$ above intertwine this complex with the complex~(\ref{eq:Operators.contact-complex}), 
that is, we have:\smallskip

- $\phi_{k+1}d_{R;k}=d_{k}\phi_{k}$ for $k=0,1,\ldots,n-1$;\smallskip

- $\psi_{n}d_{R;n}=D_{n}\phi_{n}$;\smallskip

- $\psi_{k+1}d_{R;k}=d_{k+1}\psi_{k}$ for $k=n,\ldots,2n$.\smallskip 

Let $\theta'$ be another contact form annihilating $H$ and let $J'$ be an almost complex structure on $H$ which is calibrated with respect to 
$\theta'$. Then, assigning the superscript $'$ to objects associated to $\theta'$ and $J'$, we see from the above discussion that there exist 
vector bundle isomorphisms $\phi_{k}:\Lambda_{1}^{k}\rightarrow \Lambda_{1}^{k'}$, $k=0,\ldots,n$, and $\psi_{k}:\Lambda_{2}^{k}\rightarrow 
\Lambda_{2}^{k'}$, $k=n,\ldots,2n$, intertwining the contact complexes associated to $(\theta,J)$ and $(\theta',J')$. Thus, up to the 
intertwining by vector bundle isomorphisms, the contact complex depend only on $H=\ker \theta$. 

The contact Laplacian is defined as follows. In degree $k\neq n$ this is the differential operator 
$\Delta_{R;k}:C^{\infty}(M,\Lambda^{k})\rightarrow C^{\infty}(M,\Lambda^{k})$ such that
\begin{equation}
    \Delta_{R;k}=\left\{
    \begin{array}{ll}
        (n-k)d_{R;k-1}d^{*}_{R;k}+(n-k+1) d^{*}_{R;k+1}d_{R;k},& k=0,\ldots,n-1,\\
         (k-n-1)d_{R;k-1}d^{*}_{R;k}+(k-n) d^{*}_{R;k+1}d_{R;k},& k=n+1,\ldots,2n.
         \label{eq:Operators.contact-Laplacian1}
    \end{array}\right.
\end{equation}
% \begin{gather}
%     \Delta_{R;k}=
% (n-k)d_{R;k-1}d^{*}_{R;k}+(n-k+1) d^{*}_{R;k+1}d_{R;k}, \quad k=0,\ldots,n-1,\\
%          (k-n-1)d_{R;k-1}d^{*}_{R;k}+(k-n) d^{*}_{R;k+1}d_{R;k},\quad k=n+1,\ldots,2n.
% %          \label{eq:Operators.contact-Laplacian1}
% \end{gather}
For $k=n$ we have the differential operators $\Delta_{R;nj}:C^{\infty}(M,\Lambda_{j}^{n})\rightarrow C^{\infty}(M,\Lambda^{n}_{j})$, $j=1,2$, given by the formulas, 
\begin{equation}
    \Delta_{R;n1}= (d_{R;n-1}d^{*}_{R;n})^{2}+d_{R;n}^{*}d_{R;n}, \quad   \Delta_{R;n2}=d_{R;n}d_{R;n}^{*}+  (d^{*}_{R;n+1}d_{R;n}).
         \label{eq:Operators.contact-Laplacian2}
\end{equation}

Observe that outside middle degree $k\neq n$ the operator $\Delta_{R;k}$ is a differential operator order $2$, 
whereas $\Delta_{Rn1}$ and $\Delta_{Rn2}$ are differential operators of 
order $4$. Moreover, Rumin~\cite{Ru:FDVC} proved that in every degree the contact Laplacian is maximal hypoelliptic in the sense of~\cite{HN:HMOPCV}. 

Alternatively, we can show that in every degree the contact Laplacian admits a parametrix in the Heisenberg calculus and recover its hypoellipticity from this fact 
(see~\cite{JK:OKTGSU} and Section~\ref{sec:Examples}).

\chapter{Intrinsic Approach to the Heisenberg Calculus} 
 \label{chap:Heisenberg-calculus}
 This chapter is organized as follows. In Section~\ref{sec:PsiHDO} we give  a detailed review 
of the main definitions and properties of the Heisenberg calculus, following mostly the point of view 
of~\cite{BG:CHM}. 

In Section~\ref{sec:principal-symbol} we define intrinsic notions of principal symbol and model operator for the Heisenberg calculus and 
check their main properties. 

In Section~\ref{sec:hypoellipticity} we establish an invertibility criterion for the principal symbol in terms of the so 
called Rockland condition. 

After this we explain in more details the invertibility criterion in some specific cases. 
First, we deal with general sublaplacians in Section~\ref{sec:sublaplacian},  
for which the results of~\cite{BG:CHM} even yield an explicit formula for the inverse of the principal symbol. Second, 
we deal more specifically with the main differential operators on Heisenberg manifolds in Section~\ref{sec:Examples}. 

\section{Heisenberg calculus}
\label{sec:PsiHDO}
The Heisenberg calculus is the relevant pseudodifferential tool to study hypoelliptic operators  on Heisenberg 
manifolds. It was independently invented by Beals-Greiner~\cite{BG:CHM} and Taylor~\cite{Ta:NCMA}, extending previous 
works of Boutet de Monvel~\cite{BdM:HODCRPDO}, Folland-Stein~\cite{FS:EDdbarbCAHG} and Dynin~(\cite{Dy:POHG}, \cite{Dy:APOHSC}) 
 (see also~\cite{BdMGH:POPDCM}, \cite{CGGP:POGD}, \cite{EMM:HAITH}, \cite{Gr:HPOPCD}, \cite{Ho:CHPODC}, \cite{RS:HDONG}). 

 The idea in the Heisenberg calculus, which is due to Elias Stein, is to have a pseudodifferential calculus on a Heisenberg manifold $(M,H)$ which is modeled at 
 any point $a\in M$ by the calculus of left-invariant pseudodifferential operators on the tangent group $G_{a}M$. 

 \subsection{Left-invariant pseudodifferential operators}
 Let $(M^{d+1},H)$ be a Heisenberg manifold and let $G=G_{a}M$ be the tangent group of $M$ at a point $a \in M$. We shall now recall the main facts about 
left-invariant pseudodifferential operators on $G$ (see also~\cite{BG:CHM}, \cite{CGGP:POGD}, \cite{Ta:NCMA}). 

 Recall that for any finite dimensional  real vector space $E$ the Schwartz class $\cS(E)$ is a Fr\'echet space and the Fourier transform is the continuous 
 isomorphism of $\cS(E)$ onto $\cS(E^{*})$ given by
 \begin{equation}
     \hat{f}(\xi)=\int_{E}e^{i\acou{\xi}{x}}f(x)dx, \qquad f \in \cS(E), \quad \xi \in E^{*},
 \end{equation}
 where $dx$ denotes the Lebesgue measure of $E$. 
 \begin{definition}
     $\cS_{0}(E)$ is the closed subspace of $\cS(E)$ consisting of $f \in \cS(E)$ such that for any differential operator $P$ on $E^{*}$ we have $(P\hat{f})(0)=0$. 
 \end{definition}
 
 Since $G$ has the same underlying set as that of its Lie algebra $\fg=\fg_{x}M$ we can let $\cS(G)$ and $\cS_{0}(G)$ denote the Fr\'echet spaces 
 $\cS(E)$ and $\cS_{0}(E)$ associated to the underlying linear space $E$ of $\fg$ 
 (notice that the Lebesgue measure of $E$ coincides with the Haar measure of $G$ since $G$ is nilpotent). 
 
 Next, for  $\lambda \in \R$ and $\xi=\xi_{0}+\xi'$ in $\fg^{*}= (T_{a}^{*}M/H^{*}_{a})\oplus H_{a}$ we let 
  \begin{equation}
     \lambda.\xi=  \lambda.(\xi_{0}+\xi')=\lambda^{2}\xi_{0}+\lambda \xi'.
      \label{eq:PsiHDO.Heisenberg-dilation-fg*}
 \end{equation}

 \begin{definition}
     $S_{m}(\fg^{*})$, $m \in \C$, is the space of functions $p \in C^{\infty}(\fg^{*}\setminus 0)$ which are homogeneous of degree $m$ in the sense 
     that, for any $\lambda>0$, we have 
\begin{equation}
         p(\lambda.\xi)=\lambda^{m}p(\xi) \qquad \forall \xi\in \fg^{*}\setminus 0. 
\end{equation}
    In addition $S_{m}(\fg^{*})$ is endowed with the Fr\'echet space topology coming from that of $C^{\infty}(\fg^{*}\setminus 0)$. 
 \end{definition}

 Note that the image $\hat{\cS}_{0}(G)$ of $\cS(G)$ under the Fourier transform consists of functions $v\in \cS(\fg^{*})$ such that, given any norm 
 $|.|$ on $G$,  near $\xi=0$ we have $|g(\xi)|=\op{O}(|\xi|^{N})$ for any integer $N\geq 0$. Thus, any $p \in S_{m}(\fg^{*})$ defines an element of
 $\hat{\cS}_{0}(\fg^{*})'$ by letting 
 \begin{equation}
     \acou{p}{g}= \int_{\fg^{*}}p(\xi)g(\xi)d\xi, \qquad g \in \hat{\cS}_{0}(\fg^{*}).
 \end{equation}
 This allows us to define the inverse Fourier transform of $p$ as the element $\check{p}\in \cS_{0}(G)'$ such that 
 \begin{equation}
     \acou{\check{p}}{f}=\acou{p}{\check{f}} \qquad \forall f \in \cS_{0}(G).
       \label{eq:PsiHDO.inverse-Fourier-transform-symbol}
 \end{equation}

\begin{proposition}[\cite{BG:CHM}, \cite{CGGP:POGD}]\label{prop:PsiHDO.convolution-symbols-group}
    1) For any $p \in S_{m}(\fg^{*})$ the left-convolution operator by $\check{p}$, i.e.,
\begin{equation}
   \check{p}*f(x):=\acou{\check{p}(y)}{f(x.y^{-1})}, \qquad f \in \cS_{0}(G), 
     \label{eq:PsiHDO.convolution-operator}
\end{equation}
gives rise to a continuous endomorphism of $\cS_{0}(G)$.\smallskip 

   2) There is a continuous bilinear product,
 \begin{equation}
     *:S_{m_{1}}(\fg^{*})\times S_{m_{2}}(\fg^{*}) \longrightarrow S_{m_{1}+m_{2}}(\fg^{*}),
 \end{equation} 
such that, for any  $p_{1}\in S_{m_{1}}(\fg^{*})$ and any $p_{2}\in S_{m_{2}}(\fg^{*})$, the composition of the left-convolution operators by 
$\check{p_{1}}$ and $\check{p_{2}}$ is the left-convolution operator by $(p_{1}*p_{2})^{\vee}$, that is, we have
\begin{equation}
    \check{p_{1}}*(\check{p_{2}}*f)=(p_{1}*p_{2})^{\vee}*f  \qquad \forall f \in \cS_{0}(G) .
\end{equation}
\end{proposition}

Let us also mention that if $p \in S_{m}(\fg^{*})$ then the convolution operator $Pu=\check{p}*f$ is a pseudodifferential operator. 
Indeed, let $X_{0}(a),\ldots,X_{d}(a)$ be a (linear) basis of $\fg$ so that $X_{0}(a)$ is in $T_{a}M/H_{a}$ and 
$X_{1}(a),\ldots,X_{d}(a)$ span $H_{a}$. For $j=0,\ldots,d$ let $X_{j}^{a}$ be the left-invariant vector fields on $G$ such that 
$X^{w}_{j|_{x=0}}=X_{j}(a)$. The basis 
$X_{0}(a),\ldots,X_{d}(a)$ yields a linear isomorphism $\fg\simeq \Rd$, hence a global chart of $G$. In this chart $p$ is a 
homogeneous symbol on $\Rd\setminus 0$ with respect to the dilations 
\begin{equation}
    \lambda.x=(\lambda^{2}x_{0},\lambda x_{1},\ldots,\lambda x_{d}), \qquad x\in \Rd, \quad \lambda>0.
     \label{eq:PsiHDO.Heisenberg-dilations-Rd}
\end{equation}

Similarly, each vector field $\frac{1}{i}X_{j}^{a}$, $j=0,\ldots,d$, corresponds to a vector field on $\Rd$ whose symbol is denoted 
$\sigma_{j}^{a}(x,\xi)$. Then, setting $\sigma=(\sigma_{0},\ldots,\sigma_{d})$, it can be shown that in the above chart the operator $P$ is 
given by
\begin{equation}
    Pf(x)=\int_{\Rd} e^{ix.\xi}p(\sigma^{a}(x,\xi))\hat{f}(\xi), \qquad f \in \cS_{0}(\Rd).
     \label{eq:PsiHDO.PsiDO-convolution}
\end{equation}
In other words $P$ is the pseudodifferential operator $p(-iX^{a}):=p(\sigma^{a}(x,D))$ acting on $\cS_{0}(\Rd)$.

\subsection{$\mathbf{\Psi_{H}}$DO's on an open subset of $\Rd$} 
Let $U$ be an open subset of $\Rd$ together with a hyperplane bundle $H \subset TU$ and a $H$-frame $X_{0},X_{1},\ldots,X_{d}$ of $TU$. 
Then the class of \psivdos\ on $U$  is a class of pseudodifferential operators modelled on that of homogeneous convolution operators on the fibers of $GU$.

\begin{definition} $S_{m}(\URd)$, $m\in\C$, is the space of functions 
    $p(x,\xi)$ in $C^{\infty}(U\times\Rdo)$ which are homogeneous of degree $m$ with respect to $\xi$, that is, 
    \begin{equation}
        p(x,\lambda.\xi)=\lambda^m p(x,\xi) \qquad \text{for any $\lambda>0$},
    \end{equation}
    where $\lambda.\xi$ is defined as in~(\ref{eq:PsiHDO.Heisenberg-dilations-Rd}).
\end{definition}

Observe that the homogeneity of $p\in S_{m}(\URd)$ implies that, for any compact $K \subset U$, we have the estimates
\begin{equation}
     | \partial^\alpha_{x}\partial^\beta_{\xi}p(x,\xi)|\leq C_{K\alpha\beta}\|\xi\|^{\Re m-\brak\beta}, \qquad x\in K, \quad \xi \neq 0,
     \label{eq:PsiVDO.estimates-homogeneous-symbols}
\end{equation}
where $ \|\xi\|=(|\xi_{0}|^{2}+|\xi_{1}|^{4}+\ldots +|\xi_{d}|^{4})^{1/4}$ and $\brak\alpha = 2\alpha_{0}+\alpha_{1}+\ldots+ 
\alpha_{d}$. 

\begin{definition}$S^m(\URd)$,  $m\in\C$, consists of symbols  $p\in C^{\infty}(\URd)$ with
an asymptotic expansion $ p \sim \sum_{j\geq 0} p_{m-j}$, $p_{k}\in S_{k}(\URd)$, in the sense that, for any integer $N$ and 
for any compact $K \subset U$, we have
\begin{equation}
    | \partial^\alpha_{x}\partial^\beta_{\xi}(p-\sum_{j<N}p_{m-j})(x,\xi)| \leq 
    C_{\alpha\beta NK}\|\xi\|^{\Re m-\brak\beta -N}, \quad  x\in K, \ \|\xi \| \geq 1.
    \label{eq:PsiVDO.asymptotic-expansion-symbols}
\end{equation}
\end{definition}

Next, for $j=0,\ldots,d$ let  $\sigma_{j}(x,\xi)$ denote the symbol of $\frac{1}{i}X_{j}$ (in the 
classical sense) and set  $\sigma=(\sigma_{0},\ldots,\sigma_{d})$. For any $p \in S^{m}(\URd)$ it can be shown that the symbol 
$p_{\sigma}(x,\xi):=p(x,\sigma(x,\xi))$ is in the H\"ormander class of symbols of type $(\frac{1}{2}, \frac{1}{2})$ (see~\cite[Prop.~10.22]{BG:CHM}). 
Therefore, we define a continuous linear operator from $C^{\infty}_{c}(U)$ to $C^{\infty}(U)$ by letting 
    \begin{equation}
          p(x,-iX)f(x)= (2\pi)^{-(d+1)} \int e^{ix.\xi} p(x,\sigma(x,\xi))\hat{f}(\xi)d\xi,
    \qquad f\in C^{\infty}_{c}(U).
    \end{equation}

    In the sequel we let $\Psi^{-\infty}(U)$ denotes the class of smoothing operators, i.e.,~ the class of operators $P:C^{\infty}_{c}(U)\rightarrow 
    C^{\infty}(U)$ given by a smooth kernel. 
    
\begin{definition}
   $\pvdo^{m}(U)$, $m\in \C$, consists of operators $P:C^{\infty}_{c}(U)\rightarrow C^{\infty}(U)$ of the form
\begin{equation}
         P= p(x,-iX)+R,
\end{equation}
 with $p$ in $S^{m}(\URd)$, called the symbol of $P$, and $R$ in $\Psi^{-\infty}(U)$.
\end{definition}

    The above definition of the symbol of $P$ differs from that of~\cite{BG:CHM}, since 
    there the authors defined it to be $p_{\sigma}(x,\xi)=p(x,\sigma(x,\xi))$. Note also that $p$ is unique 
    modulo $S^{-\infty}(\URd)$. 

\begin{lemma}\label{lem:PsiHDO.asymptotic-completeness}
   For $j=0,1,\ldots$ let $p_{m-j}\in S_{m-j}(\URd)$. Then there exists $P\in \pvdo^{m}(U)$ with symbol $p\sim \sum_{j\geq 0}p_{m-j}$. Moreover, the 
   operator $P$ is unique modulo smoothing operators. 
\end{lemma}
   
    The class $\pvdo^{m}(U)$ does not depend on the choice of the $H$-frame $X_{0}, \ldots, X_{d}$ (see~\cite[Prop.~10.46]{BG:CHM}). Moreover, 
since it is contained in the class 
of \psidos\ of type $(\frac{1}{2},\frac{1}{2})$ we obtain:
 
\begin{proposition}\label{prop:PsiHDO.Sobolev-regularity}
    Let $P$ be a \psivdo\ of order $m$ on $U$.\smallskip 
    
    1)  The operator $P$ extends to a continuous linear mapping from 
    $\cE'(U)$ to $\cD'(U)$ and has a distribution kernel which is smooth off the diagonal.\smallskip 
    
    3) Let $k=\Re m$ if $\Re m\geq 0$ and $k=\frac{1}{2}\Re m$ otherwise. Then  $P$ extends  
    to a continuous mapping from $L^{2}_{s,\op{comp}}(U)$ to $L^{2}_{s-k,\op{loc}}(U)$ for any $s$ in $\R$. 
\end{proposition}

\subsection{Composition of $\mathbf{\Psi_{H}}$DO's}
Recall that there is no symbolic calculus for \psidos\ of type $(\frac{1}{2},\frac{1}{2})$ since the product of two such \psidos\ needs not be 
again a \psido\ of type $(\frac{1}{2},\frac{1}{2})$. However, the fact that the \psivdos\ are modelled on left-invariant pseudodifferential operators 
allows us  to construct a symbolic calculus for \psivdos.

First, for $j=0,\ldots,d$ let $X_{j}^{(x)}$ be the leading homogeneous part of $X_{j}$ in privileged coordinates centered at $x$ 
defined according to~(\ref{eq:Heisenberg.X0u}) and~(\ref{eq:Heisenberg.Xju}). These vectors span a nilpotent Lie algebra of left-invariant vector fields
on a nilpotent graded Lie group $G^{x}$ which corresponds 
to $G_{x}U$ by pulling back the latter from the Heisenberg coordinates at $x$ to the privileged coordinates at $x$. 

As alluded to above the product law of $G^{(x)}$ defines a convolution product for symbols, 
\begin{equation}
    *^{(x)}: S_{m_{1}}(\Rd) \times S_{m_{2}}(\Rd) \longrightarrow S_{m_{1}+m_{2}}(\Rd).
     \label{eq:PsiHDO.convolution-symbol-pointwise}
\end{equation}
such that, with the notations of~(\ref{eq:PsiHDO.PsiDO-convolution}), on $\cL(\cS_{0}(\Rd))$ we have  
\begin{equation}
    p_{1}(-iX^{(x)})p_{2}(-iX^{(x)})=(p_{1}*^{(x)}p_{2})(-iX^{(x)}) \qquad \forall p_{j}\in S_{m_{j}}(\Rd).
\end{equation}

As it turns out the product $*^{(x)}$ depends smoothly on $x$ (see~\cite[Prop.~13.33]{BG:CHM}). Therefore, we get a continuous bilinear product,
\begin{gather}
    *: S_{m_{1}}(\URd) \times S_{m_{2}}(\URd) 
        \rightarrow S_{m_{1}+m_{2}} (\URd),\\   
        p_{1}*p_{2}(x,\xi)=(p_{1}(x,.)*^{(x)}p_{2}(x,.))(\xi), \qquad p_{j}\in S_{m_{j}}(\URd).
 \label{eq:PsiHDO.convolution-symbols-URd}
\end{gather}

\begin{proposition}[{\cite[Thm.~14.7]{BG:CHM}}] \label{prop:PsiHDO.composition}
    For $j=1,2$ let $P_{j}\in \pvdo^{m_{j}}(U)$ have  symbol $p_{j}\sim \sum_{k\geq 0} p_{j,m_{j}-k}$ and assume  
    that one of these operators is properly supported. Then the operator  
$P=P_{1}P_{2}$ is a \psivdo\ of order $m_{1}+m_{2}$ and has symbol  $p\sim \sum_{k\geq 0} p_{m_{1}+m_{2}-k}$, with  
\begin{equation}
     p_{m_{1}+m_{2}-k} = \sum_{k_{1}+k_{2}\leq k} \sum_{\alpha,\beta,\gamma,\delta}^{(k-k_{1}-k_{2})}
            h_{\alpha\beta\gamma\delta}  (D_{\xi}^\delta p_{1,m_{1}-k_{1}})* (\xi^\gamma 
            \partial_{x}^\alpha \partial_{\xi}^\beta p_{2,m_{2}-k_{2}}),    
\end{equation}
where $\underset{\alpha\beta\gamma\delta}{\overset{\scriptstyle{(l)}}{\sum}}$ denotes the sum over all the indices such that 
$|\alpha|+|\beta| \leq \brak\beta -\brak\gamma+\brak\delta = l$ and $|\beta|=|\gamma|$, and the functions 
$h_{\alpha\beta\gamma\delta}(x)$'s are  polynomials in the derivatives of the coefficients of 
the vector fields $X_{0},\ldots,X_{d}$.
\end{proposition}

% \begin{remark}
% It follows from~(\ref{eq:PsiHDO.convolution-symbols-URd}) that for any $x \in U$ the $x$-symbol $p_{1,m_{1}}*p_{2,m_{2}}(x,.)$ 
% depends only on $p_{m_{1}}(x,.)$ and $p_{m_{2}}(x,.)$. However, the value of  $(p_{1,m_{1}}*p_{2,m_{2}})(x,\xi)$ at $(x,\xi)\in \URd$ depends on all 
% the $\eta$-values $p_{m_{1}}(x,\eta)$ and $p_{m_{2}}(x,\eta)$ as $\eta$ ranges over $\Rd$. Thus we may localize the product of Heisenberg symbols with 
% respect to $x$, but not respect to $(x,\xi)$, that is, the product of \psivdos\ is local, but is not microlocal. 
% \end{remark}

\subsection{The distribution kernels of $\mathbf{\Psi_{H}}$DO's}
An important fact about \psidos\ is their characterization in terms of their distribution kernels.  

 First, we extend the notion of homogeneity of functions to distributions. For $K\in \cS'(\Rd)$ and for $\lambda >0$ we let $K_{\lambda}$ denote 
  the element of $\cS'(\Rd)$ such that
     \begin{equation}
           \acou{K_{\lambda}}{f}=\lambda^{-(d+2)} \acou{K(x)}{f(\lambda^{-1}.x)} \quad \forall f\in\cS(\Rd). 
            \label{eq:PsiHDO.homogeneity-K-m}
      \end{equation}
In the sequel we will also use the notation $K(\lambda.x)$ for denoting $K_{\lambda}(x)$. We then say that $K$ is homogeneous of degree $m$, $m\in\C$, 
when $K_{\lambda}=\lambda^m K$ for any $\lambda>0$. 

\begin{definition}\label{def:PsiHDO.regular-distributions}
  $\cS'_{\reg}(\Rd)$  consists of tempered distributions on $\Rd$ which are smooth outside the 
  origin. We equip it with the weakest  topology such that  the inclusions of $\cS'_{\reg}(\Rd)$  into $\cS'(\Rd)$ and $C^{\infty}(\Rdo)$ are continuous. 
\end{definition}

\begin{definition}
  $\cK_{m}(\URd)$, $m\in\C$, consists of distributions $K(x,y)$ in $C^\infty(U)\hotimes\cS'_{\reg}(\Rd)$ so that  for some functions $c_{\alpha}(x) \in 
   C^{\infty}(U)$,  $\brak\alpha=m$, we have  
    \begin{equation}
        K(x,\lambda.y)= \lambda^m K(x,y) + \lambda^m\log\lambda
                \sum_{\brak\alpha=m}c_{\alpha}(x)y^\alpha \qquad \text{for any $\lambda>0$}.
    \end{equation}
\end{definition}
 
The interest of considering the distribution class $\cK_{m}(\URd)$ stems from:

\begin{lemma}[{\cite[Prop.~15.24]{BG:CHM}},~{\cite[Lem.~I.4]{CM:LIFNCG}}]\label{prop:PsiHDO.Sm-Km}
   1) Any  $p\in S_{m}(\URd)$ agrees on $\URdo$ with  a distribution $\tau(x,\xi)\in C^{\infty}(U)\hotimes \cS'(\Rd)$ 
    such that $\check \tau_{\xiy}$ is in $\cK_{\hat m}(\URd)$, $\hat m=-(m+d+2)$. 
    
    2) If $K(x,y)$ belongs to $\cK_{\hat m}(\URd)$ then the restriction of 
    $\hat{K}_{\yxi }(x,\xi)$ to $\URdo$ belongs to $S_{m}(\URd)$. 
\end{lemma}

This result is a consequence of the solution to the problem of extending a homogeneous function $p\in C^{\infty}(\Rd\setminus 0)$ into a 
homogeneous distribution on $\Rd$ and of the fact that for $\tau \in \cS'(\Rd)$ we have 
\begin{equation}
    (\hat{\tau})_{\lambda} = |\lambda|^{-(d+2)}(\tau_{\lambda^{-1}})^{\wedge}\qquad \forall \lambda \in \R\setminus 0. 
     \label{eq:PsiHDO.dilations-Fourier-transform}
\end{equation}
In particular, if $\tau$ is homogeneous of degree $m$ then $\hat{\tau}$ is homogeneous of degree $-(m+d+2)$.

% The relevant class of kernels for the Heisenberg calculus is the following.
\begin{definition}\label{def:PsiHDO.K*}
$\cK^{m}(\URd)$, $m\in \C$, consists of distributions $K$ in $\cD'(\URd)$ with an asymptotic expansion 
     $K\sim \sum_{j\geq0}K_{m+j}$,  $K_{l}\in \cK_{l}(\URd)$, 
 in the sense that,  for any integer $N$, as soon as $J$ is large enough we have 
   \begin{equation}
K-\sum_{j\leq J}K_{m+j}\in  C^{N}(\URd). 
     \label{eq:PsiHDO.asymptotics-kernel}
 \end{equation} 
\end{definition}

Since under the Fourier transform the asymptotic expansion~(\ref{eq:PsiVDO.asymptotic-expansion-symbols}) 
for symbols corresponds to that for distributions in~(\ref{eq:PsiHDO.asymptotics-kernel}),  
using Lemma~\ref{prop:PsiHDO.Sm-Km} we get: 

\begin{lemma}[{\cite[pp.~133--134]{BG:CHM}}]\label{lem:PsiHDO.characterization.Km}
    Let $K \in \cD'(\URd)$. Then the following are equivalent:\smallskip
    
    (i) The distribution $K$ belongs to $\cK^{m}(\URd)$;\smallskip
   
   (ii) We can put $K$ into the form
  \begin{equation}
      K(x,y)=\check{p}_{\xiy}(x,y)+R(x,y),
  \end{equation}
 for some $p\in S^{\hat{m}}(\URd)$, $\hat{m}=-(m+d+2)$, and some $R\in 
    C^{\infty}(\URd)$.\smallskip 
    
    Moreover, if (i) and (ii) holds and $K \sim \sum_{j\geq 0} K_{m+j}$, $K_{l}\in \cK_{l}(\URd)$, 
    then we have $p\sim \sum_{j\geq 0}p_{\hat{m}-j}$ where $p_{\hat{m}-j} \in S_{\hat{m}-j}(\URd)$ is the restriction to $\URdo$ of 
    $(K_{m+j})^{\wedge}_{\yxi}$. 
\end{lemma}

Next, for $x \in U$ let $\psi_{x}$ denote the affine change to the privileged coordinates at $x$ and let us write  
$(A_{x}^{t})^{-1}\xi=\sigma(x,\xi)$ with  $A_{x}\in \op{GL}_{d+1}(\R)$. 
Since $\psi_{x}(x)=0$ and $\psi_{x*}X_{j}=\frac{\partial}{\partial 
y_{j}}$ at $y=0$ for $j=0,\ldots,d$,  one checks that $\psi_{x}(y)=A_{x}(y-x)$. 

Let $p\in S^{m}(\URd)$. As by definition we have $p(x,-iX)=p_{\sigma}(x,D_{x})$ with $p_{\sigma}(x,\xi)=p(x,\sigma(x,\xi))=p(x,(A_{x}^{t})^{-1}\xi)$,
the distribution kernel 
$k_{p(x,-iX)}(x,y)$ of $p(x,-iX)$ is represented by the oscillating integrals
\begin{equation}
%    (2\pi)^{-(d+1)} 
%    \frac{1}{(2\pi)^{d+1}}
   (2\pi)^{-(d+1)} \!\int e^{i(x-y).\xi} p(x,(A_{x}^{t})^{-1}\xi)d\xi=   %(2\pi)^{-(d+1)}
  \frac{|A_{x}|}{(2\pi)^{d+1}} \int e^{iA_{x}(x-y).\xi} p(x,\xi)d\xi .
\end{equation}
Since $\psi_{x}(y)=A_{x}(y-x)$ we deduce that 
\begin{equation}
    k_{p(x,-iX)}(x,y)=|\psi_{x}'| \check{p}_{\xiy}(x,-\psi_{x}(y)).
    \label{eq:PsiHDO.kernel-quantization-symbol-psiy}
\end{equation}
Combining this with Lemma~\ref{lem:PsiHDO.characterization.Km} then gives: 

\begin{proposition}[{\cite[Thms.~15.39, 15.49]{BG:CHM}}]\label{prop:PsiVDO.characterisation-kernel1}
 Let $P:C_{c}^\infty(U)\rightarrow C^\infty(U)$ be a continuous linear operator with distribution kernel $k_{P}(x,y)$. Then the following are 
 equivalent:\smallskip  
 
 (i) $P$ is a \psivdo\ of order $m$, $m\in \C$.\smallskip 
 
 (ii) We can write $k_{P}(x,y)$ in the form,
  \begin{equation}
     k_{P}(x,y)=|\psi_{x}'|K(x,-\psi_{x}(y)) +R(x,y) ,
      \label{eq:PsiHDO.characterization-kernel.privilegedx}
 \end{equation}
with  $K\in \cK^{\hat{m}}(\URd)$, $\hat{m}=-(m+d+2)$, and $R \in C^{\infty}(U\times U)$.\smallskip 

Furthermore, if (i) and (ii) hold and $K\sim \sum_{j\geq 0}K_{\hat{m}+j}$, 
$K_{l}\in \cK_{l}(\URd)$, then $P$ has symbol $p\sim \sum_{j\geq 0} p_{m-j}$, $p_{l}\in S_{l}(\URd)$, where 
$p_{m-j}$ is the restriction to $\URdo$ of $(K_{m+j})^{\wedge}_{\yxi}$. 
\end{proposition}

In the sequel we will need a version of Proposition~\ref{prop:PsiVDO.characterisation-kernel1} in Heisenberg coordinates.  
To this end let $\varepsilon_{x}$ denote the coordinate change to the Heisenberg coordinates at $x$ and set $\phi_{x}=\varepsilon_{x}\circ \psi_{x}^{-1}$. 
Recall that $\phi_{x}$ is  a Lie group isomorphism from $G^{(x)}$ to 
$G_{x}U$ such that $\phi_{x}(\lambda.y)=\lambda.\phi_{x}(y)$ for any $\lambda \in \R$. Moreover, using~(\ref{eq:Bundle.Extrinsic.Phiu}) 
one can check that $|\phi_{x}'|=1$ and $\phi_{x}^{-1}(y)=-\phi_{x}(-y)$. Therefore,  
from~(\ref{eq:PsiHDO.kernel-quantization-symbol-psiy}) we see that we can put $ k_{p(x,-iX)}(x,y)$ into the form
\begin{equation}
    k_{p(x,-iX)}(x,y)=|\varepsilon_{x}'| K_{P}(x,-\varepsilon_{x}(y)), 
    \label{eq:PsiHDO.kernel-quantization-symbol}
\end{equation}
where we have let 
\begin{equation}
    K_{P}(x,y)= \check{p}_{\xiy}(x,-\phi_{x}(-y))=\check{p}_{\xiy}(x,\phi_{x}^{-1}(y)). 
%     \label{eq:¥}
\end{equation}
In fact, the coordinate changes $\phi_{x}$, $x \in U$, give rise to an action on distributions on $\URd$ given by
   \begin{equation}
     K(x,y) \longrightarrow \phi_{x}^{*}K(x,y), \qquad \phi_{x}^{*}K(x,y)=K(x,\phi_{x}^{-1}(y)).
       \label{eq:PsiHD.isomorphism-cK}
   \end{equation}
  Since $\phi_{x}$ depends smoothly on $x$, this action induces a continuous linear isomorphisms of $C^{N}(\URd)$, $N \geq 0$, and $C^{\infty}(\URd)$ 
  onto themselves. As $\phi_{x}(y)$ is 
  polynomial in $y$ in such way that $\phi_{x}(0)=0$  and $\phi_{x}(\lambda.y)=\lambda.\phi_{x}(y)$ for every $\lambda \in \R$, we deduce that the above 
  action also yields a continuous linear isomorphism of $C^{\infty}(U)\hotimes \cS_{\reg}'(\Rd)$ onto itself and, for every $\lambda>0$,  we have 
   \begin{equation}
       (\phi_{x}^{*}K)(x,\lambda.y)= \phi_{x}^{*}[K(x,\lambda.y)], \qquad K \in \cD'(\URd).
        \label{eq:PsiHDO.homogeneity.phix*}
   \end{equation}

   Furthermore, as $\phi_{x}(y)$ is polynomial in $y$ we see that for every $\alpha \in \N^{d+1}$ we can write 
$\phi_{x}(y)^{\alpha}$ in the form $\phi_{x}(y)^{\alpha}=\sum_{\brak\beta=\brak \alpha}d_{\alpha\beta}(x)y^{\beta}$ with $d_{\alpha\beta}\in C^{\infty}(\URd)$.  
   It then follows that, for every $m \in \C$,  the map $K(x,y) \rightarrow \phi_{x}^{*}K(x,y)$ induces a linear isomorphisms of $\cK_{m}(\URd)$ and 
   $\cK^{m}(\URd)$ onto themselves. Combining this with~(\ref{eq:PsiHDO.kernel-quantization-symbol}) and 
   Proposition~\ref{prop:PsiVDO.characterisation-kernel1} then gives:

\begin{proposition}\label{prop:PsiVDO.characterisation-kernel2}
 Let $P:C_{c}^\infty(U)\rightarrow C^\infty(U)$ be a continuous linear operator with distribution kernel $k_{P}(x,y)$. Then the following are 
 equivalent:\smallskip  
 
 (i) $P$ is a \psivdo\ of order $m$, $m\in \C$.\smallskip 
 
 (ii) We can write $k_{P}(x,y)$ in the form,
  \begin{equation}
     k_{P}(x,y)=|\varepsilon_{x}'|K_{P}(x,-\varepsilon_{x}(y)) +R(x,y) ,
      \label{eq:PsiHDO.characterization-kernel.Heisenberg}
 \end{equation}
with $K_{P}\in \cK^{\hat{m}}(\URd)$, $\hat{m}=-(m+d+2)$, and $R \in C^{\infty}(U\times U)$. 

 Furthermore, if (i) and (ii) hold and $K_{P}\sim \sum_{j\geq 0}K_{P,\hat{m}+j}$,  
$K_{l}\in \cK_{l}(\URd)$, then $P$ has symbol $p\sim \sum_{j\geq 0} p_{m-j}$, $p_{l}\in S_{l}(\URd)$, where 
$p_{m-j}$ is the restriction to $\URdo$ of 
$[K_{P, \hat{m}+j}(x,\phi_{x}^{-1}(y))]^{\wedge}_{\yxi }$.
\end{proposition}

\begin{remark}\label{rem:PsiHDO.principal-symbol-at-x}
 Let $a\in U$. Then~(\ref{eq:PsiHDO.characterization-kernel.Heisenberg}) 
 shows that the distribution kernel of $\tilde{P}=(\varepsilon_{a})_{*}P$ at $x=0$ is 
 \begin{equation}
     k_{\tilde{P}}(0,y)=|\varepsilon_{a}'|^{-1} k_{P}(\varepsilon_{a}^{-1}(0),\varepsilon_{a}^{-1}(y))=K_{P}(a,-y).
      \label{eq:PsiHDO.characterization-kernel.Heisenberg-origin}
 \end{equation}
Moreover, as we are in Heisenberg coordinates already, we have $\psi_{0}=\varepsilon_{0}=\phi_{0}=\op{id}$. Thus, in the 
 form~(\ref{eq:PsiHDO.characterization-kernel.Heisenberg}) for $\tilde{P}$ we have $K_{\tilde{P}}(0,y)=K_{P}(a,y)$. 
Therefore, if we let $p_{m}(x,\xi)$ denote the principal symbol of $P$ and let $K_{P,\hat{m}}\in \cK_{\hat{m}}(\URd)$ denote the leading kernel of $K_{P}$,  
then by Proposition~\ref{prop:PsiVDO.characterisation-kernel2} we have 
 \begin{equation}
     p_{m}(0,\xi)=[K_{P,\hat{m}}]^{\wedge}_{\yxi }(a,\xi).
 \end{equation}
This shows that $[K_{P,\hat{m}}]^{\wedge}_{\yxi }(a,\xi)$ is the principal symbol of $P$ at $x=0$ in Heisenberg coordinates centered at $a$.
\end{remark}

\subsection{$\mathbf{\Psi_{H}}$DO's on a general Heisenberg manifold}
Let $(M^{d+1},H)$ be a Heisenberg manifold. As alluded to before the \psivdos\ on an subset of $\Rd$ are \psidos\ of type 
$(\frac{1}{2},\frac{1}{2})$. However, the latter don't make sense on a general manifold, for their class is not preserved by an arbitrary change of chart. Nevertheless, 
when dealing with \psivdos\ this issue is resolved if we restrict ourselves to changes of Heisenberg charts. Indeed, we have:

 \begin{proposition}\label{prop:PsiHDO.invariance}
     Let $U$ (resp.~$\tilde{U}$) be an open subset of $\Rd$ together with a hyperplane bundle $H\subset TU$ (resp.~$\tilde{H}\subset T\tilde{U}$) and a 
    $H$-frame of $TU$ 
    (resp.~a $\tilde{H}$-frame of $T\tilde{U}$). Let $\phi:(U,H)\rightarrow (\tilde{U},\tilde{H})$ be a Heisenberg diffeomorphism and let $\tilde{P}\in 
    \Psi_{\tilde{H}}^{m}(\tilde{U})$.\smallskip 
    
   1) The operator $P=\phi^{*}\tilde{P}$ is a \psivdo\ of order $m$ on $U$.\smallskip 
    
   2) If we write the distribution kernel of $\tilde{P}$ in the form~(\ref{eq:PsiHDO.characterization-kernel.Heisenberg}) 
   with $K_{\tilde{P}}(\tilde{x},\tilde{y})$ in $\cK^{\hat{m}}(\tilde{U}\times \Rd)$,  then the distribution kernel of 
   $P$ can be written in the form~(\ref{eq:PsiHDO.characterization-kernel.Heisenberg}) with $K_{P}(x,y)\in \cK^{\hat{m}}(\URd)$ such that 
   \begin{equation}
       K_{P}(x,y) \sim \sum_{\brak\beta\geq \frac{3}{2}\brak\alpha} \frac{1}{\alpha!\beta!} 
       a_{\alpha\beta}(x)y^{\beta}(\partial_{\tilde{y}}^{\alpha}K_{\tilde{P}})(\phi(x),\phi_{H}'(x)y),
        \label{eq:PsiHDO.asymptotic-expansion-KP}
   \end{equation}
   where we have let $a_{\alpha\beta}(x)=\partial^{\beta}_{y}[|\partial_{y}(\tilde{\varepsilon}_{\phi(x)}\circ \phi\circ \varepsilon_{x}^{-1})(y)| 
   (\tilde{\varepsilon}_{\phi(x)}\circ \phi\circ \varepsilon_{x}^{-1}(y)-\phi_{H}'(x)y)^{\alpha}]_{|_{y=0}}$ and $\tilde{\varepsilon}_{\tilde{x}}$ 
   denote the change to the Heisenberg coordinates at $\tilde{x}\in \tilde{U}$. In particular, 
   \begin{equation}
       K_{P}(x,y)=|\phi_{H}'(x)| K_{\tilde{P}}(\phi(x),\phi_{H}'(x)y) \qquad \bmod \cK^{\hat{m}+1}(\URd).
       \label{eq:PsiHDO.asymptotic-expansion-KP-principal}
   \end{equation}
\end{proposition}

\begin{remark}
    The version of the above statement in~\cite{BG:CHM} does not contain the asymptotics~(\ref{eq:PsiHDO.asymptotic-expansion-KP-principal}),  
    which will be crucial for giving a global definition of the 
    principal symbol of a \psivdo\ in the next section. For this reason a detailed proof of the above version is given in Appendix~\ref{chap.Appendix-Invariance}. 
    This proof will also be useful in Chapter~\ref{chap.HolPHDO} and~\cite{Po:CPDE1} for generalizing Proposition~\ref{prop:PsiHDO.invariance} 
    to holomorphic families of \psivdos\ and to \psivdos\ with parameter. 
\end{remark}

As a consequence of Proposition~\ref{prop:PsiHDO.invariance} we can define \psivdos\ on $M$  acting on the sections of 
a vector bundle $\cE$ over $M$.
\begin{definition}
    $\pvdo^{m}(M,\cE)$,  $m\in \C$,  consists of continuous operators $P$ from $C^{\infty}_{c}(M,\cE)$ to 
    $C^{\infty}(M,\cE)$ such that:\smallskip 
    
    (i) The distribution kernel of $P$ is smooth off the diagonal;\smallskip 
    
    (ii) For any trivialization $\tau:\cE_{|_{U}}\rightarrow U\times \C^{r}$ over a 
     local Heisenberg chart $\kappa:U \rightarrow V\subset \Rd$ the operator $\kappa_{*}\tau_{*}(P_{|_{U}})$ belongs to
     $\pvdo^{m}(V, \C^{r}):=\pvdo^{m}(V)\otimes \End \C^{r}$. 
\end{definition}

All the previous properties of \psidos\  on an open subset of $\Rd$ hold \emph{mutatis mutandis} for \psivdos\ on $M$ acting on sections of $\cE$. 

\subsection{Transposes and adjoints of $\mathbf{\Psi_{H}}$DO's}
Let us now look at the transpose and adjoints of \psivdos. First, given a Heisenberg chart $U \subset \Rd$ we have:

\begin{proposition}\label{prop:PsiHDO.transpose-chart}
    Let $P\in \pvdo^{m}(U)$. Then:\smallskip 
    
    1) The transpose operator $P^{t}$ is a \psivdo\ of order $m$ on $U$.\smallskip 
    
    2) If we write the distribution kernel of 
    $P$ in the form~(\ref{eq:PsiHDO.characterization-kernel.Heisenberg}) 
   with $K_{P}(x,y)$ in $\cK^{\hat{m}}(\URd)$ then $P^{t}$ can be written in the form~(\ref{eq:PsiHDO.characterization-kernel.Heisenberg}) with 
   $K_{P^{t}}\in  \cK^{\hat{m}}(\URd)$ such that 
   \begin{equation}
       K_{P^{t}}(x,y) \sim  \sum_{\frac{3}{2}\brak\alpha \leq \brak \beta} \sum_{|\gamma|\leq |\delta| \leq 2|\gamma|} 
       a_{\alpha\beta\gamma\delta}(x) y^{\beta+\delta}  
       (\partial^{\gamma}_{x}\partial_{y}^{\alpha}K_{P})(x,-y), 
        \label{eq:PsiHDO.transpose-expansion-kernel}
   \end{equation}
   where 
   $a_{\alpha\beta\gamma\delta}(x)=\frac{|\varepsilon_{x}^{-1}|}{\alpha!\beta!\gamma!\delta!}
   [\partial_{y}^{\beta}(|\varepsilon_{\varepsilon_{x}^{-1}(-y)}'|(y-\varepsilon_{\varepsilon_{x}^{-1}(y)}(x))^{\alpha})
   \partial_{y}^{\delta}(\varepsilon_{x}^{-1}(-y)-x)^{\gamma}](x,0)$. In particular, we have
   \begin{equation}
       K_{P^{t}}(x,y)=K_{P}(x,-y) \quad \bmod \cK^{\hat m+1}(\URd). 
        \label{eq:PsiHDO.transpose-principal-kernels}
   \end{equation}
\end{proposition}

 \begin{remark}
     The asymptotic expansion~(\ref{eq:PsiHDO.transpose-expansion-kernel}) is not stated in~\cite{BG:CHM}, 
     but we need it in order to determine the global principal symbol of the transpose of a 
     \psivdo\ (see next section). A detailed proof of Proposition~\ref{prop:PsiHDO.transpose-chart}
     can be found in Appendix~\ref{chap:Appendix-transpose}. %for sake of completeness and of further generalizations in~\cite{Po:CPDE1} and \cite{Po:CPDE1}. 
 \end{remark}

Using this result, or its version in~\cite{BG:CHM}, we obtain: 

\begin{proposition}[{\cite[Thm.~17.4]{BG:CHM}}]\label{prop:PsiHDO.transpose-adjoint}
  Let $P:C^{\infty}(M,\cE)\rightarrow C^{\infty}(M,\cE)$ be a \psivdo\ of order $m$. Then:
 
  1) The transpose operator $P^{t}:\cE'(M,\cE^{*})\rightarrow \cD'(M,\cE^{*})$ is a \psivdo\ of order $m$;\smallskip 
  
 2) If $M$ is endowed with a smooth positive density and $\cE$ with a Hermitian metric then the adjoint $P^{*}:C^{\infty}(M,\cE)\rightarrow 
 C^{\infty}(M,\cE)$ is a \psivdo\ of order $\overline{m}$. 
\end{proposition}

\section{Principal symbol and model operators.}
\label{sec:principal-symbol}
\subsection{Principal symbol and model operators.}%\label{sec:principal-symbol}
In this section we define the principal symbols and model operators of \psivdos\ and check their main properties.
% Let us now give an intrinsic definition of the principal symbol of  a \psivdo\ and of its model 
% operator at a point. 
% 
Let $\fg^{*}M$ be the dual bundle of $\fg M$ with canonical projection $\pi:\fg^{*}M \rightarrow M$.
%be the canonical projection of $\fg^{*}M$ onto $M$. 

\begin{definition}
For $m\in \C$ the space $S_{m}(\fg^{*}M,\cE)$ consists of sections $p(x,\xi)$ in $C^{\infty}(\fg^{*}M\setminus 0, \End \pi^{*}\cE)$ which are homogeneous of 
    degree $m$,  i.e.,  we have 
    \begin{equation}
        p(x,\lambda.\xi)=\lambda^{m}p(x,\xi) \qquad \forall \text{for any $\lambda>0$},
    \end{equation}
    where $\lambda.\xi$ is defined as in~(\ref{eq:PsiHDO.Heisenberg-dilation-fg*}). 
\end{definition}

Let $P\in \pvdo^{m}(M)$ and for $j=1,2$ let $\kappa_{j}$ be a Heisenberg chart with domain 
$V_{j}\subset M$ and let $\phi:U_{1}\rightarrow U_{2}$ be the corresponding transition map, where we have let $U_{j}=\kappa_{j}(V_{1}\cap V_{2})\subset \Rd$. 

Let us first assume that $\cE$ is the trivial line bundle, so that $P$ is a scalar operator. For $j=1,2$ we let 
$P_{j}:=\kappa_{j*}(P_{|_{V_{1}\cap V_{2}}})$, so that $P_{1}=\phi^{*} P_{2}$. Since $P_{j}$ belongs to $\pvdo^{m}(U_{j})$ its distribution kernel is of 
the form~(\ref{eq:PsiHDO.characterization-kernel.Heisenberg}) with $K_{P_{j}}\in \cK^{\hat m}(U_{j}\times \Rd)$. 
Moreover, by Proposition~\ref{prop:PsiHDO.invariance} we have 
\begin{equation}
    K_{P_{1}}(x,y)=|\phi_{H}'(x)| K_{P_{2}}(\phi(x),\phi_{H}'(x)y) \qquad \bmod \cK^{\hat{m}+1}(U_{1}\times\Rd).
\end{equation}
Therefore, if $K_{P_{j},\hat{m}}\in \cK_{\hat m}(U_{j}\times\Rd)$ is the leading kernel of $K_{P_{j}}$ then we get
\begin{equation}
    K_{P_{1},\hat m}(x,y)=|\phi_{H}'(x)|K_{P_{2},\hat m}(\phi(x),\phi_{H}'(x)y).
     \label{eq:PsiHDO.principal-kernel-transition}
\end{equation}

Next, for $j=1,2$ we define 
\begin{equation}
    p_{j,m}(x,\xi)=[K_{P_{j},\hat{m}}]^{\wedge}_{\yxi }(x,\xi), \qquad (x,\xi)\in U_{j}\times \Rdo. 
     \label{eq:PsiHDO.global-principal-symbol}
\end{equation}
By  Remark~\ref{rem:PsiHDO.principal-symbol-at-x} for any $a \in U_{j}$ the symbol $p_{j}(a,.)$ 
yields in Heisenberg coordinates centered at $a$ the principal symbol of $P_{j}$ at $x=0$.  Moreover, 
since $\phi_{H}'(a)$ is a linear map, from~(\ref{eq:PsiHDO.principal-kernel-transition}) we get
\begin{equation}
    p_{1,m}(x,\xi)= p_{2,m}(\phi(x),[\phi_{H}'(x)^{-1}]^{t}\xi).
\end{equation}
This shows that $p_{m}: =\kappa_{1}^{*}p_{1,m}$ is an element of $S_{m}(\fg^{*}(V_{1}\cap V_{1}))$ which is independent of the choice of the chart 
$\kappa_{1}$. Since $S_{m}(\fg^{*}M)$ is a sheaf this gives rise this uniquely defines a symbol $p_{m}(x,\xi)$ in $S_{m}(\fg^{*}M)$.

When $\cE$ is a general vector bundle, the above construction can be carried out similarly,  so that we obtain: 

\begin{theorem}\label{prop:PsiHDO.principal-symbol}
    For any $P \in \pvdo^{m}(M,\cE)$ there is a unique symbol $\sigma_{m}(P)(x,\xi)$ in $S_{m}(\fg^{*}M, \cE)$ such that, if in a local trivializing Heisenberg 
    chart $U\subset \Rd$ we let $K_{P,\hat{m}}(x,y)\in \cK_{\hat{m}}(\URd)$ be the leading kernel for the kernel $K_{P}(x,y)$  in the 
    form~(\ref{eq:PsiHDO.characterization-kernel.Heisenberg})  
    for $P$,  then we have
    \begin{equation}
       \sigma_{m}(P)(x,\xi)=[K_{P,\hat{m}}]^{\wedge}_{\yxi }(x,\xi), \qquad (x,\xi)\in U\times \Rdo. 
    \end{equation}
    
\noindent    Equivalently, for any $x_{0} \in M$ the symbol $\sigma_{m}(P)(x_{0},.)$ agrees in trivializing Heisenberg coordinates centered at $x_{0}$ with 
  the principal symbol of $P$ at $x=0$. 
\end{theorem}

\begin{definition}\label{def:PsiHDO.principal-symbol}
   For $P\in \pvdo^{m}(M,\cE)$ the symbol $\sigma_{m}(P)\in S_{m}(\fg^{*}M,\cE)$ provided by Theorem~\ref{prop:PsiHDO.principal-symbol}
     is called the principal symbol of $P$. 
\end{definition}

\begin{remark}
    Since we have two notions of principal symbol we shall distinguish between them by saying that  $\sigma_{m}(P)$ is the global 
    principal symbol of $P$ and that in a local trivializing chart the principal symbol $p_{m}$  of $P$ in the sense of~(\ref{eq:PsiVDO.asymptotic-expansion-symbols}) 
    is the local principal symbol of $P$ in this chart. 
   \end{remark}
   
    In a local Heisenberg chart $U\subset \Rd$ the global symbol $\sigma_{m}(P)$ and the local principal symbol $p_{m}$ of $P\in \pvdo^{m}(U)$  
   can be easily related to each other. Indeed, by Proposition~\ref{prop:PsiVDO.characterisation-kernel2} we have 
    \begin{equation}
       p_{m}(x,\xi)= [K_{P,\hat{m}}(x,\phi_{x}^{-1}(y))]^{\wedge}_{\yxi }(x,\xi),
    \end{equation}
    where $K_{P,\hat{m}}$ denotes the leading kernel for the kernel $K_{P}$ in the form~(\ref{eq:PsiHDO.characterization-kernel.Heisenberg}) for $P$. 
    By combining this with the definition~(\ref{eq:PsiHDO.global-principal-symbol}) 
    of $\sigma_{m}(P)$ we thus get 
    \begin{gather}
         p_{m}(x,\xi)= (\hat{\phi}_{x}^{*}\sigma_{m}(P))(x,\xi), 
           \label{eq:PsiHDO.local-global-principal-symbol} \\
       (\hat{\phi}_{x}^{*}\sigma_{m}(P))= 
       [[\sigma_{m}(P)]^{\vee}_{\xiy}(x,\phi_{x}^{-1}(y))]^{\wedge}_{\yxi }=  [\phi_{x}^{*}[\sigma_{m}(P)]^{\vee}_{\xiy}]^{\wedge}_{\yxi },
       \label{eq:PsiHDO.isomorphism-symbols}
     \end{gather}
    where $\phi_{x}^{*}$ is the isomorphism map~(\ref{eq:PsiHD.isomorphism-cK}).
    In particular, since the latter is a linear isomorphism of $\cK_{m}(\URd)$ onto itself, 
    we see that the map  $p\rightarrow \hat{\phi}_{x}^{*}p$ is a linear isomorphism of $S_{m}(\URd)$ onto itself. 

\begin{example}
    Let $X_{0},\ldots,X_{d}$ be a local $H$-frame of $TM$ near a point $a\in M$. In any Heisenberg chart associated with this frame the Heisenberg 
    symbol of $X_{j}$ is $\frac{1}{i}\xi_{j}$. In particular, this is true in Heisenberg coordinates centered at $a$. Thus the (global) principal symbol 
    of $X_{j}$ is equal to $\frac{1}{i}\xi_{j}$ in the local trivialization of $\fg^{*}M\setminus 0$ defined by the frame $X_{0},\ldots,X_{d}$. 
    More generally, for any differential $P=\sum_{\brak \alpha\leq m} a_{\alpha}(x) X^{\alpha}$ on $M$ we have
    \begin{equation}
        \sigma_{m}(P)(x,\xi)=\sum_{\brak \alpha\leq m} a_{\alpha}(x) i^{-|\alpha|}\xi^{\alpha}.
         \label{eq:Pincipal.example}
    \end{equation}
    Thus, for differential operators the global and local principal symbols agree in suitable coordinates. Alternatively, this result follows from the fact 
    that the isomorphism~(\ref{eq:PsiHD.isomorphism-cK}) induces the identity map on distributions supported in $U\times \{y=0\}$.
\end{example}    
    
\begin{proposition}\label{prop:PsiHDO.surjectivity-principal-symbol-map}
%     For every $m\in \C$ 
    The principal symbol $\sigma_{m}: \pvdo^{m}(M,\cE) \rightarrow S_{m}(\fg^{*}M,\cE)$ gives rise to a linear isomorphism
   $\pvdo^{m}(M,\cE)/\pvdo^{m-1}(M,\cE)\stackrel{\sim}{\longrightarrow}  S_{m}(\fg^{*}M,\cE)$. 
\end{proposition}
\begin{proof}
 Since the principal symbol of $P\in \pvdo^{m}(M,\cE)$ vanishes everywhere 
  if, and only if,  $P$ has order~$\leq m-1$, we see that the principal symbol map $\sigma_{m}$   induces an injective linear map  from 
  $\pvdo^{m}(M,\cE)/\pvdo^{m-1}(M,\cE)$ to $S_{m}(\fg^{*}M,\cE)$. 
  
  It remains to show that $\sigma_{m}$ is surjective. To this end consider a symbol 
  $p_{m}(x,\xi)\in S_{m}(\fg^{*}M,\cE)$ and let $(\varphi_{i})_{i\in I}$ be 
  a partition of the unity subordinated to an open covering $(U_{i})_{i\in I}$ of $M$ by
 domains of Heisenberg charts $\kappa_{i}:U_{i}\rightarrow V_{i}$ over which there are trivializations $\tau_{i}:\cE_{|_{U_{i}}}\rightarrow 
  U_{i}\times \C^{r}$. For each index $i$ let $\psi_{i}\in C^{\infty}(U_{i})$ be such that $\psi_{i}=1$ near 
  $\supp \varphi_{i}$ and set
\begin{equation}
     p_{m}^{(i)}(x,\xi)=(1-\chi(\xi))(\hat{\phi}_{i,x}^{*}\kappa_{i*}\tau_{i*}p_{m|_{\fg^{*}U_{i}\setminus 0}})(x,\xi) \in S_{m}(V_{i})\otimes \End 
  \C^{r},
\end{equation}
 where $\chi\in C^{\infty}(\Rd)$ is such that $\chi=1$ near the origin and $\hat{\phi}_{i,x}^{*}$ denotes the 
 isomorphism~(\ref{eq:PsiHDO.isomorphism-symbols})  with respect to the 
 chart $V_{i}$. Then we define a   a \psivdo\ of order $m$ by letting
  \begin{equation}
      P=\sum \varphi_{i}[\tau_{i}^{*}\kappa_{i}^{*} p_{m}^{(i)}(x,-iX)]\psi_{i}.
  \end{equation}
  
    For for every index $i$  the local principal symbol 
    of $\varphi_{i}[\tau_{i}^{*}\kappa_{i}^{*} p_{m}^{(i)}(x,-iX)]\psi_{i}$ in the chart $V_{i}$  
    is $\varphi_{i}\circ\kappa_{i}^{-1}(\hat{\phi}_{i,x}^{*}\kappa_{i*}\tau_{i*}p_{m|_{\fg^{*}U_{i}\setminus 0}})$, so 
    by~(\ref{eq:PsiHDO.local-global-principal-symbol})  its global principal is 
    $\varphi_{i}\circ\kappa_{i}^{-1}(\kappa_{i*}\tau_{i*}p_{m|_{\fg^{*}U_{i}\setminus 0}})$, which pulls back to $\varphi_{i}p_{m}$ 
    on $U_{i}$. It follows that the global principal symbol of $P$ is $\sigma_{m}(P)= \sum_{i} \varphi_{i} p_{m}=p_{m}$. 
    This proves the surjectivity of  the map $\sigma_{m}$, so the proof is now complete.
\end{proof}
  
  Next, granted the above definition of the principal symbol, we can define the model operator at a point as follows. 
  
  \begin{definition}\label{def:PsiHDO.model-operator}
    Let $P\in \pvdo^{m}(M,\cE)$ have (global) principal symbol $\sigma_{m}(P)$.  Then the model operator of $P$ at $a\in M$ is the left-invariant 
    \psivdo-operator $P^{a}$ from $\cS_{0}(G_{a}M,\cE_{a})$ to itself with symbol $\sigma_{m}(P)^{\vee}_{\xiy}(a,.)$, i.e., we have 
    \begin{equation}
        P^{a}f(x)=\acou{\sigma_{m}(P)^{\vee}_{\xiy}(a,y)}{f(x.y^{-1})}, \qquad f \in 
        \cS_{0}(G_{a}M,\cE_{a}).   
    \end{equation}
\end{definition}

Consider a local trivializing 
chart $U \subset \Rd$ near $a$ and let us relate the model operator $P^{a}$ on $G_{a}M$ to the operator $P^{(a)}=\tilde{p}_{m}^{a}(-iX^{(a)})$ on $G^{(a)}$ 
defined 
using the local principal symbol $\tilde{p}_{m}(x,\xi)$ of $P$ in this chart. Using~(\ref{eq:PsiHDO.convolution-operator}) 
and~(\ref{eq:PsiHDO.local-global-principal-symbol}) for $f\in \cS_{0}(\Rd)$ we get 
\begin{equation}
    P^{(a)}f(y)=\acou{(p_{m}^{a})^{\vee}(z)}{f(y.z^{-1})}=\acou{(\sigma_{m}(P)^{\vee}_{\xiy}(x,\phi_{a}^{-1}(y))}{f(y.z^{-1})}.
\end{equation}
Since $|\phi_{a}'|=1$ and $\phi_{a}$ is a Lie group isomorphism from $G^{(a)}$ onto $G_{a}M$ we obtain
\begin{equation}
      P^{(a)}f(y)= \acou{(\sigma_{m}(P)^{\vee}_{\xiy}(x,y)}{f\circ \phi_{a}^{-1}(y.\phi_{a}(x)^{-1})}=(\phi_{a}^{*}P^{a})f(y).
\end{equation}
In particular, we have
\begin{equation}
    P^{a}=(\phi_{a})_{*}p_{m}^{a}(-iX^{(a)}).
    \label{eq:PsiHDO.global-local-model-operator}
\end{equation}

\subsection{Composition of principal symbols and model operators}
Let us now look at the composition of principal symbols. To this end for $a \in M$ we 
let $*^{a}:S_{m_{1}}(\Rd)\times S_{m_{2}}(\Rd)\rightarrow S_{m_{1}+m_{2}}(\Rd)$ be the 
convolution product for symbols defined by the product law of $G_{a}M$ under the identification $G_{a}M\simeq \Rd$ provided by a $H$-frame 
$X_{0},\ldots,X_{d}$ of $TM$ near $a$, that is, 
\begin{equation}
    (p_{m_{1}}*^{a}p_{m_{j}})(-iX^{a})= p_{m_{1}}(-iX^{a})\circ p_{m_{2}}(-iX^{a}), \qquad p_{m_{j}}\in S_{m_{j}}(\Rd). 
     \label{eq:PsiHDO.global-convolution-symbols}
\end{equation}

Let $U \subset \Rd$ be a local trivializing Heisenberg chart
chart  near $a$ and for $j=1,2$ let $P_{j}\in \pvdo^{m_{j}}(U)$ have (global) principal symbol $\sigma_{m_{j}}(P_{j})$. 
Under the trivialization of $GU$ provided 
by the $H$-frame $X_{0},\ldots,X_{d}$ we have $P^{a}_{j}=\sigma(P_{j})(x, -iX^{a})$, so we obtain
\begin{equation}
    [\sigma_{m_{j}}(P_{j})(x,,)*^{a}\sigma_{m_{j}}(P_{j})(x,.)](-iX^{a})=P_{1}^{a}P_{2}^{a}.
     \label{eq:PsiHDO.convolution-symbol-product-model-operators}
\end{equation}

On the other hand, using~(\ref{eq:PsiHDO.local-global-principal-symbol}) and~(\ref{eq:PsiHDO.global-local-model-operator}) 
we see that $ \hat{\phi}_{a}^{*}[p_{m_{1}}*^{a}p_{m_{2}}](-iX^{a})$ is equal to
\begin{multline}
   \phi_{a}^{*}[p_{m_{1}}(-iX^{a})\circ p_{m_{2}}(-iX^{a})] = 
    \phi_{a}^{*}[p_{m_{1}}(-iX^{a})]\circ  \phi_{a}^{*}[p_{m_{2}}(-iX^{a})] \\ = 
    (\hat{\phi}_{a}^{*}p_{m_{1}})(-iX^{(a)})\circ (\hat{\phi}_{a}^{*}p_{m_{2}})(-iX^{(a)}) =  
     [(\hat{\phi}_{a}^{*}p_{m_{1}})*^{(a)}(\hat{\phi}_{a}^{*}p_{m_{2}})](-iX^{(a)}).
\end{multline}
Hence we have
\begin{equation}
    p_{m_{1}}*^{a}p_{m_{2}}=(\hat{\phi}_{a})_{*}[(\hat{\phi}_{a}^{*}p_{m_{1}})*^{(a)}(\hat{\phi}_{a}^{*}p_{m_{2}})] \qquad \forall p_{m_{j}}\in 
    S_{m_{j}}(\Rd),
    \label{eq:PsiHDO.global-local-convolution-symbols}
\end{equation}
where $(\hat{\phi}_{a})_{*}$ denotes the inverse of $\hat{\phi}_{a}^{*}$. Since $\hat{\phi}_{a}^{*}$, its inverse and $*^{(a)}$ depend smoothly on 
$a$, we deduce that that so does $*^{a}$.  Therefore, we get: 

\begin{proposition}\label{prop:Heisenberg.product-symbols-manifold}
  The group laws on the fibers of $GM$ give rise to a convolution product,
    \begin{gather}
        *:S_{m_{1}}(\fg^{*}M,\cE)\times S_{m_{2}}(\fg^{*}M,\cE) \longrightarrow S_{m_{1}+m_{2}}(\fg^{*}M,\cE),\\
%     \end{gather}
% such that for symbols $p_{m_{j}}\in S_{m_{j}}(\fg^{*}M,\cE)$, $j=1,2$, we have
%     \begin{gather}
        p_{m_{1}}*p_{m_{2}}(x,\xi)=[p_{m_{1}}(x,.)*^{x}p_{m_{2}}(x,.)](\xi), \qquad p_{m_{j}}\in S_{m_{j}}(\fg^{*}M,\cE),
    \end{gather}
 where $*^{x}$ denote the convolution product for symbols on $G_{x}M$.
\end{proposition}

Notice that~(\ref{eq:PsiHDO.global-local-convolution-symbols}) shows that, under the 
relation~(\ref{eq:PsiHDO.local-global-principal-symbol}) between local and global principal symbols, the convolution 
product~(\ref{eq:PsiHDO.global-convolution-symbols})  for 
global principal symbols corresponds to the convolution product~(\ref{eq:PsiHDO.convolution-symbols-URd}) 
for local principal symbols. Since by Proposition~\ref{prop:PsiHDO.composition} the latter yields 
the \emph{local} principal symbol of the product of two \psivdos\ in a local chart, we deduce that the convolution 
product~(\ref{eq:PsiHDO.global-convolution-symbols}) yields the 
\emph{global} principal  symbol of the product two \psivdos. 

Moreover,  by~(\ref{eq:PsiHDO.convolution-symbol-product-model-operators}) 
the global convolution product~(\ref{eq:PsiHDO.global-convolution-symbols}) corresponds to the product of model 
operators, so the model operator of a product of two \psivdos\ is equal to the product of the model operators. We have thus proved:

\begin{proposition}\label{prop:PsiHDO.composition2}
    For $j=1,2$ let $P_{j}\in \pvdo^{m_{j}}(M,\cE)$ and assume that $P_{1}$ or $P_{2}$ 
    is properly supported.\smallskip 
   
    1) We have $\sigma_{m_{1}+m_{2}}(P_{1}P_{2})=\sigma_{m_{1}}(P)*\sigma_{m_{2}}(P)$.\smallskip
    
    2) At every $a\in M$ the model operator of $P_{1}P_{2}$ is $(P_{1}P_{2})^{a}=P^{a}_{1}P_{2}^{a}$.
\end{proposition}

Finally, we look at the continuity of the above product for homogeneous symbols. To this end for each $m\in \C$ we endow $S_{m}(\fg^{*}M,\cE)$ with 
the Fr\'echet space topology inherited from that of $C^{\infty}(\fg^{*}M\setminus 0,\End \cE)$. 

\begin{proposition}\label{prop:PsiHDO.continuity-product-symbol}
    The product $*$ for homogeneous symbols 
   gives rise to a continuous bilinear map from $S_{m_{1}}(\fg^{*}M,\cE)\times S_{m_{2}}(\fg^{*}M,\cE)$ 
    to $S_{m_{1}+m_{2}}(\fg^{*}M,\cE)$. 
\end{proposition}
\begin{proof}
 Consider a sequence $(p_{k},q_{k})_{k\geq 0}$ in $S_{m_{1}}(\fg^{*}M,\cE)\times S_{m_{2}}(\fg^{*}M,\cE)$ converging to $(p,q)$ and such that 
 $p_{k}*q_{k}$ converges to $r$ in $S_{m_{1}+m_{2}}(\fg^{*}M,\cE)$. Let $a\in M$. Then $(p_{k}(a,.),q_{k}(a,.))$ converges to $(p(a,.),q(a,.))$ in 
 $S_{m_{1}}(\fg^{*}_{a}M,\cE_{a})\times S_{m_{2}}(\fg^{*}_{a}M,\cE_{a})$ and $(p_{k}*q_{k})(a,.)$ converges to $r(a,.)$ in 
 $S_{m_{1}+m_{2}}(\fg^{*}_{a}M,\cE_{a})$. 
 
On the other hand, by Proposition~\ref{prop:PsiHDO.convolution-symbols-group} 
the product $*^{a}$ gives rise to a continuous bilinear map from  
$S_{m_{1}}(\fg^{*}_{a}M,\cE_{a})\times S_{m_{2}}(\fg^{*}_{a}M,\cE_{a})$ to $S_{m_{1}+m_{2}}(\fg^{*}_{a}M,\cE_{a})$, so 
$(p_{k}*q_{k})(a,.)=p_{k}(a,.)*^{a}q_{k}(a,.)$ also converges to $p(a,.)*^{a}q(a,.)=p*q(a,.)$. Hence 
$p*q(a,.)=r(a,.)$ for any $a\in M$, that is, the symbols $p*q$ and $r$ agree. It then follows from the closed graph theorem that $*$ gives rise 
a continuous bilinear map from $S_{m_{1}}(\fg^{*}M,\cE)\times S_{m_{2}}(\fg^{*}M,\cE)$ 
to $S_{m_{1}+m_{2}}(\fg^{*}M,\cE)$.
\end{proof}

\subsection{Principal symbol of transposes and adjoints}
In this subsection we shall determine the principal symbols and the model operators of transposes and adjoints of \psivdos. 

Recall that by Proposition~\ref{prop:PsiHDO.transpose-adjoint} if $P\in \pvdo^{m}(M,\cE)$ then  
its transpose operator $P^{t}:C^{\infty}_{c}(M,\cE^{*})\rightarrow C^{\infty}(M,\cE)$ is a \psivdo\ of order $m$ and its adjoint 
$P^{*}:C^{\infty}_{c}(M,\cE)\rightarrow C^{\infty}(M,\cE)$ is a \psivdo\ of order $\overline{m}$ (assuming $M$ endowed with a positive density and 
$\cE$ with a Hermitian metric in order to define the adjoint). 

\begin{proposition}\label{prop:PsiHDO.transpose-global}
  Let $P \in \pvdo^{m}(M,\cE)$ have principal symbol $\sigma_{m}(P)$. Then:\smallskip 
  
 1) The principal symbol of $P^{t}$ is $\sigma_{m}(P^{t})(x,\xi)= \sigma_{m}(x,-\xi)^{t}\in S_{m}(\fg^{*}M,\cE^{*})$;\smallskip 
  
  2) If $P^{a}$ is the model operator of $P$ at $a$ then the model operator of $P^{t}$ at $a$ is the transpose operator 
  $(P^{a})^{t}: \cS_{0}(G_{x}M,\cE_{x}^{*})\rightarrow \cS_{0}(G_{x}M,\cE_{x}^{*})$.
\end{proposition}
\begin{proof}
   Let us first assume that $\cE$ is the trivial line bundle and that $P$ is a scalar operator. In a local Heisenberg chart $U \subset \Rd$ 
   we can write the distribution kernels of $P$ and $P^{t}$ in the form~(\ref{eq:PsiHDO.characterization-kernel.Heisenberg}) with distributions 
   $K_{P}$ and $K_{P^{t}}$ in  $\cK^{\hat{m}}(\URd)$. Let $K_{P,\hat{m}}$  and $K^{t}_{P^{t},\hat{m}}$ denote the principal parts of 
   $K_{P}$ and $K_{P^{t}}$ respectively. Then the principal symbols of $P$ and $P^{t}$ are
  $\sigma_{m}(P)(x,\xi)=(K_{P,\hat{m}})^{\wedge}_{\yxi}(x,\xi)$ and $\sigma_{m}(P^{t})(x,\xi)=(K_{P^{t},\hat{m}})^{\wedge}_{\yxi}(x,\xi)$ respectively.
   Since~(\ref{eq:PsiHDO.transpose-principal-kernels}) implies that $K_{P^{t},\hat{m}}(x,y)=K_{P,\hat{m}}(x,-y)$ and the Fourier transform commutes 
   with the multiplication by $-1$ we get
   \begin{equation}
       \sigma_{m}(P^{t})(x,\xi)=\sigma_{m}(P)(x,-\xi). 
   \end{equation}
   
  Next, for $a \in U$ let $p\in S_{m}(G_{a})$ and let $P$ be the left-invariant \psivdo\ with symbol $p$. Then the transpose $P^{t}$ is 
  such that, for $f$ and $g$ in $S_{0}(G_{a}U)$, we have 
   \begin{multline}
       \acou{P^{t}f}{g}=\acou{f}{Pv}= \acou{1}{f(x)(Pg)(x)}=\acou{1}{f(x)\acou{\check{p}(y)}{g(x.y^{-1})}}\\
       = \acou{1\otimes \check{p}(x,y)}{f(x)g(x.y^{-1})}.
   \end{multline}
   Therefore, using the change of variable $(x,y)\rightarrow (x.y^{-1},y^{-1})$ and the fact that $y^{-1}=-y$ we get
   \begin{equation}
        \acou{P^{t}f}{g}=  \acou{1\otimes \check{p}(x,-y)}{f(x)g(x.y^{-1})}=\acou{1}{f(x)\acou{\check{p}(-y)}{g(x.y^{-1})}}.
       \label{eq:PsiHDO.transpose-model-operator-transpose}
   \end{equation}
    Since $\check{p}(-y)=\check{p}^{t}(y)$ with $p^{t}(\xi)=p(-\xi)$, we obtain
    \begin{equation}
        \acou{P^{t}f}{g}= \acou{1}{f(x)\acou{\check{p}^{t}(y)}{g(x.y^{-1})}}=\acou{(p^{t}*f)(x)}{g(x)}.
    \end{equation}
    Thus $P^{t}$ is the left-convolution operator with symbol $p^{t}(\xi)=p(-\xi)$. 
    
    Now, since the model operator $(P^{t})^{a}$ is the 
    left-invariant \psivdo\ with symbol $\sigma_{m}(P^{t})(a,\xi)=\sigma_{m}(P)(a,-\xi)$ we see that it  agrees 
    with the transpose $(P^{a})^{t}$. 
    
    In the general case, when $\cE$ is not the trivial bundle, we can similarly show that $P^{t}$ is a \psivdo\ of order $m$ with principal symbol 
    $\sigma_{m}(P^{t})(x,\xi)=\sigma_{m}(P)(x,-\xi)^{t}$ and such that at every point $a\in M$ its model operator at $a$ is the transpose 
    $(P^{a})^{t}$.
\end{proof}

Assume now that $M$ is endowed with a positive density and $\cE$ with a Hermitian metric respectively and let $L^{2}(M,\cE)$ be the associated 
$L^{2}$-Hilbert space. 

\begin{proposition}\label{prop:PsiHDO.adjoint-manifold}
    Let $P \in \pvdo^{m}(M,\cE)$ have principal symbol $\sigma_{m}(P)$. Then:\smallskip

    1) The principal symbol of $P^{*}$ is $\sigma_{\bar{m}}(P^{*})(x,\xi)=\sigma_{m}(P)(x,\xi)^{*}$.\smallskip 
  
  2) If $P^{x}$ denotes the model operator of $P$ at $x\in M$ then the model operator of $P^{*}$ at $x$ is  
  the adjoint  $(P^{x})^{*}$ of $P^{x}$. 
 \end{proposition} 
 \begin{proof}
    Let us first assume that $\cE$ is the trivial line bundle, so that $P$ is a scalar operator. Moreover, since the above statements are local ones, it 
    is enough to prove them in a local Heisenberg chart $U \subset \Rd$ and we may assume  that $P$ is a \psivdo\ on $U$. 
    
    Let  $\overline{P}:C^{\infty}_{c}(U)\rightarrow C^{\infty}(U)$ be the conjugate operator of $P$, so that 
    $\overline{P}u=\overline{P(\overline{f})}$ for any  $f\in C^{\infty}_{c}(U)$. 
    By Proposition~\ref{prop:PsiVDO.characterisation-kernel2} 
    the distribution kernel of $P$ of the 
    form~(\ref{eq:PsiHDO.characterization-kernel.Heisenberg}) with $K_{P}(x,y)$ in $\cK_{\hat{m}}(\URd)$, so the kernel of $\overline{P}$ takes the form
    \begin{equation}
        k_{\overline{P}}(x,y)=\overline{k_{P}(x,y)}=|\varepsilon_{x}'|\overline{K_{P}(x,y)} \quad \bmod C^{\infty}(U\times U).
    \end{equation}
    
    Since the conjugation of distribution $K(x,y)\rightarrow \overline{K(x,y)}$ induces an anti-linear isomorphism from $\cK^{\hat{m}}(\URd)$ onto 
    $\cK^{\hat{\bar m}}(\URd)$, it follows from Proposition~\ref{prop:PsiVDO.characterisation-kernel2}  
    that $\overline{P}$ is a \psivdo\ of order $\hat{m}$ and its kernel can be 
    put into the form~(\ref{eq:PsiHDO.characterization-kernel.Heisenberg}) 
    with $K_{\overline{P}}(x,y)=\overline{K_{P}(x,y)}$. In particular, if $K_{P,\hat{m}}\in \cK_{\hat{m}}(\URd)$ denotes 
    the leading kernel of $K_{P}$ then the leading kernel of $K_{\overline{P}}$ is $ \overline{K_{P,\hat{m}}}$. Thus $\overline{P}$ has principal symbol
    \begin{equation}
        \sigma_{\bar{m}}(\overline{P})(x,\xi)=[\overline{K_{P,\hat{m}}}]^{\wedge}_{\xiy}(x,\xi)=\overline{[(K_{P,\hat{m}})^{\wedge}_{\xiy}(x,-\xi)}= 
        \overline{\sigma_{m}(x,-\xi)}.   
    \end{equation}
    
    Moreover, we have $\sigma_{\bar m}(\overline{P})^{\vee}_{\xiy}(x,y)= \overline{\sigma_{m}(P)^{\vee}_{\xiy}(x,y)}$. Therefore,  
    for any $f$ in $\cS_{0}(G_{a}U)$ the function $(\overline{P})^{a}f(x)$ is equal to
%     for any $a \in U$ the model operator of , since
    \begin{equation}
        \acou{\overline{\sigma_{m}(P)^{\vee}_{\xiy}(x,y)}}{f(x.y^{-1})}= \overline{\acou{\sigma_{m}(P)^{\vee}_{\xiy}(x,y)}{\overline{f(x.y^{-1}}}} 
        = \overline{P^{a}}f(x).
    \end{equation}
Hence $(\overline{P})^{a}$ agrees with $\overline{P^{a}}$. 
 
    Combining all this with Proposition~\ref{prop:PsiHDO.transpose-chart} and 
    Proposition~\ref{prop:PsiHDO.transpose-global} we see that $\overline{P}^{t}$ is a \psivdo\ of order $\overline{m}$ such that:\smallskip 
    
    - If we put the kernel of $\overline{P}^{t}$ into the form~(\ref{eq:PsiHDO.characterization-kernel.Heisenberg}) with a distribution 
    $K_{\overline{P}^{t}}(x,y)$ in $\cK_{\hat{\bar{m}}}(\URd)$, then the leading kernel of $K_{\overline{P}^{t}}$ is 
    $K_{\overline{P}^{t},\hat{\bar{m}}}=\overline{K_{P,\hat{m}}(x,-y)}$;\smallskip 
   
   - The global principal symbol of $\overline{P}^{t}$ is 
   $\sigma_{\bar{m}}(\overline{P}^{t})=\overline{\sigma_{m}(P^{t})(x,\xi)}=\overline{\sigma_{\bar m}(P)(x,\xi)}$;\smallskip
   
   - The model operator at $a \in U$ of $\overline{P}^{t}$ is $(\overline{P}^{t})^{a}=\overline{P^{a}}^{t}=(P^{a})^{*}$.\smallskip 
   
    Now, let $d\rho(x)=\rho(x)dx$ be the smooth positive density on $U$ coming from that of $M$. Then the adjoint $P^{*}: C^{\infty}_{c}(U)\rightarrow 
    C^{\infty}(U)$ of $P$ with respect to $d\rho$ is such that%, for functions $f$ and $g$ in  $C^{\infty}_{c}(U)$, we have
    \begin{equation}
        \int_{f} \overline{Pf(x)}g(x)\rho(x)dx = \int_{U}\overline{f(x)}P^{*}g(x)\rho(x)dx, \qquad f,g\in C^{\infty}_{c}(U).
    \end{equation}
    Thus $P^{*}=\rho^{-1}\overline{P}^{t}\rho$, which shows that $P^{*}$ is a \psivdo\ of order $\bar{m}$. Moreover, as in the proof of 
    Proposition~\ref{prop:PsiHDO.transpose-chart} in Appendix~\ref{chap:Appendix-transpose}, 
    we can prove that the kernel of $P^{*}$ can be put into the form~(\ref{eq:PsiHDO.characterization-kernel.Heisenberg}) with 
    $K_{P^{*}}(x,y)\in \cK^{\overline{m}}(\URd)$ equal to
       \begin{equation}
       \rho(x)^{-1}K_{\overline{P}^{t}}(x,y)\rho(\varepsilon_{x}^{-1}(y)) \sim \sum_{\alpha} \frac{1}{\alpha!} \rho(x)^{-1}
       \partial_{y}(\rho(\varepsilon_{x}^{-1}(y))_{|_{y=0}}K_{\overline{P}^{t}}(x,y).
    \end{equation}
In particular,  $K_{P^{*}}(x,y) $ and $K_{\overline{P}^{t}}(x,y)$ agree modulo $\cK^{\hat{\bar{m}}+1}(\URd)$, hence have same leading kernels. It then 
follows that $P^{*}$ and $\overline{P}^{t}$ have same principal symbol and same model operator at a point $a\in U$, that is, $\sigma_{\bar 
m}(P^{*})(a,\xi)=\overline{\sigma_{m}(P)(x,\xi)}$ and $(P^{*})^{a}=(P^{a})^{*}$. 

Finally, assume that $\cE$ is an arbitrary vector bundle of rank $r$, so that the restriction of $P$ to $U$ is given by a matrix $P=(P_{ij})$ of \psivdos\ of order $r$. 
Let $h(x)\in C^{\infty}(U,GL_{r}(\C))$, $h(x)^{*}=h(x)$, be the Hermitian metric on $U\times 
\C^{r}$ coming from that of $\cE$. Then the adjoint of $P$ with respect to this Hermitian metric is $P^{*}=\rho^{-1}h^{-1}\overline{P}^{t}h\rho$. Therefore, 
in the same way as in the scalar case we can prove that $P^{*}$ has principal symbol 
$h(x)^{-1}\overline{\sigma_{m}(P)(x,\xi)}^{t}h(x)=\sigma_{m}(P)(x,\xi)^{*}$ and its model operator at any point $a\in U$ is 
$h(a)^{-1}\overline{P^{a}}^{t}h(a)=(P^{a})^{*}$. 
\end{proof}
 
\section{Hypoellipticity and Rockland condition}
\label{sec:hypoellipticity}

In this section we define a Rockland condition for \psivdos\ and relate it to the invertibility of the principal symbol to get hypoellipticity 
criterions. 

\subsection{Principal symbol and Parametrices}
By~\cite[Sect.~18]{BG:CHM} in a local Heisenberg chart the invertibility of the local principal symbol of a \psivdo\ 
is equivalent to the existence of a \psivdo-parametrix. Using the 
global principal symbol we can give a global reformulation of this result as follows.

\begin{proposition}\label{thm:PsiHDO.hypoellipticity}
   Let $P:C^{\infty}_{c}(M,\cE) \rightarrow C^{\infty}(M,\cE)$ be a \psivdo\ of order $m$. The following are equivalent:\smallskip 

   1) The principal symbol $\sigma_{m}(P)$ of $P$ is invertible with respect to the convolution product for homogeneous 
    symbols;\smallskip 
    
    2) $P$ admits a parametrix $Q$ in $\pvdo^{-m}(M,\cE)$, so that~$PQ=QP=1  \bmod \psinf(M,\cE)$.
  \end{proposition}
 \begin{proof} 
     First, it follows from Proposition~\ref{prop:PsiHDO.composition2} that 2) implies 1). 
     Conversely, the formula~(\ref{eq:PsiHDO.global-local-convolution-symbols}) shows that  in a local trivializing Heisenberg chart 
     the invertibility of the global principal $\sigma_{m}(P)$ is equivalent to that of the local 
     principal symbol. Once the latter is granted then, as shown in~\cite[p.~142]{BG:CHM},  Lemma~\ref{lem:PsiHDO.asymptotic-completeness} 
     and Proposition~\ref{prop:PsiHDO.composition} allows us to carry out  in a local trivializing Heisenberg chart  
     the standard parametrix construction to obtain a parametrix for $P$ as a \psivdo\ of order $-m$. 
    A standard partition of the unity argument then allows us to construct a parametrix for $P$ in $\pvdo^{-m}(M,\cE)$.  
\end{proof}

When a \psivdo\ has an invertible principal symbol the Sobolev regularity properties of its parametrices allows us to get:  

\begin{proposition}[{\cite[p.~142]{BG:CHM}}]\label{prop:hypoellipticity-parametrix}
      Let $P:C^{\infty}_{c}(M,\cE) \rightarrow C^{\infty}(M,\cE)$ be a \psivdo\ of order $m$ with $\Re m\geq m$ and such that its
      principal symbol is invertible. Then $P$ is hypoelliptic with gain of $\frac{1}{2}\Re m$ derivatives, i.e.,  setting $k=\frac{1}{2}\Re m$, for any $a \in M$,  any 
   $u \in \cE'(M,\cE)$  and any $s \in \R$, we have 
   \begin{equation}
       \text{$Pu$ is $L^{2}_{s}$ near $a$} \ \Longrightarrow \  \text{$u$ is $L^{2}_{s+k}$ near $a$}.
         \label{eq:PsiHDO.local-hypoellipticity}
   \end{equation}
   In particular, if $M$ is compact then, for any reals $s$ and $s'$, we have the hypoelliptic estimate,
   \begin{equation}
       \|f\|_{L^{2}_{s+k}}\leq C_{s}(\|Pf\|_{L^{2}_{s}}+\|f\|_{L^{2}_{s'}}), \qquad f \in C^{\infty}(M,\cE).
        \label{eq:PsiHDO.subellipticity.subelliptic-estimates}
   \end{equation}
\end{proposition}

\begin{remark}
 It used to be customary to call the above property hypoellipticity with \emph{gain} of $k$ derivatives (see, e.g., \cite{BG:CHM}).  
 We have followed here the recent terminology of~\cite{Ko:SHLD},
 where is constructed an example of sum of squares which is hypoelliptic but, instead of gaining of derivatives as in~\cite{Ho:HSODE}, it gains derivatives.  
\end{remark}

\begin{remark}%\label{rem:hypoellipticity-parametrix}
 We can give sharper regularity results for \psivdos\ in terms of suitably weighted Sobolev spaces (see~\cite{FS:EDdbarbCAHG} and 
    Section~\ref{sec.Sobolev}). When $P$ is a differential operator and the Levi form is non-vanishing these results correspond to the maximal hypoellipticity 
    of $P$ as in~\cite{HN:HMOPCV}.
\end{remark}

\begin{remark}\label{rem:hypoellipticity-parametrix}
 As it follows from the proof in~\cite[p.~142]{BG:CHM} in order to have the hypoelliptic properties~(\ref{eq:PsiHDO.local-hypoellipticity}) 
    and~(\ref{eq:PsiHDO.subellipticity.subelliptic-estimates}) it is enough to have a left-parametrix for $Q$ rather than a two-sided parametrix. 
    Therefore, Proposition remains valid when we assume the principal symbol of $P$ to be left-invertible only. 
\end{remark}

% Given a $H$-frame $X_{0},\ldots,X_{d}$ and a multi-order $\alpha$ we let $X^{\alpha}
% \begin{definition}
%     Let $P\in \pvdo^{m}(M,\cE)$, $m \in \N$. We say that $P$ is hypoelliptic maximal when for any $a \in M$, given a $H$-frame $X_{0},..,X_{d}$ near 
%     $a$, for any $u \in \cE'(M,\cE)$ and any $s \in \R$, we have 
%     \begin{equation}
%        \text{$Pu$ is $L^{2}_{s}$ near $a$} \ \Longrightarrow \  \text{$X^{\alpha}u$ is $L^{2}_{s}$ near $a$ for any $\alpha\in \N^{d}$ such that 
%        $\brak \alpha\leq m$},
% %         \label{eq:}
%     \end{equation}
%     where we have let $X^{\alpha}=X_{0}^{\alpha_{0}}\ldots X_{d}^{\alpha_{d}}$. 
% \end{definition}
% 
% It is easy to check that maximal hypoellipticity implies hypoellipticity with loss of $2[\frac{m}{2}]$-derivatives. Moreover, at least 
% 
% 
% When $M$ is compact, combining this with 
% the compactness of the embedding of  $H^{k/2}(M,\cE)$ into $L^{2}(M,\cE)$ we get: 
% \begin{proposition}\label{prop:PsiHDO.spectrum}
%     Suppose $M$ compact and assume that $P$ has an invertible principal symbol and a 
%     spectrum different from $\C$. Then: \smallskip 
% 
%     1) The spectrum of $P$ consists of isolated eigenvalues with finite multiplicities.\smallskip 
% 
%     2) For any $\lambda\in \op{Sp}P$ the eigenspace $\ker (P-\lambda)$ is a finite dimensional subspace of $C^{\infty}(M,\cE)$.  
% \end{proposition}

\subsection{Rockland condition}
Assume that $M$ is endowed with a positive density and 
$\cE$ with a Hermitian metric and let $P:C^{\infty}_{c}(M,\cE)\rightarrow C^{\infty}(M,\cE)$ be a \psivdo\ of order $m$.  
Let $P^{a}$ be the model operator of $P$ at 
a point $a\in M$ 
and let  $\pi: G\rightarrow \cH_{\pi}$  be a (nontrivial) unitary representation of $G=G_{a}M$. We define the symbol 
$\pi_{P^{a}}$ as follows (see also~\cite{Ro:HHGRTC}, \cite{Gl:SSGMCAINGHG}, \cite{CGGP:POGD}).

Let $\cH_{\pi}^{0}(\cE_{a})$ be the subspace of $\cH_{\pi}(\cE_{a}):=\cH_{\pi}\otimes \cE_{a}$ spanned by the vectors, 
\begin{equation}
     \pi_{f}\xi=\int_{G}(\pi_{x}\otimes 1_{\cE_{a}})(\xi \otimes f(x)) dx, %\qquad \xi \in \cH_{\pi}(\cE_{a}), f \in % \cS_{0}(G), \quad .
\end{equation}
where $\xi$ ranges over  $\cH_{\pi}$ and $f$ over $\cS_{0}(G,\cE_{a})=\cS_{0}(G)\otimes \cE_{a}$.
Then $\pi_{P^{a}}$ is the (unbounded) operator of $\cH_{\pi}(\cE_{a})$ with domain 
$\cH_{\pi}^{0}(\cE_{a})$ such that
\begin{equation}
    \pi_{P^{a}}(\pi_{f}\xi)=\pi_{P^{a}f}\xi \qquad \forall f\in \cS_{0}(G,\cE_{a})\quad \forall \xi \in \cH_{\pi}.
\end{equation}
One can check that $\pi_{P^{a*}}$ is the adjoint of $\pi_{P^{a}}$ on $\cH_{\pi}^{0}$, hence 
is densely defined. Thus $\pi_{P^{a}}$ is closeable and we can let $\overline{\pi_{P^{a}}}$ denote its closure. 

In the sequel we let $C^{\infty}_{\pi}(\cE_{a})=C^{\infty}_{\pi}\otimes \cE_{a}$, where $C^{\infty}_{\pi}\subset \cH_{\pi}$ denotes the space of 
smooth vectors of $\pi$, i.e., the subspace of vectors $\xi \in \cH_{\pi}$ so that $x \rightarrow \pi(x)\xi$ is smooth from $G$ to $\cH_{\pi}$.

\begin{proposition}[\cite{CGGP:POGD}]\label{PsiHDO.properties-symbol-representation}
    1) The domain of $\overline{\pi_{P^{a}}}$ always contains $C^{\infty}_{\pi}(\cE_{a})$.\smallskip 
   
   2)  If $\Re m \leq 0$ then the operator $\overline{\pi_{P^{a}}}$ is bounded.\smallskip 
   
   3) We have $\overline{(\pi_{P^{a}})^{*}}=(\overline{\pi_{P^{a}}})^{*}$.\smallskip 
   
   4) If $P_{1}$ and $P_{2}$ are \psidos\ on $M$ then $\overline{\pi_{(P_{1}P_{2})^{a}}}=\overline{\pi_{P^{a}_{1}}}\, \overline{\pi_{P^{a}_{2}}}$.  
\end{proposition}

\begin{remark}
   If $\cE_{a}=\C$ and $P^{a}$ is a differentiable operator then, as it is left-invariant, $P^{a}$ belongs to the enveloping algebra 
$\cU(\fg)$ of the Lie algebra $\fg=\fg_{a}M$ of $G$. In this case $\overline{\pi_{P^{a}}}$ coincides on $C^{\infty}_{\pi}$ with the operator 
$d\pi(P^{a})$, where $d\pi$ is the representation of $\cU(\fg)$ induced by $\pi$.  
\end{remark}

\begin{definition}
    We say that $P$ satisfies the Rockland condition at $a$ if for any nontrivial unitary irreducible representation 
$\pi$ of $G_{a}M$ the operator $\overline{\pi_{P^{a}}}$ is injective on $C^{\infty}_{\pi}(\cE_{a})$.
\end{definition}

Set $2n=\rk \cL_{a}$. Under the identification $G=G_{a}M\simeq \bH^{2n+1}\times \R^{d-2n}$ given by Proposition~\ref{prop:Bundle.intrinsic.fiber-structure} 
there are left-invariant vector fields $X_{0},\ldots,X_{d}$ on $G$ 
such that $X_{0},\ldots,X_{2n}$ are given by~(\ref{eq:Examples.Heisenberg-left-invariant-basis}) 
and $X_{k}=\frac{\partial}{\partial x_{k}}$ for 
$k\geq 2n+1$. Then, up to unitary 
equivalence, the nontrivial irreducible representations of $G$  are of two types:\smallskip

(i) Infinite dimensional representations $\pi^{\lambda,\xi}:G\rightarrow L^{2}(\R^{n})$ parametrized by $\lambda\in \R\setminus0$ and 
$\xi=(\xi_{2n+1},\ldots,\xi_{2n})$ such that 
\begin{gather}
    d\pi^{\lambda,\xi}(X_{0})=i\lambda|\lambda|, \qquad d\pi^{\lambda,\xi}(X_{k})=i\lambda \xi_{k}, \quad k=2n+1,\ldots,d,\\
    d\pi^{\lambda,\xi}(X_{j})=|\lambda|\frac{\partial}{\partial \xi_{j}}, 
    \qquad d\pi^{\lambda,\xi}(X_{n+j})=i\lambda \xi_{j}, \quad j=1,\ldots,n.
%     \label{eq:}
\end{gather}
Moreover, in this case we have $C^{\infty}(\pi^{\pm,\xi})=\cS(\R^{n})$.

% for $j=1,\ldots,n$ and $k=2n+1,\ldots,d$.\smallskip 
% and given by 
% \begin{equation}
%     (\pi^{(\lambda,\xi)}(x_{0},\ldots,x_{d})f)(\xi_{1},\ldots,\xi_{n})=e^{i(\lambda x_{0}+x_{1}.(\xi_{1}+ \frac{1}{2}x_{n+1})+
%     \ldots+x_{n}.(\xi_{n}+\frac{1}{2}x_{2n})+x_{2n+1}.\xi_{2n+1}+\ldots+x_{d}.\xi_{d})}f(x_{1}+\xi_{n+1},\ldots,x_{n}+\xi_{2n}),
% %     \label{eq:}
% \end{equation}

(ii) One dimensional representations $\pi^{\xi}:G\rightarrow \C$ indexed by $\xi=(\xi_{1},\ldots,\xi_{d})$ in $\R^{d}\setminus 0$ such that
\begin{equation}
     d\pi^{\xi}(X_{0})=0, \qquad d\pi^{\xi}(X_{j})=i\xi_{j}, \qquad j=1,\ldots,d.
%     \label{eq:¥}
\end{equation}
% In this case $C^{\infty}(\pi^{\xi})=\C$. 
 
In particular, if $P=p_{m}(-iX)$ with $p \in S_{m}(G)$ then the homogeneity of the symbol $p$ implies that we have 
 $\pi^{\lambda,\xi}_{P}=|\lambda|^{m}\pi^{\pm,\xi}_{P}$, where 
 $\overline{\pi^{\pm,\xi}_{P}}=\overline{\pi^{\pm 1,\xi}_{P}}$ accordingly with the sign of $\lambda$, while 
 for the representations in (ii) we have $\overline{\pi^{\xi}_{P}}=\pi^{\xi}_{P}=p_{m}(0,\xi)$. Therefore, we obtain: 
 
 \begin{proposition}\label{prop:Rockland.reduction}
    Let $p_{m}\in S_{m}(G_{a}M)$. Then the Rockland condition for $P=p_{m}(-iX^{a})$ is satisfied if, and only if, the following two conditions hold:\smallskip 

    (i) The operators $\overline{\pi^{\pm,\xi}_{P}}$, $\xi \in \R^{d-2n}$, are injective on $\cS(\R^{n})$;\smallskip 
 
   (ii) The restriction of $p_{m}$ to $\{0\}\times (\R^{d}\setminus 0)\simeq H_{a}^{*}\setminus 0$ is pointwise invertible.
 \end{proposition}
% \begin{remark}
%    When $G_{a}M=\bH^{2n+1}$ the conditions~(i) and (ii) are essentially those considered by Taylor~\cite{Ta:NCMA}. 
% \end{remark}

\subsection{Parametrices and Rockland condition}
The aim of this subsection is to show that the Rockland condition is enough to insure us the invertibility of 
the principal symbol and the existence of a parametrix in the Heisenberg calculus. 

First, we deal with zero'th order \psivdos. In this case, we have: 

\begin{theorem}\label{thm:Chap3.Rockland-Parametrix-order0}
  Let $P:C^{\infty}_{c}(M,\cE)\rightarrow C^{\infty}(M,\cE)$ be a zero'th order \psivdo. Then the following are equivalent:\smallskip
  
  (i) The principal symbol of $P$ is invertible;\smallskip
  
 (ii) For every point $a \in M$ the model operator $P^{a}$ is invertible on $L^{2}(G_{a}M,\cE_{a})$.\smallskip
 
 (iii) $P$ and $P^{t}$ satisfy the Rockland condition at every point $a\in M$.\smallskip
 
\noindent Furthermore, if $M$ is endowed with smooth density~$>0$ and $\cE$ with a Hermitian metric, then in (iii) we can replace the Rockland condition for 
$P^{t}$ by that for $P^{*}$. 
% 
% In any case when the above conditions hold $P$ admits a parametrix in $\pvdo^{0}(M,\cE)$.
\end{theorem}
\begin{proof}
When  $GM$ is a (trivial) fiber bundle of Lie groups, i.e., the Levi form of $(M,H)$ has constant rank, 
the theorem can be deduced from the results of~\cite[Sect.~5]{CGGP:POGD}, which are based on an idea due to Christ~\cite{Ch:ISASIO} 
(see~\cite{Po:FCSAHOHM1}). As we shall now see elaborating on the same idea allows us to deal with the general case as well. 

First, assume that the principal symbol of $P$ is invertible, so that $P$ admits a parametrix $Q\in 
 \pvdo^{-m}(M,\cE)$, which without any loss of generality may be assumed to be properly property. Then,
 for any $a \in M$ the operators $Q^{a}P^{a}$ and $(Q^{a})^{t}(P^{a})^{t}$ are equal to $1$ on $\cS_{0}(G_{a}M,\cE_{a})$ 
 and $\cS_{0}(G_{a}M,\cE_{a}^{*})$ respectively. Therefore, it follows from Proposition~\ref{PsiHDO.properties-symbol-representation} 
 that for any nontrivial irreducible unitary representation $\pi$ of $G_{a}$ the operators 
 $\overline{\pi_{P^{a}}}$ and $\overline{\pi_{(P^{a})^{t}}}$ are injective on $C^{\infty}(\pi)$, i.e., $P$ and $P^{t}$ satisfy the Rockland condition 
 at every point of $M$. Hence (i) implies (iii). 

Second, by a result of G\l owacki~\cite[Thm.~4.3]{Gl:RCNDO2} for any $a\in M$ the Rockland condition for $P^{a}$ is equivalent to the 
 left-invertibility of $P^{a}$ on $L^{2}(G_{a}M,\cE_{a})$. The same is true for $(P^{t})^{a}=(P^{a})^{t}$, so $P^{a}$ has a two-sided inverse on 
 $L^{2}(G_{a}M,\cE_{a})$  if, and only if, $P^{a}$ and $(P^{t})^{a}$ satisfy the Rockland condition. Therefore, we see that the conditions (ii) 
 and (iii) are equivalent. 
 
 Notice also that if $M$ is endowed with smooth density~$>0$ and $\cE$ is endowed with a Hermitian metric then the above arguments can be carried through 
 without any changes if we replace the transpose of $P$ by the adjoint $P^{*}$. Therefore, in (iii) we may replace the Rockland condition for $P^{t}$ by that 
 for $P^{*}$. 
 
 It remains now to prove that (ii) implies (i). Observe that it is enough to proceed locally in a trivializing Heisenberg chart $U\subset \Rd$ with 
 $H$-frame $X_{0},\ldots,X_{d}$. In fact, we may further assume that $\cE$ is 
 a trivial line bundle since, as we sill see later, the arguments below can be carried out \emph{verbatim} for systems of \psivdos. Therefore, we are 
 reduced to prove: 

 \begin{proposition}\label{lem:Chap3.inverting-symbol}
     Let $p\in S_{0}(\URd)$ be such that for any $a \in U$ the operator $P^{a}=p(a,-iX^{a})$ is invertible on $L^{2}(\Rd)$. Then $p$ is invertible in 
     $S_{0}(\URd)$.
% %     i.e., there exists $q\in 
% %     S_{0}(\URd)$ such that $p*q=q*p=1$.
 \end{proposition}

The proof will follow from a series of lemmas. 
% of Proposition~\ref{lem:Chap3.inverting-symbol} will be the consequence of several intermediate lemmas. 
% The proof of this proposition follows from an idea due to Christ~(\cite{Ch:}, \cite{CGGP:POGD}) and is based on the following two lemmas
In the sequel we let $S=\{x\in \Rd; \|x\|=1\}$ and we endow its with its induced Riemannian metric. Let $K\in \cS'(\Rd)$ be 
homogeneous of degree $-(d+2)$ and such that $K_{|_{S^{1}}}$ is in $L^{2}(S)$. Then $K$ is of the form, 
\begin{equation}
    K=\op{pv}(K)+c(K)\delta_{0}, \qquad c(K) \in \C,
     \label{eq:Chap3.decomposition-homogeneous-distribution}
\end{equation}
where $c(K)$ is a complex constant and $\op{pv}(K)$ is the principal value distribution, 
\begin{equation}
    \acou{\op{pv}(K)}{f}=\lim_{\epsilon\rightarrow 0^{+}} \int_{|x|>\epsilon}K(x)f(x)dx, \qquad f \in \cS(\Rd).
%     \label{eq:}
\end{equation}
The decomposition~(\ref{eq:Chap3.decomposition-homogeneous-distribution}) 
is unique and, in particular, $K$ is uniquely determined by its restriction to $S$ and the constant $c(K)$ 
(see~\cite[Prop.~6.13]{FS:HSHG}, \cite[Lem.~2.4]{Ch:ORISIO}). 

\begin{definition}
    $\cA_{k}$, $k \in \N$, consists of distributions $K\in \cS'(\Rd)$ that are homogeneous of degree $-(d+2)$ and such that 
    $(\partial^{\alpha}K)_{|_{S}}\in L^{2}(S)$ for $\brak\alpha \leq k$. 
\end{definition}

We turn  $\cA_{k}$ into a Banach space by endowing it with the norm,
\begin{equation}
    |K|_{k}=|c(K)|+\sum_{\brak\alpha \leq k}\|(\partial^{\alpha}K)_{|_{S}}\|_{L^{2}(S)}, \qquad K\in \cA_{k}. 
%     \label{eq:¥}
\end{equation}
Notice that $\cK_{-(d+2)}(\Rd)$ is contained in all the spaces $\cA_{k}$, $k \in \N$, and the norms $|.|_{k}$'s give rise to a system 
of semi-norms on $\cK_{-(d+2)}(\Rd)$ whose corresponding topology is the weakest topology making the maps $K\rightarrow c(K)$ and $K\rightarrow 
K_{|_{S}}$ be continuous from $\cK_{-(d+2)}(\Rd)$ to $\C$ and $C^{\infty}(S)$ respectively.

In the sequel for $a \in U$ and $K \in \cA_{k}$ we let $P_{K}^{a}$ denote the convolution operator $P_{K}^{a}f=K*^{a}f$, $f\in 
\cS_{0}(\Rd)$. Then we have: 

% \begin{lemma} \label{lem:Chap3.L2-boundedness1}
% Let $a \in U$. Then for any $K\in  \cB_{1}$  the operator $P_{K}^{a}$ extends to a bounded operator on $L^{2}(\Rd)$. Furthermore, 
% there exists $C(a)>0$ such that for any $K\in  \cB_{1}$ we have
%     \begin{equation}
%         \|P_{K}^{a}\| \leq C(a)|K|_{\cB_{1}},
% %         \label{eq:¥}
%     \end{equation}
%     where $\|.\|$ denotes the operator norm on $\cL(L^{2}(\Rd))$.
%     %and the constant $C(a)$ depends continuously on $a$.
% \end{lemma}
% 
% In fact, it follows from the proof in \cite[Thm.~6.19]{FS:HSHG} that the constant $C(a)$ above can be chosen independently of $a$ provided that $a$ 
% remains in a compact set of $U$. Combining this with the continuous embedding of $\cA_{k}$ into $\cB_{1}$ for $k\geq d+3$ then gives: 

\begin{lemma}[{\cite[Thm.~1]{KS:IOSSG}}, {\cite[Thm.~6.19]{FS:HSHG}}, \cite{Ch:ISASIO}]\label{lem:Chap3.Rockland.L2-boundedness} 
Let $k$ be an integer greater than or equal to $d+3$. Then:\smallskip

    (i) For any $a \in U$ and any $K\in \cA_{k}$  the operator $ P_{K}^{a}$ extends to a 
bounded operator on $L^{2}(\Rd)$.\smallskip

 (ii) For any compact $L\subset U$ there exists $C_{Lk}>0$ such that for any $a \in L$ and any $K\in \cA_{k}$ we have 
\begin{equation}
     \|P_{K}^{a}\| \leq C_{Lk}|K|_{k}.
     \label{eq:Rockland.Operator-norm-PK-Ak}
\end{equation}
\end{lemma}
\begin{proof}
Let $\cB_{1}$ be the space of distributions $K\in \cS'(\Rd)$ that are homogeneous of degree $-(d+2)$ and such that $K_{|_{S}}$ is $C^{1}$. 
This becomes a Banach space when endowed with the norm, 
\begin{equation}
    |K|_{\cB_{1}}=|c(K)|+\|K_{|_{S}}\|_{1}, \qquad K\in \cB_{1},
%     \label{eq:¥}
\end{equation}
where $\|K_{|_{S}}\|_{1}$ is some Banach norm on $C^{1}(S)$. Then it follows from a theorem of Knapp-Stein~\cite[Thm.~1]{KS:IOSSG} (see 
also~\cite[Thm.~6.19]{FS:HSHG}) that, given $a \in U$, for any 
$K\in \cB_{1}$ the operator $P_{K^{a}}^{a}$ extends to a bounded operator on $L^{2}(\Rd)$ and there exists a constant $C_{a}>0$ such that for 
any $K\in  \cB_{1}$ we have
    \begin{equation}
        \|P_{K}^{a}\| \leq C(a)|K|_{\cB_{1}}. 
%         \label{eq:¥}
    \end{equation}
Furthermore, it follows from the proof in~\cite{FS:HSHG} that the constant $C_{a}$ can be chosen independently of $a$ provided that $a$ 
remains in a compact set of $U$.

On the other hand,  if $K\in \cA_{k}$ with $k\geq d+3$ then 
$(\partial^{\alpha}K)_{|_{S}}$ is in $L^{2}(S)$ for $|\alpha| \leq \frac{d+3}{2}$. Therefore, it follows from the Sobolev embedding theorem that $\cA_{k}$ embeds 
continuously into $\cB_{1}$. Combining this with the first part of the proof then gives the lemma.
\end{proof}

\begin{remark}
The $L^{2}$-boundedness of $P_{K}^{a}$ above is actually true for any 
$K\in \cA_{0}$ (see~\cite{Ch:ISASIO}), but the uniform dependence with respect to $a$ is more difficult to keep track in~\cite{Ch:ISASIO} than 
in~\cite{FS:HSHG}. This is not really relevant in the sequel, because in order to obtain Theorem~\ref{thm:Chap3.Rockland-Parametrix-order0}
it is enough to prove Lemma~\ref{lem:Chap3.Rockland.L2-boundedness} for $k$ large enough.
\end{remark}
     
% The proof in~\cite{Ch:ISASIO} is more involved than 
% in~\cite{FS:HSHG} and in this the uniform dependence with respect to $a$ is somewhat less apparent than in the latter, 
% This is irrelevant in the sequel since we need only the dependence to be uniform for $k$ large enough.

\begin{lemma}[{\cite[Lem.~5.7]{CGGP:POGD}}] \label{lem:Chap3.Rockland.norm-product-Ak} Let $k$ be an integer~$\geq d+3$. Then:\smallskip
    
    (i) For any $a \in U$ the convolution $*^{a}$ induces a bilinear product on $\cA_{k}$.\smallskip
   
   (ii) For any compact $L\subset U$ there exists $C_{Lk}>0$ such that, for any $a \in L$ and any $K_{1}$ and $K_{2}$ in $\cA_{k}$, we have
\begin{equation}
    |K_{1}*^{a}K_{2}|_{k}\leq C_{Lk}(\|P_{K_{1}}^{a}\||K_{2}|_{k}+|K_{1}|_{k}\|P_{K_{2}}^{a}\|).
    \label{eq:Rockland.norm-product-Ak}
\end{equation}
\end{lemma}
\begin{remark}
    The fact that the constant in~(\ref{eq:Rockland.norm-product-Ak}) can be chosen independently of $a$ when $a$ remains in a compact set of $U$ is not 
    explicitly stated in~\cite{CGGP:POGD}, but this follows from its proof, noticing that the proof of the Lemma~2.10 of~\cite{Ch:ORISIO} 
    shows that in the Lemma~5.8 of \cite{CGGP:POGD} the constant $C$ can be chosen independently of $a$ when $a$ stays in a compact set of $U$.
\end{remark}

As we will see below it is essential that the estimate~(\ref{eq:Rockland.norm-product-Ak}) 
involves also the operator norm and not just the norm of $\cA_{k}$. This trick is 
initially due to Christ~\cite{Ch:ISASIO}.  Note also that 
Lemma~\ref{lem:Chap3.Rockland.norm-product-Ak} allows us to endow $L^{\infty}_{\loc}(U,\cA_{k})$ with the convolution product, 
\begin{equation}
    (K_{1}^{a})_{a\in U}*(K_{2}^{a})_{a\in U}= (K_{1}^{a}*^{a}K_{2}^{a})_{a\in U}, \qquad (K_{j}^{a})_{a\in U}\in L^{\infty}_{\loc}(U,\cA_{k}).
%     \label{eq:}
\end{equation}
In particular, the constant family $(\delta_{0})_{a \in U}$ is a unit for this product.

\begin{lemma}\label{lem:Chap3.Rockland.inverse}
    Let $k \in \N$, $k\geq d+3$, let $L\subset U$ be compact and consider a family $K=(K^{a})_{a \in L}$ in $L^{\infty}(L,\cA_{k})$  
    such that $\sup_{a \in L}\|P_{K^{a}}^{a}\|<1$.  Then $\delta_{0}- K$ is invertible in $L^{\infty}(L,\cA_{k})$. 
\end{lemma}
\begin{proof}
    Let $C_{Lk}$ be the sharpest constant in the estimate~(\ref{eq:Rockland.norm-product-Ak}) and let us endow $L^{\infty}(L,\cA_{k})$ with the Banach norm, 
    \begin{equation}
        |K|_{k,L}= \frac{1}{C_{Lk}}\sup_{a \in L}|K^{a}|_{k}, \qquad K=(K^{a})_{a\in U}\in L^{\infty}(L,\cA_{k}).
%         \label{eq:¥}
    \end{equation}
    Then for any $K_{1}$ and $K_{2}$ in $L^{\infty}(L,\cA_{k})$ we have 
    \begin{equation}
        |K_{1}*K_{2}|_{k,L}\leq (|K_{1}|_{k,L}\|P_{K_{2}}\|_{L}+\|P_{K_{1}}\|_{L}|K_{2}|_{k,L}). 
         \label{eq:Rockland.norm-product-AkL}
    \end{equation}
    where we have let $\|P_{K_{j}}\|=\sup_{a \in L}\|P_{K^{a}_{j}}^{a}\|$. 
    
    For $j=0,1,2,\ldots$ let $K^{(j)}$ be the $j$'th power of $K$ with respect to the product $*$ on $L^{\infty}(L,\cA_{k})$. 
    From~(\ref{eq:Rockland.norm-product-AkL}) we deduce 
    by induction that for $j=1,2,\ldots$ we have 
    \begin{equation}
        |K^{(j)}|_{k,L}\leq (j-1)|K|_{k,L}\|P_{K}\|_{L}^{j-1}.
%         \label{eq:}
    \end{equation}
   It follows that $\limsup  |K^{(j)}|_{k,L}^{1/j}=\|P_{K}\|_{L}<1$, so that the series $\sum_{j=0}^{\infty} K^{(j)}$ converges normally in 
   $L^{\infty}(L,\cA_{k})$ to the inverse of $\delta_{0}-K$. Hence the result.
\end{proof}

With all this preparation we are ready to prove Proposition~\ref{lem:Chap3.inverting-symbol}. 
\begin{proof}[Proof of Proposition~\ref{lem:Chap3.inverting-symbol}]
      Let $p\in S_{0}(\URd)$ be such that for any $a \in U$ the operator $P^{a}=p(a,-iX^{a})$ is invertible on $L^{2}(\Rd)$. Then for any $a \in U$ the operators 
  $(P^{a})^{*}P^{a}$ and $P^{a}(P^{a})^{*}$ are invertible on $L^{2}(\Rd)$ and have respective symbols $\overline{p}*p$ and $p*\overline{p}$ in 
  $S_{0}(\URd)$. Moreover, if $\overline{p}*p$ and $p*\overline{p}$ admits respective inverses $q_{1}$ and $q_{2}$ in $S_{0}(\URd)$, then 
  $q_{1}*\overline{p}$ and $\overline{p}*q_{2}$ are respectively a left-inverse and a right-inverse for $p$. Therefore, it is enough to prove that 
  $\overline{p}*p$ and $p*\overline{p}$ are invertible. Incidentally, we may assume that $P^{a}$ is a positive operator on $L^{2}(\Rd)$ for any $a \in 
  U$. 
  
 On the other hand, recall that for any $a \in U$ the convolution product $*^{a}$ on $\cK_{-(d+2)}(\URd)$ corresponds under the Fourier transform to 
 the product $*^{a}$ on $S_{0}(\Rd)$. Since the latter depends smoothly on $a$ and by~\cite[Prop.~15.30]{BG:CHM} the Fourier transform is a topological 
 isomorphism 
 from $\cK_{-(d+2)}(\Rd)$ onto $S_{0}(\Rd)$ it follows that the convolution product $*^{a}$ depends smoothly on $a$, i.e., gives rise to smooth family 
 of bilinear maps from $\cK_{-(d+2)}(\URd)\times \cK_{-(d+2)}(\URd)$ to $\cK_{-(d+2)}(\Rd)$. Therefore, we get 
 a convolution product on $\cK_{-(d+2)}(\URd)$ by letting
 \begin{equation}
     K_{1}*K_{2}(x,.)=K_{1}(x,,)*^{x}K_{2}(x,.), \qquad K_{j}\in \cK_{-(d+2)}.
%      \label{eq:}
 \end{equation}
 
   Bearing all this in mind, define $K(x,y)=\check{p}_{\xiy}(x,y)$. This is an element of $\cK_{-(d+2)}(\URd)$ which we can write as a smooth family 
   $(K^{a})_{a\in U}=(K(a,.))_{a\in U}$ with values in $\cK_{(d+2)}(\Rd)$.
  
  \begin{claim}
      The family $(K^{a})_{a\in U}$ is invertible in $L^{\infty}_{\loc}(U,\cK_{-(d+2)}(\Rd))$. 
  \end{claim}
  \begin{proof}[Proof of the claim]
      It is enough to show that for any integer $k\geq d+3$ and any compact $L\subset U$ the family $(K^{a})_{a\in L}$ is invertible in 
      $L^{\infty}_{\loc}(L,\cA_{k})$. 
      
      By assumption for any $a\in L$ the operator $P_{K^{a}}=P^{a}$ is an  invertible positive operator on $L^{2}(\Rd)$.  Since
     by  Lemma~\ref{lem:Chap3.Rockland.L2-boundedness} the family $(P_{K^{a}}^{a})_{a \in L}$ is bounded in 
     $\cL(L^{2}(\Rd)$, the same is true for the family $((P_{K^{a}}^{a})^{-1})_{a\in U}$.  
   
     Furthermore, since $P_{K^{a}}=P^{a}$ is positive its spectrum is contained in the interval 
      $[\|(P^{a})^{-1}\|^{-1},\|P^{a}\|]$. Therefore, there exist constants $c_{1}$ and $c_{2}$ with $0<c_{1}<c_{2}$ such that for any $a\in L$ the 
      spectrum of $P_{K^{a}}^{a}$ is contained in the interval $[c_{1},c_{2}]$. 
      Without any loss of generality we may assume $c_{2}<1$. Then the 
      spectrum of $1-P_{K^{a}}^{a}$ is contained in $[1-c_{2},1-c_{1}]$ for any $a\in L$, hence we get 
      \begin{equation}
          \sup_{a \in L}\|1-P_{K^{a}}^{a}\|\leq 1-c_{1}<1. 
%           \label{eq:}
      \end{equation}
      Since $1-P_{K^{a}}^{a}=P_{\delta_{0}-K^{a}}^{a}$ it follows from Lemma~\ref{lem:Chap3.Rockland.inverse} 
      that $(K^{a})_{a\in L}$ is invertible in $L^{\infty}(L,\cA_{k})$. 
      The proof is therefore complete.  
\end{proof}
 
Let $J=(J^{a})_{a\in U}\in L^{\infty}_{\loc}(U,\cK_{-(d+2)}(\Rd))$ be the inverse of $K$. The next step is to prove: 

\begin{claim}
    The family $(J^{a})_{a\in U}$ belongs to $C^{\infty}(U,\cK_{-(d+2)})$.
\end{claim}
\begin{proof}
    For any $a$ and $b$ in $U$ we have the equality 
    \begin{equation}
        J^{b}-J^{a}=-J^{a}*^{a}(K^{b}-K^{a})*^{b}J^{b}.
%         \label{eq:¥}
    \end{equation}
    We know that $\lim_{b \rightarrow a}(K^{b}-K^{a})=0$ in $\cK_{-(d+2)}(\Rd)$ and that the convolution product $*^{a}$ on $\cK_{-(d+2)}$ depends 
    smoothly on $a$. As near $a$ the family $(J^{a})_{a\in U}$ is bounded in 
    $\cK_{-(d+2)}(\Rd)$ it follows that we have $\lim_{b\rightarrow a}J^{b}=J^{a}$ in 
    $\cK_{-(d+2)}(\Rd)$. Hence the family $(J_{a})_{a \in U}$ depends continuously on $a$. 
    
    Now, if we let $e_{1},\ldots,e_{d+1}$ be the canonical basis of $\Rd$ then we have 
%     and fact that $(J_{a})_{a \in U}$ is 
%     continuous family and $(K^{a})_{a\in U}$ is a smooth family that in $\cK_{-(d+2)}(\Rd)$ 
%     we have 
    \begin{multline}
        \lim_{t \rightarrow 0}\frac{1}{t}(J^{a+te_{j}}-J^{a})=-\lim_{t\rightarrow 
        0}J^{a}*^{a}\frac{1}{t}(K^{a+te_{j}}-K^{a})*^{a+te_{j}}J^{a+te_{j}}\\ = 
        -J^{a}*^{a}\partial_{a_{j}}K^{a}*^{a}J^{a}.
%         \label{eq:¥}
    \end{multline}
    Thus $\partial_{a_{j}}J^{a}$ exists and is equal to $ -J^{a}*^{a}\partial_{a_{j}}K^{a}*^{a}J^{a}$. The latter depends continuously on $a$, so 
    $(J^{a})_{a\in U}$ is a $C^{1}$-family with values on $\cK_{-(d+2)}(\Rd)$. An induction then shows that it is actually of class $C^{k}$ for every 
    integer $k\geq 1$, hence is a smooth family. The claim is thus proved.  
\end{proof}

We now can complete the proof of Proposition~\ref{lem:Chap3.inverting-symbol}. For $(x,\xi)\in \URdo$ let $q(x,\xi)=\hat{J}_{\yxi}(x,\xi)$. Then $q$ belongs to 
$S_{0}(\URd)$ and we have $(q*p)^{\vee}_{\xiy}=J*K=\delta_{0}$, 
hence $q*p=1$. Similarly, we have $p*q=1$, so $q$ is an inverse for $p$ in $S_{0}(\URd)$. The lemma is thus proved.   
% of Proposition~\ref{lem:Chap3.inverting-symbol}
% is therefore achieved. 
\end{proof}
 
All this proves that in Theorem~\ref{thm:Chap3.Rockland-Parametrix-order0} the condition (ii) implies (i) when $\cE$ is the trivial line bundle 
over $M$. 
In fact, all the Lemmas~\ref{lem:Chap3.Rockland.L2-boundedness}--\ref{lem:Chap3.Rockland.inverse} 
are true for systems as well. More precisely, for Lemma~\ref{lem:Chap3.Rockland.L2-boundedness}
this follows from the fact that the proof of~\cite[Thm.~6.19]{FS:HSHG} can be carried out \emph{verbatim} 
for systems, while the extension to systems of Lemmas~\ref{lem:Chap3.Rockland.norm-product-Ak} and~\ref{lem:Chap3.Rockland.inverse} is immediate. 
Henceforth, the arguments in the proof of Proposition~\ref{lem:Chap3.inverting-symbol} remain valid
for systems. It then follows that for general bundles too~(ii) implies (i). 
The proof of Theorem~\ref{thm:Chap3.Rockland-Parametrix-order0}  is thus achieved. 
\end{proof}
% 
% \begin{remark}\label{rem:Chap3.Rockland-Parametrix-order0}
% If $M$ is endowed with smooth density~$>0$ and $\cE$ with a Hermitian metric then in (ii) we can replace the Rockland condition for $P^{t}$ by that for $P^{*}$, 
% as the arguments in the proof of Theorem~\ref{thm:Chap3.Rockland-Parametrix-order0} remain valid when we replace $P^{t}$ by $P^{*}$. 
% In particular, if $P$ is selfadjoint 
% then the validity at every point of the Rockland condition for $P$ only is equivalent to the invertibility of the principal symbol of $P$.
% \end{remark}

\begin{theorem}\label{thm:Chap3.Rockland-Parametrix-order+}
Let $P:C^{\infty}_{c}(M,\cE)\rightarrow C^{\infty}(M,\cE)$ be a \psivdo\ of integer order $m\geq 1$. Then the following are equivalent:\smallskip
  
  (i) The principal symbol of $P$ is invertible;\smallskip
  
  (ii) $P$ and $P^{t}$ satisfy the Rockland condition at every point $a\in M$.\smallskip
   
\noindent Furthermore, when $M$ is endowed with smooth density~$>0$ and $\cE$ with a Hermitian metric then in (ii) we can replace the Rockland condition for 
$P^{t}$ by that for $P^{*}$. In any case, when (i) or (ii) holds the operator $P$ admits a parametrix in $\pvdo^{-m}(M,\cE)$.
\end{theorem}
\begin{proof}
 First, in the same way as in the proof of Theorem~\ref{thm:Chap3.Rockland-Parametrix-order0} 
 we can show that if the principal symbol of $P$ is invertible then $P$ and $P^{t}$ satisfy 
 the Rockland condition at every point $a\in M$.  
    
Conversely, assume that $P$ and $P^{t}$ satisfy the Rockland condition at every point. Without any loss 
 of generality we may assume that $P$ and $P^{t}$ are properly supported. 
 In order to show that the principal symbol of $P$ is invertible it is enough to proceed locally in a trivializing Heisenberg chart $U\subset \Rd$ with 
 $H$-frame $X_{0},\ldots,X_{d}$. As a consequence we may assume that $\cE$ is the trivial line bundle and 
 $P$ is a scalar \psivdo, since the arguments below can be carried out \emph{verbatim} for systems of \psivdos. 
 
 On $U$ consider the sublaplacian $\Delta=-(X_{1}^{2}+\ldots+X_{d}^{2})+X_{0}$. By~\cite[Thm.~18.4]{BG:CHM} 
 the principal symbol of $\Delta$ is invertible, so $\Delta$ admits a parametrix $Q$ in $\pvdo^{-2}(U)$, which may be assumed to be properly 
 supported. 
 Then for proving 
 that the principal symbol of $P$ is left-invertible it is enough to check the invertibility of that of $Q^{m}P^{t}P$, for if $q_{0}(x,\xi)$ is the inverse of 
 $\sigma_{0}(Q^{m}P^{t}P)$ then $q*\sigma_{-m}(Q^{m}P^{t})$ is a left-inverse for $\sigma_{m}(P)$. Similarly, to show that $\sigma_{m}(P)$ is 
 right-invertible it is sufficient to prove that the principal symbol of $PP^{t}Q^{m}$ is invertible. 
 
 Next, as the operators 
 $P^{a}$, $(P^{t})^{a}=(P^{a})^{t}$  and $(Q^{m})^{a}=(Q^{a})^{m}$  satisfy the Rockland condition at every point $a\in U$, the same is true 
 for the operators $(Q^{m}P^{t}P)^{a}=(Q^{m})^{a}(P^{t})^{a}P^{a}$ and $(PP^{t}Q^{m})^{a}=P^{a}(P^{t})^{a}(Q^{m})^{a}$, that is,  
 $Q^{m}P^{t}P$ and 
 $PP^{t}Q^{m}$ satisfy the Rockland condition at every point. Furthermore, as $Q^{t}$ has an invertible principal symbol it satisfies the Rockland 
 condition at every point. Therefore, by arguing as above we see that the operators $(Q^{m}P^{t}P)=P^{t}P(Q^{t})^{m}$ and 
 $(PP^{t}Q^{m})^{t}=(Q^{t})^{m}PP^{t}$ satisfy the Rockland condition at every point. 
 
 Now, $Q^{m}P^{t}P$ and $PP^{t}Q^{m}$ both have order 0, so it follows from Theorem~\ref{thm:Chap3.Rockland-Parametrix-order0} 
 that their principal symbol are invertible. As alluded to 
 above this implies that the principal symbol of $P$ is invertible. 
 
 Finally, when $M$ is endowed with a smooth density~$>0$ and $\cE$ with a Hermitian metric, the above arguments remain valid 
 if we replace $P^{t}$ by the adjoint $P^{*}$, so that in (ii) we may replace the Rockland condition for $P^{t}$ by that for $P^{*}$. 
%  and we know that in this case (ii) and (iii) are equivalent, so at every point $a \in U$ 
%  the model operators  $(Q^{m}P^{t}P)^{a}$ and $(PP^{t}Q^{m})^{a}$ are invertible on $L^{2}(G_{a}U)$. Therefore, in order to prove that 
%  the principal symbol of $P$ is invertible it is sufficient to prove: 
\end{proof}

\begin{remark}%\label{rem:Chap3.Rockland-Parametrix}
In Chapter~\ref{chap.complex-powers} we will extend Theorem~\ref{thm:Chap3.Rockland-Parametrix-order+} to \psivdos\ with non-integer 
orders. 
% $P$ is an arbitrary complex number. 
\end{remark}

% \begin{remark}\label{rem:Chap3.Rockland-Parametrix-order+}
% If , as in Remark~\ref{rem:Chap3.Rockland-Parametrix-order0},   
% in (ii) we can replace the Rockland condition for $P^{t}$ by that for $P^{*}$. In particular, if $P$ is selfadjoint 
% then the validity at every point of the Rockland condition for $P$ only is equivalent to the invertibility of the principal symbol of $P$.
% \end{remark}

Let us now mention few applications of Theorems~\ref{thm:Chap3.Rockland-Parametrix-order0} 
and~\ref{thm:Chap3.Rockland-Parametrix-order+}. First, we have the following hypoellipticity criterion: 

\begin{proposition}\label{prop:Heisenberg.Rockland-hypoellipticity}
  Let $P:C^{\infty}_{c}(M,\cE)\rightarrow C^{\infty}(M,\cE)$ be a \psivdo\ of integer order $m\geq 0$ such that $P$ satisfies the Rockland condition 
  at every point. Then $P$ is hypoelliptic with gain of $\frac{m}{2}$-derivatives in the sense of~(\ref{eq:PsiHDO.local-hypoellipticity}). 
\end{proposition}
\begin{proof}
Without any loss of generality we may assume that $M$ is endowed with smooth density~$>0$ and $\cE$ with a Hermitian metric 
and we may also assume that $P$ is properly supported. As explained in Remark~\ref{rem:hypoellipticity-parametrix} 
a sufficient condition for $P$ to be hypoelliptic with gain of 
$\frac{m}{2}$-derivatives it is that its principal 
 symbol is left-invertible. To this end it is enough to show that the principal symbol of $P^{*}P$ is invertible, because if $q$ is an inverse for 
 $\sigma_{2m}(P^{*}P)=\sigma_{m}(P)^{*}*\sigma_{m}(P)$ then $q*\sigma_{m}(P)^{*}$ is a left inverse for $\sigma_{m}(P)$. 

 On the other hand, by Theorems~\ref{thm:Chap3.Rockland-Parametrix-order0} and~\ref{thm:Chap3.Rockland-Parametrix-order+} 
 the principal symbol of $P^{*}P$ is invertible if, and only if, $P^{*}P$ 
 satisfies the Rockland condition at every point. Observe that for any $a\in M$ and any nontrivial irreducible representation of $G_{a}M$ the 
 operators $\pi_{(P^{*}P)^{a}}=\pi_{P^{a}}^{*}\pi_{P^{a}}$ and $\pi_{P^{a}}$ have same kernels, so the Rockland conditions at $a$ for $P^{*}P$ and $P$  
 are equivalent. 
 
Combining all this, we see that if $P$ satisfies the Rockland condition at every point then $P$ is hypoelliptic with gain of $\frac{m}{2}$-derivatives.
\end{proof}

Next, even though the representation theory of the tangent group may vary from point to point, 
the Rockland condition and the invertibility of the principal symbol are open properties. Indeed, we have: 

\begin{proposition}\label{prop:Chap3.Rockland.open}
    Let $P:C^{\infty}_{c}(M,\cE)\rightarrow C^{\infty}(M,\cE)$ be a \psivdo\ of integer order $m$ with principal symbol $p_{m}(x,\xi)$
    and let $a\in M$.\smallskip
    
    1) If $P$ satisfies the Rockland condition at $a$ then there exists an open neighborhood $V$ of $a$ such that $P$ satisfies the Rockland 
    condition at every point of $V$.\smallskip
    
    2) If $p_{m}(a,\xi)$ is invertible in $S_{m}(\fg^{*}_{a}M,\cE_{a})$ then there exists an open neighborhood $V$ of $a$ such that 
    $p_{m|_{V}}$ is invertible on $S_{m}(\fg^{*}V,\cE)$.
\end{proposition}
\begin{proof}
  Let us first prove 2). To achieve this it is enough to proceed in a trivializing Heisenberg chart $U$ near $a$ and we may assume that $\cE$ is the 
  trivial line bundle since the proof for systems follows along similar lines. Furthermore, arguing as in the proof of 
  Theorem~\ref{thm:Chap3.Rockland-Parametrix-order+}  allows us to reduce the proof to the case $m=0$. 
 
  Assume now that $P$ has order $0$.  Then by Theorem~\ref{lem:Chap3.Rockland.L2-boundedness} 
  the invertibility of $p_{0}(a,.)$ in $S_{0}(\Rd)$ with respect to the product $*^{a}$ is equivalent to the 
  invertibility of the model operator $P^{a}$ on $L^{2}(\Rd)$. By Lemma~\ref{lem:Chap3.Rockland.L2-boundedness}
  the $L^{2}$-extension of $P^{x}$ depends 
  continuously on $x$, so there exists an open neighborhood $V$ of $a$ such that $P^{x}$ is  invertible on $L^{2}(\Rd)$ 
  for any $x\in V$. Then Theorem~\ref{thm:Chap3.Rockland-Parametrix-order0} insures us that $p_{0|_{V}}$ is invertible in $S_{0}(V\times \Rd)$. Hence 
  the result.  
  
  Let us now deduce 1) from 2). Assume 
  that $P$ satisfies the Rockland condition at $a$ and let us endow $M$ with a density~$>0$ and $\cE$ with a Hermitian metric. Then $P^{*}P$ 
    satisfies the Rockland condition at $a$, so by Theorem~\ref{thm:Chap3.Rockland-Parametrix-order+} its principal symbol is invertible at $a$. Therefore, 
    the principal symbol of $P$ is left-invertible on $V$ and along similar lines as that at the beginning  of Theorem~\ref{thm:Chap3.Rockland-Parametrix-order0} 
    we can show that $P$ satisfies the 
    Rockland condition at every point of $V$.
 \end{proof}

Finally, we look at families of  invertible symbols parametrized by a smooth manifold $B$. As in Proposition~\ref{prop:PsiHDO.continuity-product-symbol} 
each space $S_{m}(\fg^{*}M,\cE)$ is 
endowed with the topology inherited from that of $C^{\infty}(\fg^{*}M\setminus 0, \End \cE)$. 
In addition, if $p\in S_{m}(\fg^{*}M,\cE)$ is invertible then for $k\in \Z$ we let $p^{(k)}$ denote the $k$'th power of $p$ 
with respect to the product for homogeneous symbols, e.g., $p^{(-1)}$ is the inverse of $p$. 
% 
% for any $m \in \C$ we endow $S_{m}(\fg^{*}M,\cE)$  This 
% allows us to consider smooth family with values in $S_{m}(\fg^{*}M,\cE)$ and parametrized by some base manifold $A$. Moreover, it follow the product $*$ for 
% homogeneous symbols is continuous with respect to this topology,
% 
% 
% Then we have: 
% 
\begin{proposition}\label{prop:Chap3.Rockland.family}
    Let $(p_{\nu})_{\nu \in B}$ be a smooth family of invertible symbols in $S_{m}(\fg^{*}M,\cE)$. Then $(p_{\nu}^{(-1)})_{\nu \in B}$ is a 
    smooth family with 
values in $S_{-m}(\fg^{*}M,\cE)$. 
%     
%     $m \in \C$. Then:\smallskip
%     
% 1) $S_{m}(\fg^{*}M,\cE)^{\times}$ is an open subset of $S_{m}(\fg^{*}M,\cE)$.\smallskip
% 
% 2) The map $p\rightarrow p^{(-1)}$ is continuous from $S_{m}(\fg^{*}M,\cE)^{\times}$ to $S_{-m}(\fg^{*}M,\cE)^{\times}$. In fact, 
\end{proposition}
\begin{proof}
First, as in the proofs of Theorem~\ref{thm:Chap3.Rockland-Parametrix-order0} and Proposition~\ref{prop:Chap3.Rockland.open} it is enough to prove the result in a 
trivializing Heisenberg chart $U\subset \Rd$ and we may assume that $\cE$ is the trivial line bundle. 
    
Next, we can reduce the proof to the case $m=0$ as follows. Let $X_{0},\ldots,X_{d}$ be a $H$-frame on $U$ and let $p_{2}$ be the principal 
symbol in the Heisenberg sense of the sublaplacian $-(X_{1}^{2}+\ldots+X_{d}^{2})-X_{0}$, so that by~\cite[Thm.~8.4]{BG:CHM}  $p_{2}$ is an invertible 
symbol.  Then we have $p_{\nu}^{(-1)}=p_{2}^{(-1)}*p_{\nu}*(p_{2}^{(-m)}*p_{\nu}^{(2)})^{(-1)}$. Since by Proposition~\ref{prop:PsiHDO.continuity-product-symbol}  
is bilinear continuous, hence smooth, it follows that in order to prove that $p_{\nu}^{(-1)}$ depends smoothly on $\nu$ it is enough to do it for 
$(p_{2}^{(-m)}*p_{\nu}^{(2)})^{(-1)}$. As the latter is the inverse of a symbol in $S_{0}(\URd)$ we see that it is sufficient to prove the 
proposition in the case $m=0$. 

Suppose now that $p_{\nu}$ belongs to $S_{0}(\URd)$. We shall use here the same notation as that of the proof of 
Theorem~\ref{thm:Chap3.Rockland-Parametrix-order0}   and regard 
$(p_{\nu}(a,.))_{(\nu,a)\in B\times U}$  as a smooth family with values in 
$S_{0}(\Rd)$ parametrized by $B\times U$. For $(\nu,a)$ in $B\times U$ let $K^{a}_{\nu}=(p_{\nu})^{\vee}_{\xiy}(a,.)$. 
Since by~\cite[Prop.~15.30]{BG:CHM} the inverse Fourier transform is a topological isomorphism from $S_{0}(\Rd)$ onto $\cK_{-(d+2)}(\Rd)$ we see that 
$(K^{a}_{\nu})_{(\nu,a)\in B\times U}$ is a smooth family with values in $\cK_{-(d+2)}(\URd)$. 

Since $p_{\nu}$ is an invertible symbol for 
every $\nu$ we see that for every $\nu\in B$ and every $a\in U$ the model operator $P_{K_{\nu}^{a}}^{a}$ is invertible on $L^{2}(\Rd)$. Therefore, a 
simple modification of the arguments at the end of the proof of the Theorem~\ref{thm:Chap3.Rockland-Parametrix-order0} 
shows that  there exists a smooth family $(J^{a}_{\nu})_{(\nu,a)\in 
B\times U}$ with values in $K_{-(d+2)}(\Rd)$ such that $J^{a}_{\nu}*^{a}K^{a}_{\nu}=K^{a}_{\nu}*^{a}J^{a}_{\nu}=\delta_{0}$. Clearly, we have 
$p_{\nu}^{(-1)}(a,.)=(J^{a}_{\nu})^{\wedge}$, so $(p_{\nu}^{(-1)}(a,.))_{(\nu,a)\in B\times U}$ is a smooth family with values in $S_{0}(\Rd)$, that 
is, $(p_{\nu}^{(-1)})$ is a smooth family with values in $S_{0}(\URd)$. The proof is thus achieved. 
% and for $a\in M$ we let $P_{\alpha}^{a}$ be a the left-invariant \psivdo\ on $G_{a}M$ with symbol 
% $p_{\alpha}(a,.)$.  It follows from Lemma~\ref{} that $(P_{\alpha}^{a})_{a\in M}$ is a continuous section of the bundle $\sqcup_{a \in 
% M}L^{2}(G_{a}M,\cE)$. 
\end{proof}

\begin{remark}
 For $m \in \C$ let $S_{m}(\fg^{*}M,\cE)^{\times}$ denote the set of invertible elements of  $S_{m}(\fg^{*}M,\cE)$. Then 
with some additional work it is possible to show that $S_{m}(\fg^{*}M,\cE)^{\times}$ is an open subset of $S_{m}(\fg^{*}M,\cE)$ and that the 
the map $p\rightarrow p^{(-1)}$ is continuous from $S_{m}(\fg^{*}M,\cE)^{\times}$ to $S_{-m}(\fg^{*}M,\cE)^{\times}$ and, in fact, is infinitely many 
times differentiable. 
\end{remark}

\section{Invertibility criteria for sublaplacians}
\label{sec:sublaplacian}
In this section we focus on sublaplacians, which yield several important examples of operators on Heisenberg manifolds. The scalar case was 
dealt with in~\cite{BG:CHM}, but the results were not extended to sublaplacians acting on sections of vector bundles. 
These extensions are necessary in 
order to deal with sublaplacians acting on forms such as the Kohn Laplacian or the horizontal sublaplacian (see next section). 

In this section, after having explained the scalar case from the point of view of this memoir, we extend the results to the non-scalar case.   
% In particular, this will allow us to complete the treatment of the Kohn Laplacian in~\cite{BG:CHM} (see Remark~\ref{rem:Examples.Boxb} ahead). 
% this section we focus on sublaplacians, which allows us to deal with many examples. First, we reformulate the known results for scalar 
% sublaplacians, then we explain how to extend these results This extension was not carried out 
% in~\cite{BG:CHM}, but this is a necessary step in order 
% 
% examples of operators on Heisenberg manifolds are sublaplacians (see next section). In this section we explain how the results of the previous 
% sections can be used to recover the known hypoellipticity criterions for sublaplacians.  
% In this section we review from the point of view of this paper criterions insuring us the invertibility of the principal symbol for the main geometric 
% operators on Heisenberg manifolds: H\"ormander's sum of squares, Kohn Laplacian, horizontal sublaplacian and its conformal powes, contact Laplacian. 
% All these operators but the last one these operators are sublaplacians or powers of sublaplacians, so we start by dealing with general sublaplacians 
% on Heisenberg manifolds. 
% 
% We will use the following definition for sublaplacians.

Recall that a differential operator $\Delta:C^{\infty}(M,\cE)\rightarrow C^{\infty}(M,\cE)$ is a sublaplacian when, near any point $a\in M$, we can put 
$\Delta$ in the form,
%  where which $\Delta$ takes the form
 \begin{equation}
    \Delta=-(X_{1}^{2}+\ldots+X_{d}^{2})- i\mu(x) X_{0}+ \op{O}_{H}(1),
%     \sum_{j=1}^{d} b{j}(x)X_{j}+ c(x), \quad \mu(x)=(\mu_{1}(x),\ldots,\mu_{r}(x)).
    \label{eq:Heisenberg.sublaplacian.bundle}
\end{equation}
where $X_{0},X_{1},\ldots,X_{d}$ is a local $H$-frame of $TM$, the coefficient  $\mu(x)$  is a local section of $\End \cE$ 
and the notation $\op{O}_{H}(1)$  means a differential operator of Heisenberg order~$\leq 1$.
% $\mu_{1}(x),\ldots,\mu_{d}(x)$ are smooth complex-valued functions and $b_{1}(x),\ldots,b_{d}(x)$ and $c(x)$ are smooth functions with 
% values in $M_{r}(\C)$.
% \end{definition}

Let us look at the Rockland condition for a sublaplacian $\Delta:C^{\infty}(M)\rightarrow C^{\infty}(M)$ acting on functions. Let $a\in M$ and let 
$X_{0},X_{1},\ldots,X_{d}$ be a local $H$-frame of $TM$ so that near $a$ we can write 
\begin{equation}
    \Delta=-\sum_{j=1}^{d} X_{j}^{2} - i\mu(x) X_{0}+ \op{O}_{H}(1),
     \label{eq:Heisenberg.sublaplacian.scalar}
\end{equation}
where $\mu(x)$ is a smooth function near $a$. Using~(\ref{eq:Pincipal.example}) we see that the principal symbol of $\Delta$ is 
\begin{equation}
    \sigma_{2}(\Delta)(x,\xi)=|\xi'|^{2}+\mu(a)\xi_{0}, \qquad \xi'=(\xi_{1},\ldots,\xi_{d}).
%     \label{eq:¥}
\end{equation}
In particular we have $ \sigma_{2}(\Delta)(x,0,\xi')=|\xi'|^{2}>0$ for $\xi'\neq 0$, which shows that the condition (i) of 
Proposition~\ref{prop:Rockland.reduction}  is 
always satisfied.

Let $L(x)=(L_{jk}(x))$ be the matrix of the Levi form $\cL$ with respect to the $H$-frame $X_{0},\ldots,X_{d}$, so that for $j,k=1,\ldots,d$ we have
\begin{equation}
    \cL(X_{j},X_{k})=[X_{j},X_{k}]=L_{jk}X_{0} \quad \bmod H.
%     \nonum
\end{equation}
Equivalently, if we let $g(x)$ be the metric on $H$ making orthonormal the frame $X_{1},\ldots,X_{d}$, then for any sections $X$ and $Y$ of $H$ we have 
\begin{equation}
    \cL(X,Y)=g(x)(L(x)X,Y)X_{0} \quad \bmod H.
%     \label{eq:¥}
\end{equation}

The matrix $L(x)$ is antisymmetric, so up to an orthogonal change of frame of $H$, which does not affect the form~(\ref{eq:Heisenberg.sublaplacian.scalar}), 
we may assume that $L(a)$ is in the normal form, 
\begin{equation}
   \qquad L(a)=\left( 
   \begin{array}{ccc}
       0 & D & 0  \\
       -D & 0 & 0 \\
       0 & 0 & 0
   \end{array}\right), \qquad D=\op{diag}(\lambda_{1},\ldots,\lambda_{n}), \quad \lambda_{j}>0,
     \label{eq:Sublaplacian.Levi-form.normal-form}
\end{equation}
so that $\pm i\lambda_{1},\ldots,\pm i\lambda_{2n}, 0,\ldots,0$ are the eigenvalues of $L(a)$ counted with multiplicity.  Then 
the model vector fields $X_{0}^{a},\ldots,X_{d}^{a}$ are: 
\begin{gather}
    X_{0}^{a}=\frac{\partial}{\partial x_{0}}, \qquad X_{k}^{a}=\frac{\partial}{\partial x_{k}}, \qquad k=2n+1,..,d,\\
    X_{j}^{a}=\frac{\partial}{\partial x_{j}}-\frac{1}{2}\lambda_{j}x_{n+j}\frac{\partial}{\partial x_{0}}, \quad 
    X_{n+j}^{a}=\frac{\partial}{\partial x_{j}}+\frac{1}{2}\lambda_{j}x_{j}\frac{\partial}{\partial x_{0}}, \ j=1,\ldots,n. 
\end{gather}
In terms of these vector fields the model operator of $\Delta$ at $a$ is 
\begin{equation}
    \Delta^{a}=-[(X_{1}^{a})^{2}+\ldots+(X_{1}^{a})^{2}]-i\mu(a)X_{0}^{a}.
    \label{eq:Sublaplacian.model-operator}
\end{equation}

Next, under the isomorphism $\phi:\bH^{2n+1}\times \R^{d-2n}\rightarrow G_{a}M$ given by 
\begin{equation}
    \phi(x_{0},\ldots,x_{d})=(x_{0},\lambda_{1}^{\frac{1}{2}}x_{1},\ldots,\lambda_{n}^{\frac{1}{2}}x_{n},\lambda_{1}^{\frac{1}{2}}x_{n+1}, \ldots
    \lambda_{n}^{\frac{1}{2}}x_{2n},x_{2n+1},\ldots,x_{d}), 
%     \label{eq:¥}
\end{equation}
the representations $\pi^{\pm,\xi}=\pi^{\pm 1,\xi}$, $\xi \in\{0\}^{2n}\times \R^{d-2n}$, become the representations of $G_{a}M$ such that
\begin{gather}
    d\pi^{\pm,\xi}(X_{0})=\pm i, \qquad d\pi^{\pm,\xi}(X_{k})=\pm i \xi_{k}, \quad k=2n+1,\ldots,d,\\
    d\pi^{\pm,\xi}(X_{j})=\lambda_{j}^{\frac{1}{2}}\frac{\partial}{\partial \xi_{j}}, 
    \qquad d\pi^{\pm,\xi}(X_{n+j})=\pm i\lambda_{j}^{\frac{1}{2}} \xi_{j}, \quad j=1,\ldots,n,\\
       \pi_{\Delta^{a}}^{\pm,\xi}= d\pi^{\pm,\xi}(\Delta^{a})= 
    \sum_{j=1}^{n}\lambda_{j}(-\partial_{\xi_{j}}^{2}+\xi_{j}^{2}) \pm(\xi_{2n+1}^{2}+\ldots+\xi_{d}^{2}+\mu(a)).
\end{gather}

The spectrum of the harmonic oscillator $\sum_{j=1}^{n}\lambda_{j}(-\partial_{\xi_{j}}^{2}+\xi_{j}^{2})$ is $\sum_{j=1}^{n}\lambda_{j}(1+2\N)$ and 
 all its eigenvectors belong to $\cS(\R^{n})$. Thus, the operator $ \pi_{\Delta^{a}}^{\pm,\xi}$ is injective on $\cS(\R^{n})$ if, and only if, 
  $\xi_{2n+1}^{2}+\ldots+\xi_{d}^{2}+\mu(a)$ is not $\pm \sum_{j=1}^{n}\lambda_{j}(1+2\N)$. This occurs for any $\xi \in\{0\}^{2n}\times 
  \R^{d-2n}$ if, and only if, the following condition holds
%   \begin{equation}
%  
%        \label{eq:Sublaplacian.condition.scalar}
%   \end{equation}
%  where the set $\Lambda^{a}$ is defined by 
 \begin{gather}
         \text{$\mu(a)$ is not in the singular set $\Lambda_{a}$}, 
         \label{eq:Sublaplacian.condition.scalar} \\
         \Lambda_{a}=(-\infty, -\frac12 \Tra |L(a)|]\cup [\frac12 \Tra 
  |L(a)|,\infty) \qquad \text{if $2n<d$},\label{eq:Sublaplacian.singular-set1}\\
   \Lambda_{a}=\{\pm(\frac12 \Tra |L(a)|+2\sum_{1\leq j \leq n}\alpha_{j}|\lambda_{j}|); \alpha_{j}\in \N^{d}\}\qquad \text{if $2n=d$}.
   \label{eq:Sublaplacian.singular-set2}
 \end{gather}
In particular, the condition~(ii) of Proposition~\ref{prop:Rockland.reduction} 
is equivalent to~(\ref{eq:Sublaplacian.condition.scalar}). Since the condition (i) is always satisfied, it follows that the 
Rockland condition for $\Delta$ is equivalent to~(\ref{eq:Sublaplacian.condition.scalar}). 

Notice also that, independently of the equivalence with the Rockland condition, 
the condition~(\ref{eq:Sublaplacian.condition.scalar}) does not depend on the choice of the $H$-frame, because 
as $\Lambda_{a}$ depends only on the eigenvalues of $L(a)$ which scale in the same way as $\mu(a)$ 
under a change of $H$-frame preserving the form~(\ref{eq:Heisenberg.sublaplacian.scalar}). 

On the other hand, since the transpose $(\Delta^{a})^{t}=(\Delta^{t})^{a}$ is given by the formula~(\ref{eq:Sublaplacian.model-operator}) 
 with $\mu(a)$ replaced by $-\mu(a)$, which has no 
effect on~(\ref{eq:Sublaplacian.condition.scalar}), we see that the Rockland condition for $(\Delta^{t})^{a}$ too is equivalent 
 to~(\ref{eq:Sublaplacian.condition.scalar}). Therefore, we have obtained: 

\begin{proposition}\label{prop:Sublaplacian.invertibility.scalar}
The Rockland conditions for $\Delta^{t}$ and $\Delta$ at $a$ are both equivalent to~(\ref{eq:Sublaplacian.condition.scalar}).
\end{proposition}

In particular, we see that if the principal symbol of $\Delta$ is invertible then the condition~(\ref{eq:Sublaplacian.condition.scalar}) holds at every 
point. As shown by Beals-Greiner the converse is true as well. The key result is the following.

\begin{proposition}[{\cite[Sect.~5]{BG:CHM}}]\label{prop:Sublaplacian.inverse.scalar}
    Let $U\subset \Rd$ be a Heisenberg chart near $a$ and define
    \begin{equation}
        \Omega=\{(\mu,x)\in \C \times U; \ \mu \not \in \Lambda_{x}\}.
%         \label{eq:}
    \end{equation}
     Then $\Omega$ is an 
    open set and there exists $q_{\mu}(x,\xi)\in C^{\infty}(\Omega, S_{-2}(\Rd))$ such that:\smallskip
    
   (i) $q_{\mu}(x,\xi)$ is analytic with respect to $\mu$;\smallskip
%     
%    (ii) We have $q_{\mu}(x,\lambda.\xi)=\lambda^{-2}q_{\mu}(x,\xi)$ for any $\lambda>0$.\smallskip
    
   (ii) For any $(\mu,x)\in \Omega$ the symbol $q_{\mu}(x,.)$ inverts $|\xi'|^{2}+i\mu\xi_{0}$ on $G_{x}U$, that is,
    \begin{equation}
        q_{\mu}(x,.)*^{x}(|\xi'|^{2}+i\mu \xi_{0})=(|\xi'|^{2}+i\mu \xi_{0})*^{x}q_{\mu}(x,.)=1.
%         \label{eq:¥}
    \end{equation}
More precisely, $q_{\mu}(x,\xi)$ is obtained from the analytic continuation of the function,
\begin{gather*}
    q_{\mu}(x,\xi) =\int_{0}^{\infty}e^{-t\mu \xi_{0}}G(x,\xi,t)dt, \qquad |\Re \mu|<\frac{1}{2}\Tr |L(x)|,\\
    G(x,\xi,t)=\det{}^{-\frac{1}{2}}[\cosh (t|\xi_{0}||L(x)|)] \exp[-t\acou{\frac{\tanh (t|\xi_{0}||L(x)|)}{t|\xi_{0}||L(x)|}\xi'}{\xi'}].
     \label{eq:Rockland-Sublaplacian.heat-kernel-inverse}
\end{gather*}
%   where we have let $\xi'=(\xi_{1},\ldots,\xi_{d})$.   
\end{proposition}

 This implies that if the condition~(\ref{eq:Sublaplacian.condition.scalar}) is satisfied at every point $x \in U$ then we get an inverse $q_{-2}\in 
S_{-2}(\URd)$ for $\sigma_{2}(\Delta)(x,\xi)=|\xi'|^{2}+i\mu(x)\xi_{0}$ on $\URd$ by letting 
$q_{-2}(x,\xi)=q_{\mu(x)}(x,\xi)$  for $(x,\xi)\in \URdo$.
% %     \label{eq:¥}
% \end{equation}
It thus follows that if~(\ref{eq:Sublaplacian.condition.scalar}) 
holds at every point of $M$ then the principal symbol of $\Delta$ is invertible near any point of $M$, hence admits 
an inverse in $S_{-2}(\fg^{*}M)$. Therefore, we get: 

\begin{proposition}\label{prop:Sublaplacian.Rockland.scalar}
A sublaplacian $\Delta:C^{\infty}(M)\rightarrow C^{\infty}(M)$ has an invertible principal symbol if, and only if, it satisfies the 
condition~(\ref{eq:Sublaplacian.condition.scalar}) at every point.  
\end{proposition}

Let us now deal with a sublaplacian $\Delta:C^{\infty}(M,\cE)\rightarrow C^{\infty}(M,\cE)$ acting on the sections of 
the vector bundle $\cE$.

Let $a \in M$ and let $X_{0},\ldots,X_{d}$ be a local $H$-frame near $a$ such that
\begin{equation}
    \Delta=-\sum_{j=1}^{d} X_{j}^{2} - i\mu(x) X_{0}+ \op{O}_{H}(1),
     \label{eq:Heisenberg.sublaplacian.bundle2}
\end{equation}
where $\mu(x)$ is a smooth local section of $\End \cE$. 

In a suitable basis of $\cE_{a}$ the matrix  of $\mu(a)$ is in  triangular form,
\begin{equation}
    \mu(a)=\left( 
    \begin{array}{ccc}
        \mu_{1}(a) & * & *  \\
        0 & \ddots & *  \\
        0 & 0 & \mu_{r}(a)
    \end{array}
\right).    
    %     \label{eq:}
\end{equation}
where $\mu_{1}(a),\ldots,\mu_{r}(a)$ denote the eigenvalues of $\mu(a)$ counted with multiplicity. 
Therefore, the model operator of $\Delta$ at $a$ is of the form, 
\begin{equation}
    \Delta^{a}= \left( 
    \begin{array}{ccc}
        \Delta_{1}^{a} & * & *  \\
        0 & \ddots & *  \\
        0 & 0 &\Delta_{r}^{a}
    \end{array}
\right), \qquad \Delta_{j}^{a}=-[(X_{1}^{a})^{2}+\ldots+(X_{1}^{a})^{2}]-i\mu_{j}(a)X_{0}^{a}.  
%     \label{eq:¥}
\end{equation}
It follows that $\Delta^{a}$ satisfies the Rockland condition if, and only if, so does each sublaplacian $\Delta_{j}^{a}$, $j=1,\ldots,r$. Using 
Proposition~\ref{prop:Sublaplacian.Rockland.scalar}
we then deduce that the Rockland condition $\Delta^{a}$ is equivalent to the condition, 
\begin{equation}
    \op{Sp}\mu(a) \cap \Lambda_{a}=\emptyset.
     \label{eq:Sublaplacian.condition}
\end{equation}

Notice that the same is true for the transpose $(\Delta^{a})^{t}$. Moreover, the condition~(\ref{eq:Sublaplacian.condition}) 
is independent of the choice of the basis of $\cE_{a}$ or of the $H$-frame since the condition involves $\mu(a)$ only though its eigenvalues of $\mu(a)$ 
and the latter scale in the same way as that of $L(a)$ under a change of $H$-frame preserving the form~(\ref{eq:Heisenberg.sublaplacian.bundle2}). 

Next, concerning the invertibility of the principal symbol of $\Delta$ the following extension of 
Proposition~\ref{prop:Sublaplacian.inverse.scalar} holds. 
% Moreover, as the condition involves $\mu(a)$ only via its spectrum it does dependthe choice of the 
% trivialization near $a$. Moreover, as the condition~(\ref{eq:Sublaplacian.condition.scalar}) this condition does not depend on the choice of the Heisenberg coordinates 
% either, hence makes sense intrinsically. 

\begin{proposition}\label{prop:Sublaplacian.inverse.system}
    Let $U\subset \Rd$ be a trivializing Heisenberg chart near $a$ and define
    \begin{equation}
        \Omega=\{(\mu,x)\in M_{r}(\C) \times U; \ \mu \not \in \Lambda_{x}\}. 
%     \label{eq:}
\end{equation}
Then $\Omega$ is an open set and there exists $q_{\mu}(x,\xi)\in C^{\infty}(\Omega,  S_{-2}(\Rd, \C^{r}))$ so that:\smallskip
    
    (i) $q_{\mu}(x,\xi)$ is analytic with respect to $\mu$;\smallskip
%     
%     (ii) We have $q_{\mu}(x,\lambda.\xi)=\lambda^{-2}q_{\mu}(x,\xi)$ for any $\lambda>0$.\smallskip
    
    (ii) For any $(\mu,x)\in \Omega$ the symbol $q_{\mu}(x,.)$ inverts $|\xi'|^{2}+i\mu\xi_{0}$ on $G_{x}U$, that is,
    \begin{equation}
        q_{\mu}(x,.)*^{x}(|\xi'|^{2}+i\mu \xi_{0})=(|\xi'|^{2}+i\mu \xi_{0})*^{x}q_{\mu}(x,.)=1.
%         \label{eq:¥}
    \end{equation}
\end{proposition}
\begin{proof}
  It is enough to prove that near point $(\mu_{0},x_{0})\in \Omega$ there exists an open neighborhood $\Omega'$ contained in $\Omega$ and 
  a function $q_{\mu}(x,\xi)\in C^{\infty}(\Omega', S_{-2}(\Rd, \C^{r}))$ satisfying the properties (i) and (ii) on $\Omega'$. 
  
   To this end observe that since $\op{Sp}\mu \subset \overline{D}(0,\|\mu\|)$ for any $\mu \in M_{r}(\C)$, we see that if we let  
 $K=\overline{B(0,\|\mu_{0}\|+1)}$ then any  $\mu \in M_{r}(\C)$ close enough to $\mu_{0}$ has its spectrum contained in $K$. 
 In addition, let $\delta\in (0, \frac{1}{2}\op{dist}(\op{Sp}\mu_{0},\Lambda_{x_{0}}))$ and set $V_{1}=\op{Sp} \mu_{0}+D(0,\delta)$ and 
 $V_{2}=\Lambda_{x_{0}}+D(0,\delta)$, so that $V_{1}$ and $V_{2}$ are disjoint open subsets of $\C$ containing $\Sp \mu_{0}$ and $\Lambda_{x_{0}}$ 
 respectively. 
 
 Notice that for any $\mu$ close enough to $\mu_{0}$ we have $\op{Sp} \mu \subset V_{1}$. Otherwise there 
 exists a sequence $(\mu_{k})_{k \geq 1} \subset M_{r}(\C)$ converging to $\mu_{0}$ and a sequence of eigenvalues 
 $(\lambda_{k})_{k \geq 1} \subset K$, $\lambda_{k }\in \Sp \mu_{k}$, such that $\lambda_{k} \not \in V_{1}$ for any $k\geq 1$. 
 Since the sequence $(\lambda_{k})_{k 
 \geq 1}$ is bounded,  we may assume that it converges to some $\lambda \not \in V_{1}$. Necessarily $\lambda$ is an eigenvalue of $\mu_{0}$, which 
 contradicts  the fact that $\lambda \not \in V_{1}$. Thus there exists $\eta_{1}>0$ so that for any $\mu \in B(\mu_{0},\eta_{1})$ 
 we have $\op{Sp} \mu \subset V_{1}$.
 
 Similarly, there exists $\eta_{2}>0$ so that for any $x \in B(x_{0},\eta_{2})$ we have $\Sp |L(x)|$ is contained in $\Sp |L(x_{0})|+D(0,\delta)$, which implies 
 $\Lambda_{x}\subset \Lambda_{x_{0}}+D(0,\delta)=V_{2}$. Therefore, the open set $\Omega'=B(\mu_{0},\eta_{1})\times B(x_{0},\eta_{2})$ is such that 
 for any $(\mu,x)\in \Omega'$ we have $\Sp \mu \cap \Lambda_{x}\subset V_{1}\cap V_{2}=\emptyset$, i.e., $\Omega'$ is an open neighborhood of 
 $(\mu_{0},x_{0})$ in $\Omega$.%Hence $\Omega$ is open.
 
Next, let $\Gamma$ be a smooth curve of index 1 such that the bounded connected component of $\C \setminus \Gamma$ contains $V_{1}$ and its unbounded 
component contains $V_{2}$. Then we define an element of $\Hol(B(\mu_{0},\eta_{1}))\hotimes C^{\infty}(B(x_{0},\eta_{2})\times \Rdo)$ by letting 
\begin{equation}
    q_{\mu}(x,\xi) =\frac{1}{2i\pi}\int_{\Gamma}q_{\gamma}(x,\xi)(\gamma-\mu)d\gamma, \qquad (\mu,x,\xi)\in \Omega'\times \Rdo.
%     B(\mu_{0},\eta_{1})\times     B(x_{0},\eta_{2})
%     \label{eq:}
\end{equation}
This function is homogeneous of degree $-2$ with respect to $\xi$ and for any $(\mu,x)\in \Omega'$ we have
\begin{multline}
    q_{\mu}(x,.)*^{x}(|\xi'|^{2}+i\mu\xi_{0})=\frac{1}{2i\pi}\int_{\Gamma}q_{\gamma}(x,.)*^{x}(|\xi'|^{2}+i\mu\xi_{0})(\gamma-\mu)^{-1}d\gamma\\
   = \frac{1}{2i\pi}\int_{\Gamma}[(\gamma-\mu)^{-1}-iq_{\gamma}(x,.)*\xi_{0}]d\gamma=1.
\end{multline}
Similarly, we have $(|\xi'|^{2}+i\mu\xi_{0})*^{x}q_{\mu}(x,.)=1$. Thus $q_{\mu}(x,\xi)$ is an element of $C^{\infty}(\Omega', S_{-2}(\Rd, \C^{r}))$ 
satisfying the properties (i) and (ii) on $\Omega'$. 
The proof is therefore complete. 
% Since the same is true for $\Omega'$ sufficiently small open neighborhood of any other element of $\Omega$ we deduce that there exists a 
% function $q_{\mu}(x,\xi)\in C^{\infty}(\Omega\times \Rdo)$ satisfying the properties (i), (ii) and (iii).
\end{proof}

In the same way as Proposition~\ref{prop:Sublaplacian.inverse.scalar} 
in the scalar case, Proposition~\ref{prop:Sublaplacian.inverse.system} implies that when the condition~(\ref{eq:Sublaplacian.condition}) holds 
everywhere the principal symbol of $\Delta$ admits an inverse in $S_{-2}(\fg^{*}M,\cE)$. We have thus proved:
% from Proposition~\ref{prop:Sublaplacian.inverse.system} that if the condition~(\ref{eq:Sublaplacian.condition}) is satisfied for every $a \in U$ 
% then the principal symbol $\sigma_{2}(x,\xi)=|\xi'|^{2}+i\mu(x)\xi_{0}$ 
% of $\Delta$ is invertible with inverse $q_{-2}\in S_{-2}(\URd, M_{r}(\C))$ given by $q_{-2}(x,\xi)=q_{\mu(x)}(x,\xi)$. 
% Summarizing all this we have proved: 

\begin{proposition}\label{prop:Sublaplacian.Rockland-bundle}
  1) At every point $a\in M$ the Rockland conditions for $\Delta$ and $\Delta^{t}$ are equivalent to~(\ref{eq:Sublaplacian.condition}).\smallskip
    
    2) The principal symbol of $\Delta$ is invertible if, and only if, the condition~(\ref{eq:Sublaplacian.condition}) holds everywhere. 
    Moreover, when the latter occurs $\Delta$ admits a parametrix in $\pvdo^{-2}(M,\cE)$ and is hypoelliptic with gain of 1 derivative.
\end{proposition}

\section{Invertibility criteria for the main differential operators}\label{sec:Examples}
In this section we explain how the previous results of this monograph can be used to deal with the hypoellipticity for the main geometric operators on 
Heisenberg manifolds: H\"ormander's sum of squares, Kohn Laplacian, horizontal sublaplacian and  contact Laplacian. In particular, we complete the treatment 
in~\cite{BG:CHM} of the Kohn Laplacian and we establish a criterion for the invertibility of the horizontal sublaplacian.

% All these operators but the last one these operators are sublaplacians or powers of sublaplacians, so we start by dealing with general sublaplacians 
% on Heisenberg manifolds. 

\subsection{H\"ormander's sum of squares.}
 Let $(M^{d+1},H)$ be a Heisenberg manifold and let $\Delta:C^{\infty}(M,\cE)\rightarrow 
C^{\infty}(M,\cE)$ be a generalized sum of squares of the form~(\ref{eq:Operators.generalized-sum-of-squares}), that is,
\begin{equation}
    \Delta=-(\nabla_{X_{1}}^{2}+\ldots+\nabla_{X_{m}}^{2}) +\op{O}_{H}(1),
%      \label{eq:Operators.generalized-sum-of-squares2}
\end{equation}
where $\nabla$ is a connection on $\cE$ and the vector $X_{1},\ldots,X_{m}$ span $H$. Then, in the local form~(\ref{eq:Heisenberg.sublaplacian.bundle}) 
for $\Delta$ the matrix $\mu(x)$ vanishes, so~(\ref{eq:Sublaplacian.condition}) holds at a point $a\in M$ if, and only if, the Levi form does not vanish at~$a$. 
Combining this with Proposition~\ref{prop:Sublaplacian.Rockland-bundle} then gives: 

\begin{proposition}
 1) At a point $x\in M$ the operators $\Delta$ and $\Delta^{t}$ satisfies the Rockland condition if, and only if, the Levi form $\cL$ does not vanish 
 at $x$.\smallskip
 
 2) The principal symbol of $\Delta$ is invertible if, and only if, the Levi form is non-vanishing. In particular, when the latter occurs $\Delta$ 
 admits a parametrix in $\pvdo^{-2}(M,\cE)$ and is hypoelliptic with gain of one derivative.
\end{proposition}

In particular, since the nonvanishing of the Levi form is equivalent to the  
bracket condition $H+[H,H]=TM$, we see that, the special case of Heisenberg manifolds, we recover the hypoellipticity result of~\cite{Ho:HSODE} for 
sums of squares.

\subsection{Kohn Laplacian}
Let $M^{2n+1}$ be an orientable CR manifold with CR tangent bundle $T_{1,0}\subset T_{\C}M$, so that $H=\Re(T_{1,0}\oplus T_{0,1})$ yields a Heisenberg 
structure on $M$, and let $\cN$ be a line bundle supplement of $H$ in $TM$. Assuming that $T_{\C}M$ endowed with a Hermitian metric commuting with 
complex conjugation and making orthogonal the splitting $T_{\C}M=T_{1,0}\oplus T_{0,1}\oplus (\cN \otimes \C)$, we let 
$\Boxbpq:C^{\infty}(M,\Lambda^{p,q})\rightarrow C^{\infty}(M,\Lambda^{p,q})$ be the Kohn Laplacian acting on $(p,q)$-forms.

Let $\theta$ be a global nonvanishing real $1$-form anihilating $H$ with Levi form $L_{\theta}$. 
Recall that the condition $Y(q)$ at a point $x$ requires to have 
\begin{equation}
    q\not \in \{\kappa(x),\ldots,\kappa(x)+n-r(x)\}\cup \{r(x)-\kappa(x),\ldots,n-\kappa(x)\},
%      \label{eq:Operators.Y(q)-condition}
\end{equation}
where $\kappa(x)$ denotes the number of negative eigenvalues of $L_{\theta}$ at $x$ and $r(x)$ its rank. 

It is shown in~\cite[Sect.~21]{BG:CHM} that 
at every point $a\in M$ the condition $Y(q)$ is equivalent to the 
condition~(\ref{eq:Sublaplacian.condition}). Therefore, from Proposition~\ref{prop:Sublaplacian.Rockland-bundle} we get:  

\begin{proposition}\label{prop:Examples.Boxb}
%     Let $\Box_{b}:C^{\infty}(M,\Lambda^{p,q})\rightarrow C^{\infty}(M,\Lambda^{p,q})$ be the Kohn Laplacian acting on $(p,q)$-forms.\smallskip 
%    
   1) At a point $x\in M$ the Rockland condition for $\Boxbpq$ is equivalent to the condition~$Y(q)$.\smallskip
   
   2) The principal symbol of $\Boxbpq$ is invertible if, and only if, the condition $Y(q)$ is satisfied at every point. In particular, when the latter 
   occurs $\Boxbpq$ admits a parametrix in $\pvdo^{-2}(M,\Lambda^{p,q})$ and is hypoelliptic with gain of one derivative.
\end{proposition}

The proof of the second part above is not quite complete in~\cite[Sect.~21]{BG:CHM}. In fact, Beals-Greiner claimed that diagonalizing the 
    leading part of the Kohn Laplacian allows us make use of the criterion from Proposition~\ref{prop:Sublaplacian.Rockland.scalar} for \emph{scalar} sublaplacians. 
    This fact is definitely true in case of a 
    Levi Metric (see~\cite{FS:EDdbarbCAHG}) or even a smoothly diagonalizable Levi form, but it fails in general. Indeed,  for 
    the Kohn Laplacian the eigenvalues of the 
    matrix $\mu(x)$ in~(\ref{eq:Heisenberg.sublaplacian.bundle}) with respect to an orthonormal $H$-frame of $TM$ are given in terms of eigenvalues of the Levi 
    form (see~Eq.~(21.31) in~\cite{BG:CHM}), but the latter need not depend smoothly on $x$. 
%     (unless it can be diagonalized  %(unless the metric on $T_{\C}M$ is a Levi metric). 
%    diagonalization  of 
%     the leading part of $\Box_{b}$ it is possible to apply Proposition~\ref{prop:Sublaplacian.invertibility.scalar}, 
%     which deals with the scalar case only. However, this not quite true because in an eigenframe 
%     the coefficient $\mu(x)$ in the form~(\ref{eq:Heisenberg.sublaplacian.bundle}) 
%     of $\Box_{b}$ depend on the eigenvalues of the Levi form, which need not be smooth functions
  
    This shows that in order to deal  with the Kohn Laplacian acting on forms, and more generally with sublaplacians acting on sections of a vector 
    bundle, we really need to use Proposition~\ref{prop:Sublaplacian.Rockland-bundle}, as we can cannot in general reduce the study to the scalar 
    case. 
% \end{remark}

% \begin{remark}

\subsection{The horizontal sublaplacian}
Let $(M^{d+1},H)$ be a Heisenberg manifold endowed with a Riemannian metric, 
let $\Lambda^{*}_{\C}H^{*}=\oplus_{k=0}^{d}\Lambda^{k}_{\C}H^{*}$ 
be the (complexified) bundle of horizontal forms and let 
$\Delta_{b;k}:C^{\infty}(M,\Lambda^{k}_{\C}H^{*})\rightarrow C^{\infty}(M,\Lambda^{k}_{\C}H^{*})$ be the associated horizontal sublaplacian on 
horizontal forms of degree $k$ as defined in~(\ref{eq:Operators.sublaplacian}). 

We shall now express the condition~(\ref{eq:Sublaplacian.condition}) in terms of the more geometric condition $X(k)$ below.

\begin{definition}
    For $x\in M$ let $2r(x)$ be the rank of the Levi form $\cL$ at $x$. Then we say that $\cL$ satisfies the condition 
$X(k)$ at $x$ when we have
\begin{equation}
    k\not \in\{r(x),r(x)+1,\ldots,d-r(x)\}.
    \label{eq:Operators.X(k)-condition}
\end{equation}
\end{definition}

For instance, the condition $X(0)$ is satisfied if, and only if, the Levi form does not vanish. Also, 
if $M^{2n+1}$ is a contact manifold or a nondegenerate CR manifold then the Levi form is everywhere nondegenerate, so $r(x)=2n$ and 
the $X(k)$-condition becomes $k\neq n$. In any case, we have:

\begin{proposition}\label{prop:Examples.horizontal-sublaplacian}
%     Let $\Delta_{b}:C^{\infty}(M,\Lambda^{k}_{\C}H^{*})\rightarrow C^{\infty}(M,\Lambda^{k}_{\C}H^{*})$ be the horizontal  sublaplacian acting on 
%     horizontal forms of degree $k$. 
%    
   1) At a point $x\in M$ the Rockland condition for $\Delta_{b;k}$ is equivalent to the condition~$X(k)$.\smallskip
   
   2) The principal symbol of $\Delta_{b;k}$ is invertible if, and only if, the condition $X(k)$ is satisfied at every point. In particular, when the latter 
   occurs $\Delta_{b;k}$ admits a parametrix in $\pvdo^{-2}(M,\Lambda^{k}_{\C}H^{*})$ and is hypoelliptic with gain of one derivative.
\end{proposition}
\begin{proof}
    First, thanks to Proposition~\ref{prop:Sublaplacian.Rockland-bundle}
    we only have to check that for $k=0,\ldots,d$ at any point $a$  the condition~(\ref{eq:Sublaplacian.condition}) for 
    $\Delta_{b;k}$ is equivalent to the condition $X(k)$. 
    
    Next, let $U\subset \Rd$ be a Heisenberg chart around $a$ together with an orthonormal $H$-frame $X_{0},X_{1},\ldots X_{d}$ of $TU$. 
    Let $g$ be the Riemannian metric of $M$. Then on $U$ we can write the Levi form $\cL$ in the form, 
    \begin{equation}
        \cL(X,Y)=[X,Y]=\acou{L(x)X}{Y}X_{0} \quad \bmod H, 
%         \label{eq:}
    \end{equation}
for some antisymmetric section $L(x)$ of $\End_{\R}H$. In particular, if for $j,k=1,\ldots,d$ we let $L_{jk}=\acou{LX_{j}}{X_{k}}$ then 
we have 
\begin{equation}
    [X_{j},X_{k}]=L_{jk}X_{0} \qquad \bmod H. 
     \label{eq:Examples.Levi-form.coefficients}
\end{equation}

Let $2n$  be the rank of $L(a)$. Since the condition~(\ref{eq:Sublaplacian.condition}) for $\Delta_{b;k}$ at $a$ is independent of the choice of the Heisenberg chart, 
we may assume that $U$ is chosen in such way that at $x=a$ we have $g(a)=1$ and $L(a)$ is in the normal form, 
\begin{equation}
   \qquad L(a)=\left( 
   \begin{array}{ccc}
       0 & D & 0  \\
       -D & 0 & 0  \\
       0 & 0 & 0
   \end{array}\right), \qquad D=\op{diag}(\lambda_{1},\ldots,\lambda_{n}), \quad \lambda_{j}>0,
     \label{eq:Examples.Levi-form.normal-form}
\end{equation}
so that $\pm i\lambda_{1},\ldots,\pm i\lambda_{n}$ are the nonzero eigenvalues of $L(a)$ counted with multiplicity.

Let $\omega^{1},\ldots,\omega^{n}$ be the coframe of $H^{*}$ dual to $X_{1},\ldots,X_{d}$. For a 1-form $\omega$ we let  $\varepsilon(\omega)$  denote the 
exterior product  
and $\iota(\omega)$ denote the interior product with $\omega$, that is, the contraction with the vector fields dual to $\omega$. 
For an ordered subset $J=\{j_{1},\ldots,j_{k}\}\subset 
\{1,\ldots,d\}$, so that $j_{1}<\ldots<j_{d}$, we let $\omega^{J}=\omega^{j_{1}}\wedge \ldots \wedge \omega^{j_{k}}$ (we make the convention 
that $\omega^{\emptyset}=1$). Then the forms $\omega^{J}$'s give rise to an orthonormal frame of $\Lambda^{*}_{\C}H^{*}$ over $U$. With respect to this frame 
we have
% the operators $d_{b}$ and $d_{b}^{*}$ are of the form, 
\begin{equation}
    d_{b}=\sum_{j=1}^{d}\varepsilon(\omega^{j})X_{j} \qquad \text{and} \qquad d_{b}=-\sum_{l=1}^{d}\iota(\omega^{l})X_{l}+\op{O}_{H}(1).
%     \label{eq:}
\end{equation}
Therefore, modulo first order terms we have 
\begin{multline}
 \Delta_{b}=d_{b}^{*}d_{b}+d_{b}d_{b}^{*}=-\sum_{j,l=1}^{d}[\varepsilon(\omega^{j})\iota(\omega^{l})X_{j}X_{l}+ 
    \iota(\omega^{l})\varepsilon(\omega^{j})X_{l}X_{j}] = \\
  -\frac{1}{2}\sum_{j,l=1}^{d}[ (\varepsilon(\omega^{j})\iota(\omega^{l})+\iota(\omega^{l})\varepsilon(\omega^{j}))(X_{j}X_{l}+X_{l}X_{j})
   + (\varepsilon(\omega^{j})\iota(\omega^{l})-\iota(\omega^{l})\varepsilon(\omega^{j}))[X_{j},X_{l}]].
    \end{multline}
Combining this with~(\ref{eq:Examples.Levi-form.coefficients}) and the relations,
 \begin{equation}
   \varepsilon(\omega^{j})\iota(\omega^{l})+\iota(\omega^{l})\varepsilon(\omega^{j}) 
   =  \delta_{jl}, \qquad j,l=1,\ldots,d,
% %     \label{eq:¥}
 \end{equation}
 we then obtain
\begin{equation}
    \Delta_{b}=-\sum_{j=1}^{d}X_{j}^{2}-i\mu(x)X_{0}+\op{O}_{H}(1), \qquad 
    \mu(x)= \frac{1}{i}\sum_{j,l=1}^{d}\varepsilon(\omega^{j})\iota(\omega^{l})L_{jl}.
%     \label{eq:}
\end{equation}
In particular, thanks to~(\ref{eq:Examples.Levi-form.normal-form}) at $x=a$ we have
\begin{equation}
    \mu(a)=\frac{1}{i}\sum_{j=1}^{n}(\varepsilon(\omega^{j})\iota(\omega^{n+j})-\varepsilon(\omega^{n+j})\iota(\omega^{j}))\lambda_{j}.
%     \label{eq:¥}
\end{equation}

For $j=1,\ldots,n$ define $\theta^{j}=\frac{1}{\sqrt{2}}(\omega^{j}+i\omega^{n+j})$ and 
$\theta^{\bar{j}}=\frac{1}{\sqrt{2}}(\omega^{j}-i\omega^{n+j})$. 
Then $\frac{1}{i}(\varepsilon(\omega^{j})\iota(\omega^{n+j})-\varepsilon(\omega^{n+j})\iota(\omega^{j}))$ is equal to
\begin{multline}
% \frac{1}{i}(\varepsilon(\omega^{j})\iota(\omega^{n+j})-\varepsilon(\omega^{n+j})\iota(\omega^{j}))= \\
\frac{-1}{2}[(\varepsilon(\theta^{j})+\varepsilon(\theta^{\bar{j}}))(\iota(\theta^{\bar{j}})-\iota(\theta^{j}))- 
(\varepsilon(\theta^{j})-\varepsilon(\theta^{\bar{j}}))(\iota(\theta^{\bar{j}})+\iota(\theta^{j}))]\\
 =  \varepsilon(\theta^{j})\iota(\theta^{j})-\varepsilon(\theta^{\bar{j}})\iota(\theta^{\bar{j}}).
%  
%  = (\varepsilon(\omega^{j})+i\varepsilon(\omega^{n+j}))(\iota(\omega^{j})-i\iota(\omega^{n+j}))\\
%     = \varepsilon(\omega^{j})\iota(\omega^{j})+\varepsilon(\omega^{n+j})\iota(\omega^{n+j})+
%     i(\varepsilon(\omega^{j})\iota(\omega^{n+j})-\varepsilon(\omega^{n+j})\iota(\omega^{j})).
% %         \label{eq:¥}
\end{multline}
% Similarly, we have
% \begin{equation}
%      2 = \varepsilon(\omega^{j})\iota(\omega^{j})+\varepsilon(\omega^{n+j})\iota(\omega^{n+j})-
%     i(\varepsilon(\omega^{j})\iota(\omega^{n+j})-\varepsilon(\omega^{n+j})\iota(\omega^{j})).
% %     \label{eq:¥}
% \end{equation}
Therefore, we obtain 
% $\varepsilon(\theta^{j})\iota(\theta^{j})-\varepsilon(\theta^{\bar{j}})\iota(\theta^{\bar{j}})=i(\varepsilon(\omega^{j})\iota(\omega^{n+j})- 
% \varepsilon(\omega^{n+j})\iota(\omega^{j}))$, from which we get
\begin{equation}
   \mu(a)=\sum_{j=1}^{n}(\varepsilon(\theta^{j})\iota(\theta^{j})-\varepsilon(\theta^{\bar{j}})\iota(\theta^{\bar{j}}))\lambda_{j}.
    \label{eq:Examples.mu(a)}
\end{equation}

For ordered subsets $J= \{j_{1},\ldots,j_{p}\}$ and $\overline{K}=\{k_{1},\ldots,k_{q}\}$ of $\{1,\ldots,n\}$ we let 
$   \theta^{J}=\theta^{j_{1}}\wedge \ldots \wedge \theta^{j_{p}}$ and $\theta^{\bar{K}}=\theta^{\bar{k}_{1}}\wedge 
    \ldots\wedge\theta^{\bar{k}_{q}}$.
% \begin{equation}
%     \theta^{J}=\theta^{j_{1}}\wedge \ldots \wedge \theta^{j_{p}}, \qquad \theta^{\bar{K}}=\theta^{\bar{k}_{1}}\wedge 
%     \ldots\wedge\theta^{\bar{k}_{q}},
% %     \label{eq:¥}
% \end{equation}
Then the forms $\theta^{J}\wedge \theta^{\bar{K}}\wedge \omega^{L}$ give rise to an orthonormal frame of $\Lambda_{\C}^{*}H^{*}$ as $J$ and $K$ range 
over all the ordered subsets of $\{1,\ldots,n\}$ and $L$ over all the ordered subsets of $\{2n+1,\ldots,d\}$. Moreover, for $j=1,\ldots,n$ we have 
\begin{gather}
    \varepsilon(\theta^{j})\iota(\theta^{j})(\theta^{J}\wedge \theta^{\bar{K}}\wedge \omega^{L} )= 
    \left\{ 
    \begin{array}{cl}
       \theta^{J}\wedge \theta^{\bar{K}}\wedge \omega^{L}   & \text{if $j\in J$}, \\
       0 & \text{if $j\not \in J$},
%         0 & \text{if $J=\emptyset$},
    \end{array}\right. \\
    \varepsilon(\theta^{\bar{j}})\iota(\theta^{\bar{j}})(\theta^{J}\wedge \theta^{\bar{K}}\wedge \omega^{L} )= 
    \left\{ 
    \begin{array}{cl}
       \theta^{J}\wedge \theta^{\bar{K}}\wedge \omega^{L}   & \text{if $j\in K$}, \\
        0 & \text{if $j\not \in K$}.
%         0 & \text{if $K=\emptyset$}.
    \end{array}\right.
%     \label{eq:¥}
\end{gather}
% and
% \begin{equation}
% %     \label{eq:¥}
% \end{equation}
Combining this with~(\ref{eq:Examples.mu(a)}) then gives
\begin{equation}
    \mu(a)(\theta^{J}\wedge \theta^{\bar{K}}\wedge \omega^{L}) = \mu_{J,\bar{K}}(a)\theta^{J}\wedge \theta^{\bar{K}}\wedge \omega^{L}, \quad 
        \mu_{J,\bar{K}}(a)=\sum_{j\in J} \lambda_{j}-\sum_{j\in K}\lambda_{j}.
     \label{eq:Sublaplacian.eigenvalues-mu(a)} 
\end{equation}
This shows that $\mu(a)$ diagonalizes in the basis of $\Lambda^{*}_{\C}H^{*}_{a}$ provided by the forms of $\theta^{J}\wedge \theta^{\bar{K}}\wedge 
\omega^{L}$ with eigenvalues given by the numbers $\mu_{J,\bar{K}}(a)$. In particular, for $k=0,\ldots,d$ we have 
\begin{equation}
    \op{Sp}\mu(a)_{|_{\Lambda^{k}H^{*}}}=\{\mu_{J,\bar{K}}; \ |J|+|K|\leq k\}. 
%     \label{eq:¥}
\end{equation}

Note that we always have $|\mu_{J,K}|\leq \sum_{j=1}^{n}\lambda_{j}$ with equality if, and only if, one the subsets $J$ or $K$ is empty and the other is 
$\{1,\ldots,n\}$, which occurs for eigenvectors in the subspace spanned by the forms $\theta^{1}\wedge \ldots\theta^{n}\wedge \omega^{L}$ and 
$\theta^{\bar{1}}\wedge \ldots\theta^{\bar{n}}\wedge \omega^{L}$  as $L$ ranges over all the subsets of $\{2n+1,\ldots,d\}$. 

Since $\lambda_{1},\ldots,\lambda_{n}$ are the eigenvalues of $|L(a)|$, each of them counted twice, if follows that the condition~(\ref{eq:Sublaplacian.condition}) 
for $\Delta_{b;k}$ reduces to 
$\pm \sum_{j=1}^{n}\lambda_{j}\not \in \op{Sp}\mu(a)_{|_{\Lambda^{k}H^{*}}}$. 
This latter condition is satisfied if, and only if, the space $\Lambda^{k}_{\C}H^{*}_{a}$ does contain any of the forms 
$\theta^{1}\wedge \ldots\theta^{n}\wedge \omega^{L}$ and $\theta^{\bar{1}}\wedge \ldots\theta^{\bar{n}}\wedge \omega^{L}$ with $L$ subset of 
$\{2n+1,\ldots,d\}$. Therefore, the sublaplacian $\Delta_{b;k}$ satisfies~(\ref{eq:Sublaplacian.condition}) at $a$ if, 
and only if, the integer $k$ is not between $n$ and $n+d-2n=d-n$, that is, if, and only if, the condition $X(k)$ holds at $a$. 
The proof is thus achieved. 
\end{proof}

Suppose now that $M$ is an orientable CR manifold of dimension $2n+1$ with Heisenberg structure $H=\Re(T_{1,0}\oplus T_{0,1})$ and 
let $\theta$ be a global  nonvanishing  section of 
$TM/H$ with associated Levi form $L_{\theta}$. Assume in addition that $T_{\C}M$ is endowed with a Hermitian metric compatible with its CR structure. 
Then it follows from~(\ref{eq:Operators.Tanaka-Kohn}) 
that $\Delta_{b;p,q}$ preserves the bidegree. In this we shall now refine the condition $X(k)$ in each degree $(p,q)$ in 
terms of the following condition. 

\begin{definition}
 For $x\in M$ let $r(x)$ and $\kappa(x)$ respectively denote the rank and the number of negative eigenvalues of the Levi form $L_{\theta}$  at $x$. 
 Then the condition $X(p,q)$ is satisfied at $x$ when we have 
  \begin{equation}
     \{ (p,q),(q,p)\}\cap \{(\kappa(x)+j,r(x)-\kappa(x)+k);\ j,k=0,\ldots, n-r(x)\} =\emptyset.
%        \label{eq:¥}
  \end{equation}
\end{definition}

In particular, when $M$ is $\kappa$-strictly pseudoconvex the condition $X(p,q)$ reduces to $(p,q)\neq (\kappa,n-\kappa)$ 
  and $(p,q)\neq (n-\kappa,\kappa)$. In any case we have: 
  
\begin{proposition}\label{prop:Examples.Deltabpq}
  1) At any point $x\in M$ the Rockland condition for $\Delta_{b;p,q}$ is equivalent to the condition $X(p,q)$.\smallskip
  
  2) The principal symbol of $\Delta_{b;p,q}$ is invertible if, and only if, the condition $X(p,q)$ holds at every point. In particular, 
  when the latter occurs $\Delta_{b;p,q}$ admits a parametrix in $\pvdo^{-2}(M,\Lambda^{p,q})$ and is hypoelliptic with gain of one 
  derivative.
\end{proposition}
\begin{proof}
    As in the proof of Proposition~\ref{prop:Examples.horizontal-sublaplacian}
    thanks to Proposition~\ref{prop:Sublaplacian.Rockland-bundle} the proof reduces to checking 
    that at any point $a \in M$  the condition~(\ref{eq:Sublaplacian.condition}) for 
    $\Delta_{b;p,q}$ is equivalent to the condition $X(p,q)$. 
    
    Near $a$ let $X_{0}$ be a real vector field such that $\theta(X_{0})=1$ and let $Z_{1},\ldots,Z_{n}$ be an orthonormal frame of $T_{1,0}$. Since 
    the Hermitian form $h$ commutes with complex conjugation $\overline{Z_{1}},\ldots,\overline{Z_{n}}$ is an orthonormal frame of $T_{0,1}$. In 
    addition, we write the Levi form $L_{\theta}$ in the form, 
    \begin{equation}
        L_{\theta}(Z,W)=h(L^{c}(x)Z,W), \qquad Z,W\ \text{sections of $T_{1,0}$},
%         \label{eq:¥}
    \end{equation}
    where $L^{c}(x)$ is a Hermitian section of $\End T_{1,0}$. In particular, if for $j,k=1,\ldots,n$ we let $L_{jk}^{c}(x)=h(LZ_{j},Z_{k})$ then we have 
    \begin{equation}
        [Z_{j},\overline{Z_{k}}]=-iL_{jk}^{c}(x)X_{0} \quad \bmod H.
%         \label{eq:}
    \end{equation}
    
    Let $r$ be the rank of $L_{\theta}$ at $a$ and let $\kappa$ be its number of negative eigenvalues. As in the proof of 
    Proposition~\ref{prop:Examples.horizontal-sublaplacian}  we may 
    assume that $L^{c}(a)$ is of the form, 
    \begin{equation}
       L^{c}(a)= \op{diag}(\lambda_{1},\ldots,\lambda_{n}),
%         \label{eq:¥}
    \end{equation}
    where $\lambda_{1},\ldots,\lambda_{n}$ are the eigenvalues of $L^{c}(a)$ ordered in such way that $\lambda_{j}>0$ for $j\leq r-\kappa$ and 
    $\lambda_{j}>0$ for $r-\kappa+1\leq j\leq r$, while $\lambda_{j}=0$ for $j\geq r+1$. 
    
    Let $\theta^{1},\ldots,\theta^{n}$ (resp.~$\theta^{\overline{1}},\ldots,\theta^{\overline{n}}$) be the orthonormal coframe of $\Lambda^{1,0}$ 
    (resp.~$\Lambda^{0,1}$) dual to $Z_{1},\ldots,Z_{n}$ (resp.~$\overline{Z_{1}},\ldots,\overline{Z_{n}}$). For any ordered subsets 
    $J= \{j_{1},\ldots,j_{p}\}$ and $\overline{K}=\{k_{1},\ldots,k_{q}\}$ of $\{1,\ldots,n\}$ we let 
$\theta^{J,\overline{K}}=\theta^{j_{1}}\wedge \ldots \wedge \theta^{j_{p}}\wedge \theta^{\bar{k}_{1}}\wedge 
    \ldots\wedge\theta^{\bar{k}_{q}}$. 
%     \label{eq:¥}
% \end{equation}
Then  the forms $\theta^{J,\overline{K}}$ give rise to an orthonormal coframe of $\Lambda^{*}_{\C}H^{*}$. 

As shown in~\cite{BG:CHM} near $a$ the operator $\Box_{b}$ has the form, 
\begin{gather}
    \Box_{b}=-\frac{1}{2}\sum_{1\leq j \leq n}(Z_{j}\overline{Z_{j}}+\overline{Z_{j}}Z_{j})-i\nu(x) X_{0}+\op{O}_{H}(1),\\
    \nu(x)=\sum_{1\leq j,k \leq n}\iota(\theta^{\overline{j}})\varepsilon(\theta^{\overline{k}})L_{jk}^{c}(x)-\frac{1}{2}\sum_{j=1}^{n}L_{jj}^{c}(x).
%     \label{eq:}
\end{gather}
Hence  $\overline{\Box}_{b}=-\frac{1}{2}\sum_{j=1}^{n}(Z_{j}\overline{Z_{j}}+\overline{Z_{j}}Z_{j})+i\overline{\nu}(x)X_{0}+\op{O}_{H}(1)$. 
% the conjugate operator $\overline{\Box}_{b}$ is of the form, 
% \begin{equation}
%       
% %     \label{eq:}
% \end{equation}
Thus,
\begin{equation}
    \Delta_{b}=\Box_{b}+\overline{\Box}_{b}=-\frac{1}{2}\sum_{j=1}^{n}(Z_{j}\overline{Z_{j}}+\overline{Z_{j}}Z_{j})-i\mu(x) X_{0}+\op{O}_{H}(1),
    \label{eq:Examples.Deltabpq}
 \end{equation}
where $  \mu(x)=\nu(x)-\overline{\nu}(x)=\sum_{j,k=1}^{n}(\iota(\theta^{\overline{j}})\varepsilon(\theta^{\overline{k}})L_{jk}^{c}(x)- 
   \iota(\theta^{j})\varepsilon(\theta^{k})\overline{L_{jk}^{c}}(x))$. In particular, for $x=a$ we obtain
\begin{equation}
    \mu(a)=\sum_{j=1}^{n}(\iota(\theta^{\overline{j}})\varepsilon(\theta^{\overline{j}}) - \iota(\theta^{j})\varepsilon(\theta^{j}))\lambda_{j}= 
    \sum_{j=1}^{n}(\varepsilon(\theta^{j})\iota(\theta^{j})- \varepsilon(\theta^{\overline{j}})\iota(\theta^{\overline{j}}))\lambda_{j}.
%     \label{eq:}
\end{equation}
Therefore, as in~(\ref{eq:Sublaplacian.eigenvalues-mu(a)}) 
we have $\mu(a)(\theta^{J,\overline{K}})=(\sum_{j\in J}\lambda_{j}- \sum_{k\in\overline{ K}}\lambda_{j}) 
\theta^{J,\overline{K}}$. Thus, 
\begin{equation}
    \Sp \mu(a)_{|_{\Lambda^{p,q}}}=\{ \sum_{j\in J}\lambda_{j}-\sum_{k\in \overline{K}}\lambda_{j}; \ |J|=p, |\overline{K}|=q\}. 
%     \label{eq:}
\end{equation}

For $j=1,\ldots,n$ let $X_{j}=\frac{1}{\sqrt{2}}(Z_{j}+\overline{Z_{j}})$ and $X_{n+j}=\frac{1}{i\sqrt{2}}(\overline{Z_{j}}-Z_{j})$. Then 
$X_{1},\ldots,X_{2n}$ is an orthonormal frame of $H$ and we have 
\begin{equation}
    \Delta_{b}=-(X_{1}^{2}+\ldots+X_{2n}^{2})-i\mu X_{0}+\op{O}_{H}(1).
%     \label{eq:¥}
\end{equation}

Let $L(x)=(L_{jk}(x))_{0\leq jk \leq 2n}$ be such that 
\begin{equation}
    [X_{j},X_{k}]=L_{jk}X_{0} \quad \bmod H.
%     \label{eq:}
\end{equation}
Since the integrability of $T_{1,0}$ implies that $[Z_{j},Z_{k}]=[\overline{Z_{j}},\overline{Z_{k}}]=0\bmod H$ one can check that $L(x)$ is of the 
form, 
\begin{equation}
    L(x)=\left( 
    \begin{array}{cc}
        0 & -\Re L^{c}(x)  \\
        \Re L^{c}(x) & 0
    \end{array}\right).
%     \label{eq:¥}
\end{equation}
This implies that $\frac{1}{2}\Tra |L(x)|=\Tra |L^{c}(x)|$. Thus, 
\begin{equation}
    \frac{1}{2}\Tra |L(a)|=\sum_{j=1}^{n}|\lambda_{j}|=\sum_{j=1}^{r}|\lambda_{j}|.
%     \label{eq:¥}
\end{equation}

Since we always have $|\sum_{j\in J}\lambda_{j}-\sum_{k\in  \overline{K}}\lambda_{j}|\leq \sum_{j=1}^{r}|\lambda_{j}|$, 
in the same way as in the proof of Proposition~\ref{prop:Examples.horizontal-sublaplacian} the 
condition~(\ref{eq:Sublaplacian.condition}) for $\Delta_{b;p,q}$ at $a$ becomes 
\begin{equation}
    \sum_{j=1}^{r}|\lambda_{j}|>|\sum_{j\in J}\lambda_{j}-\sum_{k\in \overline{K}}\lambda_{k}|,
%     \label{eq:¥}
\end{equation}
for any ordered subsets $J$ and $\overline{K}$ such that $|J|=p$ and $|\overline{K}|=q$. 

In fact, we have $ \sum_{j=1}^{r}|\lambda_{j}|=|\sum_{j\in 
J}\lambda_{j}-\sum_{k\in \overline{K}}\lambda_{k}|$ if, and only if, either $\overline{K}^{c}\subset \{1,\ldots,r-\kappa\}\subset J$ and $J^{c}\subset 
\{r-\kappa+1,\ldots,r\}\subset \overline{K}$, or $J^{c}\subset \{1,\ldots,r-\kappa\}\subset \overline{K}$ and $\overline{K}^{c}\subset 
\{r-\kappa+1,\ldots,r\}\subset J$. This is possible for  $J$ and $\overline{K}$ such that $|J|=p$ and $|\overline{K}|=q$ if, and only if, either $p$ 
or $q$ if of the form $r-\kappa+j$ with $0\leq j\leq n-r$ and the other is of the form $\kappa+k$ with $0\leq k\leq n-r$, that is,  if, and only if, the 
condition $X(p,q)$ fails at $a$. Therefore, the condition~(\ref{eq:Sublaplacian.condition}) 
at $a$ for $\Delta_{b;p,q}$  is equivalent to the condition $X(p,q)$ at $a$. The proof is therefore complete. 
\end{proof}

\subsection{Gover-Graham operators}  
 Let $M^{2n+1}$ be a strictly pseudoconvex CR manifold equipped with a pseudohermitian strcuture $\theta$, i.e., a contact form $\theta$ anihilating 
 $H=\Re (T_{1,0}\oplus T_{0,1})$ and such that the associated Levi form on $T_{1,0}$ is positive definite. We endow $T_{\C}M$ with the associate Levi 
 metric and for $k=1,\ldots,n+1, n+2,n+4,\ldots$ we let $\boxdot_{\theta}^{(k)}: C^{\infty}(M)\rightarrow 
 C^{\infty}(M)$ be the Gover-Graham of order $k$. In addition, we let $X_{0}$ denote the Reeb  vector field of $\theta$, so that $\imath_{X_{0}} 
 \theta=1$ and $\imath_{X_{0}}d\theta=0$ and we let $\Delta_{b;0}$ denote the horizontal sublaplacian of $M$ acting on functions. 
 
 When $k=1$ the operator $\boxdot_{\theta}^{(1)}$ agrees with the conformal sublaplacian of Jerison-Lee~\cite{JL:YPCRM}, so its principal 
 symbols agrees with that of  $\Delta_{b;0}$ and is therefore invertible. In general, we have: 
 
 \begin{proposition}\label{prop:Examples.Gover-Graham}
     1)  The operator $\boxdot_{\theta}^{(k)}$ is equal to
       \begin{equation}
      (\Delta_{b;0}+i(k-1)X_{0})(\Delta_{b;0}+i(k-3)X_{0}) \cdots (\Delta_{b;0}-i(k-1)X_{0}) +\op{O}_{H}(2k-1). 
          \label{eq:Examples.Gover-Graham-principal-part}
     \end{equation}
%     For $k=1,\ldots,n+1, n+2,n+4,\ldots$ we have
%      \begin{equation}
%    \boxdot_{\theta}^{(k)}=  \prod_{j=0}^{k-1} (\Delta_{b;0}+i(k-1-2j)X_{0})+\op{O}_{H}(2k-1). 
% %          \label{eq:¥}
%      \end{equation}
\indent  2) Unless for the value $k=n+1$ the principal symbol of $ \boxdot_{\theta}^{(k)}$ is invertible.
 \end{proposition}
\begin{proof}
    Let $Z_{1},\ldots,Z_{n}$ be a local orthonormal frame for $T_{1,0}$ and for $w \in \C$ define 
    $\Box_{\omega}=\overline{Z_{1}}Z_{1}+\ldots+\overline{Z_{n}}Z_{n}+iwX_{0}$. The operator $\boxdot_{\theta}^{(k)}$ corresponds to the operators 
    $P_{w,w'}$ and $\mathcal{P}_{w,w'}$ of~\cite{GG:CRIPSL} with $w=w'=\frac{1}{2}(k-1-n)$ under the canonical trivializations of the density  
    bundles $\cE(w,w)=|\Lambda^{n,n}|^{-w/n+2}$, $w \in \R$. Therefore, 
    as explained in~\cite[pp.~15, 25]{GG:CRIPSL}, the operator $\boxdot_{\theta}^{(k)}$ is equal to
    \begin{equation}
      (-2\Box_{-\frac{1}{2}(n+k-1)})(-2\Box_{1-\frac{1}{2}(n+k-1)}) \cdots (-2\Box_{\frac{1}{2}(k-1-n)}) +\op{O}_{H}(2k-1). 
         \label{eq:Examples.Gover-Graham-principal-part-1}
    \end{equation}
    
    Next, notice that for $j=1,\ldots,n$ we have 
    \begin{equation}
        2\overline{Z_{j}}Z_{j}=\overline{Z_{j}}Z_{j} +Z_{j}\overline{Z_{j}}-[Z_{j},\overline{Z_{j}}]= \overline{Z_{j}}Z_{j} 
        +Z_{j}\overline{Z_{j}}+iX_{0}+\op{O}_{H}(1).
%         \label{eq:¥}
    \end{equation}
    In view of the formula~(\ref{eq:Examples.Deltabpq}) for $\Delta_{b;0}$ it follows that $ -2\Box_{w}$ is equal to
   \begin{equation}
      -\sum_{j=1}^{n}(\overline{Z_{j}}Z_{j} +Z_{j}\overline{Z_{j}})-i(2w+n)X_{0}+\op{O}_{H}(1)
       = \Delta_{b;0}-i(2w+n)X_{0}+\op{O}_{H}(1). 
%        \label{eq:¥}
   \end{equation}
   Combining this with~(\ref{eq:Examples.Gover-Graham-principal-part-1}) then yields the formula~(\ref{eq:Examples.Gover-Graham-principal-part}). 
   
   Now, in the same way as in the proof of Proposition~\ref{prop:Examples.Deltabpq} 
   we can show that the sublaplacian $\Delta_{b;0}+i\alpha T$, $\alpha \in \C$, satisfies the 
   condition~(\ref{eq:Sublaplacian.condition}) iff $\alpha$ is not in the singular set $\pm(n+2\N)$. 
   Thanks to~(\ref{eq:Examples.Gover-Graham-principal-part}) we see that:\smallskip 
   
   (i) For $k=1,\ldots,n$ the principal term in~(\ref{eq:Examples.Gover-Graham-principal-part}) is the product of sublaplacians $\Delta_{b;0}+i\alpha 
   X_{0}$ 
   with $|\alpha|\leq 
   k-1<n$.\smallskip
   
   (ii) For $k=n+1$ the principal term in~(\ref{eq:Examples.Gover-Graham-principal-part}) contains the factors $\Delta_{b;0}\pm inX_{0}$, whose principal symbols 
   are not invertible.\smallskip
   
   (iii) For $k=n+2,n+4,\ldots$ the integers $k-1$ and $n$ have opposite parities, so  the principal term in~(\ref{eq:Examples.Gover-Graham-principal-part}) is 
   the product of sublaplacians 
   $\Delta_{b;0}+i\alpha X_{0}$ with integers $\alpha$ which are not not in the singular set $\pm(n+2\N)$, since their parity is the opposite 
   to that of $n$.\smallskip
   
  \noindent Therefore, unless for $k=n+1$ the principal symbol of $ \boxdot_{\theta}^{(k)}$ appears as the product of invertible symbols, hence is 
  invertible.  
\end{proof}

\subsection{Contact Laplacian}
Let  $(M^{2n+1},H)$ be an orientable contact manifold, let $\theta$ be a contact form on $H$ with  Reeb 
vector field $X_{0}$ and let $J$ be a calibrated almost complex structure on $H$,  so that we can endow 
$M$ with the Riemannian metric $g_{\theta, J}=d\theta(.,J.)+\theta^{2}$.

Consider the contact complex of $M$, 
\begin{equation}
    C^{\infty}(M)\stackrel{d_{R;0}}{\rightarrow}
    \ldots C^{\infty}(M,\Lambda^{n}_{1})\stackrel{D_{R;n}}{\rightarrow} C^{\infty}(M,\Lambda^{n}_{2}) 
    \ldots \stackrel{d_{R;2n-1}}{\rightarrow} C^{\infty}(M,\Lambda^{2n}). 
%      \label{eq:Operators.contact-complex}
\end{equation}
Let $\Delta_{R;k}: C^{\infty}(M,\Lambda^{k})\rightarrow C^{\infty}(M,\Lambda^{k})$ be the contact Laplacian in degree $k\neq n$ and let 
$\Delta_{R;nj}:C^{\infty}(M,\Lambda^{n}_{j})\rightarrow C^{\infty}(M,\Lambda^{n}_{j})$, $j=1,2$, be the contact Laplacians in degree $n$ as defined 
in~(\ref{eq:Operators.contact-Laplacian1})--(\ref{eq:Operators.contact-Laplacian1}).

The almost complex structure $J$ of $H$ defines a bigrading on $\Lambda^{*}_{\C}H^{*}$. More precisely, we have an 
   orthogonal splitting $H\otimes \C=T_{1,0}\oplus T_{0,1}$ with $T_{1,0}=\ker(J+i)$ and $T_{0,1}=\overline{T_{1,0}}=\ker(J-i)$.
   Therefore, if we consider the subbundles $\Lambda^{1,0}=T_{1,0}^{*}$ and $\Lambda^{0,1}=T_{0,1}^{*}$ of $H^{*}\otimes \C\subset T^{*}_{\C}M$ 
   then we have the orthogonal decomposition $  \Lambda^{*}_{\C}H^{*}=\bigoplus_{0\leq p,q\leq n}\Lambda^{p,q}$ with 
   $\Lambda^{p,q}=(\Lambda^{1,0})^{p}\wedge (\Lambda^{0,1})^{q}$. We then get a bigrading on $\Lambda_{j}^{n}$, $j=1,2$, by letting 
 \begin{equation}
     \Lambda_{j}^{n}=\bigoplus_{p+q=n}\Lambda^{p,q}_{j}, \qquad \Lambda^{p,q}_{j}= \Lambda^{n}_{j}\cap \Lambda^{p,q}.
 \end{equation}

As shown by Rumin~\cite[Prop.~7]{Ru:FDVC} there exist first order differential operators $P_{k}$, $k=0,..,n-1,n+1,\ldots,2n$ and second order differential operators 
$P_{p,q}^{(j)}$, $j=1,2$, $p+q=n$, such that 
\begin{gather}
    (n-k+2)\Delta_{R;k}= (n-k)(n-k+1)\Delta_{b;k} +P_{k}^{*}P_{k}, \qquad \text{$k=0,..,n-1$},
    \label{eq:Operators.equalities-DeltaR-Deltab.1}\\
  (k-n+2)\Delta_{R;k}= (k-n)(k-n+1)\Delta_{b;k} +P_{k}^{*}P_{k}, \quad \text{$k=n+1,\ldots,2n$},\\
    4\Delta_{R;nj}=\Delta_{b;n}^{2}+(P_{p,q}^{(j)})^{*}P_{p,q}^{(j)} \qquad \text{on $\Lambda^{p,q}_{\pm}$ with $\sup(p,q)\geq 1$},\\
    \Delta_{R;nj}=(\Delta_{b;n}+iX_{0})^{2}  \qquad \text{on $\Lambda^{n,0}_{j}$},\\
    \Delta_{R;nj}=(\Delta_{b;n}-iX_{0})^{2} \qquad \text{on $\Lambda^{0,n}_{j}$}.
    \label{eq:Operators.equalities-DeltaR-Deltab.5}
\end{gather}
These formulas enabled him to prove: 

\begin{proposition}[{\cite[p.~300]{Ru:FDVC}}]\label{prop:Examples.Rockland-contact-Laplacian}
    The operators $\Delta_{R;k}$, $k\neq n$, and $\Delta_{R;nj}$, $j=1,2$, satisfy the Rockland condition at every point. 
\end{proposition}
% \begin{proof}
%     
% \end{proof}

Proposition~\ref{prop:Examples.Rockland-contact-Laplacian} 
allowed Rumin to apply results of Helffer-Nourrigat~\cite{HN:HMOPCV} for proving the maximal hypoellipticity of $\Delta_{R}$ 
in every degree. In particular, unlike the Kohn Laplacian and the horizontal sublaplacian, the contact Laplacian is hypoelliptic in every bidegree. 

Alternatively, we may combine Proposition~\ref{prop:Examples.Rockland-contact-Laplacian}  with 
Theorem~\ref{thm:Chap3.Rockland-Parametrix-order+} to get: 

\begin{proposition} 1) The contact Laplacian $\Delta_{R;k}$, $k\neq n$, has an invertible 
principal symbol of degree $-2$, hence admits a parametrix in $\pvdo^{-2}(M,\Lambda^{k})$ and is 
hypoelliptic with gain of one derivative.\smallskip

2) The contact Laplacian $\Delta_{R;nj}$, $j=1,2$,  has an invertible 
principal symbol of degree $-4$, hence admits a parametrix in $\pvdo^{-4}(M,\Lambda^{n}_{j})$ and 
is hypoelliptic with gain of two derivatives.
\end{proposition}

%%%%%%%%%%%%%%%%%%%%%%%%%%%%%%%%%%
%%%%%%%%%%% Chap4.tex %%%%%%%%%%%%%%%%%%
%%%%%%%%%%%%%%%%%%%%%%%%%%%%%%%%%%%

\chapter{Holomorphic families of $\mathbf{\Psi_{H}}$DOs}
\label{chap.HolPHDO}

In this chapter we define  holomorphic families of \psivdos\ and check their main properties. To this end we make use of 
an "almost homogeneous'' approach to the Heisenberg calculus which we describe in the first section.

 \section{Almost homogeneous approach to the Heisenberg calculus} 
In this section we explain how the \psidos\ can be described in terms of symbols and 
kernels which are almost homogeneous, in the sense that there are homogeneous modulo infinite order terms.

\begin{definition}\label{def:HolPHDO.almost-homgenous-symbols}
A symbol $q(x,\xi)\in C^{\infty}(\URd)$ is almost homogeneous of degree $m$, $m \in \C$, when we have
\begin{equation}
     q(x,\lambda.\xi)-\lambda^m q(x,\xi)\in S^{-\infty}(\URd)\qquad \text{for any $\lambda>0$}.
\end{equation}
The space of almost homogeneous symbols of degree~$m$ is denoted $S_{\ah}^m(\URd)$.
\end{definition}\begin{lemma}[{\cite[Prop.~12.72]{BG:CHM}}]\label{lem:HolPHDO.almost-homogeneity}
Let  $q(x,\xi) \in C^{\infty}(\URd)$. Then  the following are equivalent:\smallskip 

(i) $q(x,\xi)$ is almost homogeneous of degree $m$;\smallskip 

(ii)  $q(x,\xi)$ is in $S^{m}(\URd)$ and we have $q\sim p_{m}$ with $p_{m}\in S_{m}(\URd)$, i.e., the only nonzero term in the asymptotic 
expansion~(\ref{eq:PsiVDO.asymptotic-expansion-symbols}) for $q$ is $p_{m}$. 
\end{lemma}
 
Granted this we shall now prove: 
\begin{lemma}\label{lem:HolPHDO.almost-homogeneous-characterization}
    Let $p \in C^{\infty}(\URd)$. Then we have equivalence:\smallskip
    
    (i) $p$ belongs to $S^{m}(\URd)$.\smallskip 
    
    (ii) For $j=0,1,..$ there exists $q_{m-j}\in S_{\ah}^{m-j}(\URd)$ such that $p\sim \sum_{j\geq 0} q_{m-j}$.
\end{lemma}
\begin{proof}
    Suppose that $p\sim \sum_{j\geq 0} q_{m-j}$ with $q_{m-j}\in S_{\ah}^{m-j}(\URd)$. 
    By Lemma~\ref{lem:HolPHDO.almost-homogeneity} there exists $p_{m-j}\in S_{m-j}(\URd)$ such that $q_{m-j}\sim p_{m-j}$. Then 
    we have $p\sim \sum_{j\geq 0}p_{m-j}$ and so  $p$ belongs to $S^{m}(\URd)$. Hence (ii) implies (i).
        
    Conversely, assume that $p$ is in $S^{m}(\URd)$.  Then we have $p\sim \sum_{j\geq 0} 
    p_{m-j}$ with $p_{m-j}\in S_{m-j}(\URd)$. Let $\varphi \in C^{\infty}_{c}(\Rd)$ be such that $\varphi(\xi)=1$ near $\xi=1$ and $\varphi(\xi)=0$ for 
    $\|\xi\|\leq 1$. For $j=0,1,..$ set 
    $q_{m-j}(x,\xi)=(1-\varphi(\xi))p_{m-j}(x,\xi)$. For any $t>0$ the symbol $q_{m-j}(x,t.\xi)-t^{m-j}q_{m-j}(x,\xi)$ is equal to
    \begin{equation}
        (\varphi(t.\xi)-\varphi(\xi))p_{m-j}(x,\xi) \in S^{-\infty}(\URd),
%         \label{eq:}
    \end{equation}
     so $q_{m-j}$ belongs to 
    $S^{m-j}_{\ah}(\URd)$. Moreover, as $q_{m-j}(x,\xi)=p_{m-j}(x,\xi)$ for $\|\xi\|>1$ we see that $p\sim \sum_{j\geq 0} 
    \tilde{p}_{m-j}$. Hence (i) implies (ii). 
\end{proof}

The almost homogeneous symbols have been considered in~\cite[Sect.~12]{BG:CHM} already. In the sequel it will be important to have a 
"dual'' notion of almost homogeneity for distributions as follows. \begin{definition}
    The space $\cD'_{\reg}(\Rd)$ consists of the distributions on $\Rd$ that are smooth outside the 
    origin.  It is  endowed with
    the weakest topology that makes continuous the inclusions of  $\cD'_{\reg}(\Rd)$ into 
    $\cD'(\Rd)$ and  $C^{\infty}(\Rdo)$.
\end{definition}
\begin{definition}\label{def:HolPHDO.almost-homogeneous-kernels}
 A distribution $K(x,y)\in C^\infty(U)\hotimes \cD'_{\reg}(\Rd)$  is said to be almost homogeneous of degree $m$, $m \in \C$, when
    \begin{equation}
         K(x, \lambda.y)-\lambda^{m}K(x,y)\in C^{\infty}(\URd) \quad \text{for any  $\lambda>0$}. 
     \end{equation}
We let $\cK_{\ah}^m(\URd)$ denote the space of almost homogeneous distributions of degree $m$.
 \end{definition}

\begin{proposition}%[See also~{\cite[pp.~18-21]{Ta:NCMA}}]
    \label{prop:HolPHDO.characterization-almost-homogeneous-kernels}
 Let $K(x,y)\in C^\infty(U)\hotimes \cD'_{\reg}(\Rd)$. Then the following are equivalent:\smallskip 
 
(i) $K(x,y)$ belongs to $\cK_{\ah}^m(\URd)$.\smallskip

(ii) We can put $K(x,y)$ into the form,
\begin{equation}
    K(x,y)=K_{m}(x,y)+R(x,y),
\end{equation}
for some $K_{m}\in \cK_{m}(\URd)$ and $R\in C^{\infty}(\URd)$.\smallskip

(iii) We can put $K(x,y)$ into the form,
\begin{equation}
    K(x,y)=\check{p}_{\xiy}(x,y)+R(x,y),
\end{equation}
for some $p\in S_{\ah}^{\hat{m}}(\URd)$, $\hat{m}=-(m+d+2)$, and $R\in C^{\infty}(\URd)$. 
\end{proposition}
\begin{proof}
    First, if $K_{m}\in \cK_{m}(\URd)$ then~(\ref{eq:PsiHDO.homogeneity-K-m}) 
    implies that, for any  $\lambda>0$, the distribution $K(x, \lambda.y)-\lambda^{m}K(x,y)$ is in $C^{\infty}(\URd)$. 
    Thus (ii) implies~(i). 
    
    Second, let $p\in S_{\ah}^{\hat{m}}(\URd)$. By Lemma~\ref{lem:HolPHDO.almost-homogeneity} there is $p_{\hat{m}}\in S_{m}(\URd)$ such that 
    $p\sim p_{\hat{m}}$.  
    Thanks to Lemma~\ref{prop:PsiHDO.Sm-Km} we extend $p_{\hat{m}}$ into a distribution $\tau(x,\xi)$  in  
 $C^{\infty}(U) \hotimes \cS'(\Rd)$ such that $\check{\tau}_{\xiy}(x,y)$ is in $\cK_{m}(\URd)$. Let  
 $\varphi \in C^{\infty}_{c}(\Rd)$ be such that $\varphi =1$ near the origin. Then we can write 
 \begin{equation}
     p=\tau+\varphi(p-\tau) +(1-\varphi)(p-p_{\hat{m}}).
 \end{equation}
 Here $\varphi(\xi)(p(x,\xi)-\tau(x,\xi))$ belongs to $C^{\infty}(U)\otimes \cD'(\Rd)$ and is supported on a fixed compact set with respect to $\xi$, 
 so $[ \varphi(p-\tau)]^{\vee}_{\xiy}$ is smooth. 
 Moreover, as $p\sim p_{\hat{m}}$ both $(1-\varphi)(p-p_{m})$ and  $[(1-\varphi)(p-p_{m})]^{\vee}_{\xiy}$ are in $S^{-\infty}(\URd)$. It then follows  
 that $\check{p}_{\xiy}$ coincides with $\check{\tau}_{\xiy}$ up to an element of $C^{\infty}(\URd)$.  Since $\check{\tau}_{\xiy}(x,y)$ is in 
 $\cK_{m}(\URd)$ we deduce from this that (iii) implies~(ii).
 
 To complete the proof it remains to show that (i) implies (iii). Assume that $K(x,y)$ belongs to $\cK_{\ah}^m(\URd)$. 
 Let  $\varphi(y) \in C^{\infty}_{c}(\Rd)$ be  
 so that $\varphi(y) =1$ near $y=0$ and set $p=(\varphi K)^{\wedge}_{\yxi}(x,y)$. Then $p$ is a smooth and has slow growth with respect to $\xi$. 
 Moreover $\check{p}_{\xiy}(x,y)$ differs from 
 $K(x,y)$ by the smooth function $(1-\varphi(y))K(x,y)$. 
 
 Next, using~(\ref{eq:PsiHDO.dilations-Fourier-transform}) we see that, for any $\lambda>0$, the function 
 $p(x,\lambda.\xi)-\lambda^{\hat m}p(x,\xi)$ is equal to
\begin{multline}
    \lambda^{-(d+2)}[\varphi(\lambda^{-1}.y)K(x,\lambda^{-1}.y)-\lambda^{-m}\varphi(y) K(x,y)]^{\wedge}_{\yxi}\\ 
    =  \lambda^{-(d+2)}[(\varphi(\lambda^{-1}.y)-\varphi(y))K(x,\lambda^{-1}.y)+
    \varphi(y) (K(x,\lambda^{-1}.y)-\lambda^{-m}K(x,y))]^{\wedge}_{\yxi}.
 \end{multline}
Notice that $(\varphi(\lambda^{-1}.y)-\varphi(y))K(x,\lambda^{-1}.y)+
    \varphi(y) (K(x,\lambda^{-1}.y)-\lambda^{-m}K(x,y))$ belongs to $C^{\infty}(\URd)$ and is  compactly supported with respect to 
$y$, so it belongs to $S^{-\infty}(\URd)$. Since the latter is also true for Fourier transform with respect to $y$ we see that 
 $ p(x,\lambda.\xi)-\lambda^{\hat 
m}p(x,\xi) $ is in $S^{-\infty}(\URd)$ for any $\lambda>0$, that is, the symbol $p$ is almost homogeneous of degree $\hat{m}$. It then follows that 
the distribution $K$ satisfies~(iii). This proves that (i) implies (iii). The proof is thus achieved.
\end{proof}
\section{Holomorphic families of $\mathbf{\Psi}_{H}$DO's} From now on we let $\Omega$ denote an open subset of $\C$. 
\begin{definition}\label{def:complex.symbols}
A family $(p_{z})_{z\in\Omega}\subset S^*(\URd)$ is holomorphic when:  \smallskip  
     
     (i) The order $m(z)$ of $p_{z}$ depends analytically on $z$; \smallskip 

     (ii) For any $(x,\xi)\in \URd$ the function $z\rightarrow p_{z}(x,\xi)$ is holomorphic on $\Omega$; \smallskip 
     
     (iii) The bounds of the asymptotic expansion~(\ref{eq:PsiVDO.asymptotic-expansion-symbols}) 
     for $p_{z}$ are locally uniform with respect to $z$, i.e., we have 
     $p_{z} \sim \sum_{j\geq 0} p_{z, m(z)-j}$, 
         $p_{z,m(z)-j}\in S_{m(z)-j}(\URd)$, and for
         any integer $N$ and for any compacts $K\subset U$ and $L\subset \Omega$ we have 
\begin{equation}
   | \partial_{x}^\alpha\partial_{\xi}^\beta (p_{z}-\sum_{j<N}  
            p_{z,m(z)-j})(x,\xi)| \leq C_{NKL\alpha\beta} \|\xi\|^{\Re m(z)-N-\brak\beta}, \quad 
    \label{eq:complex.symbols.asymptotic-expansion}
\end{equation}
for $(x,z)\in K\times L$ and $\|\xi\|\geq  1$. 
% $(p_{z})_{z\in\Omega}\subset S^*(\URd)$ that are holomorphic. 

We let $\Hol(\Omega,S^*(\URd))$ denote the set of holomorphic families with  values in $S^*(\URd)$.
\end{definition}

\begin{remark}\label{rem:HolPHDO.analyticity-homogeneous-symbols}
    If $(p_{z})_{z\in \Omega}$ is a holomorphic family of symbols then the homogeneous symbols 
    $p_{z,m(z)-j}$ depend   
    analytically on $z$. Indeed, for $\xi\neq 0$ we have  
    \begin{equation}
        p_{z,m(z)}(x,\xi)=\lim_{\lambda \rightarrow \infty} \lambda^{-m(z)}  
            p_{z}(x,\lambda.\xi).
    \end{equation}
     Since  the above axioms imply that the family $(\lambda^{-m(z)}  p_{z}(x,\lambda.\xi) )_{\lambda \geq 1}$ is bounded in the Fr\'echet-Montel space $\Hol(\Omega, 
C^\infty(\URdo)$ the convergence actually holds in $\Hol(\Omega, C^\infty(\URdo)$. Hence $p_{z,m(z)}$ depends 
analytically on $z$. Moreover, as for $j=1,2,\ldots$ and for $\xi\neq 0$ we also have 
    \begin{equation}
        p_{j,z}(x,\xi)=\lim_{\lambda \rightarrow \infty} \lambda^{j-m(z)} 
                     (p_{z}(x,\lambda.\xi) -\sum_{l< j}\lambda^{m(z)-l}p_{z,m(z)-l}(x,\xi)),
    \end{equation}
an easy induction shows that all the symbols $p_{z,m_{z}-j}$ depend analytically on $z$.
\end{remark}

Recall that $\psinf(U)=\cL(\cE'(U),C^{\infty}(U))$ is naturally a Fr\'echet space which is isomorphic to $C^{\infty}(U\times U)$ by the Schwartz's kernel 
theorem. Therefore holomorphic families of smoothing operators make sense and 
we may define holomorphic families of \psivdos\ as follows.

\begin{definition}\label{def:HolPHDO.holomorphic-families}
 A family $(P_{z})_{z\in \Omega}\subset \pvdo^{m}(U)$ is holomorphic when it can be put into the form 
\begin{equation}
     P_{z} = p_{z}(x,-iX) + R_{z} \qquad z \in \Omega, 
\end{equation}
 with $(p_{z})_{z\in \Omega}\in \Hol(\Omega, S^{*}(\URd))$ and $(R_{z})_{z\in \Omega}\in \Hol(\Omega, \Psi^{-\infty}(U))$. 

  We let $\Hol(\Omega,\pvdo^{*}(U))$ denote the set of holomorphic families of \psivdos.
\end{definition}

For technical sake it will be useful to consider the symbol class below.

\begin{definition}
The space $\Svb^k(\URd)$, $k \in \R$, consists of functions $p(x,\xi)$ in $C^{\infty}(\URd)$ such that, for any 
   compact $K\subset U$,  we have 
   \begin{equation}
     |\partial_{x}^\alpha\partial^\beta_{\xi}p(x,\xi)| \leq 
    C_{K\alpha\beta}(1+\|\xi\|)^{k-\brak\beta}, \qquad (x,\xi)\in K\times\Rd.  
    \label{eq:HolPHDO.SvbU.estimates}
   \end{equation} 
 Its space topology  is defined by the family of seminorms given by sharpest constants $C_{K\alpha\beta}$'s 
 in the estimates~(\ref{eq:HolPHDO.SvbU.estimates}). 
\end{definition}

Note that the estimates~(\ref{eq:PsiVDO.asymptotic-expansion-symbols})  imply that $S^{m}(\URd)$, $m \in \C$, 
is contained in $\Svb^k(\URd)$ for any $k \geq \Re m$.  

\begin{proposition}\label{prop:Complexpowers.operators.properties}
    Let $(P_{z})_{z\in \Omega}$ be a holomorphic family of \psivdo's. Then:\smallskip 
    
    1) $(P_{z})_{z\in \Omega}$ gives rise to holomorphic families with values in 
        $\cL(C_{c} ^\infty(U),C^\infty(U))$ and $\cL(\cE'(U)),\cD'(U))$. \smallskip 
	
        2) Off the diagonal of $U\times U$ the distribution kernel of $P_{z}$ is represented by a holomorphic family of smooth  functions. 
    \end{proposition}
\begin{proof}
    Without any loss of generality we may suppose that $P_{z}=p_{z}(x,-iX)$, with $(p_{z})_{z\in \Omega}$ in 
    $\Hol(\Omega, S^{*}(\URd))$. Moreover, shrinking $\Omega$ if necessary,  we may also assume that the degree $m_{z}$ of 
    $p_{z}$ stays bounded, as much so   $(p_{z})_{z\in \Omega}$ is contained in $\Svb^k(\URd)$ for some real $k\geq 0$.  
For $j=1,\ldots,d$ let $\sigma_{j}(x,\xi)$ denote the 
classical symbol of $-iX_{j}$ and set $\sigma=(\sigma_{0},\ldots,\sigma_{d})$. 
Then the proof of~\cite[Prop.~10.22]{BG:CHM} shows that the map $p(x,\xi) \rightarrow p^{\sigma}(x,\xi):=p(x,\sigma(x,\xi))$ is 
continuous from 
$\Svb^{k}(\URd)$ to $S^{k}_{\frac{1}{2},\frac{1}{2}}(\URd)$. Thus, the family $(p^{\sigma}_{z})_{z \in \Omega}$ 
belongs to $\Hol(\Omega,S^{k}_{\frac{1}{2},\frac{1}{2}}(\URd))$.

Next, it follows from the proof of~\cite[Thm.~2.2]{Ho:PDOHE} that:\smallskip   

(i) The  quantization map $q \rightarrow q(x,D)$ induces continuous $\C$-linear maps from $S^{k}_{\frac{1}{2},\frac{1}{2}}(\URd)$ 
 to $\cL(C_{c} ^\infty(U),C^\infty(U))$ and to $\cL(\cE'(U)),\cD'(U))$; 
\smallskip 

(ii) The linear map $q(x,\xi) \rightarrow \check{q}_{\xiy}(x,y)$ is continuous from $S^{k}_{\frac{1}{2},\frac{1}{2}}(\URd)$ to $C^{\infty}(U)\otimes 
\cD'_{\reg}(\Rd)$, in such way that for any $q\in S^{k}_{\frac{1}{2},\frac{1}{2}}(\URd)$ the distribution kernel $\check{q}_{\xiy}(x,x-y)$ 
of $q(x,D)$ is represented off the diagonal by a smooth function depending continuously on $q$.\smallskip 

As  a continuous $\C$-linear map is  analytic it follows that, on the 
one hand,  $(p_{z}(x,-iX))_{z\in \Omega}$ gives rise to  holomorphic families with values in  $\cL(C_{c} ^\infty(U),C^\infty(U))$ and 
$\cL(\cE'(U)),\cD'(U))$ and, 
 on the other hand,  the distribution kernel of $P_{z}$ is represented off the diagonal by a holomorphic family of smooth  functions. The proposition 
 is thus proved.
\end{proof}

\begin{definition}
    Let $(P_{z})_{z \in \Omega}\subset \cL(C^{\infty}_{c}(U),C^{\infty}(U))$ and for $z \in \Omega$ let $k_{P_{z}}(x,y)$ denote the distribution kernel of 
    $P_{z}$. Then the family $(P_{z})_{z \in \Omega}$ is said to be uniformly properly supported when, for any 
compact $K\subset U$, there exist compacts $L_{1}\subset U$ and $L_{2}\subset K$ such that for any $z\in \Omega$ we have
\begin{equation}
    \supp k_{P_{z}}(x,y)\cap (U\times K) \subset L_{1} \quad \text{and} \quad \supp k_{P_{z}}(x,y)\cap (K\times U) \subset L_{2}.
     \label{eq:HolPHDO.uniformly-properly-supported}
\end{equation}
\end{definition}

Bearing in mind the above definition we have:
\begin{proposition}\label{prop:Complexpowers.operators.properties2}
Let $(P_{z})_{z\in \Omega}$ be a holomorphic family of \psivdos.  \smallskip 
                   
 1) We can write $P_{z}$ in the form $P_{z}=Q_{z}+R_{z}$ with  $(Q_{z})_{z\in \Omega}\in \Hol(\Omega, \pvdo^{*}(U)$ uniformly properly supported 
 and $(R_{z})_{z\in \Omega}\in \Hol(\Omega, \Psi^{-\infty}(U)$.\smallskip 

2) If the family $(P_{z})_{z \in \Omega}$ is uniformly properly supported then it gives rise to 
 holomorphic families of continuous endomorphisms of $C_{c}^\infty(U)$ and $C^\infty(U)$ and of $\cE'(U)$ and $\cD'(U)$.
\end{proposition}
\begin{proof}
    Let $(\varphi_{i})_{i \geq 0} \subset C^{\infty}_{c}(U)$ be a partition of the unity subordinated to a locally finite covering 
    $(U_{i})_{i \geq 0}$ of $U$ by 
    relatively compact open subsets. For each $i\geq 0$ let $\psi_{i}\subset C^{\infty}_{c}(U)_{i}$ be such that 
    $\psi_{i}=1$ near $\supp \varphi_{i}$ and set $\chi(x,y)=\sum \varphi_{i}(x)\psi_{i}(y)$. Then $\chi$ is a smooth function on $U\times U$ which is 
    properly supported and such that $\chi(x,y)=1$ near the diagonal of $U\times U$. 
    
    For $z \in \Omega$ let $k_{P_{z}}(x,y)$ denote the distribution kernel of $P_{z}$ and let $Q_{z}$ and $R_{z}$ be the elements of 
    $\cL(C^{\infty}_{c}(U), C^{\infty}(U))$ with respective distribution kernels,
     \begin{equation}
        k_{Q_{z}}(x,y)=\chi(x,y)k_{P_{z}}(x,y) \quad \text{and} \quad k_{R_{z}}(x,y)=(1-\chi(x,y))k_{P_{z}}(x,y).
    \end{equation}
    Notice that since $\chi$ is properly supported the family $(Q_{z})_{z\in \Omega}$ is uniformly properly supported.
    As  by Proposition~\ref{prop:Complexpowers.operators.properties} the distribution $k_{P_{z}}(x,y)$ is represented outside the 
        diagonal of $U\times U$ by a holomorphic family of smooth  functions, we see that $(k_{R_{z}}(x,y))$ is a holomorphic family of smooth kernels, 
        i.e.,~$(R_{z})_{z\in \Omega}$ is a holomorphic family of smoothing operators. Since $Q_{z}=P_{z}-R_{z}$ it follows that $(Q_{z})_{z\in \Omega}$ 
        is a holomorphic family of \psivdos. Hence the first assertion.
        
       Assume now that $(P_{z})_{z\in \Omega}$ is uniformly properly supported. Thanks to Proposition~\ref{prop:Complexpowers.operators.properties}
        we already know that  $(P_{z})_{z \in \Omega}$ gives rise to holomorphic families with 
        values in $\cL(C^{\infty}_{c}(U),C^{\infty}(U))$ and $\cL(\cE'(U),\cD'(U))$. Let $K$ be a compact subset of $U$. 
        Then~(\ref{eq:HolPHDO.uniformly-properly-supported})  implies that  
        there exists a compact $L\subset U$ such that for every $z \in \Omega$ the operator $P_{z}$ maps 
        $C_{K}^{\infty}(U)$ to $C^{\infty}_{L}(U)$ and $\cE'_{K}(U)$ to $\cE'_{L}(U)$, in such way that  
        $(P_{z})_{z\in \Omega}$ gives rise to holomorphic families with values in $\cL(C_{K}^{\infty}(U),C^{\infty}_{L}(U))$ and $\cL(\cE'_{K}(U),\cE'_{L}(U))$. 
        In view of the definitions of the topologies of $C^{\infty}_{c}(U)$ and $\cE'(U)$ as the inductive limit topologies of $C^{\infty}_{K}(U)$ and 
        $\cE_{K}(U)$ as $K$ ranges over compacts of $U$, this shows that the family 
        $(P_{z})_{z\in \Omega}$ gives rise to elements of $\Hol(\Omega,C_{c}^\infty(U))$ and $\Hol(\Omega, \cE'(U))$.
       
        Next, let $(\varphi_{i})_{i\geq 0}\subset C^{\infty}_{c}(U)$ be a partition of unity. For each index $ i$ let $K_{i}$ be a compact 
        neighborhood of $\supp \varphi_{i}$. Then~(\ref{eq:HolPHDO.uniformly-properly-supported}) 
        implies that there exists a compact $L_{i}\subset U$ such that $\supp 
        k_{P_{z}}(x,y)\cap(K_{i}\times U)\subset L_{i}$ for every $z \in \Omega$. Let $\psi_{i}\in C^{\infty}(U)$ be such that $\psi_{i}=1$ near 
        $K_{i}$. Then we have 
        \begin{equation}
            P_{z}=\sum_{i\geq 0}\varphi_{i}P_{z}=\sum_{i\geq 0}\varphi_{i}P_{z}\psi_{i}.
             \label{eq:HolPHDO.uniformly-properly-supported-locally-finite-sum}
        \end{equation}
        Since each family $(\varphi_{i}P_{z}\psi_{i})_{z\in \Omega}$ is holomorphic with values in $\cL(C^{\infty}(U))$ and $\cL(\cD'(U))$ and the 
        sums are locally finite this shows that $(P_{z})_{z\in \Omega}$ gives rise to elements of $\Hol(\Omega,C^\infty(U))$ and $\Hol(\Omega,\cD'(U))$. 
\end{proof}

\section{Composition of holomorphic families of \psivdos}
Let us now look at the analyticity of the composition of \psivdos. To this 
end we need to deal with holomorphic families of almost homogeneous symbols as follows. 

\begin{definition}
    A family $(q_{z})_{z\in \Omega}\subset S_{\ah}^{*}(\URd)$ is holomorphic when:\smallskip
    
    (i) The degree $m(z)$ of $q_{z}$ is a holomorphic function on $\Omega$;\smallskip
    
    (ii) The family $(q_{z})_{z\in \Omega}$ belongs to $\Hol(\Omega,  C^{\infty}(\URd))$;\smallskip 
    
    (iii) $(q_{z}(x,t.\xi) -t^{m(z)}q_{z}(x,\xi))_{z\in \Omega}$ is in $\Hol(\Omega, S^{-\infty}(\URd))$ for any $t>0$.\smallskip
    
   \noindent We let $\Hol(\Omega,  S_{\ah}^{*}(\URd))$ denote the set of holomorphic families of almost homogeneous symbols.
\end{definition}

\begin{lemma}\label{lem:HolPHDO.characterization-almost-homogeneous-symbols}
 Let $(q_{z})_{z\in\Omega}\in \Hol(\Omega, C^{\infty}(\URd))$. Then  the following are equivalent:\smallskip 
 
 (i) The family $(q_{z})_{z\in\Omega}$ is in $ \Hol(\Omega, S_{\ah}^{*}(\URd))$ and has degree 
  degree $m(z)$;\smallskip

 (ii) The family $(q_{z})_{z\in \Omega}$ lies in $\Hol(\Omega, S^{*}(\URd))$ and, in the sense of~(\ref{eq:complex.symbols.asymptotic-expansion}), 
 we have $q_{z} \sim p_{z}$ where, for 
 every $z \in \Omega$, the symbol $p_{z}$ belongs to $S_{m(z)}(\URd)$.
\end{lemma}
\begin{proof}
    Assume that  $(q_{z})_{z\in \Omega}$ is in $\Hol(\Omega, S^{*}(\URd))$ and, in the sense of~(\ref{eq:complex.symbols.asymptotic-expansion}), 
      we have $q_{z} \sim p_{z}$ where, for 
 every $z \in \Omega$, the symbol $p_{z}$ belongs to $S_{m(z)}(\URd)$. Then the order $m(z)$ of $q_{z}$ is a holomorphic function on $\Omega$ and,
for any compact subset $K \subset U$, any integer $N$ and any open $\Omega' \subsubset \Omega$, we have
     \begin{equation}
          | \partial_{x}^\alpha\partial^\beta_{\xi}(q_{z}-p_{z})(x,\xi)\leq C_{NK\Omega'\alpha\beta}\|\xi\|^{-N}, 
         \label{eq:HolPHDO.asymptotics-symbol-single-term}
     \end{equation}
for $x\in K$, $\|\xi\|\geq 1$ and $z \in \Omega'$. Therefore, for any $t>0$, the family 
$\{q_{z}(x,t.\xi)-t^{m(z)}q_{z}(x,\xi)\}_{z\in \Omega}=\{(q_{z}(x,t.\xi)-q_{z}(x,t.\xi))-t^{m(z)}(q_{z}(x,\xi)-q_{z}(x,\xi))\}_{z\in \Omega}$  is 
contained in $\Hol(\Omega, S^{-\infty}(\URd))$. Hence $(q_{z})_{z\in \Omega}$ belongs to $ \Hol(\Omega, S_{\ah}^{*}(\URd))$.
         
  Conversely, suppose that $(q_{z})_{z\in\Omega}$ is contained in $ \Hol(\Omega, S_{\ah}^{*}(\URd))$ and has 
  degree $m(z)$. Then, for any $t>0$, any compact 
  $K \subset U$, any integer $N$ and any open $\Omega'\subsubset \Omega$, we have   
  \begin{equation}   
      |\partial_{x}^\alpha\partial^\beta_{\xi}(q_{z}(x,t.\xi) -  t^{m(z)} q_{z}(x,\xi))| \leq 
     C_{tNK\Omega'\alpha\beta}  (1+\|\xi\|)^{-{N}},  
      \label{eq:HolPHDO.almost-homogeneity}
  \end{equation}
  for $(x,\xi)\in K \times\Rd$ and $z\in \Omega'$.
  Then replacing $\xi$ by $s.\xi$, 
 $s>0$, in~(\ref{eq:HolPHDO.almost-homogeneity}) shows that when $N\geq \sup_{z\in\Omega'}\Re \hat{m}(z)$ we have 
 \begin{multline}
     | \partial_{x}^\alpha\partial^\beta_{\xi}(s^{m(z)}q_{z}(x,st.\xi) - (st)^{m(z)}q_{z}(x,s.\xi))|  \leq 
        C_{tNK\Omega'\alpha\beta}   s^{\Re m(z)-N}\|\xi\|^{-{N}}\\ \leq  C_{tNK\Omega'\alpha\beta}   s^{-1}\|\xi\|^{-{N}}.
     \label{eq:HolPHDO.twisted-almost-homogeneity}
 \end{multline}
 for $(x,\xi)\in K \times\Rdo$ and $z\in \Omega'$. 
 
 Next, for $ k\in\N$ let $q_{z,k}(x,\xi)=(2^{k})^{-m}q_{z}(x,2^{k}.\xi)$. Then, for any compact 
$K \subset U$, any open $\Omega'\subsubset \Omega$ and any integer $N \geq \sup_{z\in\Omega'}\Re \hat{m}(z)$, we have
 \begin{equation}
     | \partial_{x}^\alpha\partial^\beta_{\xi}(q_{z,k+1}(x,\xi) - q_{z,k}(x,\xi))|\leq C_{2NK\Omega'\alpha\beta} 2^{-k}\|\xi\|^{-{N}},
 \end{equation}
for $(x,\xi)\in K \times\Rdo$ and $z\in \Omega'$.
This shows that the series 
 $\sum_{k \geq 0} (q_{z,k+1}-q_{z,k})$ is convergent in $\Hol(\Omega,C^\infty(\URdo))$. Hence the sequence $(q_{z,k})_{k\geq 0}$ 
 converges in $\Hol(\Omega,C^\infty(\URdo))$  to some  family $(p_{z})_{z\in \Omega}$. 
 In fact, taking $s=2^{k}$ in~(\ref{eq:HolPHDO.twisted-almost-homogeneity})  
  and letting $k\rightarrow \infty$ with $t$ fixed 
 shows that $q_{z}$ is homogeneous of degree $m(z)$ with respect to the $\xi$-variable. Moreover, for any compact 
$K \subset U$, any open $\Omega'\subsubset \Omega$ and any integer $N \geq \sup_{z\in\Omega'}\Re \hat{m}(z)$, we have 
\begin{equation}
     | \partial_{x}^\alpha\partial^\beta_{\xi}(q_{z}- p_{z})(x,\xi)| \leq \sum  |\partial_{x}^\alpha\partial^\beta_{\xi}(q_{z,k+1}-q_{z,k})(x,\xi)|  
     \leq C_{2NK\Omega'\alpha\beta} \|\xi\|^{-{N}}, 
\end{equation}
for $(x,\xi)\in K \times(\Rdo)$ and $z\in \Omega'$, i.e.~we have $q_{z}\sim p_{z}$ in the sense of~(\ref{eq:complex.symbols.asymptotic-expansion}). 
\end{proof}

Using Lemma~\ref{lem:HolPHDO.characterization-almost-homogeneous-symbols} and arguing as in the proof of 
Lemma~\ref{lem:HolPHDO.almost-homogeneous-characterization} we get the following 
characterization of holomorphic families of symbols. \begin{lemma}\label{lem:complex.symbols.criterion-almost-homogeneity}
    Let $(p_{z})_{z\in \Omega} \in \Hol(\Omega,C^{\infty}(\URd))$. Then 
    we have equivalence:\smallskip 
    
    (i) The family $(p_{z})_{z\in \Omega}$ is in $\Hol(\Omega, S^{*}(\URd))$ and has order $m(z)$.\smallskip 
    
    (ii) For $j=0,1,\ldots$ there exists $(q_{j,z})_{z\in \Omega}\in \Hol(\Omega, S_{\ah}(\URd))$ almost homogeneous of degree 
    $m(z)-j$ so that we have $p_{z}\sim \sum_{j \geq 0} q_{j,z}$ in the sense of~(\ref{eq:complex.symbols.asymptotic-expansion}).
\end{lemma}

Next, it is shown in~\cite[Sect.~12]{BG:CHM} that, as for homogeneous symbols 
in~(\ref{eq:PsiHDO.convolution-symbol-pointwise})--(\ref{eq:PsiHDO.convolution-symbols-URd}), there is a continuous bilinear product, 
\begin{equation}
    *:\Svb^{k_{1}}(\URd) \times \Svb^{k_{2}}(\URd) \longrightarrow \Svb^{k_{1}+k_{2}}(\URd),  
\end{equation}
which is homogeneous in the sense that, for any $\lambda \in \R$, we have 
\begin{equation}
    (p_{1}*p_{2})_{\lambda}=p_{1,\lambda}*p_{2,\lambda}, \qquad p_{j}\in \Svb^{k_{j}}(\URd).
     \label{eq:HolPHDO.homogeneity-*product}
\end{equation}
This product is related to the product of homogeneous symbols as follows.

\begin{lemma}[{\cite[Sect.~13]{BG:CHM}}]
   For $j=1,2$ let $p_{j}\in S^{m_{j}}(\URd)$ have principal symbol 
   $p_{m_{j}}\in S_{m_{j}}(\URd)$. Then $p_{1}*p_{1}$ lies in $S^{m_{1}+m_{2}}(\URd)$ 
   and has principal symbol $p_{m_{1}}*p_{m_{2}}$. 
\end{lemma}

Furthermore, this product is holomorphic, for we have: 

\begin{lemma}\label{lem:HolPHDO.*convolution} 
   For $j=1,2$ let $(p_{j,z})_{z \in \Omega} \subset S^{*}(\URd)$ be a holomorphic family of symbols. 
   Then $(p_{1,z}*p_{2,z})_{z \in \Omega}$ is a holomorphic family of symbols as well.  \end{lemma}
\begin{proof}
    For $j=1,2$ let $m_{j}(z)$ be the order of $p_{j,z}$. Since $m_{1}(z)$ and $m_{2}(z)$ are holomorphic functions, possibly by 
    shrinking $\Omega$, we may assume that we have $\sup_{z\in \Omega}m_{j}(z)\leq k<\infty$, so that the families $(p_{j,z})_{z\in \Omega}$ are in 
    $\Hol(\Omega, S_{\|}^{k}(\URd))$. 
    As $*$ gives rise to  a continuous $\C$-bilinear map from  
$S_{\|}^{k}(\URd)\times S_{\|}^{k}(\URd)$ to $S_{\|}^{2k}(\URd)$ we 
see that $p_{1,z}*p_{2,z}$ is in $\Hol(\Omega, S_{\|}^{2k}(\URd))$, hence is in $\Hol(\Omega, C^{\infty}(\URd))$.

Now, assume that $p_{j,z}$, $j=1,2$, is almost homogeneous of degree $m_{j}(z)$. Then 
using~(\ref{eq:HolPHDO.homogeneity-*product}) 
we see that for $\lambda>0$ the symbol $(p_{1,z}*p_{2,z})_{\lambda} -\lambda^{m_{1}(z)+m_{2}(z)} p_{1,z}*p_{2,z}$ is equal to
\begin{multline}
  p_{1,z,\lambda}*p_{1,z,\lambda} -\lambda^{m_{1}(z)+m_{2}(z)} p_{1,z}*p_{2,z}\\ 
 = (p_{1,z,\lambda}-\lambda^{m_{1}(z)}p_{1,z})*p_{2,z}+\lambda^{m_{1}(z)} p_{1,z,\lambda}*(p_{2,z,\lambda}-\lambda^{m_{2}(z)}p_{2,z}).
\end{multline}
As $(p_{1,z,\lambda}-\lambda^{m_{1}(z)}p_{1,z})_{z\in \Omega}$ and $(p_{2,z,\lambda}-\lambda^{m_{2}(z)}p_{2,z})_{z\in \Omega}$ are holomorphic 
families with values in $S^{-\infty}(\URd)$, 
combining this with the analyticity of $*$ on $\Svb^{*}(\URd)$ shows that, for any $\lambda>0$, the family   
$(p_{1,z}*p_{2,z})_{\lambda} -\lambda^{m_{1}(z)+m_{2}(z)} p_{1,z}*p_{2,z}$ belongs to 
$\Hol(\Omega,S^{-\infty}(\URd))$. Then Lemma~\ref{lem:complex.symbols.criterion-almost-homogeneity}  
 implies that $p_{1,z}*p_{2,z}$ is a holomorphic family of almost homogeneous 
of symbols of degree $m_{1}(z)+m_{2}(z)$. 

In general, by Lemma~\ref{lem:complex.symbols.criterion-almost-homogeneity} we have
$p_{j,z} \sim \sum_{l\geq 0} p_{j,z,l}$, with $(p_{j,z,l})_{z\in 
\Omega}$ in $\Hol(\Omega, S^{*}_{\ah}(\URd))$ of degree 
$m_{j}(z)-l$ and $\sim$ taken in the sense of~(\ref{eq:complex.symbols.asymptotic-expansion}). 
In particular, for any integer $N$  we have 
$ p_{j,z}$ is equal to $\sum_{l<N} p_{j,z,l}$ modulo a family in $\Hol(\Omega, S_{\|}^{k-N}(\URd))$. 
Thus, 
\begin{equation}
    p_{1,z}*p_{2,z}=\sum_{l+p<N} p_{1,z,l}*p_{2,z,p} \quad \bmod \Hol(\Omega, S_{\|}^{2k-N}(\URd)).
     \label{eq:HolPHDO.analyticity-*convolution-symbols}
\end{equation}
As  explained above $(p_{1,z,l}*p_{2,z,p})_{z\in \Omega}$ is a holomorphic family of almost homogeneous 
of symbols of degree $m_{1}(z)+m_{2}(z)-l-p$. It then follows from 
Lemma~\ref{lem:complex.symbols.criterion-almost-homogeneity} that $(p_{1,z}*p_{2,z})_{z\in \Omega}$ 
is a holomorphic family of symbols.
\end{proof}

We are now ready to prove: 
\begin{proposition}\label{prop:HolPHDO.composition}
    For $j=1,2$ let $(P_{j,z})_{z\in \Omega}$ be in 
    $\Hol(\Omega,\pvdo^{*}(U))$ and suppose that at least one the families $(P_{1,z})_{z 
    \in \Omega}$ or $(P_{2,z})_{z \in \Omega}$ is uniformly properly supported. 
    Then the family $(P_{1,z}P_{2,z})_{z  \in \Omega}$ is a holomorphic family of \psivdos. 
\end{proposition}
\begin{proof}
By assumption $(P_{1,z})_{z\in \Omega}$ or $(P_{2,z})_{z\in \Omega}$ 
is uniformly properly supported, 
so by Proposition~\ref{prop:Complexpowers.operators.properties2} it gives rise to holomorphic families with values in $\cL(C^{\infty}(U))$ and 
$\cL(\cE'(U))$.  Moreover, Proposition~\ref{prop:Complexpowers.operators.properties2} 
tells us that the other family at least coincides with a uniformly properly supported holomorphic family of \psivdos\ 
 up to a holomorphic family of smoothing operators. 
 It thus follows that $(P_{1,z}P_{2,z})_{z\in \Omega}$ is the product of two uniformly properly supported holomorphic families of 
 \psivdos\ up to a holomorphic families of smoothing operators. 

As a consequence of this we may assume that the families 
$(P_{1,z})_{z\in \Omega}$ 
and $(P_{2,z})_{z\in \Omega}$ are both uniformly properly supported. Thanks 
to~(\ref{eq:HolPHDO.uniformly-properly-supported-locally-finite-sum})   
this allows us to write 
\begin{equation}
    P_{1,z}P_{2,z}=\sum_{i \geq 0} \varphi_{i} P_{1,z}\psi_{i}P_{2,z},
    \label{eq:PsiHDO.composition-locally-finite-sum}
\end{equation}
where $(\varphi_{i})_{i \geq 0}\subset C^{\infty}_{c}(U)$ and $(\psi_{i})_{i\geq 0} \subset  
   C^{\infty}_{c}(U)$ are locally finite families 
  such that $(\varphi_{i})$ is a partition of the unity and $\psi_{i}=1$ near $\supp\varphi_{i}$. 
  
  Next, for $j=1,2$ let us write $P_{j,z}=p_{j,z}(x,-iX)+R_{j,z}$, with  $(p_{j,z})_{z\in \Omega}$ in 
  $\Hol(\Omega, S^{*}(\URd))$ and 
    $(R_{j,z})_{z\in \Omega}$ in $\Hol(\Omega, \psinf(U))$. Since by Proposition~\ref{prop:Complexpowers.operators.properties}
    each family $(p_{j,z}(x,-iX))_{z\in \Omega}$ is holomorphic with values in 
    $\cL(C_{c}^{\infty}(U), C^{\infty}(U))$ and $\cL(\cE'(U),\cD'(U))$ 
    using~(\ref{eq:PsiHDO.composition-locally-finite-sum}) we see that 
    \begin{equation}
        P_{1,z}P_{2,z}=\sum \varphi_{i} 
        p_{1,z}(x,-iX)\psi_{i}p_{2,z}(x,-iX) \quad \bmod 
        \Hol(\Omega, \Psi^{-\infty}(U)).
    \end{equation}
    
    At this stage we make appeal to: 
    
    \begin{lemma}[{\cite[Prop.~14.45]{BG:CHM}}] \label{lem:HolPHDO.proof-composition}
 For $j=1,2$ let $p_{j}\in \SvbU{k_{j}}$ and let $\psi \in 
 C_{c}^\infty(U)$. Then:
   \begin{equation}
     p_{1}(x,-iX) \psi p_{2}(x,-iX)= p_{1}\#_{\psi} p_{2}(x,-iX), 
 \end{equation}
 where $ \#_{\psi}$ is a continuous bilinear map from 
 $\SvbU{k_{1}}\times\SvbU{k_{2}}$ to $\SvbU{k_{1}+k_{2}}$. Moreover, for any integer $N\geq 1$ we have  
   \begin{equation}
      p_{1}\#_{\psi}p_{2} = \sum_{j<N}
        \sum_{\alpha\beta\gamma\delta}^{(j)} h_{\alpha\beta\gamma\delta} \psi (D_{\xi}^\delta p_{1})* (\xi^\gamma 
            \partial_{x}^\alpha \partial_{\xi}^\beta p_{2}) + R_{N, \psi}(p_{1},p_{2}),
  \end{equation}
  where the notation is the same as that of 
  Proposition~\ref{prop:PsiHDO.composition} and $R_{N, \psi}$ is a 
  continuous bilinear map from  $\SvbU{k_{1}}\times\SvbU{k_{2}}$ to $\SvbU{k_{1}+k_{2}-N}$.
\end{lemma}
\begin{remark}
 The continuity content of Lemma~\ref{lem:HolPHDO.proof-composition} is not explicitly stated in Proposition~14.45 of \cite{BG:CHM}, but it follow 
 from its proof or from a standard use of the closed graph theorem.   
\end{remark}

Now, thanks to Lemma~\ref{lem:HolPHDO.proof-composition} we have
\begin{equation}
    P_{1,z}P_{2,z}=p_{z}(x,-iX)+R_{z}, \qquad p_{z}= \sum \varphi_{i}p_{1,z}\#_{\psi_{i}}p_{2,z}. 
     \label{eq:HolPHDO.composition-PsiHDO}
\end{equation}
Furthermore, possibly by shrinking $\Omega$, we may assume that there is a real $k$ such that 
for $j=1,2$ we have $\Re \ord p_{j,z}\leq k$ for any $z \in \Omega$. Then 
the continuity contents of Lemma~\ref{lem:HolPHDO.proof-composition} 
imply that $(p_{z})_{z\in \Omega}$ belongs to $\Hol(\Omega, \Svb^{2k}(\URd))$ 
and for any integer $N\geq 1$ we can write
\begin{equation}
       p_{z} =  \sum_{r<N}  q_{r,z} +R_{N,z}, \quad 
        q_{j,z}= \overset{(r-s-t)}{\underset{\alpha\beta\gamma\delta}{\sum}}
h_{\alpha\beta\gamma\delta} (D_{\xi}^\delta p_{1,z})* (\xi^\gamma 
            \partial_{x}^\alpha \partial_{\xi}^\beta p_{2,z}),
\end{equation}
where $(R_{N,z})_{z\in \Omega}:=(\sum_{i\geq 0}\varphi_{i} R_{N, 
\psi_{i}}(p_{1,z},p_{2,z}))_{z\in \Omega}$ is in $\Hol(\Omega,S_{\|}^{2k-N}(\URd))$. 
Thanks to Lemma~\ref{lem:HolPHDO.*convolution} 
            the family $(q_{r,z})_{z\in \Omega}$ is in $\Hol(\Omega, S^{*}(\URd))$  
            and has order $m_{1}(z)+m_{2}(z)-r$. Therefore, we have $p_{z} \sim \sum_{j\geq 0} q_{j,z}$ 
            in the sense of~(\ref{eq:complex.symbols.asymptotic-expansion}), 
            which by Lemma~\ref{lem:HolPHDO.*convolution} implies that $(p_{z})_{z\in \Omega}$ 
            belongs to $\Hol(\Omega, S^{*}(\URd))$.
           Combining this with~(\ref{eq:HolPHDO.composition-PsiHDO}) 
           then shows that $(P_{1,z}P_{2,z})_{z \in \Omega}$ is a holomorphic family of \psivdos. \end{proof}

\section{Kernel characterization of holomorphic families of $\mathbf{\Psi}$DO's}
We shall now give a characterization of holomorphic families of \psivdos\ in terms of holomorphic families with values in 
$\cK^{*}(\URd)$. Since the 
latter is defined in terms of asymptotic expansions involving distributions in 
$\cK_{*}(\URd)$ a difficulty occurs, since for a family $(K_{z})_{z\in 
\Omega}\subset \cK_{*}(\URd)$ logarithmic singularities may appear as 
the order of $K_{z}$ crosses non-negative integer values. 

This issue is resolved by making use of holomorphic families of almost homogeneous distributions as follows.
\begin{definition}
   A family $(K_{z})_{z \in \Omega}\subset \cK^{*}_{\ah}(\URd)$ is holomorphic when:\smallskip 
   
   (i) The degree $m(z)$ of $K_{z}$ is a holomorphic function on $\Omega$;\smallskip 
   
   (ii) The family $(K_{z})_{z\in \Omega}$ belongs to $\Hol(\Omega,C^{\infty}(U)\otimes \cD'_{\reg}(\Rd))$;\smallskip 
  
   (iii) For any $\lambda >0$ the family $\{K_{z}(x,\lambda.y)-\lambda^{m(z)}K_{z}(x,y)\}_{z \in \Omega}$ is a holomorphic family with values in 
   $C^{\infty}(\URd)$.\smallskip 
   
We let $\Hol(\Omega,\cK^{*}_{\ah}(\URd))$ denote the set of holomorphic 
families of almost homogeneous distributions.
\end{definition}

\begin{lemma}\label{lem:HolPHDO.almost-homogeneous-kernel-symbol}
    Let $(K_{z})_{z\in \Omega}\in \Hol(\Omega, C^{\infty}(U)\otimes 
    \cD'_{\reg}(\Rd))$.
    Then we have equivalence:\smallskip 
    
    (i) The family $(K_{z})_{z \in \Omega}$ belongs to $\Hol(\Omega, \cK_{\ah}^{*}(\URd))$ and has degree $m(z)$.\smallskip 
    
    (ii) We can put $(K_{z})_{z \in \Omega}$ into the form,
    \begin{equation}
        K_{z}(x,y)=\check{p}_{z,\xiy}(x,y) +R_{z}(x,y), \qquad z\in \Omega,
    \end{equation}
    with $(p_{z})_{z \in \Omega}$ in 
    $\Hol(\Omega,S^{*}_{\ah}(\URd))$ of degree $\hat{m}(z):=-(m(z)+d+2)$ and  $(R_{z})_{z\in \Omega}$ in $\Hol(\Omega, 
    C^{\infty}(\URd))$.   
\end{lemma}
\begin{proof}
 Assume that $(K_{z})_{z \in \Omega}$ belongs to $\Hol(\Omega, \cK_{\ah}^{*}(\URd))$ and has degree 
 $m(z)$.  Let $\varphi \in 
 C^{\infty}_{c}(\URd)$ be such that $\varphi(y)=1$ near $y=0$ and define $p_{z}=(\varphi(y)K_{z}(x,y))_{\yxi}$. 
 As $(\varphi(y)K_{z}(x,y))_{z\in \Omega}$ is in $\Hol(\Omega, C^{\infty}(U)\otimes\cE'_{L}(\Rd))$ with $L=\supp \varphi$, 
 we see that $(p_{z})_{z\in \Omega}$ belongs to $\Hol(\Omega,C^{\infty}(\URd))$. Moreover, since 
 $((1-\varphi(y))K_{z}(x,y))_{z\in \Omega}$ is in $\Hol(\Omega, C^{\infty}(\URd))$ we also see that $K_{z}(x,y)=\check{p}_{z,\xiy}(x,y)$  
%  Thus  $K_{z}(x,y)=\check{p}_{z,\xiy}(x,y)  +(1-\varphi(y))K_{z}(x,y) = \check{p}_{z,\xiy}(x,y)$ 
 modulo $\Hol(\Omega, C^{\infty}(\URd))$. 
%  \end{equation}

On the other hand, it follows from~(\ref{eq:PsiHDO.dilations-Fourier-transform}) that for any $\lambda>0$ the symbol 
$\lambda^{-(d+2)}(p_{z}(x,\lambda.\xi)-\lambda^{\hat{m}(z)}p_{z}(x,\xi))$ is 
 equal to 
\begin{equation}
         [(\varphi(\lambda^{-1}.y)-\varphi(y))K_{z}(x,\lambda^{-1}.y)+
        \varphi(y) (K_{z}(x,\lambda^{-1}.y)-\lambda^{-m}K_{z}(x,y))]_{\yxi}^{\wedge}.
%     \label{eq:}
\end{equation}
This is  the Fourier transform with respect to $y$ of a 
    holomorphic family with values in $C_{c}^{\infty}(\URd)$, so it belongs to $ \Hol(\Omega, S^{-\infty}(\URd))$. 
% 
%     \begin{multline}
%     =\\ \lambda^{-(d+2)}
%      [(\varphi(\lambda^{-1}.y)-\varphi(y))K_{z}(x,\lambda^{-1}.y)+
%     \varphi(y) (K_{z}(x,\lambda^{-1}.y)-\lambda^{-m}K_{z}(x,y))]_{\yxi}^{\wedge}(x,\xi).
%  \end{multline}
% Since the r.h.s.~above iswe see that $(p_{z}(x,\lambda.\xi)-\lambda^{\hat{m}(z)}p_{z}(x,\xi))_{z \in \Omega}$ is in 
%     $\Hol(\Omega, S^{-\infty}(\URd))$ for any $\lambda>0$. 
  
    Combining all this with Lemma~\ref{lem:HolPHDO.characterization-almost-homogeneous-symbols}
    shows that 
    $(p_{z})_{z\in \Omega}$ is a family of almost homogeneous symbols of degree $\hat{m}(z)$.     
 
    Conversely, let $(p_{z})_{z \in \Omega}\in  \Hol(\Omega, S^{*}_{\ah}(\URd))$ have degree 
    $\hat{m}(z)$. Possibly by shrinking $\Omega$ we may assume that we have $\sup_{z \in \Omega}\Re \hat{m}(z)\leq k<\infty$. 
    Then the family $(p_{z})_{z \in \Omega}$ belongs to $\Hol(\Omega, \Svb^{k}(\URd))$, 
    hence to $\Hol(\Omega, S_{\frac{1}{2}\frac{1}{2}}^{k}(\URd))$. As mentioned in  the proof of 
    Proposition~\ref{prop:Complexpowers.operators.properties}, the map $q(x,\xi) \rightarrow \check{q}_{\xiy}(x,y)$ is analytic from 
    $S_{\frac{1}{2}\frac{1}{2}}^{k}(\URd)$ to $C^{\infty}(U)\otimes \cD'_{\reg}(\Rd)$, so 
    $(\check{p}_{z,\xiy})_{z\in \Omega}$ is a holomorphic family with values in $C^{\infty}(U)\hotimes \cD'_{\reg}(\Rd)$. 
    
    Next, using~(\ref{eq:PsiHDO.dilations-Fourier-transform}) we see that, for any $\lambda>0$, we have
    \begin{equation}
        \check{p}_{z,\xiy}(x,\lambda.y)-\lambda^{m(z)}\check{p}_{z,\xiy}(x,y)= 
        [p_{z}(x,\lambda^{-1}.\xi)-\lambda^{-\hat{m}(z)}p_{z}(x, \xi)]_{\xiy}^{\wedge}.
         \label{eq:HolPHDO.almost-homogeneity-checkp}
    \end{equation}
    Since by Lemma~\ref{lem:HolPHDO.characterization-almost-homogeneous-symbols}
    the r.h.s.~of~(\ref{eq:HolPHDO.almost-homogeneity-checkp}) is the inverse Fourier transform with respect to $\xi$ of an element 
    of $\Hol(\Omega,S^{-\infty}(\URd))$, we see that the family
    $\{\check{p}_{z,\xiy}(x,\lambda.y)-\lambda^{m(z)}\check{p}_{z,\xiy}(x,y)\}_{z \in \Omega}$ is contained in $\Hol(\Omega, C^{\infty}(\URd))$. 
    It then follows that $(\check{p}_{z,\xiy})_{z \in \Omega}$ is a holomorphic family of almost homogeneous distributions of degree $m(z)$. 
\end{proof}

\begin{definition}\label{def:HolPHDO.kernels-families}
   A family $(K_{z})_{z\in \Omega}\subset \cK^{*}(\URd)$ is holomorphic when:\smallskip 
    
    (i) The order $m_{z}$ of $K_{z}$ is a 
    holomorphic function of $z$;\smallskip 
    
   (ii) For $j=0,1,..$ there exists $(K_{j,z})\in \Hol (\Omega, \cK_{\ah}^{*}(\URd))$ of degree 
   $m(z)+j$ such that $K_{z}\sim \sum_{j \geq 0} K_{j,z}$  in the sense that, for any open 
    $\Omega' \subsubset \Omega$ and any integer $N$, as soon as $J$ is large enough we have 
 \begin{equation}
         K_{z}- \sum_{j\leq J} K_{z,m_{z}+j}\in \Hol(\Omega',C^{N}(\URd)).
                      \label{eq:HolPHDO.kernel-asymptotic}
 \end{equation}
\end{definition}\begin{proposition}\label{prop:HolPHDO.characterization-Hol-cK*}
   For a family  $(K_{z})_{z\in \Omega}\subset \cK^{*}(\URd)$ the following are equivalent:\smallskip
   
   (i) The family   $(K_{z})_{z\in \Omega}$ is holomorphic and has order $m(z)$.\smallskip
   
   (ii) We can put $(K_{z})_{z\in \Omega}$ into the form, 
   \begin{equation}
       K_{z}(x,y)=(p_{z})^{\vee}_{\xiy}(x,y)+R_{z}(x,y),
       \label{eq:HolPHDO.kernel-symbol}
   \end{equation}
  for some family $(p_{z})_{z \in \Omega}\in \Hol(\Omega,S^{*}(\URd))$ of order $\hat{m}(z):=-(m(z)+d+2)$ and 
  some family $(R_{z})_{z\in \Omega}\in \Hol(\Omega, C^{\infty}(\URd))$.   
   \end{proposition}
 \begin{proof}
     Assume that $(K_{z})_{z\in \Omega}$ is in
     $\Hol(\Omega, \cK^{*}(\URd))$. Let $\varphi\in C^{\infty}(\Rd)$ be such that 
    $\varphi(y)=1$ near $y=0$ and for $z\in \Omega$ set $p_{z}=(\varphi(y)K_{z}(x,y))^{\wedge}_{\yxi}$. Since $(K_{z})_{z\in \Omega}$ lies in 
    $\Hol(\Omega, C^{\infty}(U)\hotimes \cD_{\reg}'(\Rd))$ we see that $(p_{z})_{z\in \Omega}$ belongs to $\Hol(\Omega, C^{\infty}(\URd))$ and we 
    have 
      \begin{multline}
       K_{z}(x,y)=(p_{z})^{\vee}_{\xiy}(x,y)+(1-\varphi(y))K_{z}(x,y)\\ =(p_{z})^{\vee}_{\xiy}(x,y) \quad \bmod \Hol(\Omega,C^{\infty}(U\times U)).
       \label{eq:HolPHDO.kernel-symbol2}
   \end{multline}
    
   Let us write $K_{z}\sim \sum_{j\geq 0}K_{j,z}$ with 
  $(K_{j,z})_{z \in \Omega}\in \Hol(\Omega, \cK_{\ah}^{*}(\URd))$ of 
  degree $\hat{m}(z)+j$ and $\sim$ taken in the sense of~(\ref{eq:HolPHDO.kernel-asymptotic}) and for $j=0,1,\ldots$ define  
  $p_{j,z}=(\varphi(y)K_{j,z})_{z \in \Omega}$. Then arguing as in the proof of Lemma~\ref{lem:HolPHDO.almost-homogeneous-kernel-symbol} shows 
  that $(p_{j,z})_{z\in \Omega}$ is a family in $\Hol(\Omega, 
  S^{*}_{\ah}(\URd))$ of degree $\hat{m}(z)-j$. 
  
  Next, in the sense of~(\ref{eq:HolPHDO.kernel-asymptotic}) we have  $(p_{z})^{\vee}_{\xiy}(x,y)\sim \sum_{j \geq 0}(p_{j,z})^{\vee}_{\xiy}(x,y)$.
 Under the Fourier transform with respect to $y$ this shows  that, 
 for any compact $L\subset U$, any integer $N$ and any open $\Omega' \subsubset \Omega$, as soon as $J$ is large enough 
 we have estimates,
     \begin{equation}
         |\partial_{x}^{\alpha}\partial_{\xi}^{\beta}(p_{z}- \sum_{j\leq J} (\varphi(y) K_{j,z}(x,y))_{\yxi}^{\wedge})(x,\xi)| \leq 
         C_{\Omega'NJL\alpha\beta}(1+|\xi|^{2})^{-[N/2]},
          \label{eq:HolPHDO.kernel-characterization.asymptotic-expansion-symbol}
     \end{equation}
 for $(x,\xi)\in L\times\Rd$ and $z \in \Omega'$. Hence $p_{z}\sim \sum_{j\geq 0}p_{j,z}$ in the sense 
 of~(\ref{eq:complex.symbols.asymptotic-expansion}).  It thus follows that
 $(p_{z})_{z\in\Omega}$ is in $\Hol(\Omega, S^{*}(\URd))$ and has order $\hat{m}(z)$, so using~(\ref{eq:HolPHDO.kernel-symbol2}) we see that the family 
 $(K_{z})_{z\in 
 \Omega}$ is of the form~(\ref{eq:HolPHDO.kernel-symbol}). 
 
      Conversely, assume that $K_{z}(x,y)=(p_{z})^{\vee}_{\xiy}(x,y)+R_{z}(x,y)$ with 
      $(p_{z})_{z \in \Omega}$ in $\Hol(\Omega,S^{*}(\URd))$ of order $\hat{m}(z)$ and 
  $(R_{z})_{z\in \Omega}$ in $\Hol(\Omega, C^{\infty}(\URd))$. 
      By Lemma~\ref{lem:complex.symbols.criterion-almost-homogeneity} 
    we have $p_{z}\sim \sum_{j\geq 0}p_{j,z}$ with $(p_{j,z})_{z\in \Omega}\in 
    \Hol(\Omega,S^{*}_{\ah}(\URd))$ of degree $m(z)-j$ and $\sim$ taken in the sense 
    of~(\ref{eq:complex.symbols.asymptotic-expansion}). Thus, under 
    the inverse Fourier transform with respect to $\xi$, we have $(p_{z})^{\vee}_{\xiy}\sim \sum_{j\geq 
    0}(p_{j,z})^{\vee}_{\xiy}$ in the sense of~(\ref{eq:HolPHDO.kernel-asymptotic}). 
    As Lemma~\ref{lem:HolPHDO.almost-homogeneous-kernel-symbol} tells us that  
    $((p_{j,z})^{\vee}_{\xiy})_{z \in \Omega}$ is a holomorphic family of almost 
    homogeneous distributions of degree $\hat{m}(z)+j$, it follows that 
    $((p_{z})^{\vee}_{\xiy})_{z\in \Omega}$ belongs to $\Hol(\Omega,\cK^{*}(\URd))$ and has order $\hat{m}(z)$. Since $(K_{z})_{z \in \Omega}$ agrees with 
    $((p_{z})^{\vee}_{\xiy})_{z\in \Omega}$ up to an element of $\Hol(\Omega,C^{\infty}(U\times U))$, we see that $(K_{z})_{z \in \Omega}$ is a 
    holomorphic family with values in $\cK^{*}(\URd)$.
\end{proof}
   We are now ready to prove the kernel characterization of holomorphic families of \psivdos. As before for $x\in U$ we let $\psi_{x}$ and 
   $\varepsilon_{x}$ respectively denote the coordinate changes to the privileged coordinates and Heisenberg coordinates at $x$.  

\begin{proposition}\label{prop:HolPHDO.kernel.characterization}
 Let $(P_{z})_{z\in \Omega} \in \Hol(\Omega,\cL(C_{c}^{\infty}(U),C^{\infty}(U)))$ have distribution kernel $k_{P_{z}}(x,y)$.  Then the following are 
 equivalent:\smallskip 
 
 (i) The family $(P_{z})_{z\in \Omega}$ is a holomorphic family of \psivdos\ of order $m(z)$.\smallskip 
 
 (ii) We can put $k_{P_{z}}(x,y)$ in the form, 
%  
%  There exist families  
%  such that the distribution kernel $k_{P_{z}}(x,y)$ of $P_{z}$ is of the form
  \begin{equation}
     k_{P_{z}}(x,y) = |\psi_{x}'| P_{z}(x,-\varepsilon_{x}(y)) + R_{z}(x,y),
     \label{eq:HolPHDO.kernel.characterization-psix}
 \end{equation}
with $(K_{z})_{z\in \Omega}$ in $\Hol(\Omega,\cK^{*}(\URd))$ of order $\hat{m}(z):=-(m(z)+d+2)$ 
  and $(R_{z})_{z\in \Omega}$ in $\Hol(\Omega,C^{\infty}(U\times U))$.\smallskip   

 (iii) We can put $k_{P_{z}}(x,y)$ in the form,
%  There exist $(K_{P_{z}})_{z\in \Omega}\in \Hol(\Omega,\cK^{*}(\URd))$ of order $\hat{m}(z):=-(m(z)+d+2)$ 
%   and $(R_{z})_{z\in \Omega}\in \Hol(\Omega,C^{\infty}(U\times U))$ 
%  such that the distribution kernel $k_{P_{z}}(x,y)$ of $P_{z}$ is of the form
  \begin{equation}
     k_{P_{z}}(x,y) = |\varepsilon_{x}'| P_{z}(x,-\varepsilon_{x}(y)) + R_{z}(x,y).
     \label{eq:HolPHDO.kernel.characterization}
 \end{equation}
with $(K_{P_{z}})_{z\in \Omega}$ in $\Hol(\Omega,\cK^{*}(\URd))$ of order $\hat{m}(z):=-(m(z)+d+2)$ 
  and $(R_{z})_{z\in \Omega}$ in $\Hol(\Omega,C^{\infty}(U\times U))$.\smallskip   
 \end{proposition}
\begin{proof}
   First, it follows from~(\ref{eq:PsiHDO.kernel-quantization-symbol-psiy}) and 
   Proposition~\ref{prop:HolPHDO.characterization-Hol-cK*} that (i) and (ii) are equivalent. 
   
   Next, for $x\in U$ let $\phi_{x}$ denote the transition map from the privileged coordinates at $x$ to the Heisenberg coordinates at $x$. 
Recall that this gives rise to an action on distributions on $\URd$ given by
   \begin{equation}
     K(x,y) \longrightarrow \phi_{x}^{*}K(x,y), \qquad \phi_{x}^{*}K(x,y)=K(x,\phi_{x}^{-1}(y)).
        \label{eq:HolPHDO.action-phix*}
   \end{equation}

   Since $\phi_{x}$ depends smoothly on $x$, this action gives rise to continuous linear isomorphisms of $C^{N}(\URd)$, $N \geq 0$, and $C^{\infty}(\URd)$ 
  onto themselves, hence to analytic isomorphisms. Moreover, since $\phi_{x}(0)=0$  this also yields an analytic isomorphism of $C^{\infty}(U)\hotimes 
  \cD_{\reg}'(\Rd)$ onto itself. Combining this with the homogeneity property~(\ref{eq:PsiHDO.homogeneity.phix*}) 
  we then deduce that~(\ref{eq:HolPHDO.action-phix*}) induces linear isomorphisms of 
  $\Hol(\Omega,\cK^{*}_{\ah}(\URd))$ and  $\Hol(\Omega,\cK^{*}(\URd))$ onto themselves. Together 
  with~(\ref{eq:PsiHDO.kernel-quantization-symbol}) this shows that the statements 
  (ii) and (iii) are equivalent.  
\end{proof}

\section{Holomorphic families of $\mathbf{\Psi}$DO's on a general Heisenberg manifold}
 Let us now define holomorphic families of \psivdos\ a general Heisenberg manifold. To this end we will need the following lemma.

\begin{lemma}\label{lem:HolPHDO.invariance}
    Let $(K_{z})_{z\in \Omega}\in \Hol(\Omega, \cK*(\URd))$ and assume there exists an integer $N$ such that $\inf_{z\in \Omega}\Re \ord K_{z}\geq 
    2N$.   Then the family $(K_{z})_{z\in \Omega}$ is contained in $\Hol(\Omega, C^{N}(\URd))$. 
\end{lemma}
\begin{proof}
    Thanks to Proposition~\ref{prop:HolPHDO.kernel.characterization} we may assume that $K_{z}$ is of the form 
    $K_{z}(x,y)=\check{p}_{z,\xiy}(x,y)$ with $(p_{z})_{z\in \Omega}$ in $\Hol(\Omega, 
    S^{*}(\URd))$.  Since by assumption we have $-(\Re \ord p_{z}+d+2)= \Re \ord K_{z}\geq 2N$ we see that $(p_{z})_{z\in \Omega}$ is contained in $\Hol(\Omega, 
    \Svb^{-(2N+d+2)}(\URd))$. Since the map 
         $p\rightarrow \check{p}_{\xiy}$ is continuous from $S_{||}^{-(2N+d+2)}(\URd)$ to $C^{\infty}(U)\hotimes C^{N}(\URd)$ 
         (see~Lemma~\ref{lem:Appendix.Heisenberg.invariance} in Appendix~\ref{chap.Appendix-Invariance}), 
    it follows that $(K_{z})_{z \in \Omega}$ lies in $\Hol(\Omega, C^{N}(\URd))$. 
\end{proof}

% \begin{lemma}\label{lem:Appendix.Heisenberg.invariance}
%   1)  Assume $k<-(d+2)$ and set $N=[-\frac{k+d+2}2]$. Then the map 
%          $p\rightarrow \check{p}_{\xiy}$ is continuous from $S_{||}^{k}(\URd)$ to $C^{\infty}(U)\hotimes C^{N}(\URd)$.  \smallskip 
%          
%          2) For $\Re m>0$ we have $\cK^{m}(\URd)\subset C^{\infty}(U)\hotimes C^{[\frac{\Re m}{2}]}(\Rd)$.  
% \end{lemma}
% \begin{proof}  
%      Let $\alpha$ be a multi-order such that $|\alpha|\leq N$. Then we have $\brak\alpha \leq 2|\alpha|\leq 
% -(k+d+2)$, so the multiplication by $\xi^{\alpha}$ maps continuously $S_{||}^{k}(\URd)$ to $C^{\infty}(U)\hotimes 
% L^{1}(\Rd)$. Composing it with the inverse Fourier transform with respect to $\xi$ then shows that the map 
% $p\rightarrow \partial_{y}^{\alpha}\check{p}_{\xiy}$ is continuous from $S_{||}^{k}(\URd)$ to $C^{\infty}(U)\hotimes 
% C^{0}(\Rd)$. Thus the map $p\rightarrow \check{p}_{\xiy}$  is continuous from $S_{||}^{k}(\URd)$ to 
% $C^{\infty}(U)\hotimes C^{N}(\Rd)$.\smallskip  
%      
%      2) Let $K \in \cK^{m}(\URd)$. Then by Lemma~\ref{lem:PsiHDO.characterization.Km} 
%      there exists $p\in S^{\hat{m}}(\URd)$, $\hat{m}=-(m+d+2)$, such that $K(x,y)$ is equal to $\check{p}_{\xiy}(x,y)$ modulo 
%      a smooth function. Since $S^{\hat m}(\URd) \subset S_{||}^{\Re \hat m}(\URd)$ it follows from the first part that $K$ is in 
%      $C^{\infty}(U)\hotimes C^{[\frac{\Re m}{2}]}(\Rd)$. Thus $\cK^{m}(\URd)$ is contained in 
%      $C^{\infty}(U)\hotimes C^{[\frac{\Re m}{2}]}(\Rd)$. 
%  \end{proof}

\begin{proposition}\label{prop:HolPHDO.invariance}
Let $\tilde{U}$ be an open subset of $\Rd$ together with a hyperplane bundle $\tilde{H}\subset T\tilde{U}$ and a $\tilde{H}$-frame of $T\tilde{U}$ 
and let $\phi:(U,H)\rightarrow (\tilde{U},\tilde{H})$ be a Heisenberg diffeomorphism. Then for any family 
$(\tilde{P}_{z})_{z \in \Omega}\in \Hol(\Omega, \Psi_{\tilde{H}}^{*}(\tilde{U}))$ the family 
$(P_{z})_{z \in \Omega}:=(\phi^{*}\tilde{P}_{z})_{z \in \Omega}$ is contained in $\Hol(\Omega, \Psi_{H}^{*}(U))$.
\end{proposition}
\begin{proof}
For $x \in U$ and $\tilde{x}\in \tilde{U}$ let $\varepsilon_{x}$ and $\tilde{\varepsilon}_{\tilde{x}}$ 
denote the coordinate changes to the Heisenberg coordinates at $x$ and $\tilde{x}$ respectively. Then
by Proposition~\ref{prop:HolPHDO.kernel.characterization} the distribution kernel $k_{\tilde{P}_{z}}(\tilde{x},\tilde{y})$ of $\tilde{P}_{z}$ 
is of the form
  \begin{equation}
      k_{\tilde{P}_{z}}(\tilde{x},\tilde{y}) = |\tilde{\varepsilon}_{\tilde{x}}'|K_{\tilde{P}_{z}}(\tilde{x},-\tilde{\varepsilon}_{\tilde{x}}(\tilde{y})) + 
      \tilde{R}_{z}(\tilde{x},\tilde{y}),
  \end{equation}
with $(K_{\tilde{P}_{z}})_{z\in \Omega}$ in  $\Hol(\Omega, \cK^{*}(\tilde{U}\times \tilde{U}))$ and $(\tilde{R}_{z})_{z\in \Omega}$ in $\Hol(\Omega, 
C^{\infty}(\tilde{U}\times\Rd))$. Then the proof of Proposition~\ref{prop:PsiHDO.invariance} in Appendix~\ref{chap.Appendix-Invariance} 
shows that the distribution kernel 
$k_{P_{z}}(x,y)$ of $P_{z}$ takes the form, 
\begin{multline}
    k_{P_{z}}(x,y)= |\varepsilon_{x}'| K_{P_{z}}(x,-\varepsilon_{x}(y)) \\
    + (1-\chi(x,-\varepsilon_{x}(y))) |\tilde{\varepsilon}_{\phi(x)}'| 
    K_{\tilde{P}_{z}}(\phi(x),-\tilde{\varepsilon}_{\phi(x)}(\phi(y))) +\tilde{R}_{z}(\phi(x),\phi(y)),
\end{multline}
where we have let
\begin{gather}
    K_{P_{z}}(x,y)= \chi(x,y) |\partial_{y}\Phi(x,y)| K_{\tilde{P}_{z}}(\phi(x),\Phi(x,y)),\\  \Phi(x,y)=-\tilde{\varepsilon}_{\phi(x)}\circ \phi\circ 
    \varepsilon_{x}^{-1}(-y),
    \label{eq:HolPHDO.kernel-characterization.KPz}
\end{gather}
and the function $\chi(x,y)\in C^{\infty}(\URd)$ has a supported contained in the open $\cU=\{ (x,y)\in \URd; \  \varepsilon_{x}^{-1}(-y)\in U\}$, is properly 
supported with respect to $x$ and is such that $\chi(x,y)=1$ near $U\times \{0\}$. In particular, we have
\begin{equation}
    k_{P_{z}}(x,y)= |\varepsilon_{x}'| K_{P_{z}}(x,-\varepsilon_{x}(y)) \quad \bmod \Hol(\Omega,C^{\infty}(U\times U)).
      \label{eq:HolPHDO.invariance.kPz}
\end{equation}

Let us now prove that $(K_{P_{z}})_{z \in \Omega}$ belongs to $\Hol(\Omega, \cK^{*}(\URd))$. To this end, possibly by shrinking 
$\Omega$, we may assume that $\inf_{z \in \Omega} \Re \hat{m}(z)\geq \mu >-\infty$.
Moreover, the proof of Proposition~\ref{prop:PsiHDO.invariance} in Appendix~\ref{chap.Appendix-Invariance} also shows that for any integer $N$  we have 
\begin{gather}
    K_{P_{z}}(x,y)= \sum_{\brak \alpha<N} \sum_{\frac32\brak\alpha \leq \brak\beta < \frac{3}{2}N} 
    K_{\alpha\beta, z}(x,y)+ \sum_{j=1}^{3} R_{N, z}^{(j)}(x,y),
    \label{eq:HolPHDO.invariance-KPz}\\
    K_{\alpha\beta, z}(x,y)=a_{\alpha\beta}(x) y^\beta 
    (\partial^\alpha_{\tilde{y}}K_{\tilde{P}_{z}})(\phi(x),\phi_{H}'(x)y),
\end{gather}
where the smooth functions $a_{\alpha\beta}(x)$ are as in Proposition~\ref{prop:PsiHDO.invariance} and the remainder terms 
   $R_{N,z}^{(j)}(x,y)$, $j=1,2,3$, take the forms:\smallskip 
   
     - $R_{N, z}^{(2)}(x,y)= \sum_{\brak \alpha<N} \sum_{\brak\beta \dot{=} \frac{3}{2}N} r_{M\alpha}(x,y) y^{\beta} 
    (\partial^\alpha_{\tilde{y}}K_{\tilde{P}_{z}})(\phi(x),\phi_{H}'(x)y)$ for some functions $r_{\alpha\beta}(x,y)$ in 
    $C^{\infty}(\URd)$;\smallskip
    
    - $R_{N, z}^{(3)}(x,y)= \sum_{\brak\alpha=N}\sum_{\brak\beta\dot{=}\frac{3}{2}N} \int_{0}^{1} 
    r_{\alpha\beta}(t,x,y)(\tilde{y}^{\beta}\partial^\alpha_{\tilde{y}}K_{\tilde{P}_{z}})(\phi(x),t\Phi(x,y)+(1-t)\phi_{H}'(x)y)dt$, 
    for some functions $r_{\alpha\beta}(t,x,y)$ in $C^{\infty}([0,1]\times \URd)$;\smallskip
    
  - $R_{N, z}^{(3)}(x,y)= \sum_{\brak \alpha<N} \sum_{\frac32\brak\alpha \leq \brak\beta < \frac{3}{2}N} (1-\chi(x,y))K_{\alpha\beta, z}(x,y)$.\smallskip 
    
Observe that the map $\Phi(x,y)=(\phi(x),\phi'_{H}(x)y)$ is a smooth diffeomorphism from $\URd$ onto $\tilde{U}\times \Rd$ such that $\Phi(x,0)=(\phi(x),0)$ 
and for any $\lambda \in \R$ we have $\Phi(x,\lambda.y)=(\phi(x),\lambda.\phi'_{H}(x)y)$, so along similar lines as that of the 
proof of Proposition~\ref{prop:HolPHDO.kernel.characterization} 
we can prove that the map  
\begin{equation}
 \cD'(\URd) \ni K(x,y)\longrightarrow K(\phi(x),\phi'_{H}(x)y) \in \cD'(\tilde{U}\times \Rd)
 \end{equation}
gives rise to a linear map from $\Hol(\Omega, \cK^{*}(\URd))$ to $\Hol(\Omega, \cK^{*}(\tilde{U}\times \Rd))$ preserving the order. Therefore 
$(K_{\alpha\beta, z}(x,y))_{z \in \Omega}$ is in $\Hol(\Omega, \cK^{*}(\tilde{U}\times \Rd))$ and has order 
$\hat{m}(z)+\brak\beta-\brak\alpha$.  Incidentally $(R_{N, z}^{(3)})_{z \in \Omega}$ 
belongs to $\Hol(\Omega, C^{\infty}(\URd))$. 
% _{z \in 
% \Omega}\in \Hol(\Omega, \cK^{*}(\URd))
    
On the other hand, if $\frac32\brak\alpha \leq \brak\beta  \dot{=} \frac{3}{2}N$ then the order 
$\hat{m}_{\alpha\beta}(z)=\hat{m}(z)+\brak\beta-\brak\alpha$ of $\tilde{y}^{\beta}\partial^\alpha_{\tilde{y}}K_{\tilde{P}_{z}}$ is such that   
$ \Re\hat{m}_{\alpha\beta}(z)\geq   \Re \hat{m}(z)+\frac{1}{3}\brak\beta\geq \mu +\frac{N}{2}$. Therefore, 
Lemma~\ref{lem:HolPHDO.invariance} implies that, for any integer $J$, as soon as $N$ is  large enough 
$(\tilde{y}^{\beta}\partial^\alpha_{\tilde{y}}K_{\tilde{P}_{z}})_{z \in \Omega}$ is in $\Hol(\Omega, C^{J}(\tilde{U}\times \Rd))$, so that 
remainder terms $(R_{N, z}^{(2)})_{z \in \Omega}$ and $(R_{N, z}^{(3)})_{z \in \Omega}$ are in $\Hol(\Omega, C^{J}(\URd))$.
   
All this shows that  we have 
$K_{P_{z}}(x,y)\sim \sum_{\frac32\brak\alpha \leq \brak\beta} K_{\alpha\beta, z}(x,y)$ in the sense of~(\ref{eq:HolPHDO.kernel-asymptotic}), 
which implies that $(K_{P_{z}})_{z \in \Omega}$ belongs to $\Hol(\Omega, \cK^{*}(\URd))$. Combining this 
    with~(\ref{eq:HolPHDO.invariance.kPz}) and Proposition~\ref{prop:HolPHDO.kernel.characterization}
    then shows that $(P_{z})_{z\in \Omega}$ is a holomorphic family of \psivdos.
\end{proof}

Now, let $(M^{d+1},H)$ be a Heisenberg manifold and let $\cE$ be a smooth vector bundle over $M$. Then Proposition~\ref{prop:HolPHDO.invariance}
allows us to define 
holomorphic families with values in $\Psi^{*}_{H}(M,\cE)$ as follows.
\begin{definition}
 A  family $(P_{z})_{z\in \Omega}\subset \pvdo^{*}(M,\cE)$ is holomorphic when:\smallskip 
 
     (i) The order $m(z)$ of $P_{z}$ is a holomorphic function of $z$; \smallskip 
 
     (ii)  For $\varphi$ and  $\psi$ in $C^\infty_{c}(M)$ with disjoint supports $(\varphi P_{z}\psi)_{z\in \Omega}$ is 
     a holomorphic family of smoothing operators, i.e., is given by a holomorphic family of smooth distribution kernels); \smallskip 
 
     (iii) For any trivialization $\tau:\cE_{|_{U}}\rightarrow U\times \C^{r}$ over a 
     local Heisenberg chart $\kappa:U \rightarrow V\subset \Rd$ the family $(\kappa_{*}\tau_{*}(P_{z|_{U}}))_{z\in\Omega}$ belongs to 
     $\Hol(\Omega, \pvdo^{*}(V,\C^{r})):=\Hol(\Omega, \pvdo^{*}(V))\otimes \End \C^{r}$.
\end{definition}
All the preceding properties of holomorphic families of \psivdos\ on an open subset of $\Rd$ 
hold \emph{verbatim} for holomorphic families with values in $\Psi_{H}^{*}(M,\cE)$. Moreover, we have: 

\begin{proposition}
    The principal symbol map $\sigma_{*}:\pvdo^{*}(M,\cE) \rightarrow S_{*}(\fg^{*}M,\cE)$ is analytic, in the sense that for any holomorphic family 
    $(P_{z})_{z\in \Omega}\subset \pvdo^{*}(M,\cE)$ the family of symbols $(\sigma_{*}(P_{z}))_{z\in\Omega}$ is in 
    $\Hol(\Omega, C^{\infty}(\fg^{*}M\setminus 0,\End \cE))$. 
\end{proposition}
\begin{proof}
   Let $(P_{z})_{z\in \Omega}\subset \pvdo^{*}(M,\cE)$ be a holomorphic family of \psivdos\ of order $m(z)$ and let us show that 
   the family of symbols $(\sigma_{*}(P_{z}))_{z\in\Omega}$ belongs to $\Hol(\Omega, C^{\infty}(\fg^{*}M\setminus 0,\End \cE))$. 
   Since this a purely local issue we may as well assume that  $(P_{z})_{z\in \Omega}$ is a holomorphic family of scalar \psivdos\ on a local trivializing 
   Heisenberg chart $U\subset \Rd$. 
   
   By Proposition~\ref{prop:HolPHDO.kernel.characterization} we can put the distribution kernel of $P_{z}$ into the form, 
   \begin{equation}
         k_{P_{z}}(x,y) = |\varepsilon_{x}'| K_{P_{z}}(x,-\varepsilon_{x}(y)) + R_{z}(x,y),
   \end{equation}
   with $(K_{P_{z}})_{z\in \Omega}\in \Hol(\Omega,\cK^{*}(\URd))$ of order $\hat{m}(z)=-(m(z)+d+2)$ and 
   $(R_{z})_{z\in \Omega}\in \Hol(\Omega,C^{\infty}(U\times U))$. Let $\varphi \in C^{\infty}(\Rd)$ be such that $\varphi(y)=1$ near $y=0$ and let 
   $p_{z}=(\varphi(y)K_{P_{z}}(x,y))^{\wedge}_{\xiy}$. Then the proof of Proposition~\ref{prop:HolPHDO.characterization-Hol-cK*}
   shows that $(p_{z})_{z\in \Omega}$ is a holomorphic family 
   of symbols. Moreover, we have 
   \begin{multline}
       K_{P_{z}}(x,y)= (p_{z})^{\vee}_{\xiy}(x,y)+(1-\varphi(y))K_{P_{z}}(x,y)\\ =  (p_{z})^{\vee}_{\xiy}(x,y) \ \bmod \Hol(\Omega,C^{\infty}(U\times U)).
        \label{eq:HolPHDO.KPz-local-symbol}
   \end{multline}
   
   Let $z\in \Omega$ and let $K_{\hat{m}(z)}  \in \cK_{\hat{m}(z)}(\URd)$ be the principal kernel of $K_{P_{z}}$. 
   Then~(\ref{eq:HolPHDO.KPz-local-symbol}) and  
   Proposition~\ref{lem:PsiHDO.characterization.Km} show 
   that the leading symbol of $p_{z}$ is the restriction to $\URdo$ of $(K_{\hat{m}(z)})^{\wedge}_{\yxi}$. Since the latter is equal to $\sigma_{m(z)}(P_{z})$, 
   we see that the leading symbol of $p_{z}$ is just $\sigma_{m(z)}(P_{z})$.
   Since $(p_{z})_{z\in \Omega}$ is a holomorphic family of symbols it then follows from Remark~\ref{rem:HolPHDO.analyticity-homogeneous-symbols} 
   that the family $(\sigma_{m(z)}(P_{z}))_{z \in \Omega}$ 
   belongs to $\Hol(\Omega, C^{\infty}(\URdo))$. The proof is thus achieved.
\end{proof}

\section{Transposes and adjoints of holomorphic families of \psivdos}
Let us now look at the analyticity and anti-analyticity of transposes and adjoints of holomorphic family of \psivdos. 

\begin{proposition}\label{prop:HolPHDO.transpose}
    Let $(P_{z})_{z\in \Omega} \in \Hol(\Omega, \pvdo^{*}(M,\cE))$. Then the transpose family $(P_{z}^{t})\subset 
    \pvdo^{*}(M,\cE^{*})$ is a holomorphic family of \psivdos.  
\end{proposition}
\begin{proof}
    For $z \in \Omega$ let $k_{P_{z}}(x,y)$ denote the distribution kernel of $P_{z}$. The distribution kernel of $P_{z}^{t}$ is 
    $k_{P_{z}}(x,y)=k_{P_{z}}(y,x)^{t}$, hence is represented outside the diagonal by a holomorphic family of smooth kernels. Therefore, we need only to prove 
    the statement for a holomorphic 
    family of scalar \psivdos\ on a Heisenberg chart $U\subset \Rd$, as we shall now suppose that the family $(P_{z})_{z\in \Omega}$ is. In addition, there 
    is no loss of generality in assuming that the order $m(z)$ of $P_{z}$ is such that there exists $\mu\in \R$ so that $\Re \hat{m}(z)\geq \mu$ for any 
    $z \in \Omega$. 
    
    Next, thanks to Proposition~\ref{prop:HolPHDO.kernel.characterization}
    the kernel of $P_{z}$ is of the form
    \begin{equation}
        k_{P_{z}}=|\varepsilon_{x}'|K_{P_{z}}(x,-\varepsilon_{x}(y))+R_{z}(x,y),
    \end{equation}
    with $(K_{P_{z}})_{z \in \Omega}$ in $\Hol(\Omega,\cK^{*}(\URd))$ and $(R_{z})_{z \in \Omega}$ in $\Hol(\Omega, C^{\infty}(U\times U))$. 
    Then the proof of 
    Proposition~\ref{prop:PsiHDO.transpose-chart} in Appendix~\ref{chap:Appendix-transpose} shows that $k_{P_{z}}(x,y)$ is equal to
    \begin{equation}
         |\varepsilon_{x}'|K_{P_{z}^{t}}(x,-\varepsilon_{x}(y)) +  
        (1-\chi(x,-\varepsilon_{x}(y)))|\varepsilon_{y}'|K_{P_{z}}(y,-\varepsilon_{y}(x))+R_{z}(y,x),
    \end{equation}
where we have let 
    \begin{equation}
        K_{P_{z}^{t}}(x,y)= \chi(x,y)|\varepsilon_{x}'|^{-1}|\varepsilon_{y}'|K_{P_{z}}(\varepsilon_{x}^{-1}(-y), 
        -\varepsilon_{\varepsilon_{x}^{-1}(-y)}(x)),  
    \end{equation}
and the function  $\chi(x,y)\in C^{\infty}(\URd)$ has a support contained in the open subset $\cU=\{(x,y)\in \URd; \ \varepsilon_{x}^{-1}(-y)\in 
U\}$,  is properly supported with respect $x$ and is equal to $1$ near $U\times \{0\}$. 
    In particular, we have 
    \begin{equation}
        k_{P_{z}^{t}}(x,y)= |\varepsilon_{x}'|K_{P_{z}^{t}}(x,-\varepsilon_{x}(y)) \qquad \bmod \Hol(\Omega,C^{\infty}(U\times U)).
         \label{eq:HolPHDO.kernel-transpose}
    \end{equation}

    Moreover, it follows from the proof of Proposition~\ref{prop:PsiHDO.transpose-chart} in Appendix~\ref{chap:Appendix-transpose} 
    that for any integer $N$ we can write
    \begin{gather}
        K_{P_{z}^{t}}(x,y)= \sum_{\alpha,\beta,\gamma,\delta}^{(N)} K_{\alpha\beta\gamma\delta,z}+ \sum_{j=1}^{4}R_{N,z}(x,y),\\ 
     K_{\alpha\beta\gamma\delta,z}=  a_{\alpha\beta\gamma\delta}(x) y^{\beta+\delta}  
       (\partial^{\gamma}_{x}\partial_{y}^{\alpha}K_{P_{z}})(x,-y),
\end{gather}
where the smooth functions $a_{\alpha\beta\gamma\delta}(x)$ are as in Proposition~\ref{prop:PsiHDO.transpose-chart}, the summation goes over all the multi-orders 
$\alpha$, $\beta$, $\gamma$ and $\delta$ such that $\brak\alpha<N$, $\frac{3}{2}\brak\alpha \leq 
\brak \beta < \frac{3}{2}N$ and $|\gamma|\leq |\delta| \leq 2|\gamma|<2N$ and the remainder terms $R_{N,z}^{(j)}(x,y)$ take the forms:\smallskip 

- $R_{N,z}^{(1)}=\sum_{\brak\alpha=N}\sum_{\brak \beta\dot{=}\frac{3}{2}\brak \alpha} 
   \frac{|\varepsilon_{y}'|}{ |\varepsilon_{x}'|}\int_{0}^{1}r_{N\alpha\beta}(t,x,y)
    (y^{\beta}\partial_{y}^{\alpha}K_{P_{z}})(\varepsilon_{x}^{-1}(-y), \Phi_{t}(x,y))$, where 
the functions $r_{N\alpha\beta}(t,x,y)$ are in $C^{\infty}([0,1]\times \URd)$, the equality $k\dot{=}\frac{3}{2}l$ means that $k$ is equal to $\frac{3}{2}l$ if 
$\frac{3}{2}l$ is integer and  to $\frac{3}{2}l+\frac{1}{2}$ otherwise, and we have let 
$\Phi_{t}(x,y)=-y+t(y-\varepsilon_{\varepsilon_{x}^{-1}(-y)}(x))$;\smallskip

- $R_{N,z}^{(2)}(x,y)= \sum_{\brak \beta \dot{=}\frac{3}{2}N}r_{N\alpha}(x,y)y^{\beta} (\partial_{y}^{\alpha}K_{P_{z}})(\varepsilon_{x}^{-1}(-y),-y)$,  
with $r_{N\alpha}(x,y)$ in $C^{\infty}(\URd)$;\smallskip

- $R_{N,z}^{(3)}(x,y)=\sum_{|\gamma|=N}\sum_{N\leq |\delta|\leq 2N} a_{\alpha\beta\gamma\delta}(x)y^{\beta+\delta} \int_{0}^{1}(1-t)^{N-1} 
  (\partial^{\gamma}_{x}\partial_{y}^{\alpha}K_{P_{z}}) (x+t(\varepsilon_{x}^{-1}(-y)-x), -y)dt$, with $a_{\alpha\beta\gamma\delta}(x)\in 
  C^{\infty}(U)$;\smallskip

  - $R_{N,z}^{(4)}(x,y)=\sum_{\alpha,\beta,\gamma,\delta}^{(N)} (1-\chi(x,y))K_{\alpha\beta\gamma\delta,z}(x,y)$.\smallskip

Each family $(K_{\alpha\beta\gamma\delta,z})_{z \in \Omega}$ belongs to $\Hol(\Omega,\cK^{*}(\URd))$. Moreover, the remainder term
 $R_{N,z}^{(4)}$ belongs to $\Hol(\Omega, C^{\infty}(\URd))$ and, along similar lines as that of the proof of 
Proposition~\ref{prop:HolPHDO.invariance}, we can show that for any integer $J$ the other remainder terms $(R_{N,z}^{(j)})_{z\in \Omega}$ are 
in $\Hol(\Omega, C^{J}(\URd))$ as soon $N$ is large enough. Therefore, we have 
$K_{P_{z}^{t}} \sim  \sum_{\frac{3}{2}\brak\alpha \leq \brak \beta} \sum_{|\gamma|\leq |\delta| \leq 2|\gamma|}  K_{\alpha\beta\gamma\delta,z}$
in the sense of~(\ref{eq:HolPHDO.kernel-asymptotic}),  which means that $K_{P_{t}^{z}}$ belongs to $\Hol(\Omega,\cK^{*}(\URd))$. Combining this 
with~(\ref{eq:HolPHDO.kernel-transpose}) and Proposition~\ref{prop:HolPHDO.kernel.characterization}  
then shows that $(P_{z}^{t})_{z\in \Omega}$ is a holomorphic family of \psivdos. 
\end{proof}

Assume now that $M$ is endowed with a density~$>0$ and $\cE$ with a Hermitian metric. Then Proposition~\ref{prop:HolPHDO.transpose} 
allows us to carried out the proof of Proposition~\ref{prop:PsiHDO.adjoint-manifold} in the setting of holomorphic families, so that we get: 

\begin{proposition}\label{prop:HolPHDO.adjoint}
    Let $(P_{z})_{z\in \Omega} \subset \pvdo^{*}(M,\cE)$ be a holomorphic family of \psivdos. Then the family 
    $(P_{z}^{*})_{z \in \Omega}\subset \pvdo^{*}(M,\cE^{*})$ is an anti-holomorphic family of \psivdos, in the sense 
    that $(P_{\overline{z}}^{*})_{z \in \Omega}$ is a holomorphic family of \psivdos.
\end{proposition}

 %%%%%%%%%%%%%%%%%%%%%%%%%%%%%%%%%%%%%
%%%%%%%%%%%%%%% Chap5.tex %%%%%%%%%%%%%%%%
%%%%%%%%%%%%%%%%%%%%%%%%%%%%%%%%%%%%%
 
 \chapter{Heat Equation and Complex Powers of Hypoelliptic Operators}
\label{chap.complex-powers}

In this chapter we deal with complex powers of hypoelliptic differential operators in connection with the heat equation. Due to the lack of 
microlocality of the Heisenberg we cannot carry out in the Heisenberg setting the standard approach of Seeley~\cite{Se:CPEO} to the complex powers of 
elliptic operators. Instead we rely on the pseudodifferential representation of the heat kernel of~\cite{BGS:HECRM}, which is especially suitable for 
dealing with positive differential operators.

In Section~\ref{sec.Volterra-PsiHDO-calculus} we recall the pseudodifferential representation of the heat kernel of an hypoelliptic operator in terms of the 
Volterra-Heisenberg calculus of~\cite{BGS:HECRM} and we extend to this setting the intrinsic approach of 
Chapter~\ref{chap:Heisenberg-calculus}. We then specialize the results to sublaplacians and integer powers of sublaplacians in 
Section~\ref{sec.Heat-complex-powers.sublaplacians}.   

In Section~\ref{sec.powers1}  we make use of the framework of Section~\ref{sec.Volterra-PsiHDO-calculus} 
to prove that the complex powers of a positive differential operators form a 
holomorphic family of \psivdos, provided that the principal symbol of the corresponding heat operator is invertible in the Volterra-Heisenberg calculus. 

The other two sections are devoted to applications of the above result. First, in Section~\ref{sec:Rockland-heat} 
we make use of it to extend Theorem~\ref{thm:Chap3.Rockland-Parametrix-order+} to 
non-integer \psivdos\ and to show that the invertibility of the principal symbol of the heat operator is implied by the Rockland condition 
when the bracket condition $H+[H,H]=TM$ holds. 

Second, in  Section~\ref{sec.Sobolev} we construct the weighted Sobolev spaces $W_{H}^{s}(M,\cE)$ 
and check their main properties. In particular, we prove that they yield sharp regularity results for \psivdos. 

% extended to the Heisenberg 
% setting the results of~\cite{Gr:AEHE}. 
% % This approach to complex powers of (hypo)elliptic operators was introduced in~\cite{Po:PhD} and has 
% % been used subsequently in~\cite{MMS:FAI} and~\cite{Me:SPLLB} in the context of projective pseudodifferential operators on Azamaya bundles. 
% 
% This section is divided into two subsections. In the first one we recall the  pseudodifferential representation of the heat kernel of a hypoelliptic 
% operator of~\cite{BGS:HECRM}, and in the second one we deal with the complex powers of a positive hypoelliptic operator. 

Throughout this chapter we let $(M^{d+1},H)$ be a compact Heisenberg manifold equipped with a (smooth) density~$>0$ and let $\cE$ be a Hermitian vector 
bundle over $M$ of rank~$r$.

\section{Pseudodifferential representation of the heat kernel}
\label{sec.Volterra-PsiHDO-calculus}
In this section we recall the pseudodifferential representation of the heat kernel of a hypoelliptic operator of~\cite{BGS:HECRM}, which extends to 
the Heisenberg setting the approach to the heat kernel asymptotics of~\cite{Gr:AEHE}. 
% On the way we complete the treatment in~\cite{BGS:HECRM} of the heat kernel of the Kohn 
% Laplacian. 
% subelliptic sublaplacian. 
% We will make use of the results of~\cite{Gr:AEHE} in the next section to show that the complex 
% powers of positive hypoelliptic operators give rise to a holomorphis

Let $P:C^{\infty}(M,\cE)\rightarrow C^{\infty}(M,\cE)$ be a differential operator of 
even (Heisenberg) order $v$ which is selfadjoint and bounded from below.  We also assume that the principal symbol of $P$ is an invertible principal 
symbol, which by Theorem~\ref{thm:Chap3.Rockland-Parametrix-order+}  is equivalent to say that $P$ satisfies the Rockland condition at every point. 
In particular,  $P$ is hypoelliptic with gain of 
$\frac{v}{2}$-derivatives by Theorem~\ref{thm:PsiHDO.hypoellipticity}. 

Since $P$ is bounded from below it generates on $L^{2}(M,\cE)$ a heat semigroup $e^{-tP}$, $t\geq 0$.
In fact, the hypoellipticity implies for $t>0$ the operator $e^{-tP}$ is smoothing, i.e., is given a  smooth kernel $k_{t}(x,y)$  in 
$C^\infty(M\times M, \cE\boxtimes(\cE \otimes|\Lambda|(M)))$, where $|\Lambda|(M)$ denotes the bundle of densities on $M$. 
% and for any integer $k\geq 1$ the operators 
% $P^{k}e^{-tP}$ and $e^{-tP}P^{k}$ coincide and are bounded. Thus $e^{-tP}$ maps continuously 
% $L^{2}(M,\cE)$ to $\cap_{k\geq 0} \dom P^{k}$, which is just $C^{\infty}(M,\cE)$ since $P$ is hypoelliptic with gain of $\frac{v}{2}$-derivatives.  
% Moreover, as $e^{-tP}$ is selfadjoint it defines by duality a continuous map from $\cD'(M,\cE)$ to $L^{2}(M,\cE)$. Since
% $e^{-tP}=e^{-tP/2}e^{-tP/2}$ it follows that for $t>0$ the operator $e^{-tP}$ extends to a continuous map from 
% $\cD'(M,\cE)$ to $C^{\infty}(M,\cE)$, i.e.,     
%  $e^{-tP}$ is a smoothing operator and so is given by 

On the other hand, the heat semi-group allows us to invert the heat equation. Indeed, the operator given by 
\begin{equation}
    Q_{0}f(x,t)=\int_{0}^\infty e^{-sP} f(x,t-s)dt, \qquad f \in C^\infty_{c}(M\times \R, \cE), 
     \label{eq:volterra.inverse-heat-operator}
\end{equation}
maps continuously $C^\infty_{c}(M\times \R, \cE)$ into $C^{0}(\R, L^{2}(M,\cE)) \subset \cD(M\times\R, \cE)$ and satisfies
\begin{equation}
    (P+\partial_{t})Q_{0}f = Q_{0}(P+\partial_{t})f=f \qquad \forall f \in C^\infty_{c}(M\times\R,\cE).
\end{equation} 

Notice that the operator $Q_{0}$  has the \emph{Volterra property} of~\cite{Pi:COPDTV}, i.e., it is translation invariant and satisfies the causality 
principle with respect to the time variable, or equivalently, has a distribution kernel of the form $K_{Q_{0}}(x,y,t-s)$ with $K_{Q_{0}}(x,y,t)$ supported outside 
the region $\{t<0\}$. Indeed, at the level of distribution kernels the formula~(\ref{eq:volterra.inverse-heat-operator}) implies that we have 
\begin{equation}
    K_{Q_{0}}(x,y,t) = \left\{ 
    \begin{array}{ll}
         k_{t}(x,y) & \quad \text{if $t> 0$},\\
        0 &  \quad \text{if $t<0$}. \label{eq:Powers1.kernel-inverse-heat-operator-heat-kernel}
    \end{array}\right. 
\end{equation}

The above equalities are the main motivation for using pseudodifferential operators to study the heat kernel $k_{t}(x,y)$. 
The idea is to consider a class of \psivdo's, the Volterra \psivdo's, taking into account:  \smallskip 

(i)  The aforementioned Volterra property;\smallskip 

(ii) The parabolic homogeneity of the heat operator $P+ \partial_{t}$, i.e., the homogeneity with respect 
             to the dilations of $\Rdd=\Rd\times \R$ defined by
             \begin{equation}
                  \lambda.(\xi,\tau)=(\lambda.\xi,\lambda^{v}\tau), \qquad  \lambda\neq 0. 
             \end{equation}

 In the sequel for $g\in \cS'(\R^{d+2})$  and $\lambda\neq 0$ we let $g_{\lambda}$ denote the element of $\cS'(\R^{d+2})$ such that
\begin{equation}
    \acou{g_{\lambda}( \xi,\tau)}{f(\xi,\tau)} =   |\lambda|^{-(d+2+v)}\acou{g(\xi,\tau)} {f(\lambda^{-1}\xi, 
    \lambda^{-v}\tau)}, \quad f\in \cS(\R^{d+2}).               
\end{equation}

\begin{definition}
    A distribution $ g\in \cS'(\R^{d+2})$ is parabolic 
homogeneous of degree $m$, $m\in \Z$, when we have $g_{\lambda}=\lambda^m g$ for any $\lambda \neq 0$.  
\end{definition}

Let $\C_{-}$ denote the complex halfplane $\{\Im \tau <0\}$ with closure $\overline{\C_{-}}\subset \C$. 
Then we define Volterra symbols and Volterra \psivdos\ as follows.

\begin{definition}\label{def:Volterra-symbols-Rd}
   The space $S_{\op{v},m}(\Rdd_{(v)})$, $m\in\Z$,  consists of functions $q(\xi,\tau)$ in $C^{\infty}(\Rdd\setminus 0)$ such 
    that:\smallskip
    
    (i) $q(\lambda.\xi,\lambda^{v}\tau)=\lambda^{m}q(x,\xi,\tau)$ for any $\lambda\neq 0$;\smallskip
    
    (ii) $q(\xi,\tau)$ extends to a function in $C^{0}((\Rd\times\overline{\C_{-}})\setminus 0)$ whose restriction to $\Rd\times \C_{-}$
   belongs to $C^{\infty}(\Rd)\hotimes \Hol(\C_{-})$.\smallskip
\end{definition}

We also endow $S_{\op{v},m}(\Rdd_{(v)})$ with the Fr\'echet space topology inherited from that of $C^{\infty}(\Rdd\setminus 
 0) \cap C^{0}((\Rd\times\overline{\C_{-}})\setminus 0)\cap [C^{\infty}(\Rd)\hotimes \Hol(\C_{-})]$.
%  i.e., the weakest topology making continuous the 
%  embeddings of $S_{\op{v},m}(\Rdd_{(v)})$ into these Fr\'echet spaces.
%  i.e., with the topology defined by 
% by the sytem of the seminorms on $S_{\op{v},m}(\Rdd_{(v)})$ coming from those defining the topologies of these Fr\'echet spaces. 

The interest of the above definition stems from:
\begin{lemma}[{\cite[Prop.~1.9]{BGS:HECRM}}] \label{lem:volterra.volterra-extension}
Let $q(\xi,\tau)\in S_{\op{v},m}(\Rdd_{(v)})$. Then  there exists a unique distribution $g\in \cS'(\R^{d+2})$ agreeing with $q$ 
on $\R^{d+2}\setminus 0$ such that: \smallskip 

(ii) $g$ is parabolic homogeneous of degree $m$; \smallskip 

(iii) The inverse Fourier transform $\check g(x,t)$ vanishes for $t<0$.
\end{lemma}

Let $U$ be an open subset of $\R^{d+1}$ together with a hyperplane bundle $H\subset TU$ and $H$- frame $X_{0},\ldots, 
X_{d}$ of $TU$.  

\begin{definition}\label{def:Volterra-symbols-URd}
 $S_{\op{v},m}(\URdd_{(v)})$, $m\in\Z$,  consists of functions $q(x,\xi,\tau)$ in $C^{\infty}(U\times(\R^{d+2}\setminus 0))$ such 
    that:\smallskip
    
    (i) $q(x,\lambda.\xi,\lambda^{v}\tau)=\lambda^{m}q(x,\xi,\tau)$ for any $\lambda\neq 0$;\smallskip
    
    (ii) $q(x,\xi,\tau)$ extends to an element of $C^{\infty}(U)\hotimes C^{0}((\Rd\times\overline{\C_{-}})\setminus 0)$ whose restriction to 
    $U\times \Rd\times \C_{-}$
   belongs to $C^{\infty}(\URd)\hotimes \Hol(\C_{-})$.
\end{definition}

\begin{definition}
    $S_{\op{v}}^m(\URdd_{(v)})$, $m\in\Z$,  consists of functions $q(x,\xi,\tau)$ in 
    $C^{\infty}(U\times\Rdd)$ with an asymptotic expansion  $q \sim \sum_{j\geq 0} q_{m-j}$ with $q_{m-j}$ in $S_{\op{v},m-j}(U\times\Rdd_{(v)})$ 
    and $\sim $ taken   in the sense that, for any integer $N$ and any compact $K\subset U$,  we have  
            \begin{equation}
                |\partial^{\alpha}_{x}\partial^{\beta}_{\xi} \partial^k_{\tau}(q-\sum_{j< N} 
            q_{m-j})(x,\xi,\tau) | 
                \leq C_{NK\alpha\beta k} (\|\xi\|+|\tau|^{1/v})^{m-N-\brak\beta-v k},
                          \label{eq:volterra.asymptotic-symbols}
            \end{equation}
            for  $x\in K$ and $|\xi|+|\tau|^{\frac{1}{v}}>1$.
\end{definition}

\begin{definition}\label{def:Volterra.inverse-Fourier-transform}
  Let $q(x,\xi,\tau)\in S_{\op{v},m}(U\times\Rdd_{(v)})$ and let $g$ be the distribution in  
  $C^{\infty}(U)\hotimes \cS'(\Rd)$ be the unique homogeneous extension of $q$ provided by Lemma~\ref{lem:volterra.volterra-extension}. Then we let 
  $\check{q}(x,y,t)$ denote the inverse Fourier transform of $g(x,\xi,\tau)$ with respect to the variables $(\xi,\tau)$.
\end{definition}

\begin{remark}
    The above definition makes sense since it follows from the proof of Lemma~\ref{lem:volterra.volterra-extension} 
    in \cite{BGS:HECRM} that the extension process of Lemma~\ref{lem:volterra.volterra-extension}  applied to every symbol $q(x,.,.)$, $x \in U$, 
    is smooth with respect to $x$, so really gives rise to an element of $C^\infty(U)\hotimes\cS'(\R^{d+2})$. 
 \end{remark}

\begin{definition}\label{def:volterra.PsiDO}
    $\pvhdo^m(U\times\R_{(v)})$, $m\in\Z$,  consists of continuous operators 
    $Q:C_{c}^\infty(U_{x}\times\R_{t}) \rightarrow C^\infty(U_{x}\times\R_{t})$ such that $Q$ has the Volterra property and can be put into the form
\begin{equation}
    Q=q(x,-iX,D_{t})+R, 
\end{equation}  
   with $q$ in $S^m_{\op{v}}(U\times\Rdd_{(v)})$ and $R$ in $\psinf(U\times \R)$. 
\end{definition}

\begin{remark}
    It is immediate to extend the properties of \psivdos\ on $U$ alluded to in Section~\ref{sec:PsiHDO} to Volterra \psivdos\ on $U\times \R$ except for 
    the asymptotic completeness as in Lemma~\ref{lem:PsiHDO.asymptotic-completeness}, 
    which is crucial for carrying out the standard parametrix construction. The issue is that the cut-off arguments of the classical proof 
    the asymptotic completeness of standard \psidos\ cannot be carried through in Volterra setting because we require analyticity with respect to the time 
    covariable. A proof of the asymptotic completeness of Volterra \psidos\ is given in~\cite{Pi:COPDTV}, but simpler proofs which can be carried out 
    \emph{verbatim} for Volterra \psivdos\ can be found in~\cite{Po:JAM1}.  
\end{remark}

Let $a \in U$. Then, as for Heisenberg symbols, the convolution on the groups $G^{(a)}\times \R$, gives rise to a continuous bilinear product, 
\begin{equation}
    *^{(a)}:S_{\op{v},m_{1}}(\Rdd_{(v)})\times S_{\op{v},m_{2}}(\Rdd_{(v)})\rightarrow S_{\op{v},m_{1}+m_{2}}(\Rdd_{(v)}).
%     \label{eq:¥}
\end{equation}
% \begin{equation}
%     *^{(a)}:S_{\op{v},m_{1}}(\Rdd_{(v)})\times S_{\op{v},m_{2}}(\Rd\times\R_{(v)})\rightarrow S_{\op{v},m_{1}+m_{2}}(\Rd\times\R_{(v)}).
% %     \label{eq:¥}
% \end{equation}
Here again $*^{a}$ depends smoothly on $a$ and so we get a bilinear product,
% $*$ from  $S_{\op{v},m_{1}}(\URd\times\R_{(v)})\times 
% S_{\op{v},m_{2}}(\URd\times\R_{(v)})$ to $S_{\op{v},m_{1}+m_{2}}(\URd\times\R_{(v)})$ such that we have 
% , for any $q_{1}\in S_{\op{v},m_{1}}(\URd\times\R_{(v)})$ and $q_{2}\in S_{\op{v},m_{2}}(\URd\times\R_{(v)})$, 
\begin{gather}
    S_{\op{v},m_{1}}(\URdd_{(v)})\times S_{\op{v},m_{2}}(\URdd_{(v)})\rightarrow S_{\op{v},m_{1}+m_{2}}(\URdd_{(v)}),\\
    q_{1}*q_{2}(x,\xi,\tau)=(q_{1}(x,.,.)*^{x}q_{2}(x,.,.))(\xi,\tau), \quad q_{j}\in S_{\op{v},m_{j}}(\URdd_{(v)}).
%     \label{eq:¥}
\end{gather}

\begin{proposition}
   For $j=1,2$ let $Q_{j} \in \pvhdo^{m_{j}}(U\times\R_{(v)})$ and assume  
    that one of these operators is properly supported with respect to the space variable $x$. Then $Q_{1}Q_{2}$ belongs to 
    $\pvhdo^{m_{j}}(U\times\R_{(v)})$ and if $q_{j} \sim \sum_{k \geq 0} q_{j,m-k}$ denote the symbol of $Q_{j}$ then $Q_{1}Q_{2}$ 
    has symbol $q\sim \sum_{k\geq 0}q_{m_{1}+m_{2}-k}$ where, using the notation of Proposition~\ref{prop:PsiHDO.composition}, we have 
    \begin{equation}
        q_{m_{1}+m_{2}-k} = \sum_{k_{1}+k_{2}\leq k} \sum_{\alpha,\beta,\gamma,\delta}^{(k-k_{1}-k_{2})}
            h_{\alpha\beta\gamma\delta}  (D_{\xi}^\delta q_{1,m_{1}-k_{1}})* (\xi^\gamma 
            \partial_{x}^\alpha \partial_{\xi}^\beta q_{2,m_{2}-k_{2}}).      
%         \label{eq:¥}
    \end{equation}
\end{proposition}

\begin{remark}\label{rem:Volterra-time-microlocality}
Since $G^{(a)}\times \R$ is Abelian with respect to the time variable, the product $*^{a}$ is merely the pointwise product with respect to 
$\tau$, e.g., we have $\tau*q=\tau*q=\tau q$ for any $q\in S_{\op{v},m}(\URdd_{(v)})$. In particular, the Volterra-Heisenberg calculus, 
while not microlocal with respect to the space variable, is to a large extent microlocal with respect to the time variable. 
\end{remark}

On the other hand, thanks to the Volterra property the kernels of \psivdos\ can be characterized as follows.

\begin{definition}
$\cK_{\op{v}, m}(U\times\Rdd_{(v)})$, $m\in \Z$, consists of distributions $K(x,y,t)$ in 
 $C^{\infty}(U)\hotimes\cS'_{\reg}(\R^{d+2})$ such that:\smallskip 
     
    \indent (i) The support of $K(x,y,t)$ is contained in $U\times\Rd\times\R_{+}$;\smallskip 
    
    \indent (ii) $K(x,\lambda.y,\lambda^{v}t)=(\op{sign} \lambda )^{d}\lambda^{m}K(x,y,t)$ for any $\lambda\in \R\setminus 0$.   
\end{definition}

\begin{definition}
 $\cK_{\op{v}}^{m}(U\times\Rdd_{(v)})$, $m\in \Z$, is the space of distributions $K(x,y,t)$ in $\cD'(U\times\R^{d+2})$ which  
 admit an asymptotic expansion $K \sim \sum_{j \geq 0} K_{m+j}$ with $K_{m+j}$ in $\cK_{\op{v}, m+j}(U\times\R^{d+2})$ and 
 $\sim$ taken in the sense of~(\ref{eq:PsiHDO.asymptotics-kernel}).  
\end{definition}

In the sequel, for $x\in U$ we let $\psi_{x}$ and $\varepsilon_{x}$ respectively denote the changes of variable to the privileged coordinates and 
to the Heisenberg coordinates at $x$. Then,  along the same lines as that of the proofs of 
Proposition~\ref{prop:PsiVDO.characterisation-kernel1} and Proposition~\ref{prop:PsiVDO.characterisation-kernel2}, 
we obtain the following characterization of Volterra \psivdos.

\begin{proposition}\label{prop:volterra.kernel-charaterization} 
    Let $Q:C_{c}^\infty(U_{x}\times\R_{t}) \rightarrow C^\infty(U_{x}\times\R_{t})$ be a continuous operator with distribution kernel $k_{Q}(x,t;y,s)$. Then 
    the following are equivalent:\smallskip 
    
    (i) The operator $Q$ belongs to $\pvhdo^{m}(U\times \R)$;\smallskip
    
    (ii) The kernel of $Q$ can be put into the form,
    \begin{equation}
        k_{Q}(x,t;y,s)=|\psi_{x}'|K(x,-\psi_{x}(y),t-s)+R(x,y,t-s),
     \label{eq:Powers1.characterization-Volterra-PsiHDO's-psix}
\end{equation}
    with $K$ in $\cK_{\op{v}}^{\hat m}(U\times\Rdd_{(v)})$, $\hat{m}=-(m+d+2+v)$,  and $R$ in $C^{\infty}(U\times 
    \Rdd)$.

(iii) The kernel of $Q$ can be put into the form,
    \begin{equation}
        k_{Q}(x,t;y,s)=|\varepsilon_{x}'|K_{Q}(x,-\varepsilon_{x}(y),t-s)+R(x,y,t-s),
     \label{eq:Powers1.characterization-Volterra-PsiHDO's}
\end{equation}
    with  $K_{Q}$ in $\cK_{\op{v}}^{\hat m}(U\times\Rdd_{(v)})$, $\hat{m}=-(m+d+2+v)$, and $R$ in $C^{\infty}(U\times 
    \Rdd)$.
\end{proposition}

An interesting consequence of Proposition~\ref{prop:volterra.kernel-charaterization} 
is the following small time asymptotics for the kernel of a Volterra \psivdo. 

\begin{proposition}[{\cite[Thm.~4.5]{BGS:HECRM}}]\label{prop:Volterra.asymptotics-kQ}
Let $Q\in \pvhdo^{m}(U\times\R_{(v)})$ have symbol $q \sim \sum_{j \geq 0} q_{m-j}$ and kernel $k_{Q}(x,y,t-s)$. 
Then as $t\rightarrow 0^{+}$ the following asymptotics 
holds in  $C^\infty(U)$, 
\begin{equation}
    k_{Q}(x,x,t) \sim t^{-\frac{2[\frac{m}2]+d+4}{v}} \sum_{j\geq 0} t^j 
    |\varepsilon_{x}'|\check{q}_{2[\frac{m}2]-2j}(x,0,1). 
    \label{eq:volterra.asymptotics-Q}
\end{equation}
\end{proposition}

Let $\tilde{U}$ be an open subset of $\Rd$ together with a hyperplane bundle $\tilde{H}\subset T\tilde{U}$ and a 
  $\tilde{H}$-frame of $T\tilde{U}$ and let $\phi:(U,H)\rightarrow (\tilde{U},\tilde{H})$ be a Heisenberg diffeomorphism.  
Then using Proposition~\ref{prop:volterra.kernel-charaterization} and arguing along similar lines that of the proof 
of Proposition~\ref{prop:PsiHDO.invariance} allows us to prove: 

\begin{proposition}\label{prop:Powers1.invariance}
   Let $\tilde{Q}\in \Psi_{\tilde{H},\op{v}}^{m}(\tilde{U}\times\R_{(v)})$ and set $Q=(\phi\oplus 1_{\R})^{*}\tilde{Q}$.\smallskip 
     
   1) The operator $Q$ %=(\phi\oplus 1_{\R})^{*}\tilde{Q}:C^{\infty}_{c}(U\times \R)\rightarrow C^{\infty}(U\times \R)$ 
   belongs to $\pvhdo^{m}(U\times\R_{(v)})$.\smallskip 
    
   2) If the distribution kernel of $\tilde{P}$ is of the form~(\ref{eq:Powers1.characterization-Volterra-PsiHDO's}) 
   with $K_{\tilde{Q}}(\tilde{x},\tilde{y},t)$ in $\cK_{\op{v}}^{\hat m}(\tilde{U}\times\Rdd_{(v)})$ then the distribution kernel of 
   $P$ can be written in the form~(\ref{eq:Powers1.characterization-Volterra-PsiHDO's}) with $K_{Q}(x,y,t)$ in 
   $\cK_{\op{v}}^{\hat m}(U\times\Rdd_{(v)})$ such that 
   \begin{equation}
       K_{Q}(x,y,t) \sim \sum_{\brak\beta\geq \frac{3}{2}\brak\alpha} \frac{1}{\alpha!\beta!} 
       a_{\alpha\beta}(x)y^{\beta}(\partial_{\tilde{y}}^{\beta}K_{\tilde{Q}})(\phi(x),\phi_{H}'(x)y,t),
        \label{eq:Powers1.asymptotic-expansion-KQ}
   \end{equation}
   where the functions $a_{\alpha\beta}(x)$ are as in~(\ref{eq:PsiHDO.asymptotic-expansion-KP}).
 \end{proposition}

This allows us to define Volterra \psivdos\ on the manifold $M\times \R$ and acting on the sections of the 
bundle $\cE$ (or rather on the sections of the pullback of $\cE$ by the projection $M\times \R\rightarrow M$, again denoted $\cE$). 

\begin{definition}
  $\pvhdo^{m}(M\times\R_{(v)},\cE)$, $m\in \Z$, consists of continuous operators $Q:C^{\infty}_{c}(M\times \R,\cE)\rightarrow C^{\infty}(M\times \R,\cE)$ such 
  that:\smallskip 
  
  (i) $Q$ has the Volterra property;\smallskip
  
  (ii) The distribution kernel of $Q$ is smooth off the diagonal of $(M\times \R)\times(M\times \R)$;\smallskip
  
  (iii) For any trivialization $\tau:\cE_{|_{U}}\rightarrow U\times \C^{r}$ of $\cE$ over a 
     local Heisenberg chart $\kappa:U \rightarrow V\subset \Rd$ the operator $(\kappa\otimes\op{id})_{*}\tau_{*}(Q_{|_{U\times \R}})$ belongs to
     $\pvhdo^{m}(V\times\R_{(v)}, \C^{r}):=\pvhdo^{m}(V\times \R_{(v)})\otimes \End \C^{r}$. 
\end{definition}

Using Proposition~\ref{prop:Powers1.invariance} 
we can define the global principal symbol of a Volterra \psivdo\ as follows. Let $\fg^{*}M$ denote the dual bundle of the Lie algebra 
bundle $\fg M$ of $M$ and consider the canonical projection 
$\pi:\fg^{*}M \times \R\rightarrow M$. 

In the sequel, 
depending on the context, $0$ 
denotes either the zero section of $\fg^{*}M$ or  the zero section of $\fg^{*}M$  
crossed with $\{0\}\subset \bar\C_{-}$. 

\begin{definition}\label{def:Volterra-symbols-fgM}
 $S_{\op{v},m}(\fg^{*}M\times \R_{(v)},\cE)$, $m\in \Z$, consists of sections $q(x,\xi,\tau)$ in $C^{\infty}((\fg^{*}M\times \R)\setminus 0, \pi^{*}\End 
 \cE)$ such that:\smallskip
 
 (i) $q(x,\lambda.\xi,\lambda^{v}\tau)=\lambda^{m}q(x,\xi,\tau)$ for any $\lambda \in \R\setminus 0$;\smallskip 
 
 (ii) $q(x,\xi,\tau)$ extends to a section of $\pi^{*}\End  \cE$ over $(\fg^{*}M\times \bar\C_{-})\setminus 0$ which is 
 smooth with respect to the base space variable and continous with 
 respect to the others and which restricts on $\fg^{*}M\times \C_{-}$ to an element of 
 $C^{\infty}(\fg^{*}M, \pi^{*}\End  \cE)\hotimes \Hol(\C_{-})$.
\end{definition}

Using~(\ref{eq:Powers1.asymptotic-expansion-KQ}) and arguing as in the proof of Proposition~\ref{prop:PsiHDO.principal-symbol} we get:

\begin{proposition}\label{prop:Powers1.principal-symbol}
    For any $Q \in \pvhdo^{m}(M\times\R_{(v)},\cE)$ there is a unique symbol $\sigma_{m}(Q)\in S_{\op{v},m}(\fg^{*}M \times\R_{(v)}, \cE)$ 
    such that, if in a local trivializing Heisenberg 
    chart $U\subset \Rd$ we let $K_{Q,\hat{m}}\in \cK_{\hat{m}}(\URd)$ be the leading kernel for the kernel $K_{Q}$  in the 
    form~(\ref{eq:Powers1.characterization-Volterra-PsiHDO's})   
    for $Q$,  then for $(x,\xi,\tau)\in  U\times [(\Rd\times \R)\setminus  0]$ we have
    \begin{equation}
       \sigma_{m}(Q)(x,\xi,\tau)=[K_{Q,\hat{m}}]^{\wedge}_{(y,t)\rightarrow (\xi,\tau)}(x,\xi,\tau) . 
    \end{equation}
    
    Equivalently, on any trivializing Heisenberg coordinates   centered at $a\in M$ the symbol $\sigma_{m}(Q)(a,.,.)$ coincides 
    with the (local) principal symbol of $Q$ at $x=0$. 
\end{proposition}

\begin{definition}\label{def:Powers1.principal-symbol}
   For $Q \in \pvhdo^{m}(M\times\R_{(v)},\cE)$ the symbol $\sigma_{m}(Q)(x,\xi,\tau)$  provided by 
   Proposition~\ref{prop:Powers1.principal-symbol} is called the (global) principal symbol of $Q$. 
\end{definition}

Extending Definition~\ref{def:Volterra.inverse-Fourier-transform} 
to Volterra-Heisenberg symbols on $M\times \R$ we can define the model operator of a Volterra \psivdo\ as follows. 
\begin{definition}
 Let $Q \in \pvhdo^{m}(M\times\R_{(v)},\cE)$ have principal symbol $\sigma_{m}(Q)$ and let $a \in M$.   Then the model operator $Q^{a}$ 
 of $Q$ at $a$ is the left-convolution 
 operator by  $\sigma_{m}(Q)^{\vee}(a,.,.)$, that is, $Q^{a}$ is the continuous endomorphism of 
 $\cS(G_{a}M\times \R,\cE_{a})$ such that, for any $f \in \cS(G_{a}M,\cE_{a})$, we have 
    \begin{equation}
        Q^{a}f(x,t)=\acou{\sigma_{m}(Q)^{\vee}(a,y,t)}{f(x.y^{-1},t-s)}. 
    \end{equation}
\end{definition}
\begin{remark}
    The model operator $Q^{a}$ can be defined as an endomorphism of $\cS(G_{a}M,\cE_{a})$, not just as an endomorphism of $\cS_{0}(G_{a}M,\cE_{a})$ as in 
    Definition~\ref{def:PsiHDO.model-operator}, because $\sigma_{m}(Q)^{\vee}_{(\xi,\tau)\rightarrow (y,t)}(a,.,.)$ makes sense as an element of 
    $\cS'(G_{a}M,\cE_{a})$. 
\end{remark}

\begin{proposition}
  The group law on the fibers of $GM\times \R$ gives rise to a convolution product $*$ from $S_{\op{v},m_{1}}(\fg^{*}M\times\R_{(v)},\cE) \times 
  S_{\op{v},m_{2}}(\fg^{*}M\times\R_{(v)},\cE)$ to $S_{\op{v},m_{1}+m_{2}}(\fg^{*}M\times\R_{(v)},\cE)$
%   
%     \begin{equation}
%         *: S_{\op{v},m_{1}}(\fg^{*}M\times\R_{(v)},\cE) \times S_{\op{v},m_{2}}(\fg^{*}M\times\R_{(v)},\cE) 
%         \rightarrow S_{\op{v},m_{1}+m_{2}}(\fg^{*}M\times\R_{(v)},\cE),
%          \label{eq:Powers1.volterra-Heisenberg-symbol-product}
%     \end{equation}
such that, for $q_{m_{j}}\in S_{\op{v},m_{j}}(\fg^{*}M\times\R_{(v)},\cE)$,we have
    \begin{gather}
        q_{1}*q_{2}(x,\xi,\tau)=[q_{1}(x,.,.)*^{x}q_{2}(x,.,.)](\xi,\tau),
        \label{eq:Powers1.volterra-Heisenberg-symbol-product}
        %\quad q_{j}\in S_{\op{v},m_{j}}(\fg^{*}M\times\R_{(v)},\cE),
        %\qquad (x,\xi,\tau)\in (\fg^{*}M\times \R)\setminus 0,
    \end{gather}
 where $*^{x}$ denote the convolution product for symbols on $G_{x}M\times \R$.
\end{proposition}
 
In a local trivializing Heisenberg chart the symbolic calculus for Volterra \psivdos\ reduces the existence of a Volterra \psivdo\ parametrix to the 
invertibility of the local and global principal symbols. Therefore, we obtain:

\begin{proposition}
    Let $Q\in \pvhdo^{m}(M\times\R_{(v)},\cE)$, $m \in \Z$. Then we have equivalence:\smallskip 
    
    (i) The principal symbol of $Q$ is invertible with respect to the product~(\ref{eq:Powers1.volterra-Heisenberg-symbol-product}) 
    of Volterra-Heisenberg symbols;\smallskip 
    
    (ii) The operator $Q$ admits a parametrix in $\pvhdo^{-m}(M\times\R_{(v)},\cE)$.
\end{proposition}

In the case of the heat operator $P+\partial_{t}$, comparing a parametrix with the inverse~(\ref{eq:volterra.inverse-heat-operator}) and 
using~(\ref{eq:Powers1.kernel-inverse-heat-operator-heat-kernel}) allows us to obtain the pseudodifferential representation of the heat kernel of $P$ 
below.

\begin{theorem}[{\cite[pp.~362--363]{BGS:HECRM}}]\label{thm:volterra.inverse} 
Suppose that the principal symbol of $P+\partial_{t}$ is an invertible Volterra-Heisenberg symbols. Then:\smallskip  

1) The heat operator  $P+\partial_{t}$ has an inverse  $(P+\partial_{t})^{-1}$ in $\pvhdo^{-v}(M\times\R_{(v)},\cE)$.\smallskip 

2) Let $K_{(P+\partial_{t})^{-1}}(x,y,t-s)$ denote the kernel of $(P+\partial_{t})^{-1}$. Then the heat kernel $k_{t}(x,y)$ of $P$ satisfies
\begin{equation}
    k_{t}(x,y)=K_{(P+\partial_{t})^{-1}}(x,y,t) \quad \text{for $t>0$.}
     \label{eq:Powers1.kernel-inverse-heat-operator-heat-kernel2}
\end{equation}
\end{theorem}
Combining this with Proposition~\ref{prop:Volterra.asymptotics-kQ} then gives the heat kernel asymptotics for $P$ in the form below.

\begin{theorem}[{\cite[Thm.~5.6]{BGS:HECRM}}]\label{thm:Powers1.heat-kernel-asymptotics}
  If the principal symbol of $P+\partial_{t}$ is an invertible Volterra-Heisenberg symbol, then as $t\rightarrow 0^{+}$  the following asymptotics holds in 
  $C^{\infty}(M,(\End \cE)\otimes|\Lambda|(M))$, 
    \begin{equation}
     k_{t}(x,x) \sim t^{-\frac{d+2}{v}} \sum t^{\frac{2j}{v}} a_{j}(P)(x), \quad  
    a_{j}(P)(x) =|\varepsilon_{x}'|(q_{-v-2j})^{\vee}%_{\xitauyt}
    (x,0,1), 
    \label{eq:Rockland-Heat.heat-kernel-asymptotics}
    \end{equation}
where the equality on the right shows how to compute $a_{j}(P)(x)$ in a local trivializing Heisenberg chart by means of the symbol 
$q_{-v-2j}(x,\xi,\tau)$ of degree $-v-2j$ of any parametrix of $P+\partial_{t}$ in $\pvhdo^{-v}(M\times\R_{(v)},\cE)$. 
\end{theorem}

% % We will show in Section~\ref{sec:Rockland-heat} that, when the bracket condition $H+[H,H]=TM$ holds, i.e., the rank Levi form of $(M,H)$ is always 
% % nonzero, the fact that $P$ satisfies the Rockland condition at every point implies that the principal symbol of $P+\partial_{t}$ is an invertible 
% % Volterra-Heisenberg symbol. 

\section{Heat equation and sublaplacians}
\label{sec.Heat-complex-powers.sublaplacians}
In this section we specialize the results of the previous sections to sublaplacians and their integer powers. In particular, 
we complete the treatment in~\cite{BGS:HECRM} of the heat kernel of a sublaplacian. 
% This will  also be needed  which we need anyway in order to deal with more general operators in Section~\ref{sec:Rockland-heat}. 

Throughout this section we let $k$ be an integer~$\geq 1$ and we set $v=2k$. In order to deal with sublaplacians 
it will be convenient to enlarge the definition of homogeneous Volterra-Heisenberg symbols as follows. 

Let $U\subset \Rd$ be a Heisenberg chart with $H$-frame $X_{0},\ldots,X_{d}$ and let $\Theta$ be an open angular sector whose closure 
contains $\R$ and which is contained in $\C\setminus i[0,\infty)$. 

\begin{definition}
For $m\in \Z$ and $N \in \N\cup \{\infty\}$ we let $S_{\vo,m}(\URd \times \Theta_{(v)}^{N})$ be the space of functions $q(x,\xi,\tau)$ on 
 $U\times (\Rdd\setminus 0)$ such that:\smallskip
 
    (i) $q(x,\xi,\tau)$ is $C^{\infty}$ on $U\times (\Rdd\setminus 0)$ and near the region $\{\tau=0\}$ 
         is $C^{\infty}$ with respect to $x$ and $\xi$ and $C^{N}$ with respect to $\tau$;\smallskip 
%    $C^{0}((\Rd \times \overline{\Theta})\setminus 0)$ such that:\smallskip 
    
    (ii) We have $q(\lambda.\xi,\lambda^{v}\tau)=\lambda^{m} q(\xi,\tau)$ for any $\lambda \in \R\setminus 0$;\smallskip
    
    (iii) $q(x,\xi,\tau)$ extends to a continuous function on  $U\times [(\Rd\times \overline{\Theta})\setminus 0]$ whose restriction to 
          $\URd  \times\Theta$ belongs to $C^{\infty}(\Rd)\hotimes \Hol (\Theta)$.%\smallskip
% \noindent In addition, we endow $S_{\vo,m}^{N}(\URd \times \Theta_{(v)})$ with the Fr\'echet space topology generated by that of 
% $C^{0}((\Rd \times \overline{\Theta})\setminus 0)$, $C^{\infty}(\URdd\setminus 0)$ and $C^{\infty}(\URd)\otimes \Hol (\Theta)$. 
\end{definition}

In particular, when $\Theta=\C_{-}$ and $N=\infty$ we recover the class $S_{\op{v},m}(\URdd_{(v)})$ defined in the previous section. 

As alluded to in Remark~\ref{rem:Volterra-time-microlocality} the Volterra-Heisenberg calculus is microlocal with the respect to the time variable. 
As we shall now see this allow us to extend the product for homogeneous Volterra-Heisenberg symbols to the symbols in the class 
$S_{\vo,m}^{N}(\URd \times \Theta_{(v)}^{N})$. 
%  In fact, this will follow along the same  lines as that of the construction of the product for homogeneous 
% Heisenberg symbols in~\cite{BG:CHM}, but we carry out the details for reader's convenience. 

\begin{definition}\label{def:Heat-Sublaplacian.symbols}
$S^{-\infty}(\URd \times \Theta^{N})$ consists of functions $q(x,\xi,\tau)$ on $\URdd$ such that:\smallskip

(i) $q(x,\xi,\tau)$ is smooth on $\URd \times (\R \setminus 0)$ and for any integer $N'$, any compact $K \subset U$ and any $c>0$ we have 
estimates, 
\begin{equation}
    |\partial_{x}^{\alpha}\partial_{\xi}^{\beta}\partial_{\tau}^{k}q(x,\xi,\tau)|\leq C_{N'Kc\alpha\beta k}(1+|\xi|+|\tau|)^{-N'}, \
    (x,\xi)\in K\times \Rd, |\tau|>c, 
%     \times (\R \setminus [-c,c]),
%     \label{eq:¥}
\end{equation}
with $\alpha$, $\beta$ and $k$ arbitrary.\smallskip

(ii) $q(x,\xi,\tau)$ is smooth with respect to $x$ and $\xi$ and is $C^{N}$ with respect to $\tau$ near $\URd \times 0$  and 
for any integer $N'$, any compact $K \subset U$ we have estimates, 
\begin{equation}
    |\partial_{x}^{\alpha}\partial_{\xi}^{\beta}\partial_{\tau}^{k}q(x,\xi,\tau)|\leq C_{N'K\alpha\beta k}(1+|\xi|+|\tau|)^{-N'},\quad 
    (x,\xi,\tau)\in K\times \Rdd,
%     \label{eq:¥}
\end{equation}
with $\alpha$ and $\beta$ arbitrary and with $k\leq N$.\smallskip

(iii) $q(x,\xi,\tau)$ extends to a continuous function on  $\URd\times \overline{\Theta}$ whose restriction to 
          $\URd  \times\Theta$ belongs to $C^{\infty}(\URd)\hotimes \Hol (\Theta)$ and such that, for any integer $N'$ and any compact $K \subset U$,
%           and any angular sector $\Theta' \subsubset \Theta$ 
          we have estimates,
\begin{equation}
    |\partial_{x}^{\alpha}\partial_{\xi}^{\beta} q(x,\xi,\tau)|\leq C_{N'K\alpha\beta}(1+|\xi|+|\tau|)^{-N'},\quad (x,\xi,\tau)\in K\times \Rd \times\overline{ \Theta},
%     \label{eq:¥}
\end{equation}
with $\alpha$ and $\beta$ arbitrary.\smallskip
\end{definition}

\begin{definition}
    $S_{\ah, \vo}^{m}(\URd \times \Theta^{N})$, $m \in \Z$, consists of functions $q(x,\xi,\tau)$ on $\URdd$ such that:\smallskip
    
    (i) $q(x,\xi,\tau)$ is smooth on $U\times (\Rdd\setminus 0)$ and near the region $U\times(\Rdo)\times 0$ 
         it is smooth with respect to $x$ and $\xi$ and is $C^{N}$ with respect to $\tau$;
         
    (ii) $q(x,\xi,\tau)$ extends to a continuous function on  $U\times [(\Rd\times \overline{\Theta})\setminus 0]$ whose restriction to 
          $\URd  \times\Theta$ belongs to $C^{\infty}(\URd)\hotimes \Hol (\Theta)$.\smallskip
     
    (iii) For any $\lambda \in \R\setminus 0$ the function $q(x,\lambda.\xi,\lambda^{v}\tau)-\lambda^{m}q(x,\xi,\tau)$ belongs to $S^{-\infty}(\URd \times 
    \Theta^{N})$. 
\end{definition}

 Along similar lines as that of the proof of Lemma~\ref{lem:HolPHDO.characterization-almost-homogeneous-symbols}  and of the proof 
 of~\cite[Prop.~3.3]{Po:JAM1} we obtain: 
 
 \begin{lemma}
1)     Let $q\in S_{\ah, \vo}^{m}(\URd \times \Theta^{N})$. Then $q$ admits a unique homogeneous part, i.e., there 
     is a unique symbol $q(x,\xi,\tau)$ in $S_{\vo,m}(\URd \times \Theta_{(v)}^{N})$ such that on $U\times (\Rdd\setminus 0)$ we have
     \begin{equation}
         q_{m}(x,\xi,\tau)=\lim_{\lambda \rightarrow 0}\lambda^{-m}q(x,\lambda.\xi,\lambda^{v}\tau).
                   \label{eq:Heat-Sublaplacian.homogeneous-part}
     \end{equation}

 2) Let $q_{m}\in S_{\vo,m}(\URd \times \Theta_{(v)}^{N})$. Then there exists $q\in S_{\ah, \vo}^{m}(\URd \times \Theta^{N})$ with homogeneous 
 part $q_{m}$. Moreover $q$ is unique up to the addition of a symbol in $S^{-\infty}(\URd \times \Theta^{N})$.
 \end{lemma}

 In particular, this lemma implies that~(\ref{eq:Heat-Sublaplacian.homogeneous-part}) gives rise to a linear isomorphism, 
 \begin{equation}
     S_{\vo,m}(\URd \times \Theta_{(v)}^{N})\simeq S_{\ah, \vo}^{m}(\URd \times \Theta^{N})/S^{-\infty}(\URd \times \Theta^{N}).
     \label{eq:Heat-Sublaplacian.identification-symbols}
 \end{equation}
 
Let $q \in S_{\ah, \vo}^{m}(\URd \times \Theta^{N})$. For any compact $K \subset U$ we have estimates,
\begin{equation}
    |\partial_{x}^{\alpha}\partial_{\xi}^{\beta}q(x,\xi,\tau)|\leq C_{K\alpha\beta}(1+\|\xi\|+|\tau|^{v}), \quad (x,\xi,\tau)K\times \Rd \times 
    \overline{\Theta},
%     \label{eq:¥}
\end{equation}
Combining this with the inequalities, 
\begin{gather}
    (1+|\tau|^{v})^{\frac{1}{2}}(1+\|\xi\|)^{\frac{1}{2}}\leq 1+\|\xi\|+|\tau|^{v} \leq (1+\|\xi\|)(1+|\tau|^{v}),
    \label{eq:Heat-Sublaplacians.inequalities}\\
    ( 1+\|\xi\|+|\tau|^{v})^{m} \leq (1+\|\xi\|)^{m}(1+|\tau|^{v})^{|m|},
\end{gather}
the latter being the Peetre's inequality, we see that we can regard $(q(.,.,\tau))_{\tau \in \overline{\Theta}}$ as a continuous family with values 
in  $S_{\|}^{m}(\URd)$ which is a $\op{O}(|\tau|^{|m|})$ in this Fr\'echet space as $|\tau|\rightarrow \infty$. 

In fact, this family is also holomorphic 
on $\Theta$,  smooth on $\R\setminus 0$ and $C^{N}$ near $\tau=0$ and its $\tau$-derivatives too are $\op{O}(|\tau|^{|m|})$ in $S_{\|}^{m}(\URd)$ as $\tau$ 
becomes large. 

Recall that the convolution products on the groups $G^{(x)}$, $x\in U$, give rise to a smooth family of bilinear products,
\begin{equation}
    *^{(x)}:S_{\|}^{m_{1}}(\Rd)\times S_{\|}^{m_{2}}(\Rd) \longrightarrow S_{\|}^{m_{1}+m_{2}}(\Rd).
%     \label{eq:¥}
\end{equation}
Therefore, if $q_{j}\in S_{\ah, \vo}^{m_{j}}(\URd \times \Theta^{N})$ then we define a family with values in 
$S_{\ah, \vo}^{m_{1}+m_{2}}(\URd \times \Theta^{N})$ by letting 
\begin{equation}
    q_{1}*q_{2}(x,\xi,\tau)=(q_{1}(x,.,\tau)*^{(x)}q_{2}(x,.,\tau))(\xi), \quad (x,\xi,\tau)\in \URd \times \overline{\Theta}.
     \label{eq:Heat-Sublaplacians.product}
\end{equation}
This family is continuous on $\overline{\Theta}$, is holomorphic on $\Theta$, is smooth on $\R\setminus 0$ and is $C^{N}$ near $\tau=0$. Moreover, along with 
its derivatives it is a $\op{O}(|\tau|^{|m|})$ in $S_{\|}^{m}(\URd)$ as $\tau$ becomes large. 

On the other hand, it also follows from the inequalities~(\ref{eq:Heat-Sublaplacians.inequalities}) 
that if $q$ is a symbol in $S^{-\infty}(\URd \times \Theta_{(v)}^{N})$ then we can regard 
$(q(.,.,\tau))_{\tau \in \overline{\Theta}}$ as a continuous family with values in $S^{-\infty}(\URd)$ which is holomorphic on $\Theta$, is smooth on $\R\setminus 0$, 
is $C^{N}$ near $\tau=0$ and which, for any integer $N'$,  together with  its derivatives is a $\op{O}(|\tau|^{-N'})$ in 
$S^{-\infty}(\URd)$  as $\tau$ becomes large. Therefore, if in~(\ref{eq:Heat-Sublaplacians.product}) 
we replace $q_{1}$ or $q_{2}$ by an element in 
$S^{-\infty}(\URd \times \Theta_{(v)}^{N})$ then the 
resulting symbol belongs to $S^{-\infty}(\URd \times \Theta_{(v)}^{N})$. 

Now, for $j=1,2$ let $q_{j}\in S_{\ah, \vo}^{m_{j}}(\URd \times \Theta^{N})$. Then thanks to~(\ref{eq:HolPHDO.homogeneity-*product}) 
for any $\lambda \in \R\setminus 0$  
the symbol $q_{1}*q_{2}(x,\lambda.\xi,\lambda^{v}\tau)- 
\lambda^{m_{1}+m_{2}}q_{1}*q_{2}(x,\xi,\tau)$ is equal to 
\begin{multline}
   [q_{1}(x,\lambda.\xi,\lambda^{v}\tau)-\lambda^{m_{1}}q_{1}(x,\xi,\tau)]*q_{2}(x,\lambda.\xi,\lambda^{v}\tau) \\ + 
   \lambda^{m_{1}}q_{1}(x,\xi,\tau)*[q_{2}(x,\lambda.\xi,\lambda^{v}\tau)-\lambda^{m_{2}}q_{2}(x,\xi,\tau)], 
\end{multline}
and so belongs to $S^{-\infty}(\URd \times \Theta_{(v)}^{N})$ by the observations above. This shows that $q_{1}*q_{2}$ belongs to 
$S_{\ah, \vo}^{m_{1}+m_{2}}(\URd \times \Theta^{N})$.  Therefore, the formula~(\ref{eq:Heat-Sublaplacians.product}) 
defines a bilinear product $*$ from 
 $S_{\ah, \vo}^{m_{1}}(\URd \times \Theta^{N})\times S_{\ah, \vo}^{m_{2}}(\URd \times \Theta^{N})$ to $S_{\ah, \vo}^{m_{1}+m_{2}}(\URd \times 
     \Theta^{N})$. 
     
Moreover, if $\tilde{q}_{j}\in S_{\ah, \vo}^{m_{j}}(\URd \times \Theta^{N})$ has the same homogeneous part as that of $q_{j}$ then thanks to the equality 
$q_{1}*q_{2}-\tilde{q}_{1}*\tilde{q}_{2}=(q_{1}-\tilde{q}_{1})*q_{1}+q_{1}*(q_{2}-\tilde{q}_{2})$ we see that $q_{1}*q_{2}$ and 
$\tilde{q}_{1}*\tilde{q}_{2}$ agree up to an element in $S^{-\infty}(\URd \times \Theta_{(v)}^{N})$ and so have same homogeneous part. Therefore, 
using the isomorphism~(\ref{eq:Heat-Sublaplacian.identification-symbols}) we get: 

\begin{lemma}\label{lem:Heat-Sublaplacian.product-symbols}
   The convolution products on the groups $G^{(a)}$, $a \in U$, give rise to a bilinear product $*$ from $S_{\vo,m_{1}}(\URd \times \Theta_{(v)}^{N})\times 
   S_{\vo,m_{2}}(\URd \times \Theta_{(v)}^{N})$ to $S_{\vo,m_{1}+m_{2}}(\URd \times \Theta_{(v)}^{N})$. 
\end{lemma}

\begin{definition}
For $m \in \Z$ and $N\in \N\cup \{\infty\}$ the space $S_{\vo, m}(\fg^{*}M\times \Theta^{N}_{(v)},\cE)$ consists of sections $q(x,\xi,\tau)$ of $\cE$ over 
$(\fg^{*}M\times \R)\setminus 0$ such that:\smallskip

    (i) $q(x,\xi,\tau)$ is $C^{\infty}$ on $\fg^{*}M\times (\R\setminus 0)$ and near the region $\fg^{*}M\times 0$ 
         it is $C^{\infty}$ with respect to $x$ and $\xi$ and $C^{N}$ with respect to $\tau$;\smallskip 
    
    (ii) We have $q(\lambda.\xi,\lambda^{v}\tau)=\lambda^{m} q(\xi,\tau)$ for any $\lambda \in \R\setminus 0$;\smallskip
    
    (iii) $q(x,\xi,\tau)$ extends to a section of $\pi^{*}\End  \cE$ over $(\fg^{*}M\times  \overline{\Theta})\setminus 0$ which is 
 smooth with respect to the base space variable and continous with 
 respect to the others and which restricts on $\fg^{*}M\times \C_{-}$ to an element of 
 $C^{\infty}(\fg^{*}M, \pi^{*}\End  \cE)\hotimes \Hol(\Theta)$.
\end{definition}

% We also define the symbol classes  $S_{\vo,m}(\URd \times \Theta_{(v)})$ and $S_{\vo,m}(\fg^{*} \times \Theta_{(v)},\cE)$ 
% as in Definitions~\ref{def:Volterra-symbols-URd} and~\ref{def:Volterra-symbols-fgM}. 
% We then get a bilinear product, 
% \begin{equation}
%     S_{\vo,m_{1}}(\fg^{*} M \times \Theta_{(v)},\cE)\times S_{\vo,m_{2}}(\fg^{*} M \times \Theta_{(v)},\cE) 
%     \longrightarrow S_{\vo,m_{1}+m_{2}}(\fg^{*} M \times \Theta_{(v)},\cE).
% %     \label{eq:¥}
% \end{equation}
For instance, if $P:C^{\infty}(M,\cE)\rightarrow C^{\infty}(M,\cE)$ is a differential operator of Heisenberg order $v$ then 
the symbol $\sigma_{v}(P)(x,\xi)+i\tau$ belongs to $S_{\vo,v}(\fg^{*} \times \Theta_{(v)}^{\infty},\cE)$ for any open angular 
sector $\Theta $ whose closure containing $\R$.    

Along the same lines as that of the proof of Proposition~\ref{prop:Heisenberg.product-symbols-manifold},  
using~(\ref{eq:PsiHDO.global-local-convolution-symbols}) and Lemma~\ref{lem:Heat-Sublaplacian.product-symbols} we obtain: 

\begin{lemma}
   The convolution products on the tangent groups $G_{a}M$, $a \in M$,  give rise to a bilinear product $*$ from 
   $S_{\vo,m_{1}}(\fg^{*}M\times \Theta^{N}_{(v)},\cE)\times 
   S_{\vo,m_{2}}(\fg^{*}M\times \Theta^{N}_{(v)},\cE)$ to $S_{\vo,m_{1}+m_{2}}(\fg^{*}M\times \Theta^{N}_{(v)},\cE)$. 
\end{lemma}

In particular,  for a symbol $q \in S_{\vo, m}(\fg^{*}M\times \Theta^{N}_{(v)},\cE)$ it makes sense to speak about its inverse in $S_{\vo, -m}(\fg^{*}M\times 
\Theta^{N}_{(v)},\cE)$.  

Next, let $U\subset \Rd$ be a trivializing Heisenberg chart together with a $H$-frame $X_{0},\ldots,X_{d}$ and 
let $L(x)=(L_{jk}(x))$ be the matrix of the Levi form 
$\cL$ with respect to this $H$-frame, so that for $j,k=1,\ldots,d$ we have
\begin{equation}
    \cL(X_{j},X_{k})=[X_{j},X_{k}]=L_{jk}X_{0} \quad \bmod H.
%     \nonum
\end{equation}
Then we have the following extension of Theorem~5.22 of~\cite{BGS:HECRM}. 

\begin{proposition}\label{prop:heat-sublaplacian-invertibility-system}
 For $\omega \in (-\frac{\pi}{2},\frac{\pi}{2})$ define 
  \begin{equation}
      \Omega_{\omega}=\{(\mu,x)\in M_{r}(\C) \setminus U; \ [(\cos \omega)^{-1}\Re (e^{i\omega}\Sp \mu)]\cap  \Lambda_{x}=\emptyset  \},
      \label{eq:Heat-Sublaplacian.Omega-w}
%       \Theta_{\omega}=\{\tau \in \C; \ \Im (e^{-i\omega}\tau) <0\},
  \end{equation}
 where the singular set $\Lambda_{x}$ is defined as in~(\ref{eq:Sublaplacian.singular-set1})--(\ref{eq:Sublaplacian.singular-set2}). 
Then $\Omega_{\omega}$ is  an open set and there exists 
   $q_{\mu}^{(\omega)}(x,\xi,\tau)\in C^{\infty}(\Omega_{\omega}, S_{\vo,-2}(\Rdd_{(2)},\C^{r}))$ such that:\smallskip
    
   (i) $q_{\mu}^{(\omega)}(x,\xi,\mu)$ depends analytically on $\mu$;\smallskip
 
   (ii) For any $(\mu,x)\in \Omega$ the symbol $q_{\mu}^{(\omega)}(x,.,.)$ inverts $|\xi'|^{2}+i\mu\xi_{0}+ie^{-i\omega}\tau$, i.e.,
    \begin{equation}
        q_{\mu}^{(\omega)}(x,.)*^{x}(|\xi'|^{2}+i\mu \xi_{0}+ie^{-i\omega}\tau)=(|\xi'|^{2}+i\mu \xi_{0}+ie^{-i\omega}\tau)*^{x}q_{\mu}^{(\omega)}(x,.)=1.
%         \label{eq:¥}
    \end{equation}
% 
% \noindent Moreover, if $(\mu,x)\in \Omega_{\omega_{1}}\cap \Omega_{\omega_{2}}$ then on 
% $\Rd\times (e^{-i\omega_{1}}\C_{-}\cap e^{-i\omega_{2}}\C_{-})$ we have 
% 
\noindent Moreover, if $(\mu,x)\in \Omega_{\omega_{1}}\cap \Omega_{\omega_{2}}$ then we have 
% and $\tau \in $, we have 
\begin{equation}
    q_{\mu}^{(\omega_{1})}(x,\xi,e^{i\omega_{1}}\tau)= q_{\mu}^{(\omega_{2})}(x,\xi,e^{i\omega_{2}}\tau),
     \label{eq:Heat-Sublaplacian.equality-inverse-w1-w2-1}
\end{equation}
for any $\xi$ in $\Rd$ and any $\tau$ in $e^{-i\omega_{1}}\C_{-}\cap e^{-i\omega_{2}}\C_{-}$
\end{proposition}
\begin{proof}
 First, in the same way as in the proof of 
 Proposition~\ref{prop:Sublaplacian.inverse.system} we can show that $\Omega_{\omega}$ is an open subset of $M_{r}(\C)\times U$. 
 
 Second, let us assume that $r=1$. For $x \in U$ define 
\begin{equation}
        \Lambda_{x}^{0}= (-\infty, -\frac12 \Tra |L(a)|]\cup [\frac12 \Tra  |L(a)|,\infty)
%     \label{eq:}
\end{equation}
    For $k=1,2,\ldots$ let $\Lambda_{x}^{k}=\Lambda_{x}^{0}$ if $\rk \cL_{x}< d$ and otherwise let 
    \begin{equation}
        \Lambda_{x}^{k}=\{\pm(\frac12 \Tra |L(x)|+2\sum_{1\leq j \leq n}\alpha_{j}|\lambda_{j}|); \alpha_{j}\in \N, \ \alpha_{j}\leq k\}.
%         \label{eq:¥}
    \end{equation}
   Then for $k=0,1,\ldots$ we let $\Omega^{k}_{\omega}$ denote the subset of $\C\times U$ defined  as in~(\ref{eq:Heat-Sublaplacian.Omega-w}) 
   with $r=1$ 
  and with $\Omega_{\omega}$ replaced by $\Omega_{\omega}^{k}$. Again $\Omega_{\omega}^{k}$ is open and 
  when $\rk \cL_{x}<d$ for every $x \in U$ we have 
   $\Omega_{\omega}=\Omega_{\omega}^{0}=\Omega_{\omega}^{k}$. 
   
   We are now going to prove: 
   
   \begin{claim}
      For $k=0,1,\ldots$ there exists 
   $q_{\mu}^{(\omega)}(x,\xi,\tau)$ in $C^{\infty}(\Omega_{\omega}^{k}, S_{\vo,-2}(\Rdd_{(2)}))$  satisfying (i) and (ii) on $\Omega_{\omega}^{k}$. 
   \end{claim}
   
   For $\omega=0$ and $k=0$ the claim follows from~\cite[Thm.~5.22]{BGS:HECRM} and in this case the symbol 
   $q_{\mu}^{(0)}\in C^{\infty}(\Omega_{0}^{0}, S_{\vo,-2}(\Rdd_{(2)}))$  is given by the formula,
    \begin{equation}
           q_{\mu}^{(0)}(x,\xi,\tau) =\int_{0}^{\infty}e^{-t\mu \xi_{0}-it\tau}G(x,\xi,t)dt, 
        \label{eq:Heat-Sublaplacian.heat-kernel-inverse0}
    \end{equation}
    where $G(x,\xi,t)$ is as in~(\ref{eq:Rockland-Sublaplacian.heat-kernel-inverse}). 
    
    Similarly, for $ \omega \neq 0$ we define a symbol $q_{\mu}^{(\omega)}(x,\xi,\tau)\in 
    C^{\infty}(\Omega_{\omega}^{0}, S_{\vo,-2}(\Rdd_{(2)}))$  satisfying (i) and (ii) by letting 
        \begin{equation}
           q_{\mu}^{(0)}(x,\xi,\tau) =\int_{e^{i\omega}(0,\infty)}e^{-t\mu \xi_{0}-ite^{-i\omega}\tau}G(x,\xi,t)dt. 
        \label{eq:Heat-Sublaplacian.heat-kernel-inverse}
    \end{equation}
 Notice also that thanks to the analyticity of $G(\xi,t)$ with respect to $t$,  if 
 $(\mu,x)$ is in $\Omega_{\omega_{1}}^{0}\cap \Omega_{\omega_{2}}^{0}$ then, for any $\xi\in \Rd$ and any $\tau \in 
e^{-i\omega_{1}}\C_{-}\cap e^{-i\omega_{2}}\C_{-}$, we have 
\begin{equation}
    q_{\mu}^{(\omega_{1})}(x,\xi,e^{i\omega_{1}}\tau)= q_{\mu}^{(\omega_{2})}(x,\xi,e^{i\omega_{2}}\tau).
     \label{eq:Heat-Sublaplacian.equality-inverse-w1-w2}
\end{equation}
   
For $k\geq 1$ the claim can be proved by making integration by parts in the integral~(\ref{eq:Heat-Sublaplacian.heat-kernel-inverse}) as 
in~\cite{BG:CHM}. Moreover, as in the case $k=0$ if  $(\mu,x)\in \Omega_{\omega_{1}}^{k}\cap \Omega_{\omega_{2}}^{k}$ then the 
equality~(\ref{eq:Heat-Sublaplacian.equality-inverse-w1-w2}) 
holds for any $\xi\in \Rd$ and any $\tau \in e^{-i\omega_{1}}\C_{-}\cap e^{-i\omega_{2}}\C_{-}$.    
    
All this shows that the lemma is true for $r=1$. The case $r\geq 2$ is then deduced from the case $r=1$ by  arguing as in the proof of 
 Proposition~\ref{prop:Sublaplacian.inverse.system}. 
\end{proof}

In the sequel given open angular sectors $\Theta_{1}$ and $\Theta_{2}$ we write $\Theta_{1} \subsubset \Theta_{2}$ to mean that 
$\overline{\Theta_{1}}\setminus 0$ is contained in $\Theta_{2}$.

\begin{proposition}\label{thm:Powers1.heat-sublaplacians}
Let $\Delta: C^{\infty}(M,\cE) \rightarrow C^{\infty}(M,\cE)$ be a selfadjoint sublaplacian which is bounded from below and satisfies 
the condition~(\ref{eq:Sublaplacian.condition}) at every point. Then the principal symbol of $\Delta+\partial_{t}$ admits an inverse in  
$S_{\vo,-2}(\fg^{*} M \times \Theta_{(2)}^{\infty},\cE)$ for any open angular sector $\Theta \subsubset \C\setminus i[0,\infty)$ containing $\R$.
% invertible Volterra-Heisenberg symbol, 
% hence Proposition~\ref{thm:volterra.inverse} and Proposition~\ref{thm:Powers1.heat-kernel-asymptotics} hold for $\Delta$.   
\end{proposition}
\begin{proof}
  It is enough to prove the proposition in a local trivializing Heisenberg chart $U\subset \Rd$ with a $H$-frame $X_{0},\ldots,X_{d}$ with respect to which 
  $\Delta$ is of the form, 
  \begin{equation}
      \Delta=-\sum_{j=1}^{d} X_{j}^{2} - i\mu(x) X_{0}+ \op{O}_{H}(1).
      \label{eq:Heat.sublaplacian.bundle2}
  \end{equation}
  
  Since $\Delta$ is selfadjoint the matrix $\mu(x)$ is selfadjoint for every $x \in U$ and so the condition~(\ref{eq:Sublaplacian.condition}) implies 
  that, with the notation of Proposition~\ref{prop:heat-sublaplacian-invertibility-system}, 
  the pair $(x,\mu(x))$ is in $\Omega_{\omega}$ for any $x \in U$ and any $\omega \in (-\frac{\pi}{2},\frac{\pi}{2})$. Therefore, thanks to 
   Proposition~\ref{prop:heat-sublaplacian-invertibility-system} we define an inverse $q_{-2}\in S_{\vo,-2}(\URdd_{(2)},\C^{r})$ for 
    the principal symbol $|\xi'|^{2}+i\mu(x)\xi_{0}+i\tau$ of $\Delta$ by letting
  \begin{equation}
      q_{-2}(x,\xi,\tau)=q_{\mu(x)}^{(0)}(x,\xi,\tau),\quad  (x,\xi,\tau)\in U\times [(\Rd \times \overline{\C_{-}})\setminus 0].
%       \label{eq:¥}
  \end{equation}
  
  For $\omega \in (-\frac{\pi}{2},\frac{\pi}{2})$ let $\Theta_{\omega}=\C_{-}\cup (e^{-i\omega}\C_{-})$. 
  Since~(\ref{eq:Heat-Sublaplacian.equality-inverse-w1-w2-1}) shows that $q_{-2}(x,\xi,\tau)=q_{\mu(x)}^{(\omega)}(x,\xi,e^{i\omega}\tau)$ when $\tau$ 
  is in $\C_{-}\cap (e^{-i\omega}\C_{-})$, we see that we can extend the definition of $q_{-2}$ to $\URd \times \Theta_{\omega}$ by letting 
  \begin{equation}
       q_{-2}(x,\xi,\tau)=q_{\mu(x)}^{(\omega)}(x,\xi,e^{i\omega}\tau),\quad  (x,\xi,\tau)\in \URd \times (e^{-i\omega}\C_{-}).
%       \label{eq:}
  \end{equation}
  This defines an inverse for $|\xi'|^{2}+i\mu(x)\xi_{0}+i\tau$ in $S_{\vo,-2}(\URd\times \Theta_{\omega(2)}^{\infty},\C^{r})$. 
  The proof is now completed by noticing that any angular sector $\Theta \subsubset \C\setminus i[0,\infty)$ containing $\R$ is of the form 
  $\Theta=\Theta_{\omega_{1}}\cup \Theta_{\omega_{2}}$ with $-\frac{\pi}{2}<\omega_{1}<0<\omega_{2}<\frac{\pi}{2}$.  %Therefore, 
\end{proof}

\begin{remark}
    It can be shown that any selfadjoint sublaplacian with an invertible principal symbol  is bounded from below (see 
    \cite{Po:PhD},~\cite{Po:CPDE1}). Therefore,  the assumption on the boundedness from below of $\Delta$ in 
    Proposition~\ref{thm:Powers1.heat-sublaplacians} in not necessary.
\end{remark}

\begin{example}
 Proposition~\ref{thm:Powers1.heat-sublaplacians} is true for the following sublaplacians:\smallskip 
  
  (a) A selfadjoint sum of squares $\Delta_{\nabla,X}=\nabla_{X_{1}}\nabla_{X_{1}}^{*}+\ldots+\nabla_{X_{1}}\nabla_{X_{1}}^{*}$, 
  where $\nabla$ is a connection on $\cE$ and the vector fields $X_{1},\ldots,X_{m}$ span $H$, under the bracket condition $H+[H,H]=TM$; \smallskip 
  
  (b) The Kohn Laplacian on a CR manifold acting on $(p,q)$-forms under the condition $Y(q)$;\smallskip 
  
  (c) The horizontal sublaplacian on a Heisenberg manifold acting on horizontal forms of degree $k$ under the condition $X(k)$;\smallskip 
   
  (d) The horizontal sublaplacian on a CR manifold acting on $(p,q)$-forms under the condition $X(p,q)$.\smallskip
 
  In particular, Proposition~\ref{thm:Powers1.heat-sublaplacians} allows us to complete the treatment of the heat kernel of the 
  Kohn Laplacian in~\cite{BGS:HECRM} because, as with the invertibility of its principal symbol in Section~\ref{sec:sublaplacian}, we really need the version for 
  systems provided by Proposition~\ref{prop:heat-sublaplacian-invertibility-system}, but not established in~\cite{BGS:HECRM}.  
\end{example}

Next, assume that the bracket condition $H+[H,H]=TM$ holds, that is, the Levi form of $(M,H)$ has positive rank everywhere, and consider the 
sum of squares, 
\begin{equation}
    \Delta_{\nabla,X}=-(\nabla_{X_{1}}^{*}\nabla_{X_{1}}+\ldots+\nabla_{X_{m}}^{*}\nabla_{X_{m}}),
%     \label{eq:¥}
\end{equation}
where $\nabla$ is a connection on $\cE$ and the vector fields $X_{1},\ldots, X_{m}$ span $H$.

\begin{proposition}\label{prop:Heat1.Delta^k}
    The principal symbol of $\Delta^{k}_{\nabla,X}+\partial_{t}$ admits an inverse in $S_{\op{v},-2k}(\fg^{*}M\times \R_{(2k)},\cE)$ and so 
 Theorems~\ref{thm:volterra.inverse} and~\ref{thm:Powers1.heat-kernel-asymptotics} hold for $\Delta^{k}_{\nabla,X}$.
\end{proposition}
\begin{proof}
  Since proving the above statement is a purely local issue and $\Delta_{\nabla,X}$ is a scalar operator modulo lower order terms, we may assume that $\cE$ is 
  the trivial line bundle and proceed on $U\times \R$ where $U$ is  a Heisenberg chart with $H$-frame $Y_{0},\ldots,Y_{d}$ with respect to  which 
  $\Delta_{\nabla,X}$ takes the form, 
  \begin{equation}
      \Delta_{\nabla,X}=-(Y_{1}^{2}+\ldots+Y_{d}^{2})+\op{O}_{H}(1).
       \label{eq:Heat.Rockland.local-form-Delta}
  \end{equation}
 In particular the principal symbol of $\Delta$ is just  $|\xi'|^{2}$. In addition, let $p_{2k}(x,\xi)$ be the local symbol of $\Delta^{k}_{\nabla,X}$ on $U$. 

  Let $\arg z$ be the continuous determination of the argument on $\C\setminus 0$ with values in $[0, 2\pi)$ and for $z \in \C\setminus 
    0$ let $z^{\frac{1}{k}}=|z|^{\frac{1}{k}}e^{\frac{i}{k}\arg z}$ and set $0^{\frac{1}{k}}=0$, so that the function $z 
    \rightarrow z^{\frac{1}{k}}$ is analytic on $\C \setminus [0,\infty)$ and is continuous at the origin. Set $\omega=e^{\frac{2i\pi}{k}}$. Then for any $z 
    \in \C$ we have the polynomial identity,
    \begin{equation}
       T^{k}-z=(T-z^{\frac{1}{k}}) (T-\omega z^{\frac{1}{k}})\ldots (T-\omega^{k-1}z^{\frac{1}{k}}). 
%         \label{eq:}
    \end{equation}
    Therefore, for $T=|\xi'|^{2}$ and $z=-i\tau$ in $S_{\vo,-2k}(\URdd_{(2k)})$ we get
    \begin{equation}
       p_{2k}+i\tau= (|\xi'|^{2}+i(i(-i\tau)^{\frac{1}{k}})*\ldots* (|\xi'|^{2}+i(i\omega^{k-1}(-i\tau)^{\frac{1}{k}})).
        \label{eq:Heat-sublaplacian.p-2k+tau}
    \end{equation}

Let $\Theta$ be the angular sector $|\arg (-i\tau)|>\frac{\pi}{2k}$. As $\Theta \subsubset \C\setminus i[0,\infty)$ and we have $H+[H,H]=TM$, so 
that the condition~(\ref{eq:Sublaplacian.condition}) for $\Delta_{\nabla,X}$ is satisfied at every point, Proposition~\ref{thm:Powers1.heat-sublaplacians} tells us 
that $|\xi'|^{2}+i\tau$ admits an inverse $q_{(-2)}(x,\xi,\tau)$ in $S_{\vo,-2}(\URd \times \Theta_{(2)})$. 

Notice that for $j=0,1,\ldots,k-1$ if $\Im \tau<0$ then $i\omega^{j}(-i\tau)^{\frac{1}{k}}$ is in $\Theta$, so for $(x,\xi,\tau)\in U\times 
[(\Rd \times \C_{-})\setminus 0]$ we can let 
\begin{equation}
    \tilde{q}_{(-2),j}(x,\xi,\tau)=q_{(-2)}(x,\xi,i\omega^{j}(-i\tau)^{\frac{1}{k}})).
%     \label{eq:¥}
\end{equation}
In the sense of Definition~\ref{def:Heat-Sublaplacian.symbols} this gives rise to a symbol in $S_{\vo,-2}(\URd \times \C_{-(2k)}^{0})$. Therefore, it follows 
from~(\ref{eq:Heat-sublaplacian.p-2k+tau}) and 
Lemma~\ref{lem:Heat-Sublaplacian.product-symbols} that we get an inverse for $ p_{2k}+i\tau$ in $S_{\vo,-v}(\URd \times \C_{-(2k)}^{0})$ by letting
% Although  $\tilde{q}_{(-2),j}(x,\xi,\tau)$ has the right homogeneity, this not does not define an element of $S_{\vo,-2}(\URdd_{(v)})$ because of the 
% lack of smoothness of $\tau^{\frac{1}{l}}$ near the origin. Yet, thanks to~(\ref) we can expect that the inverse of $p_{v}+i\tau$ should be 
\begin{equation}
    q_{(-2k)}=\tilde{q}_{(-2),0}*\ldots*\tilde{q}_{(-2),k-1}.
%     \label{eq:}
\end{equation}

Let us now show that $q_{(-2k)}$ is in $S_{\vo,-2k}(\URd \times \C_{-(2k)}^{\infty})=S_{\vo,-2k}(\URdd_{(2k)})$. 

\begin{claim}
    For $j=1,2$ let $q_{j}\in S_{\vo,m_{j}}(\URd \times \C_{-(2k)}^{0})$. Then, regarding $q_{1}$, $q_{2}$ and $q_{1}*q_{2}$ as smooth families with 
    values in $C^{\infty}(\URd)$ over $\R\setminus 0$, we have 
    \begin{equation}
        \partial_{\tau}(q_{1}*q_{2})=(\partial_{\tau}q_{1})*q_{2}+q_{1}*(\partial_{\tau}q_{2}).
         \label{eq:Heat-sublaplacian.Leibniz-formula}
    \end{equation}
\end{claim}
\begin{proof}[Proof of the claim]
 It follows from~(\ref{eq:Heat-Sublaplacians.product}) that the Leibniz formula~(\ref{eq:Heat-sublaplacian.Leibniz-formula}) 
 holds for almost homogeneous symbols in $S_{\ah, \vo}^{*}(\URd \times \C_{-(2k)}^{1})$. Since 
 $\partial_{\tau}$ and $*$ are homogeneous maps the formula remains true for homogeneous symbols in $S_{\vo,*}(\URd \times \C_{-(2k)}^{1})$. 
\end{proof}

Now, differentiating the equality $1=(p_{2k}+i\tau)*q_{(-2k)}$ using~(\ref{eq:Heat-sublaplacian.Leibniz-formula}) 
we get $ 0= (p_{2k}+i\tau)*(\partial_{\tau}q_{(-2k)}) + iq_{(-2k)}$, 
which after multiplication by $q_{(-2k)}$ gives
\begin{equation}
   \partial_{\tau}q_{(-2k)}= -iq_{(-2k)}*q_{(-2k)}. 
%     \label{eq:¥}
\end{equation}
This equality holds in $C^{\infty}(\R\setminus 0, C^{\infty}(\URd))$ but, since $q_{(-2k)}*q_{(-2k)}$ belongs to $S_{\vo,-2v}(\URd \times 
\C_{-(2k)}^{0})$, this actually shows that the symbol $q_{(-2k)}$ belongs to $S_{\vo,-2k}(\URd \times \C_{-(2k)}^{1})$. 

Finally, an induction shows that 
$q_{(-2k)}$ is contained in $S_{\vo,-2k}(\URd \times \C_{-(2k)}^{N})$ for any integer $N$, hence belongs to $S_{\vo,-2k}(\URdd_{(2k)})$. This proves 
that the principal symbol $p_{(2k)}$ of $\Delta^{k}+\partial_{t}$ has an inverse in $S_{\vo,-2k}(\URd \times \C_{-(2k)}^{1})$. The proof is thus achieved. 
\end{proof}

\section{Complex powers of hypoelliptic differential operators}
\label{sec.powers1}
In this section we show that the complex powers of a positive 
hypoelliptic differential operator, \emph{a priori} defined as unbounded operators on $L^{2}(M,\cE)$, give rise to a holomorphic family of \psivdos. 

Let $P:C^{\infty}(M,\cE)\rightarrow C^{\infty}(M,\cE)$ be a selfadjoint differential operator of 
even (Heisenberg) order $v$ such that $P$ has an invertible principal symbol and is positive, i.e.,~we have $ \acou{Pu}{u}\geq 0$ for any $u\in C^{\infty}(M,\cE)$.

Let $\Pi_{0}(P)$ be the orthogonal projection onto $\ker P$ and set $P_{0}=(1-\Pi_{0}(P))P+\Pi_{0}(P)$. Then $P_{0}$ is 
selfadjoint with spectrum contained in $[c,\infty)$ for some $c>0$. Thus by standard functional calculus, for any $s\in \C$, the power 
$P_{0}^{s}$ is a well defined unbounded operator on $L^{2}(M,\cE)$. We then define the power $P^{s}$, $s \in \C$, by letting 
\begin{equation}
    P^{s}=(1-\Pi_{0}(P))P_{0}^{s}=P_{0}^{s}-\Pi_{0}(P),
     \label{eq:Powers1.definition}
\end{equation}
so that  $P^{s}$ coincides with $P_{0}^{s}$ on $(\ker P)^{\perp}$ and is zero on $\ker 
P$.  In particular, we have $P^{0}=1-\Pi_{0}(P)$ and $P^{-1}$ is the partial inverse of $P$.

The main result of this section is the following. 

\begin{theorem}\label{thm:Powers1.main} 
Suppose that the principal symbol of $P+\partial_{t}$ admits an inverse in $S_{\vo,-v}(\fg^{*}M\times \R_{(v)},\cE)$. Then:\smallskip 
    
    (i) For any $s \in \C$ the operator $P^{s}$ defined by~(\ref{eq:Powers1.definition}) is a \psivdo\ of order $vs$;\smallskip 
    
    (ii) The family $(P^{s})_{s\in \C}$ forms a holomorphic 1-parameter group of \psivdos. 
\end{theorem}
\begin{proof}
Let us first assume that $\cE$ is the trivial line bundle over $M$, so that $P$ is a scalar operator.  For  $\Re s >0$ the function $x\rightarrow x^{-s}$ is 
bounded on $[0,\infty)$, so the operators $P_{0}^{-s}$ and 
$P^{-s}$ are bounded. Moreover, by the Mellin formula we have
\begin{equation}
     P^{-s}= (1-\Pi_{0}(P))P_{0}^{s}=\frac{1}{\Gamma(s)} \int_{0}^{\infty} t^{s}(1-\Pi_{0}(P))e^{-tP}\frac{dt}{t}. 
     \label{eq:Powers1.}
\end{equation}
This leads us to define  
\begin{equation}
    A_{s}=\int_{0}^{1}t^{s-1}e^{-tP}dt,  \qquad \Re s>0.
     \label{eq:Powers1.Ds}
\end{equation}
Then we have 
\begin{equation}
    \begin{split}
      \Gamma(s)P^{-s}-A_{s} & 
      = \int_{0}^{1} t^{s-1}\Pi_{0}(P)e^{-tP}dt + \int_{1}^{\infty} t^{s-1} (1-\Pi_{0}(P)) e^{-tP} dt,  \\  
      &=\frac{1}{2}\Pi_{0}(P)+ e^{-P/2} (\int_{0}^{\infty} (1+t)^{s-1}e^{-tP}dt)e^{-P/2}.
     \end{split}
\end{equation}
Since $\Pi_{0}(P)$ and $e^{-P/2}$ are smoothing operators and  $(\int_{0}^{\infty} (1+t)^{s-1}e^{-tP}dt)_{\Re s>0}$ is a 
holomorphic family of bounded operators on $L^{2}(M)$, we get
\begin{equation}
     (\Gamma(s)P^{-s}-A_{s})_{\Re s>0} \in \Hol(\Re s>0, \psinf(M)).
    \label{eq:Powers1.relation-As}
\end{equation}

Let us now show that $(A_{s})_{\Re s>0}$ defined by~(\ref{eq:Powers1.Ds}) is a holomorphic family of \psivdos\ such that $\ord A_{s}=-vs$. To this  end 
observe that, in terms of distribution kernels, the formula~(\ref{eq:Powers1.Ds}) means that $A_{s}$ has distribution kernel 

\begin{equation}
    k_{A_{s}}(x,y)= \int_{0}^{1}t^{s-1}k_{t}(x,y)dt.
\end{equation}
where $k_{t}(x,y)$ denotes the heat kernel of $P$.

On the other hand, since $P$ is bounded from below and the principal symbol of $P+\partial_{t}$ is an invertible Volterra-Heisenberg symbol, 
Theorem~\ref{thm:volterra.inverse} 
tells us that $P+\partial_{t}$ has an inverse $Q_{0}:=(P+\partial_{t})^{-1}$ in
$\pvhdo^{-v}(M\times\R_{(v)},\cE)$ and that the distribution kernel $K_{Q_{0}}(x,y,t-s)$ of $Q_{0}$ is related to the heat kernel of $P$ by 
\begin{equation}
    K_{Q_{0}}(x,y,t)=k_{t}(x,y) \qquad \text{for $t>0$}.
\end{equation}
Therefore, for $\Re s>0$ we have
\begin{equation}
      k_{A_{s}}(x,y)= \int_{0}^{1}t^{s-1}K_{Q_{0}}(x,y,t)dt..
      \label{eq:Powers1.kAs}
\end{equation}

Let $\varphi$ and $\psi$ be smooth functions on $M$ 
with disjoint supports. Then using~(\ref{eq:Powers1.kAs}) we see that $\varphi A_{s}\psi$ has distribution kernel 
\begin{equation}
    k_{\varphi A_{s}\psi}(x,y)= \int_{0}^{1}t^{s-1}\varphi(x)K_{Q_{0}}(x,y,t)\psi(y)dt.
     \label{eq:Powers1.phi-kDs-psi}
\end{equation}
Since the distribution kernel of a Volterra-\psivdo\ is smooth off the diagonal of $(M\times \R)\times (M\times \R)$ the distribution $K_{Q_{0}}(x,y,t)$ is smooth 
on 
the region $\{x\neq y\}\times \R$, so~(\ref{eq:Powers1.phi-kDs-psi}) defines a holomorphic family of smooth kernels. Thus, 
\begin{equation}
    (\varphi A_{s}\psi)_{\Re s>0} \in \Hol(\Re s>0,\Psi^{-\infty}(M)).
    \label{eq:Powers1.smoothing}
\end{equation}
 
Next, the following holds. 

\begin{lemma}\label{lem:Powers1.analyticity-Bs}
 Let $V\subset \Rd$ be a Heisenberg chart, let $Q\in \pvhdo^{-v}(V\times\R_{(v)})$ have distribution kernel 
$K_{Q}(x,y,t-s)$ and for $\Re s>0$ let $B_{s}:C_{c}^{\infty}(V)\rightarrow C^{\infty}C(V)$  be given by the distribution kernel,
\begin{equation}
    k_{B_{s}}(x,y)= \int_{0}^{1}t^{s-1}K_{Q}(x,y,t)dt, \quad \Re s>0.
\end{equation} 
Then $(B_{s})_{\Re s >0}$ is a holomorphic family of \psivdos\ such that $\ord B_{s}=-vs$.
\end{lemma}
\begin{proof}[Proof of the lemma]
Let $\varepsilon_{x}$ denote the change to the Heisenberg coordinates at $x$. By Proposition~\ref{prop:volterra.kernel-charaterization} 
on $V\times V\times \R$ the distribution $K_{Q}(x,y,t)$ is of the form
\begin{equation}
    K_{Q}(x,y,t)=|\epsilon_{x}'|K(x,-\varepsilon_{x}(y),t)+R(x,y,t),
\end{equation}
for some $K\in \cK_{\op{v}}^{-(d+2)}(V\times \Rdd_{(v)})$ and some $R\in C^{\infty}(U\times U\times \R)$. Let us write 
$K\sim \sum_{j\geq 0} K_{j-(d+2)}$ 
with $K_{l}\in \cK_{\op{v},l}(V\times \Rdd_{(v)})$. Thus, for any integer $N$, as soon as $J$ large enough we have
\begin{equation}
    K(x,y,t) =\sum_{j\leq J} K_{j-(d+2)}(x,y,t) + R_{NJ}(x,y,t), \quad R_{NJ}\in C^{N}(U\times\R^{d+2}).
\end{equation}
In particular, on $V\times V$ we can write
\begin{equation}
    k_{B_{s}}(x,y)=|\epsilon_{x}'|K_{s}(x,\varepsilon_{x}(y))+R_{s}(x,y), \quad K_{s}(x,y)=\int_{0}^{1} t^{s-1} K (x,y,t)dt,
    \label{eq:Powers1.kernelBs1}
\end{equation}
with $(R_{s})_{\Re s>0}$ in  $\Hol(\Re s>0, C^{\infty}(V \times V))$. 
Moreover, $K_{s}(x,y)$ is of the form 
\begin{equation}
 K_{s}=\sum_{j\leq J}K_{j,s}+R_{NJ,s}, \quad    K_{j,s}(x,y)=\int_{0}^{1} t^{s-1} K_{j-(d+2)} 
    (x,y,t)dt, 
        \label{eq:Powers1.kernelBs2}
\end{equation}
with $(R_{NJ,s})_{\Re s>0}$ contained in $\Hol(\Re s>0, C^{N}(V \times V))$. 
 
Notice that $K_{j-(d+2)}(x,y,t)$ is in $C^{\infty}(V)\hotimes \cD_{\reg}'(\Rd\times \R)$ and is parabolic homogeneous of 
degree $j-(d+2)\geq -(d+2)$. Thus the family $(K_{j,s})_{\Re s>0}$ belongs to 
$\Hol(\Re s>0,C^{\infty}(U)\hotimes \cD_{\reg}'(\Rd))$. 
Moreover, for any $\lambda>0$, the difference $K_{j,s}(x, \lambda.y)-\lambda^{vs+j-(d+2)}K_{j,s}(x,y)$ is equal to 
\begin{equation}
  \int_{1}^{\lambda^{2}}t^{s-1} K(x,y,t) dt \in 
    \Hol(\Re s>0,C^{\infty}(V\times \R^{d+2})).
\end{equation}
Hence $(K_{j,s})_{\Re s>0}$ is a holomorphic family of almost homogeneous distributions of degree $vs-(d+2)+j$. Combining this 
with~(\ref{eq:Powers1.kernelBs2})
then shows that $(K_{s})_{\Re 
s>0}$ belongs to $\Hol(\Re s>0, \cK^{*}(V\times \Rd))$ and has order $vs-(d+2)$. Therefore,
using~(\ref{eq:Powers1.kernelBs1}) and Proposition~\ref{prop:HolPHDO.kernel.characterization}
we see that $(B_{s})_{\Re s>0}$ is a holomorphic family of \psivdos\ such that $\ord B_{s}=-(\ord K_{s}+d+2)=-vs$.
\end{proof}

It follows from Lemma~\ref{lem:Powers1.analyticity-Bs} that for any  local Heisenberg 
chart $\kappa:U\rightarrow V$ the family $(\kappa_{*}A_{s|_{U}})_{\Re s>0}$ is a 
holomorphic family of \psivdos\ on $V$ of order~$-vs$. Combining this with~(\ref{eq:Powers1.smoothing}) and~(\ref{eq:Powers1.relation-As}) then 
shows that $(A_{s})_{\Re s>0}$  and $(P^{s})_{\Re s<0}$ are
holomorphic families of \psivdos\ of orders $-vs$ and $vs$ respectively.

Now,  let $s \in \C$ and let $k$ be a positive  integer such that $k> \Re s$. Then on $C^{\infty}(M,\cE)$ we have 
$P^{s}=P^{s-k}P^{k}$.
Since $P^{k}$ is a differential operator and $P^{s-k}$ is a \psivdo\ of order $v(s-k)$ this shows that 
$P^{s}$ is a \psivdo\ of order $vs$. In fact, as by Proposition~\ref{prop:HolPHDO.composition} the product of \psivdos\ is 
analytic this proves that $(P^{s})_{s \in \C}$ is a holomorphic family of \psivdos\ such that $\ord P^{s}=vs$ for every $s\in \C$.

Now,  let $s \in \C$ and let $k$ be a positive  integer such that $k> \Re s$. Then on $C^{\infty}(M,\cE)$ we have 
\begin{equation}
    P^{s}u=P^{s-k}P^{k}u \quad \text{for any $u \in C^{\infty}(M,\cE)$}.
\end{equation}
As $P^{k}$ is a differential operator and $P^{s-k}$ is a \psivdo\ of order $m(s-k)$ this proves that 
$P^{s}$ is a \psivdo\ of order $ms$. In fact, as by Proposition~\ref{prop:HolPHDO.composition} the product of \psivdos\ is 
analytic this actually shows that $(P^{s})_{s \in \C}$ is a holomorphic family of \psivdos\ such that $\ord P^{s}=vs$ for every $s\in \C$.

Finally, when $\cE$ is a general vector bundle we can similarly prove that the complex powers $P^{s}$, $s \in \C$,  forms a holomorphic family of \psivdos\ 
such that $\ord P^{s}= vs$ for any $s\in \C$. The proof is thus complete. 
\end{proof}

\begin{example}\label{ex:Powers1.sublaplacian}
Thanks to Proposition~\ref{thm:Powers1.heat-sublaplacians} Theorem~\ref{thm:Powers1.main} holds for the following sublaplacians:\smallskip 

    (a) A selfadjoint sums of squares $  \Delta_{\nabla,X}=\nabla_{X_{1}}^{*}\nabla_{X_{1}}+\ldots+\nabla_{X_{m}}^{*}\nabla_{X_{m}}$, 
    where  $X_{1},\ldots,X_{m}$ span $H$ and $\nabla$ is a connection on $\cE$, under the condition that the Levi form is nonvanishing.\smallskip 

    (b) The Kohn Laplacian  $\Boxbpq$ on a CR manifold acting on $(p,q)$-forms under the condition~$Y(q)$.\smallskip 
    
   (c) The  horizontal sublaplacian $\Delta_{b;k}$ on a Heisenberg manifold $(M^{d+1}, H)$ acting on sections of $\Lambda^{k}_{\C}H^{*}$ 
   when the condition~$X(k)$ holds everywhere.

   (d) The  horizontal sublaplacian $\Delta_{b;p,q}$ on a CR manifold acting on $(p,q)$-forms under the condition $X(p,q)$. 
%    when the condition~$X(k)$ holds everywhere.
\end{example}

In the next section we will make use of Theorem~\ref{thm:Powers1.main} to show that when the bracket condition $H+[H,H]=TM$ 
holds,   the Rockland condition insures us that the principal symbol of $P+\partial_{t}$ admits an inverse in  
$S_{\vo,-v}(\fg^{*}M\times \R_{(v)},\cE)$ (see Theorem~\ref{thm:Heat1.main}). Therefore, we obtain: 

\begin{theorem}\label{thm:Powers1.main2}
Assume that the bracket condition $H+[H,H]=TM$ holds and that $P$ satisfies the Rockland condition at every point of $M$. Then the complex powers 
$P^{s}$, $s \in \C$, of $P$ form a  holomorphic 1-parameter group of \psivdos\ such that $\ord P^{s}=vs$ $\forall s \in \C$. 
\end{theorem}

\begin{example}
Assume that $(M,H)$ is an orientable contact manifold. Then: 

- In degree $k\neq n$ the complex powers of the contact Laplacian $\Delta_{R,k}$ form 
a holomorphic 1-parameter group of \psivdos\ such that $\ord \Delta_{R,k}^{s}=2s$ $\forall s \in \C$.\smallskip

- In degree $n$ the complex powers of the contact Laplacians $\Delta_{R,nj}$, $j=1,2$, form 
holomorphic 1-parameter groups of \psivdos\ such that $\ord \Delta_{R,nj}^{s}=4s$ $\forall s \in \C$.
% 
% Theorem~\ref{thm:Powers1.main2} is valid for the contact Laplacian $\Delta_{R}$ on a compact orientable 
% contact manifold $(M^{2n+1},\theta)$. 
%    In this case $\Delta_{R}^{s}$ has order $2s$ 
%    on degree $k=0,\ldots, 2n$ with $k\neq n$ and has order $4s$ on degree $n$.
\end{example}

\begin{remark}
    The above example allows us to fill a technical gap in~\cite{JK:OKTGSU} in the proof that the complex powers of the contact Laplacian are 
    \psivdos\  (see~\cite{Po:Crelle1}). The latter result is one ingredient in~\cite{JK:OKTGSU}, among others, in a proof of the Baum-Connes conjecture 
    for $SU(n,1)$. 
\end{remark}

\begin{remark}
   Given a (stratified) graded nilpotent group $G$ it was shown by Folland that some complex powers of left invariant sum of squares are left 
   invariant \psidos. This was made via the use of the fundamental solution of the associated heat operators. These results have later been extended to more general 
   operators in~\cite{CGGP:POGD}. Using a different approach Geller~\cite{Ge:CPCOHG} dealt with general positive left invariant \psidos\ on the Heisenberg 
   group $\bH^{2n+1}$ provided that the Rockland condition is satisfied. In fact, on each tangent group $G_{a}M$ 
   these results can be recovered from Theorem~\ref{thm:Powers1.main2} by looking at the principal symbols of the complex powers.  
\end{remark}

\section{Rockland condition and the heat equation}
 \label{sec:Rockland-heat}
In this section we make use of Theorem~\ref{thm:Powers1.main} in connection with the Rockland condition. 

First, we extend Theorem~\ref{thm:Chap3.Rockland-Parametrix-order+} 
and Proposition~\ref{prop:Heisenberg.Rockland-hypoellipticity} to \psivdos\ with non-integer order as follows. 

\begin{theorem}\label{thm:Heat-Rockland.invertibility-non-integer-order}
 Assume that the bracket condition $H+[H,H]=TM$  holds and let $P:C^{\infty}_{c}(M,\cE)\rightarrow C^{\infty}(M,\cE)$ be a \psivdo\ of order $m \in 
 \C$. Then the following are equivalent:\smallskip
  
  (i) The principal symbol of $P$ is invertible;\smallskip
  
  (ii) $P$ and $P^{t}$ satisfy the Rockland condition at every point $a\in M$.\smallskip

  (iii) $P$ and $P^{*}$ satisfy the Rockland condition at every point $a\in M$.
\noindent Furthermore, when (i) or (ii) holds the operator $P$ admits a parametrix in $\pvdo^{-m}(M,\cE)$.
\end{theorem}
\begin{proof}
  Consider a sum of squares  $\Delta_{\nabla,X}=\nabla_{X_{1}}^{*}\nabla_{X_{1}}+\ldots+\nabla_{X_{m}}^{*}\nabla_{X_{m}}$,
  where $\nabla$ is a connection on $\cE$ and the vector fields $X_{1},\ldots,X_{m}$ span $H$. Since by assumption the bracket condition $H+[H,H]=TM$ 
  holds, we know from Proposition~\ref{thm:Powers1.heat-sublaplacians} that 
  the principal symbol of $\Delta_{\nabla,X}+\partial_{t}$ admits an inverse in $S_{\vo,-2}(\fg^{*}M\times \R_{(2)},\cE)$, 
  so by Theorem~\ref{thm:Powers1.main} the 
  operator $\Delta_{\nabla,X}^{-\frac{m}{2}}$ is a \psivdo\ of order $-m$ with an invertible principal symbol. Therefore, along the same lines as that 
  of the proof of Theorem~\ref{thm:Chap3.Rockland-Parametrix-order+} we can show that the conditions~(i), (ii) and (iii) are equivalent. 
\end{proof}

\begin{proposition}\label{prop:Heat-Rockland.Rockland-hypoellipticity}
 Assume that the bracket condition $H+[H,H]=TM$  holds and let $P:C^{\infty}_{c}(M,\cE)\rightarrow C^{\infty}(M,\cE)$ be a \psivdo\ of order $m \in 
 \C$ with $\Re m\geq 0$ and such that $P$ satisfies the Rockland condition at every point. Then $P$ is hypoelliptic with gain of 
 $\frac{1}{2}\Re m$ derivatives. 
\end{proposition}
\begin{proof}
    Since $P$ satisfies the Rockland condition at every point, we can argue as in the proof of Proposition~\ref{prop:Heisenberg.Rockland-hypoellipticity}, 
    using Theorem~\ref{thm:Heat-Rockland.invertibility-non-integer-order} instead of 
    Theorem~\ref{thm:Chap3.Rockland-Parametrix-order+}, to show that the principal symbol of $P$ is left-invertible. Therefore,  
    $P$ admits a left-parametrix in $\pvdo^{-m}(M,\cE)$ and it then follows from Remark~\ref{rem:hypoellipticity-parametrix} that $P$ is is hypoelliptic with gain of 
    $\frac{1}{2}\Re m$ derivatives. 
\end{proof}

Let us now show that the Rockland condition for a (positive) differential operator  $P$ implies the invertibility of the 
principal symbol of $P+\partial_{t}$. This will show that the frameworks of~\cite{BGS:HECRM} and of Theorem~\ref{thm:Powers1.main} 
apply to a large class of operators. 

% 
% consider a selfadjoint differential operato
% $P:C^{\infty}(M,\cE)\rightarrow C^{\infty}(M,\cE)$ of even Heisenberg order $v$. We shall now show that the Rockland condition for 
%  
%  In this section we shall show that the Rockland condition for $P$ is enough to insure us the invertibility of the principal symbol of $P+\partial_{t}$. 
%  This is needed for the pseudodifferential representation of the heat kernel and the invertibility of the heat operator $P+\partial_{t}$ 
%  in~\cite{BGS:HECRM} and in Section~\ref{sec.Volterra-PsiHDO-calculus}, as well as in Theorem~\ref{thm:Powers1.main} for the complex powers of 
%  $P$ to give rise to a holomorphic family of \psivdos. 

To reach aim it will be convenient to enlarge the class of Volterra \psivdos\ as follows. 
% by considering the class of \psivdos\ on $M\times \R$ defined as the Volterra 

% \psivdos\ with the Volterra property replaced by the translation invariance with respect to the time variable. 
% Given a Heisenberg chart $U \subset \Rd$ with $H$-frame $X_{0},\ldots,X_{d}$ we 
% We first extend Theorem~\ref{thm:Chap3.Rockland-Parametrix-order0} to the Volterra setting. To this end it would be convenient to introduce 
%  variant of the Volterra \psivdos. Let $U\subset \Rd$ be a Heisenberg chart with $H$-frame $X_{0},\ldots,X_{d}$. 
 
 \begin{definition}
     Let $U$ be an open subset of $\Rd$. Then:\smallskip 
     
    1) $S_{m}(\URdd_{(v)})$, $m \in \Z$, consists of functions $q(x,\xi,\tau)$ in $C^{\infty}(U\times (\Rdd\setminus 0))$ 
    such that  $q(x,\lambda.\xi,\lambda^{v}\tau)=\lambda^{m}q(x,\xi,\tau)$ for any $\lambda >0$.\smallskip
    
    2) $S^{m}(\URdd_{(v)})$, $m \in \Z$, consists of functions $q(x,\xi,\tau)$ in $C^{\infty}(U\times \Rdd)$ admitting  an 
    asymptotic expansion  $q \sim \sum_{j\geq 0} q_{m-j}$, $q_{m-j}\in S_{m-j}(\URd \times \R_{(v)})$, with $\sim$ taken in the sense 
    of~(\ref{eq:volterra.asymptotic-symbols}).
\end{definition}
 
\begin{definition}
    $\pvdo^{m}(M\times \R_{(v)},\cE)$, $m \in \Z$, consists of continuous operators $Q$ from $C^{\infty}_{c}(M\times \R, \cE)$ to $C^{\infty}_{c}(M\times 
    \R, \cE)$ such that:\smallskip
    
    (i) $Q$ is translation invariant with respect to the time variable;\smallskip
    
    (ii) The kernel of $Q$ is smooth outside the diagonal;\smallskip
    
    (iii) In any local trivializing chart $U \subset \Rd$ with $H$-frame $X_{0},\ldots,X_{d}$ we can write $Q$ in the form,
    \begin{equation}
        Q=q(x,-iX,D_{t})+R,
%         \label{eq:¥}
    \end{equation}
    with $q$ in $S^{m}(\URdd_{(v)}, \C^{r})$ and $R$ smoothing operator. 
\end{definition}

The main properties of \psivdos\ and of Volterra \psivdos\ hold \emph{mutatis mutandis} for these operators. In particular, it makes sense to define 
then on $M \times \R$ and the substitute to the class $\cK_{m}(\URd)$ and $\cK_{\vo,m}(\URdd_{(v)})$ is given by the distributions below. 

\begin{definition}
    $\cK_{m}(\URdd_{(v)})$, $m \in \Z$, consists of distributions $K(x,y,t)$ in $C^{\infty}\hotimes \cD_{\reg}(\Rdd\setminus 0)$ 
    for which there are $c_{\alpha,l}\in C^{\infty}(U)$,  $\brak\alpha+vl=m$, such that
    \begin{equation}
        K(x,\lambda.y,\lambda^{v}t)= \lambda^m K(x,y,t) + \lambda^m\log\lambda
                \sum_{\brak\alpha+vl=m}c_{\alpha}(x)y^\alpha t^{l}\qquad \forall \lambda>0.
    \end{equation}
\end{definition}

Moreover, in this context the principal symbol of a \psivdo\ in $\pvdo^{m}(M\times R_{(v)},\cE)$ makes sense as an element of the following symbol  
class. 

\begin{definition}
    $S_{m}(\fg^{*}M\times \R_{(v)},\cE)$, $m \in \Z$, consists of sections $q(x,\xi,\tau)$ in $C^{\infty}((\fg^{*}M\times \R)\setminus 0, \End 
    \pi^{*}\cE)$ such that $q(x,\lambda.\xi,\lambda^{v}\tau)=\lambda^{m}q(x,\xi,\tau)$.  
\end{definition}

This allows us to define the model operator and the Rockland condition at a point $a\in M$ in the same way as for \psivdos\ on $M$.

Moreover, the class $\pvdo^{m}(M\times \R_{(v)},\cE)$ is closed under taking adjoints, so if $Q$ is a Volterra \psivdo\ in 
$\pvhdo^{m}(M\times \R_{(v)},\cE)$ then its adjoint is not a Volterra \psivdos\ but at least it belongs to $\pvdo^{m}(M\times \R_{(v)},\cE)$. 

We can now establish the analogue of Theorem~\ref{thm:Chap3.Rockland-Parametrix-order0} in the Volterra setting. 
% 
% in terms of which 
% we can characterize the distribution kernels of \psivdos\ in the classes  $\pvdo^{m}(M\times R_{(v)},\cE)$, $m \in \Z$. 
% All the properties of \psivdos\ carry over to operators in the classes $\pvdo^{m}\times \R_{(v)})$, $m \in \Z$. In particular, the distribution kernels 
% can be characterized in terms of the following   

\begin{proposition}\label{prop:Chap4.Rockland-Parametrix-order0}
    Let $Q  \in \pvdo^{0}(M\times \R_{(v)},\cE)$. The following are equivalent:\smallskip
    
    (i) For any $a \in M$ the model operator $Q^{a}$ is invertible on $L^{2}(G_{a}M\times \R_{(v)},\cE_{a})$.\smallskip
    
    (ii) $Q$ and $Q^{*}$ satisfy the Rockland condition at every point of $M$.\smallskip
    
    (iii) The principal symbol of $Q$ is invertible in $S_{0}(\fg^{*}M\times \R_{(v)},\cE)$.\smallskip
    
\noindent Moreover, if the principal symbol of $Q$ is invertible and belongs to $S_{0}(\fg^{*}M\times \R_{(v)},\cE)$ then its inverse is in 
$S_{0}(\fg^{*}M\times \R_{(v)},\cE)$ too.
 \end{proposition}
 \begin{proof}
     The proof of the equivalence between the conditions (i), (ii) and (iii) follows along the same lines as that of Theorem~\ref{thm:Chap3.Rockland-Parametrix-order0}, 
     since the results of \cite{CGGP:POGD}, \cite{Ch:ORISIO}, \cite{Ch:ISASIO}, \cite{FS:HSHG} and \cite{KS:IOSSG} 
     used in the proof of Theorem~\ref{thm:Chap3.Rockland-Parametrix-order0} hold in this setting. 
     
     It remains to show that if the principal symbol of $Q$ is invertible and belongs to 
     $S_{0}(\fg^{*}M\times \R_{(v)},\cE)$ then its inverse is in $S_{0}(\fg^{*}M\times 
    \R_{(v)},\cE)$. In view of the proof of Lemma~\ref{lem:Chap3.Rockland.inverse} in order to reach this aim it is enough to show that, for any Heisenberg chart 
    $U \subset \Rd$, the Volterra class $\cK_{\vo,-(d+2+v)}(\URdd_{(v)})$ is a closed subalgebra of $\cK_{-(d+2+v)}(\URdd_{(v)})$.
    
    It is clear that if $K_{1}$ and $K_{2}$ are in $\cK_{\vo,-(d+2+v)}(\URdd_{(v)})$ then $K_{1}*K_{2}$ is also in $\cK_{\vo,-(d+2+v)}(\URdd_{(v)})$. 
    %so the space $\cK_{\vo,-(d+2+v)}(\URd \times \R_{(v)})$ is a subalgebra of  $\cK_{-(d+2+v)}(\URd \times \R_{(v)})$. 
    Moreover, if  a sequence $(K_{j})_{j \geq 0}\subset \cK_{\vo,-(d+2+v)}(\URdd_{(v)})$  converges to $K$ in $\cK_{-(d+2+v)}(\URdd_{(v)})$ 
    then $K_{j}$ converges to $K$ in $\cD'(\URdd)$ and so we have $\supp K \subset \cup \supp K_{j}\subset \URd\times [0, \infty)$, that 
    is, $K$ belongs to $\cK_{\vo,-(d+2+v)}(\URdd_{(v)})$. Thus $\cK_{\vo,-(d+2+v)}(\URdd_{(v)})$ is a closed subalgebra of 
    $\cK_{-(d+2+v)}(\URdd_{(v)})$ and henceforth the proof is achieved. 
 \end{proof}
 
In the sequel we will need the notion of positive symbols below. 
\begin{definition}\label{def:Heat1.positive-symbol}
 A symbol $p\in S_{m}(\fg^{*}M,\cE)$ is positive when it can be put into the form $p=\bar{q}*q$ for some symbol $q\in S_{\frac{m}{2}}(\fg^{*}M,\cE)$. 
\end{definition}

The interest of this definition stems from:
\begin{lemma}\label{lem:Rockland-heat.positivity}
  1) The principal symbol of $P$ is positive if, and only if, there exist $Q\in \pvdo^{\frac{v}{2}}(M,\cE)$ and $R\in \pvdo^{v-1}(M,\cE)$ so that 
  $P=Q^{*}Q+R$.\smallskip
   
  2) If $P$ satisfies the Rockland condition at every point and has a positive principal symbol, then the operators $P\pm \partial_{t}$ satisfy the Rockland condition at 
  every point.
% 
%   Let $P\in \pvdo^{m}(M,\cE)$ have principal symbol $\sigma_{m}(P)\in S^{m}(\fg^{*}M,\cE)$. Then we have equivalence:\smallskip
%    
%    (i) The symbol $\sigma_{m}(P)$ is positive;\smallskip
% 
%    (ii) There . 
\end{lemma}
\begin{proof}
Assume that the principal symbol $p=\sigma_{v}(P)$ of $P$ is positive, so that there exists $q_{\frac{v}{2}}\in S_{\frac{v}{2}}(\fg^{*}M,\cE)$ such that  
  $p=\overline{q_{\frac{v}{2}}}*q_{\frac{v}{2}}$. By Proposition~\ref{prop:PsiHDO.surjectivity-principal-symbol-map}   the principal 
 symbol map $\sigma_{\frac{v}{2}}:\pvdo^{\frac{v}{2}}(M,\cE)\rightarrow S^{\frac{v}{2}}(\fg^{*}M,\cE)$ is surjective, so there  
 exists $Q\in \pvdo^{\frac{v}{2}}(M,\cE)$ such that $\sigma_{\frac{v}{2}}(Q)=q_{\frac{v}{2}}$. Then by Proposition~\ref{prop:PsiHDO.composition2} 
 and~Proposition~\ref{prop:PsiHDO.adjoint-manifold}  the operator $Q^{*}Q$ has principal symbol 
 $\overline{q_{\frac{v}{2}}}*q_{\frac{v}{2}}=\sigma_{v}(P)$, hence coincides with $P$ modulo $\pvdo^{v-1}(M,\cE)$. 
  
 Conversely, if $P$ is of the form $P=Q^{*}Q+R$ with $Q$ in $\pvdo^{\frac{v}{2}}(M,\cE)$ and $R$ in $\pvdo^{v-1}(M,\cE)$ then $P$ and $Q^{*}Q$ have same 
 principal symbol. Therefore, we have $\sigma_{v}(P)=\overline{\sigma_{\frac{v}{2}}(Q)}*\sigma_{\frac{v}{2}}(Q)$, that is, $\sigma_{v}(P)$ is a positive 
 symbol.\smallskip 
 
 Finally, suppose that $P$ satisfies the Rockland condition at every point and has a positive principal symbol, i.e.,
   there exists $\tilde{p}_{\frac{v}{2}}\in S_{\frac{v}{2}}(U\times G)$ so that 
   $\sigma_{v}(P)=\overline{\tilde{p}_{\frac{v}{2}}}*\tilde{p}_{\frac{v}{2}}$. Therefore,   if we let $\tilde{P}^{a}$ be the left-convolution operator with 
   symbol $\tilde{p}_{\frac{v}{2}}(a,.)$ then by Proposition~\ref{prop:PsiHDO.composition2} and~Proposition~\ref{prop:PsiHDO.adjoint-manifold} 
   we have $P^{a}=(\tilde{P}^{a})^{*}\tilde{P}^{a}$. Thus, by 
   Proposition~\ref{PsiHDO.properties-symbol-representation} for 
   every non-trivial irreducible representation $\pi$ of $G$ we have 
   $\overline{\pi_{P^{a}}}=(\overline{\pi_{\tilde{P}^{a}}})^{*}\overline{\pi_{\tilde{P}^{a}}}$, 
   which shows that $\overline{\pi_{P}^{a}}$ is positive. We then can argue as in the proof of~\cite[Lem.~4.21]{FS:HSHG} to show that the operators 
   $P^{a}\pm \partial_{t}$ satisfy the Rockland condition on $G_{a}M\times \R_{(v)}$. Hence $P+\partial_{t}$ and $P-\partial_{t}$ 
   satisfy the Rockland condition at every point of $M$. 
\end{proof}

We are now ready to prove the main result of this section. 

\begin{theorem}\label{thm:Heat1.main}
 Assume that the  bracket condition $H+[H,H]=TM$ holds and that $P$ satisfies the Rockland condition at every point.\smallskip
 
 1) $P$ is bounded from below if, and only if, it has a positive principal symbol.\smallskip
 
 2) If $P$ has a positive principal symbol, then the principal symbol $P+\partial_{t}$ has an inverse in $S_{\vo,-v}(\fg^{*}M\times \R_{(v)},\cE)$.  
 Hence Theorems~\ref{thm:volterra.inverse} and~\ref{thm:Powers1.heat-kernel-asymptotics} hold for~$P$.  
\end{theorem}
\begin{proof}
      Let us first assume that $P$ has a positive principal symbol. Since $P$ satisfies the Rockland condition at every point,  
      Lemma~\ref{lem:Rockland-heat.positivity}  tells us 
      that the heat operator $P +\partial_{t}$ and its adjoint satisfy the Rockland condition at every point. 
      
      On the other hand, let $k=\frac{v}{2}$ and let $\Delta_{\nabla,X}:C^{\infty}(M,\cE)\rightarrow C^{\infty}(M,\cE)$ be a sum of squares of the form,
      \begin{equation}
          \Delta_{\nabla,X}=\nabla_{X_{1}}^{*}\nabla_{X_{1}}+\ldots+\nabla_{X_{m}}^{*}\nabla_{X_{m}},
%           \nonumb
      \end{equation}
      where $\nabla$ is a connection on $\cE$ and the vector fields $X_{1},\ldots,X_{m}$. Since by assumption the bracket condition $H+[H,H]=TM$ holds 
      it follows from Proposition~\ref{prop:Heat1.Delta^k} that the principal symbol of $\Delta_{\nabla,X}^{k}+\partial_{t}$ admits an inverse in 
      $S_{\vo,-v}(\fg^{*}M\times \R_{(v)},\cE)$ and that $(\Delta^{k} +\partial_{t})^{-1}$ belongs to $\pvhdo^{-v}(M\times R_{(v)},\cE)$. 
      In particular, since $(\Delta^{\frac{v}{2}} +\partial_{t})^{-1}$ has an invertible principal symbol, together with its adjoint it satisfies the 
      Rockland condition at every point. 
      
      Let $Q_{1}=(\Delta^{k}+\partial_{t})^{-1}(P+\partial_{t})$ and $Q_{2}=(P+\partial_{t})(\Delta^{k}+\partial_{t})^{-1}$. They 
      are elements of $\pvhdo^{0}(M\times R_{(v)},\cE)$ which  together with their adjoints satisfy the Rockland condition at every point, since 
      they are products of such operators.  
     We  may therefore apply Proposition~\ref{prop:Chap4.Rockland-Parametrix-order0} 
     to deduce that the principal symbols of $Q_{1}$ and $Q_{2}$ are invertible in $S_{0, \op{v}}(\fg^{*}M\times 
     \R_{(v)},\cE)$. In the same way as in the proof of Theorem~\ref{prop:Chap4.Rockland-Parametrix-order0} 
     this implies that the principal symbol of $P+\partial_{t}$ is an invertible Volterra-Heisenberg symbol.   
    
    It remains to show that $P$ is bounded from below if, and only if, its principal symbol is positive.  
    
    Suppose that $P$ is bounded from below. 
    Possibly replacing $P$ by $P+c$ with $c>0$ large enough we may assume that 
    $P$ is positive. Notice that $P^{2}$ is selfadjoint and satisfies the Rockland condition at every point, since so does $P$. Therefore, we may 
    apply the first part of the proof to deduce that the principal symbol of $P^{2}+\partial_{t}$ is an invertible 
    Volterra-Heisenberg symbol. As $P^{2}$ is positive we then may use Theorem~\ref{thm:Powers1.main} 
    to see that $(P^{2})^{\frac{1}{4}}=P^{\frac{1}{2}}$ is a \psivdo\ of order $k=\frac{v}{2}$. 
    Since $P=(P^{\frac{1}{2}})^{2}$ it then follows from Lemma~\ref{lem:Rockland-heat.positivity} that the principal symbol of $P$ is a positive symbol. 

    Conversely, suppose that $P$ has a positive principal symbol. Then thanks to Lemma~\ref{lem:Rockland-heat.positivity} the operator $P$ can be written as
    $P=Q^{*}Q+R$ with $Q\in \pvdo^{k}(M,\cE)$ and $R\in \pvdo^{v-1}(M,\cE)$.  Let $P_{1}=Q^{*}Q$. 
    Since $P$ and $P_{1}$ have same principal symbol, the operator $P_{1}$ also satisfies the Rockland condition at every point. 
    
      Next, since $P_{1}$ is positive for $\lambda<-1$ we have
    \begin{equation}
        (P_{1}-\lambda)^{-1}=\int_{0}^{\infty}e^{-tP_{1}}e^{t\lambda}dt,
    \end{equation}
    where the integral converges in $\cL(L^{2}(M,\cE))$. Let $\alpha=\frac{v-1}{v}$. Then we have
    \begin{equation}
        R(P_{1}-\lambda)^{-1}=RP^{-\alpha}_{1}R_{(\lambda)}, \qquad R_{(\lambda)}=\int_{0}^{\infty}P^{\alpha}_{1}e^{-tP_{1}}e^{t\lambda}dt.
         \label{eq:Heat1.decomposition-Rlambda}
    \end{equation}

    Since $P_{1}=Q^{*}Q$ has a positive principal symbol and satisfies the Rockland condition the first part of the proof tells us that 
    the principal symbol of $P_{1}+\partial_{t}$ is an invertible Volterra-Heisenberg symbol. This allows us to apply  
    Theorem~\ref{thm:Powers1.main2} to deduce that $P^{-\alpha}_{1}$ is a \psivdo\ of order $-v\alpha=-(v-1)$. 
    Therefore, the operator $RP^{-\alpha}_{1}$ is a \psivdo\ of order~$0$, hence is bounded on $L^{2}(M,\cE)$.  
  
    On the other hand, we have 
  \begin{equation}
      \|R_{(\lambda)}\|_{L^{2}}\leq \int_{0}^{\infty}t^{-\alpha} \|(tP_{1})^{\alpha}e^{-tP_{1}}\|e^{t\lambda}dt.
       \label{eq:Heat1.Rlambda-estimates1}
  \end{equation}
  As $\alpha \in(0,1)$  the function $x\rightarrow x^{\alpha}e^{-x}$ maps $[0,\infty)$ to $[0,1]$. Therefore, we have
  $\|(tP_{1})^{\alpha}e^{-tP_{1}}\|\leq 1$, from which we get
  \begin{equation}
       \|R_{(\lambda)}\|_{L^{2}}\leq \int_{0}^{\infty}t^{-\alpha} e^{t\lambda}dt=|\lambda|^{\alpha-1} \int_{0}^{\infty}u^{-\alpha} e^{-u}du.
       \label{eq:Heat1.Rlambda-estimates2}
  \end{equation}
  
 Now, by combining~(\ref{eq:Heat1.decomposition-Rlambda}), (\ref{eq:Heat1.Rlambda-estimates1}) and (\ref{eq:Heat1.Rlambda-estimates2}) together we obtain
  \begin{equation}
      \|R(P_{1}-\lambda)^{-1}\|_{L^{2}}\leq  \|RP^{-\alpha}_{1}\|_{L^{2}}\|R_{(\lambda)}\|_{L^{2}}\leq C_{\alpha}|\lambda|^{\alpha-1},
  \end{equation}
  where the constant $C_{\alpha}$ does not depend on $\lambda$. Since $\alpha<1$ it follows that for $\lambda$ negatively large enough we have 
  $\|R(P_{1}-\lambda)^{-1}\|_{L^{2}}\leq \frac{1}{2}$, so that $1+R(P_{1}-\lambda)^{-1}$ is invertible. 
 Since we have 
 \begin{equation}
     P-\lambda=P_{1}+R-\lambda=(P_{1}+R-\lambda)(P_{1}-\lambda)^{-1},
 \end{equation}
it follows that $P-\lambda$ has a right inverse for $\lambda$ negatively large enough. Since we can similarly show that for $\lambda$ negatively large enough 
$P-\lambda$ is left-invertible, we deduce that 
as soon as $\lambda$ is negatively large enough $P-\lambda$ admits a bounded two-sided inverse. This means
 that the spectrum of $P$ is contained in some interval $[c,\infty)$, that is, $P$ is bounded from below.  The proof is now complete.   
\end{proof}

\begin{example}
Assume that $(M,H)$ is an orientable contact manifold. Then the 2nd part of Theorem~\ref{thm:Heat1.main}  
applies to the contact Laplacian in every degree. Thus, in 
degree $k\neq n$ the principal symbol of $\Delta_{R,k}+\partial_{t}$ admits an inverse in $S_{\vo,-2}(\fg^{*}M\times \R_{(2)},\Lambda^{k})$ and in 
degree $n$ the principal symbol of $\Delta_{R,nj}+\partial_{t}$, $j=1,2$, admits an inverse in $S_{\vo,-4}(\fg^{*}M\times \R_{(4)},\Lambda^{n}_{j})$.
\end{example}

Finally, let $M^{2n+1}$ be a compact orientable strictly pseudoconvex CR manifold endowed with a pseudohermitian contact form $\theta$ and the 
associated Levi metric. For $k=1,\ldots,n+1$ and for $k=n+2,n+4,\ldots$ let $\boxdot_{\theta}^{(k)}:C^{\infty}(M)\rightarrow C^{\infty}(M)$ be the  
Gover-Graham operator of order $k$. We know that $\boxdot_{\theta}^{(k)}$ is selfadjoint and satisfies the Rockland condition at every point. 
Furthemore, we have: 

\begin{proposition}
   For $k\neq n+1$ the operator $\boxdot_{\theta}^{(k)}$ is bounded from below and the principal symbol of $\boxdot_{\theta}^{(k)}+\partial_{t}$ is 
   invertible in $S_{\vo,-2k}(\fg^{*}M\times \R_{(2k)})$. Hence Theorems~\ref{thm:volterra.inverse} and~\ref{thm:Powers1.heat-kernel-asymptotics} hold 
   for~$\boxdot_{\theta}^{(k)}$ with $v=2k$. 
\end{proposition}
\begin{proof}
  Assume $k\neq n+1$.  Since $\boxdot_{\theta}^{(k)}$ is selfadjoint and by Proposition~\ref{prop:Examples.Gover-Graham} 
  it satisfies the Rockland condition at every point, 
  thanks to Theorem~\ref{thm:Heat1.main} we are reduced to show that its principal symbol is positive. Let $\Delta_{b;0}$ be the horizontal sublaplacian on 
  functions and let $X_{0}$ be the Reeb vector field of $\theta$. Then by Proposition~\ref{prop:Examples.Gover-Graham} 
  the principal symbol of $\boxdot_{\theta}^{(k)}$ agrees 
  with that of 
  \begin{equation}
      \boxdot=    (\Delta_{b;0}+i(k-1)X_{0})(\Delta_{b;0}+i(k-3)X_{0})\cdots (\Delta_{b;0}-i(k-1)X_{0}).
      \label{eq:Heat.Gover-Graham}
  \end{equation}

  Since $k \neq n+1$ it follows from the proof of Proposition~\ref{prop:Examples.Gover-Graham} 
  that each factor $\Delta_{b;0}+i\alpha X_{0}$ in~(\ref{eq:Heat.Gover-Graham}) is a selfadjoint 
  sublaplacian satisfying the condition~(\ref{eq:Sublaplacian.condition}) at every point. It can be shown that any such sublaplacian is bounded from below 
  (see~\cite{Po:PhD}, \cite{Po:CPDE1}). Since the product of selfadjoint operators which are bounded from below is bounded from below provided that 
  the product is itself selfadjoint, we see that the operator $\frac{1}{2}(\boxdot +\boxdot^{*})$ is bounded from below. 
  
  Notice that $\frac{1}{2}(\boxdot +\boxdot^{*})$ has same principal as $\boxdot$ and $\boxdot_{\theta}^{(k)}$, hence satisfies the Rockland condition 
  at every point. Therefore, Theorem~\ref{thm:Heat1.main} tells us that the principal symbol of $\frac{1}{2}(\boxdot +\boxdot^{*})$ is positive. 
  Incidentally, the operator $\boxdot_{\theta}^{(k)}$ has a positive principal. The proof is thus complete. 
\end{proof}

\section{Weighted Sobolev Spaces}
\label{sec.Sobolev}
In this section   we construct weighted Sobolev spaces $W_{H}^{s}(M)$, $s \in \R$, which provide us with sharp regularity estimates for \psivdos. To 
this end  we assume throughout the section that the bracket condition $H+[H,H]=TM$ holds. 
% 
% extend to any real 
% parameter the weighted Sobolev spaces $S_{k}^{2}(M)$, $k\in \N$, of Folland-Stein~\cite{FS:EDdbarbCAHG} (see also~\cite{RS:HDONG}). 
% As a consequence these Sobolev spaces will provide us with sharp regularity estimates for hypoelliptic \psivdos. 

Let $X_{1},\ldots,X_{m}$ be real vector fields spanning $H$ and consider the 
positive sum of squares, 
\begin{equation}
    \Delta_{X}= X_{1}^{*}X_{1}+\ldots+X_{m}^{*}X_{m}.
\end{equation}
As mentioned in Example~\ref{ex:Powers1.sublaplacian} (a), 
since the the bracket condition $H+[H,H]=TM$ holds Theorem~\ref{thm:Powers1.main} is valid for $1+\Delta_{X}$. Thus, 
the complex powers $(1+\Delta_{X})^{s}$, $s \in \C$,  
gives rise to a holomorphic 1-parameter group  of \psivdos\ such that we have $\ord (1+\Delta_{X})^{s}=2s$ for any $s\in \C$. 

\begin{definition}
  $W_{H}^{s}(M)$, $s\in \R$, consists of all distributions $u\in \cD'(M)$ such that 
  $(1+\Delta_{X})^{\frac{s}{2}}u$ is in $L^{2}(M)$. It is endowed with the Hilbert norm given by
  \begin{equation}
      \|u\|_{W_{H}^{s}}=\|(1+\Delta_{X})^{\frac{s}{2}}u\|_{L^{2}}, \qquad u \in W_{H}^{s}(M). 
  \end{equation}
 \end{definition}

\begin{proposition}\label{prop:Sobolev.embeddings}
1) Neither $W_{H}^{s}(M)$, nor its topology, depend on the choice of the vector fields $X_{1},\ldots,X_{m}$.\smallskip 
    
    2) We have the following continuous embeddings: 
    \begin{equation}
   \begin{array}{rccccl}
       L^{2}_{s}(M)  & \hookrightarrow & W_{H}^{s}(M)&\hookrightarrow & L^{2}_{s/2}(M) & \qquad \text{if $s\geq 0$},\\
       L^{2}_{s/2}(M)  & \hookrightarrow & W_{H}^{s}(M) & \hookrightarrow & L^{2}_{s}(M) & \qquad \text{if $s< 0$}.
   \end{array}
         \label{eq:Sobolev.embeddings}
    \end{equation}
\end{proposition}
\begin{proof}
    1) Let $Y_{1},\ldots,Y_{p}$ be other vector fields spanning $H$. The operator 
    $(1+\Delta_{Y})^{s}(1+\Delta_{X})^{-s}$ is a \psivdo\ 
    of order $0$, hence is bounded on $L^{2}(M)$ by Proposition~\ref{prop:PsiHDO.Sobolev-regularity}. Therefore,  we get the
    estimates,
    \begin{multline}
        \|(1+\Delta_{Y})^{s}u\|_{L^{2}}=\|(1+\Delta_{Y})^{s}(1+\Delta_{X})^{-s}(1+\Delta_{X})^{s}u\|_{L^{2}}\\ \leq C_{XYs} 
        \|(1+\Delta_{X})^{s}u\|_{L^{2}},
    \end{multline}
     which hold for any $u \in C^{\infty}(M)$. Interchanging the roles of the $X_{j}$'s and of the $Y_{k}$'s also gives the estimates
     \begin{equation}
         \|(1+\Delta_{X})^{s}u\|_{L^{2}} \leq C_{YXs}\|(1+\Delta_{Y})^{s}u\|_{L^{2}}, \qquad u \in C^{\infty}(M). 
     \end{equation}
     Therefore, whether we use the $X_{j}$'s or the $Y_{k}$'s to define the $W_{H}^{s}(M)$ changes neither the space, nor its topology. \smallskip 
     
     2) Let $ s\in [0,\infty)$. Since $(1+\Delta_{X})^{{\frac{s}{2}}}$ is a \psivdo\ of order $s$, 
     Proposition~\ref{prop:PsiHDO.Sobolev-regularity}   
     tells us that 
     it is bounded from $L^{2}_{s}(M)$ to $L^{2}(M)$. Thus,
     \begin{equation}
         \|u\|_{W_{H}^{s}}=\|(1+\Delta_{X})^{{\frac{s}{2}}}u\|_{L^{2}}\leq C_{s} \|u\|_{L^{2}_{s}}, \qquad u \in W_{H}^{s}(M), 
     \end{equation}
     which shows that  $L^{2}_{s}(M)$  embeds continuously into $W_{H}^{s}(M)$.
     
    On the other hand, as $(1+\Delta_{X})^{-{\frac{s}{2}}}$ has order $-s$ Proposition~\ref{prop:PsiHDO.Sobolev-regularity}   
     also tells us that $(1+\Delta_{X})^{-{\frac{s}{2}}}$ is bounded from $L^{2}(M)$ to 
     $L^{2}_{s/2}(M)$. Therefore, on $W_{H}^{s}(M)$ we have the estimates
     \begin{multline}
         \|u\|_{L^{2}_{s/2}}=\|(1+\Delta_{X})^{-{\frac{s}{2}}}(1+\Delta_{X})^{{\frac{s}{2}}}u\|_{L^{2}_{s/2}}\\ \leq 
         C_{s}\|(1+\Delta_{X})^{-{\frac{s}{2}}}\|_{L^{2}}=C_{s}\|u\|_{W_{H}^{s}},
     \end{multline}
     which shows that $W_{H}^{s}(M)$ embeds continuously into $L^{2}_{s/2}(M)$.  
     
     Finally, when  $s<0$ we can similarly show that we have continuous embeddings 
     $L^{2}_{s/2}(M)  \hookrightarrow W_{H}^{s}(M)$ and $ W_{H}^{s}(M)\hookrightarrow L^{2}_{s}(M)$. 
\end{proof}

As an immediate consequence of Proposition~\ref{prop:Sobolev.embeddings} we obtain: 

\begin{proposition}\label{prop:Sobolev.Cinfty-cD'}
   The following equalities between topological spaces hold:
   \begin{equation}
       C^{\infty}(M)=\cap_{s \in \R} W_{H}^{s}(M) \qquad \text{and} \qquad \cD'(M)= \cup_{s \in \R}  W_{H}^{s}(M).
   \end{equation}
\end{proposition}

Let us now compare the weighted Sobolev spaces $W_{H}^{s}(M)$ to the weighted Sobolev spaces $S_{k}^{2}(M)$, $k=1,2,\ldots$, 
of Folland-Stein~\cite{FS:EDdbarbCAHG}. 

In we sequel we let $\N_{m}=\{1,\ldots,m\}$ and for any $I=(i_{1}, \ldots, i_{k})$ in $\N_{m}^{k}$ we set 
\begin{equation}
    X_{I}=X_{i_{1}}\ldots X_{i_{l}}.
\end{equation}

\begin{definition}[\cite{FS:EDdbarbCAHG}]
The Hilbert space $S_{k}^{2}(M)$, $k\in\N$, consists of functions $u\in L^{2}(M)$ such that 
  $(X_{I})u \in L^{2}(M)$ for any $I \in \cup_{j=1}^{k}\N_{m}^{j}$. It is endowed with the Hilbertian norm given by
  \begin{equation}
      \|u\|_{S_{k}^{2}}^{2}=\|u\|_{L^{2}}^{2} + \sum_{1\leq j\leq k} \sum_{I \in \N_{m}^{j}} 
      \|X_{I}u\|_{L^{2}}^{2}, \qquad u \in S_{k}^{2}(M). 
  \end{equation}
\end{definition}

\begin{proposition}\label{prop:Sobolev.equality-FS}
  For $k=1,2,\ldots$ the weighted Sobolev spaces $W_{H}^{k}(M)$ and $S_{k}^{2}(M)$ agree as spaces and carry the same topology.

\end{proposition}
\begin{proof}
    First, if $I \in \cup_{j=1}^{k}\N_{m}^{j}$ then 
    the differential operator $X_{I}$ is bounded from $W_{H}^{k}(M)$ to $L^{2}(M)$, so we get the estimate, 
    \begin{equation}
        \|u\|_{S^{k}}^{2}=\|u\|_{L^{2}}^{2} + \sum_{1\leq j\leq k} \sum_{I \in \N_{m}^{j}} 
      \|X_{I}u\|_{L^{2}}^{2} \leq C_{k}^{2}\|u\|_{W_{H}^{k}} \qquad u \in 
        C^{\infty}(M).
         \label{eq:Sobolev.SkWHk}
    \end{equation}
    
    On the other hand,   the vector fields $X_{j}$'s and their adjoints 
    $X_{j}^{*}$'s give rise to bounded linear maps from $S_{l+1}^{2}(M)$ to $S_{l}^{2}(M)$ for any $l\in \N$. Therefore,  for any integer $p$,  
    the operator 
    $(1+\Delta_{X})^{p}$ is bounded from $S^{2}_{l+2p}(M)$ to $S^{2}_{l}(M)$. It follows that, when  $k$ is even, we have 
    \begin{equation}
        \|u\|_{W_{H}^{k}}=\|(1+\Delta_{X})^{\frac{k}{2}}u\|_{L^{2}} \leq C_{k}\|u\|_{S^{k}}, \qquad u \in C^{\infty}(M).
    \end{equation}
    Moreover, for any $u \in C^{\infty}(M)$ we have 
    \begin{multline}
       \|(1+\Delta_{X})^{\frac{1}{2}}u\|_{L^{2}}^{2}=\acou{(1+\sum_{1\leq j\leq m}X^{*}_{j}X_{j})u}{u}_{L^{2}}\\ = 
        \|u\|_{L^{2}}^{2}+\sum_{1\leq j\leq m}\|X_{j}^{2}\|_{L^{2}}^{2} = \|u\|_{S^{2}_{1}}^{2}. 
         \label{eq:Sobolev.WHkSk}
    \end{multline}
    Thus $(1+\Delta_{X})^{\frac{1}{2}}$ is bounded from $S^{2}_{1}(M)$ to $L^{2}(M)$. Therefore,  if $k$ is odd, say 
    $k=2p+1$, then the operator 
    $(1+\Delta_{X})^{\frac{k}{2}}=(1+\Delta_{X})^{p}(1+\Delta_{X})^{\frac{1}{2}}$ 
    is bounded from $S_{k}^{2}(M)$ to $L^{2}(M)$. Hence~(\ref{eq:Sobolev.WHkSk}) is valid in the odd case as well. Together 
    with~(\ref{eq:Sobolev.SkWHk}) this implies that  $W_{H}^{k}(M)$ and $S_{k}^{2}(M)$ agree as spaces and carry the same topology. 
\end{proof}

Now, let $\cE$ be a Hermitian vector bundle over $M$. We can also define weighted Sobolev spaces of sections of $\cE$ as follows. 
Let $\nabla: C^{\infty}(M,\cE)\rightarrow C^{\infty}(M, T^{*}M\times \cE)$  be a connection on $\cE$ and define 
\begin{equation}
    \Delta_{\nabla,X}=\nabla_{X_{1}}^{*}\nabla_{X_{1}}+\ldots+\nabla_{X_{m}}^{*}\nabla_{X_{m}}.
\end{equation}
As in the scalar case,  
the complex powers $(1+\Delta_{\nabla,X})^{s}$, $s\in \C$, form
an analytic 1-parameter group of \psivdos\ such that $\ord (1+\Delta_{\nabla,X})^{s}=2s$ for any $s \in \C$.  

\begin{definition}
  $W_{H}^{s}(M,\cE)$, $s\in \R$, consists of all distributional sections $u\in \cD'(M,\cE)$ such that 
  $(1+\Delta_{\nabla,X})^{\frac{s}{2}}u\in L^{2}(M,\cE)$. It is endowed with the Hilbertian norm given by
  \begin{equation}
      \|u\|_{W_{H}^{s}}=\|(1+\Delta_{\nabla,X})^{\frac{s}{2}}u\|_{L^{2}}, \qquad u \in W_{H}^{s}(M,\cE). 
  \end{equation}
 \end{definition}

Along similar lines as that of the proof of Proposition~\ref{prop:Sobolev.embeddings} we can prove: 

\begin{proposition}
   1) As a topological space  $W_{H}^{s}(M,\cE)$, $s\in \R$, is independent of the choice of the connection $\nabla$ and 
    of the vector fields $X_{1},\ldots,X_{m}$.\smallskip 
    
    2)  We have the following continuous embeddings: 
    \begin{equation}
   \begin{array}{rccccl}
       L^{2}_{s}(M,\cE)  & \hookrightarrow & W_{H}^{s}(M,\cE) &\hookrightarrow& L^{2}_{s/2}(M,\cE) & \qquad \text{if $s\geq 0$},\\
       L^{2}_{s/2}(M,\cE)  & \hookrightarrow & W_{H}^{s}(M,\cE)& \hookrightarrow &L^{2}_{s}(M,\cE) & \qquad \text{if $s< 0$}.
   \end{array}
         \label{eq:Sobolev.embeddings2}
    \end{equation}
\end{proposition}

As a consequence we obtain the equalities of topological spaces,
\begin{equation}
       C^{\infty}(M,\cE)=\cap_{s \in \R} W_{H}^{s}(M,\cE) \qquad \text{and} \qquad \cD'(M,\cE)= \cup_{s \in \R}  W_{H}^{s}(M,\cE).
        \label{eq:Sobolev.smooth-distributions-WHs}
\end{equation}

Notice that we can also define Folland-Stein spaces $S_{k}^{2}(M,\cE)$, $k=1,2,\ldots$, as in the scalar case, by using the 
differential operators $\nabla_{X_{I}}=\nabla_{X_{i_{1}}}\ldots \nabla_{X_{i_{k}}}$, $I \in 
\cup_{j=1}^{k}\N_{m}^{j}$. Then, by arguing as in the proof  of Proposition~\ref{prop:Sobolev.equality-FS}, we can show that the spaces $W_{H}^{k}(M,\cE)$ 
and $S_{k}^{2}(M,\cE)$ agree and bear the same topology.  

Now, the Sobolev spaces $W_{H}^{s}(M,\cE)$ yield sharp regularity results for \psivdos. 

\begin{proposition}\label{prop:Sobolev.regularity-PsiHDOs1}
    Let $P:C^{\infty}(M,\cE)\rightarrow C^{\infty}(M,\cE)$ be a \psivdo\ of order $m$ and set $k=\Re m$. 
    Then, for any $s\in \R$, the operator $P$ extends to a continuous linear mapping from 
    $W_{H}^{s+k}(M,\cE)$ to $W_{H}^{s}(M,\cE)$. 
\end{proposition}
\begin{proof}
As $P_{s}=(1+\Delta_{\nabla,X})^{\frac{s}{2}}P(1+\Delta_{\nabla,X})^{-\frac{(s+k)}{2}}$ is a \psivdo\ with purely imaginary order, 
    Proposition~\ref{prop:PsiHDO.Sobolev-regularity} tells us that it gives rise to a bounded operator on $L^{2}(M,\cE)$. Therefore, we have 
    \begin{equation}
        \|Pu\|_{W^{s}_{H}}=\|P_{s}(1+\Delta_{\nabla,X})^{s+k}u\|_{L^{2}} \leq C_{s}\|u\|_{W^{s+k}_{H}}, \qquad u \in C^{\infty}(M,\cE), 
    \end{equation}
    It then follows that $P$ extends to a continuous linear mapping from $W_{H}^{s+k}(M,\cE)$ to $W_{H}^{s}(M,\cE)$. \
\end{proof}

\begin{proposition}\label{prop:Sobolev.regularity-PsiHDOs2}
  Let $P:C^{\infty}(M,\cE)\rightarrow C^{\infty}(M,\cE)$ be a \psivdo\ of order $m$ such that $P$ satisfies the Rockland condition at every point and 
  set $k=\Re m$. 
  Then for any $u \in \cD'(M,\cE)$ we have
    \begin{equation}
        Pu \in W_{H}^{s}(M,\cE) \Longrightarrow u \in W_{H}^{s+k}(M,\cE).
         \label{eq:Sobolev.hypoellipricity-WHs}
    \end{equation}
   In fact, for any $s'\in \R$ we have the estimate,
    \begin{equation}
        \|u\|_{W^{s+k}_{H}} \leq C_{ss'}(\|Pu\|_{W^{s}_{H}}+\|u\|_{W^{s'}_{H}}), \qquad u \in W^{s+k}_{H}(M,\cE).
        \label{eq:Sobolev.hypoelliticity-WHs-estimates}
    \end{equation}
\end{proposition}
\begin{proof}
 Since $P$ satisfies the Rockland condition at every point and the bracket condition $H+[H,H]=TM$ holds, it follows from the proof of 
 Proposition~\ref{prop:Heat-Rockland.Rockland-hypoellipticity} 
 that $P$ admits a left parametrix in $\pvdo^{-m}(M,\cE)$, i.e., there exist $Q$ in $\pvdo^{-m}(M,\cE)$ and $R$ in 
   $\psinf(M,\cE)$ such that $QP=1-R$. 
    Therefore, for any $u \in \cD'(M,\cE)$ we have $ u=QPu+Ru$.  
     Thanks to the first part we know that  $Q$ maps 
    $W_{H}^{s}(M,\cE)$ to $W_{H}^{s+k}(m,\cE)$. Since $R$ is smoothing, and so $Ru$ always is smooth, it follows
    that if $Pu$ is in $W_{H}(M,\cE)$ then $u$ 
    must be in $W_{H}(M,\cE)$.
    
    In fact, as $Q$ is actually bounded from 
    $W_{H}^{s}$ to $W_{H}^{s+k}$ and~(\ref{eq:Sobolev.smooth-distributions-WHs}) 
    implies that $R$ is bounded from any space $W_{H}^{s'}(M,\cE)$ to $W_{H}^{s+k}(M,\cE)$, on $ W^{s+k}_{H}(M,\cE)$ we have 
    estimates, 
    \begin{equation}
        \|u\|_{W^{s+k}_{H}}\leq \|QPu\|_{W^{s+k}_{H}} + \|Ru\|_{W^{s+k}_{H}} \leq 
        C_{ss'}(\|Pu\|_{W^{s}_{H}}+\|u\|_{W^{s'}_{H}}). 
    \end{equation}
    The proof is thus complete.
\end{proof}
\begin{remark}
    By combining Proposition~\ref{prop:Sobolev.regularity-PsiHDOs2} with the embeddings~(\ref{eq:Sobolev.embeddings}) 
    and~(\ref{eq:Sobolev.embeddings2}) we recover the Sobolev regularity of \psivdos\ as in Proposition~\ref{prop:PsiHDO.Sobolev-regularity} 
    as well as the hypoelliptic estimates~(\ref{eq:PsiHDO.subellipticity.subelliptic-estimates}).
\end{remark}

\begin{remark}
    When $\Re m$ is an integer it follows from Proposition~\ref{prop:Sobolev.equality-FS} that the estimates~(\ref{eq:Sobolev.hypoelliticity-WHs-estimates}) 
    are equivalent to maximal hypoellipticity in the sense of~\cite{HN:HMOPCV}. 
\end{remark}

On the other hand, the spaces $W_{H}^{s}(M,\cE)$ can be localized as follows. 

\begin{definition}\label{def:Sobolev.loaclization-WHs}
    We say that $u \in \cD'(M,\cE)$ is $W_{H}^{s}$ near a point $x_{0}\in M$ whenever there exists $ \varphi\in C^{\infty}(M)$ such that 
    $\varphi(x_{0})\neq 0$  and $\varphi u$ is in $W_{H}^{s}(M,\cE)$. 
\end{definition}
 
This definition depends only on the germ of $u$ at $x_{0}$ because we have:
\begin{lemma}\label{lem:Sobolev.localization-WHs}
    Let $u \in \cD'(M,\cE)$ be  $W_{H}^{s}$ near $x_{0}$. Then for any $\varphi\in C^{\infty}(M)$ with a small enough support 
    about $x_{0}$ the distribution $\varphi u$ lies in $W_{H}^{s}(M,\cE)$.
\end{lemma}
\begin{proof}
    Let $\varphi\in C^{\infty}(M)$ be such that $\varphi(x_{0})\neq 0$ and $\varphi u$ is in $W_{H}^{s}(M,\cE)$ and let $\psi\in C^{\infty}(M)$ be 
    so that $\psi(x_{0})\neq 0$ and $\supp \psi \cap \varphi^{-1}(0)\neq \emptyset $. 
    Then $\chi:=\frac{\psi}{\varphi}$ is in $C^{\infty}(M)$ and $(1+\Delta_{\nabla,X})^{\frac{s}{2}}\psi u$ is equal to
    \begin{equation}
         (1+\Delta_{\nabla,X})^{\frac{s}{2}}\chi \varphi u = 
       (1+\Delta_{\nabla,X})^{\frac{s}{2}}\chi  (1+\Delta_{\nabla,X})^{-\frac{s}{2}}.  (1+\Delta_{\nabla,X})^{\frac{s}{2}}\varphi u.
    \end{equation}
    Since $(1+\Delta_{\nabla,X})^{\frac{s}{2}}\varphi u$ is in $L^{2}(M,\cE)$ and $ (1+\Delta_{\nabla,X})^{\frac{s}{2}}\chi  
    (1+\Delta_{\nabla,X})^{-\frac{s}{2}}$ is a zero'th order \psivdo, so maps $L^{2}(M,\cE)$ to itself, it follows that $\psi u$ lies in 
    $W_{H}^{s}(M,\cE)$. Hence the lemma.
\end{proof}

We can now get a localized version of~(\ref{eq:Sobolev.hypoellipricity-WHs}).
\begin{proposition}\label{prop:Sobolev.hypoellipricity-WHs-localized}
  Let $P \in \pvdo^{m}(M,\cE)$ have an invertible principal symbol, set $k=\Re m$ and let $s\in \R$. Then for any $u \in \cD'(M,\cE)$ we have 
  \begin{equation}
      \text{$Pu$ is $W_{H}^{s}$ near $x_{0}$}\ \Longrightarrow \ \text{$u$ is $W_{H}^{s}$ near $x_{0}$}.
  \end{equation}
\end{proposition}
\begin{proof}
    Assume that $Pu$ is $W_{H}^{s}$ near $x_{0}$ and let $\varphi\in C^{\infty}(M)$ be such that $\varphi(x_{0})\neq 0$ and $\varphi u$ is in 
    $W_{H}^{s}(M,\cE)$. Thanks to Lemma~\ref{lem:Sobolev.localization-WHs} we may assume that 
 $\varphi=1$ near $x_{0}$. Let $\psi \in C^{\infty}(M)$ be such that $\psi(x_{0})\neq 0$ and $\varphi=1$ near $\supp \psi$. 
    Since the principal symbol of $P$ is invertible there exist $Q$ in $\pvdo^{-m}(M,\cE)$ and $R$ in $\psinf(M,\cE)$ such that $QP=1-R$. 
    Thus for any $u \in \cD'(M,\cE)$ we have 
    \begin{equation}
        \psi u=\psi QPu+\psi Ru= \psi Q\varphi Pu+ \psi Q(1-\varphi) Pu +\psi Ru.
    \end{equation}
   
    In the above equality $Ru$ is a smooth function since $R$ is a smoothing operator. Similarly, as  $\psi$ and $1-\varphi$ have disjoint supports, the 
    operator $\psi Q(1-\varphi) P$ is a smoothing operator and so $\psi Q(1-\varphi) Pu$ is a smooth function.  
    In addition, since $\varphi Pu$ is in $W_{H}^{s}(M,\cE)$ and $ \psi Q$ is a \psivdo\ of order $-m$, and so
    maps $W_{H}^{s}(M,\cE)$ to $W_{H}^{s+k}(M,\cE)$ , it follows 
    that $\psi u$ is in $W_{H}^{s+k}$. Hence $u$ is $W_{H}^{s+k}$ near $x_{0}$. 
\end{proof}

Next, Proposition~\ref{prop:Sobolev.regularity-PsiHDOs1} admits a generalization to holomorphic families.

\begin{proposition}\label{prop:Sobolev.regularity-PsiHDOs-families}
    Let $\Omega \subset \C$ be open and let $(P_{z})_{z\in \Omega}\in \Hol(\Omega, \Psi^{*}_{H}(M,\cE))$ be such that 
    $m:=\sup_{z \in \Omega}\Re\ord P_{z}<\infty$. Then for any $s\in \R$ the family $(P_{z})_{z\in \Omega}$ defines 
    a holomorphic family with values in $\cL(W_{H}^{s+m}(M,\cE),W_{H}^{s}(M,\cE))$. 
\end{proposition}
\begin{proof}
    Let $V\subset \Rd$ be a Heisenberg chart with $H$-frame $Y_{0},Y_{1},\ldots,Y_{d}$ and let $(Q_{z})_{z\in 
    \Omega}\in \Hol(\Omega, \pvdo^{*}(V))$ be such that $\Re \ord Q_{z}\leq 0$. 
    Then  we can write
    \begin{equation}
        Q_{z}=p_{z}(x,-iY)+R_{z},
    \end{equation}
    for some families $(p_{z})_{z\in \Omega}\in \Hol(\Omega, S^{*}(V\times \Rd))$ and $(R_{z})_{z\in \Omega}\in \Hol(\Omega, \Psi^{\infty}(V))$.
    Since $\Re \ord p_{z}=\Re \ord Q_{z}\leq 0$ we see that $(p_{z})_{z\in \Omega}$ belongs to $\Hol(\Omega,S^{0}_{\|}(\URd))$. 
    
    Next, for $j=1,\ldots,d$ let $\sigma_{j}$ denote the standard symbol of 
    $-iY_{j}$ and set $\sigma=(\sigma_{0},\ldots,\sigma_{d})$. 
    Then it follows from the proof of Proposition~\ref{prop:Complexpowers.operators.properties} that the family $(p_{z}(x,\sigma(x,\xi)))_{z\in 
    \Omega}$ lies in $\Hol(\Omega, S_{\frac{1}{2},\frac{1}{2}}(V\times \Rd))$. 
  Since by~\cite{Hw:L2BPO} the 
    quantization map $q\rightarrow q(x,D_{x})$ is continuous from $S_{\frac{1}{2}\frac{1}{2}}(\URd)$ to 
    $\cL(L^{2}_{\op{loc}}(U), L^{2}(U))$, we deduce that $(p_{z}(x,-iX))_{z\in \Omega}$ and $(Q_{z})_{z\in \Omega}$ are 
    holomorphic families with values in $\cL(L^{2}_{\op{loc}}(U), L^{2}(U))$. 
    
    It follows from the above result that any family $(Q_{z})_{z\in \Omega}\in \Hol(\Omega, \pvdo^{*}(M,\cE))$ 
    such that $\Re \ord Q_{z}\leq 0$ gives rise to a holomorphic family with values in $\cL(L^{2}(M,\cE))$. 
    This applies in particular to the family, 
    \begin{equation}
        Q_{z}^{(s)}=(1+\Delta_{\nabla,X})^{\frac{s}{2}}Q_{z}(1+\Delta_{\nabla,X})^{-\frac{m+s}{2}}, \quad z \in \Omega. 
    \end{equation}
    As $Q_{z}=(1+\Delta_{\nabla,X})^{-\frac{s}{2}}Q_{z}^{(s)}(1+\Delta_{\nabla,X})^{\frac{m+s}{2}}$ it follows that $(Q_{z})_{z\in 
    \Omega}$ gives rise to a holomorphic family with values in $\cL(W_{H}^{s+m}(M,\cE),W_{H}^{s}(M,\cE))$.
\end{proof}

Combining Proposition~\ref{prop:Sobolev.regularity-PsiHDOs-families} with Theorem~\ref{thm:Powers1.main2} 
we immediately get:

\begin{proposition}\label{prop:Sobolev.regularity-complex-powers}
   Let $P:C^{\infty}(M,\cE)\rightarrow C^{\infty}(M,\cE)$ be a positive differential operator of 
even (Heisenberg) order $v$ such that $P$ satisfies the Rockland condition at every point. 
Then, for any reals $m$ and $s$, we have
\begin{equation}
    (P^{z})_{\Re z<m} \in \Hol(\Re z<m, \cL(W_{H}^{s+mv}(M,\cE),W_{H}^{s}(M,\cE))).
\end{equation}
\end{proposition}

%%%%%%%%%%%%%%%%%%%%%%%%%%%%%%%%%%%%
%%%%%%%%%% Chap.6 tex %%%%%%%%%%%%%%%%%%%%
%%%%%%%%%%%%%%%%%%%%%%%%%%%%%%%%%%%%

\chapter{Spectral Asymptotics for Hypoelliptic Operators}
%\section{Weyl asymptotics}
\label{chap:Spectral}
In this chapter we apply the results of the previous chapters to derive spectral asymptotics for 
hypoelliptic differential operators on Heisenberg manifolds. 

In Section~\ref{sec:spectral} we get spectral asymptotics for general hypoelliptic differential operators.  We then give precise geometric expressions for 
these asymptotics. First, in  Section~\ref{sec:Spectral-CR} we deal with the Kohn Laplacian and the horizontal sublaplacian on a $\kappa$-strictly pseudoconvex 
CR manifold, as well as for the 
Gover-Graham operators in the strictly pseudoconvex case.  Second, we deal with 
the horizontal sublaplacian and the contact Laplacian on a contact manifold in Section~\ref{sec:Spectral-contact}. 

\section{Spectral asymptotics for hypoelliptic operators}
\label{sec:spectral}
In this section we explain how Theorem~\ref{thm:Heat1.main} enables us to derive spectral asymptotics for hypoelliptic operators on Heisenberg 
manifolds. Throughout this section we let $(M^{d+1},H)$ denote a compact Heisenberg manifold  equipped with a smooth density~$>0$ and 
let $\cE$ be a Hermitian vector bundle over $M$. 

Let $P:C^{\infty}(M,\cE)\rightarrow C^{\infty}(M,\cE)$ be a selfadjoint differential operator 
of even Heisenberg order $m$ which is bounded from below and such that the principal symbol of $P+\partial_{t}$ is an invertible Volterra-Heisenberg 
symbol. Recall that by Proposition~\ref{thm:Powers1.heat-sublaplacians} and Theorem~\ref{thm:Heat1.main} 
the latter occurs when at least one of the following conditions holds:\smallskip 
% 
% satisfies the Rockland condition
% 
% be a selfadjoint of  such that $P$ is bounded from below and 
% and the principal symbol of . Recall that by Proposition~\ref{thm:Powers1.heat-sublaplacians} 
% and Theorem~\ref{thm:Heat1.main} the latter condition is satisfied in the following cases:\smallskip

- $P$ is a sublaplacian and satisfies the condition~(\ref{eq:Sublaplacian.condition}) at every point.\smallskip

- The bracket condition $H+[H,H]=TM$ holds and $P$ satisfies the Rockland condition at every point.\smallskip

\noindent As an immediate consequence of Theorem~\ref{thm:Powers1.heat-kernel-asymptotics} we get: 

\begin{proposition}[{\cite[Thm.~5.6]{BGS:HECRM}}]\label{thm:Heat1.heat-trace-asymptotics}
    As $t\rightarrow 0^{+}$ we have
        \begin{equation}
         \Tr e^{-tP} \sim t^{-\frac{d+2}{m}} \sum t^{\frac{2j}{m}} A_{j}(P), \qquad  
   A_{j}(P)=\int_{M}\tr_{\cE} a_{j}(P)(x),
          \label{eq:Heat1.heat-trace-asymptotics}
     \end{equation}
     where the density $a_{j}(P)(x)$ is the coefficient of $t^{\frac{j-d+2}{m}}$ in the heat kernel 
     asymptotics~(\ref{eq:Rockland-Heat.heat-kernel-asymptotics}) for $P$.\smallskip
\end{proposition}

Let $\lambda_{0}(P)\leq \lambda_{1}(P)\leq \ldots$ denote the eigenvalues of $P$ counted with multiplicity and let $N(P;\lambda)$ be the counting 
function of $P$, that is, 
\begin{equation}
     N(P;\lambda)=\#\{k\in \N;\ \lambda_{k}(P)\leq \lambda \}, \qquad \lambda \in \R. 
\end{equation}

\begin{proposition}\label{thm:Heat1.spectral-asymptotics}
     1) We have $A_{0}(P)>0$.\smallskip 
     
     2) As $\lambda\rightarrow \infty$ we have 
\begin{equation}
    N(P;\lambda) \sim \nu_{0}(P)\lambda^{\frac{d+2}{m}}, \qquad 
    \nu_{0}(P)=\Gamma(1+\frac{d+2}{m})^{-1}A_{0}(P).
    \label{eq:Heat1.counting-function-asymptotics}
\end{equation}

  3) As $k\rightarrow \infty$ we have 
     \begin{equation}
         \lambda_{k}(P)\sim \left(\frac{k}{\nu_{0}(P)}\right)^{\frac{m}{d+2}}. 
         \label{eq:Heat1.eigenvalue-asymptotics}
     \end{equation}     
 \end{proposition}
 \begin{proof}
 First,  we  have $A_{0}(P)=\lim_{t\rightarrow 0^{+}}t^{\frac{d+2}{m}}\Tr e^{-tP}\geq 0$, so if we can show that $A_{0}(P)\neq 0$ then we get 
     $A_{0}(P)>0$ and the 
     asymptotics~(\ref{eq:Heat1.counting-function-asymptotics}) would follow from~(\ref{eq:Heat1.heat-trace-asymptotics}) by Karamata's Tauberian theorem 
     (see, e.g., \cite[Thm.~108]{Ha:DS}). This would also give~(\ref{eq:Heat1.eigenvalue-asymptotics}), since the latter is equivalent 
     to~(\ref{eq:Heat1.counting-function-asymptotics}) (see, e.g.,~\cite[Sect.~13]{Sh:POST}). Therefore, the bulk of the proof is to show that we 
     have $A_{0}(P)\neq 0$.
       
  Second, notice that there is at least one integer~$< \frac{m}{d+2}$ such that $A_{j}(P)\neq 0$. Otherwise, by~(\ref{eq:Heat1.heat-trace-asymptotics}) 
  there would exist a constant $C> 0$ such that $\Tr e^{-tP}\leq C$ for $0\leq t\leq 1$. Thus, 
    \begin{equation}
        k e^{-t\lambda_{k}(P)}\leq \sum_{j<k} e^{-t\lambda_{j}(P)}\leq \Tr e^{-tP}\leq C, \qquad 0<t<1.
    \end{equation}
   Therefore, letting $t\rightarrow 0^{+}$ would give  $k\leq C$ for every $k\in \N$, which is not possible since $P$ is unbounded.
   
   Next, when $m\geq d+2$ the only integer~$<\frac{m}{d+2}$ is $j=0$, so in this case we must have $A_{0}(P)\neq  0$. 
    
    Suppose now that  we have $m<d+2$ and $A_{0}(P)=0$. Let $\mu=\frac{d+2}{m}-j_{0}$, where $j_{0}$ is the smallest integer $j$ such that $A_{j}(P)\neq 0$. 
    Notice that since 
  $1\leq j_{0}< \frac{d+2}{m}$, we have $0<\mu\leq \frac{d+2}{m}-1$. Moreover, the asymptotics~(\ref{eq:Heat1.heat-trace-asymptotics}) becomes
   \begin{equation}
       \Tr e^{-tP}=A_{j_{0}}(P)t^{-\mu}+\op{O}(t^{1-\mu}) \qquad \text{as $t\rightarrow 0^{+}$}. 
   \end{equation}
   This implies that we have $A_{j_{0}}(P)\geq 0$. Since we have $A_{j_{0}}(P)\neq 0$ by definition of $j_{0}$, we get 
   $A_{j_{0}}(P)>0$. Therefore, as alluded to above, it follows from Karamata's Tauberian theorem that as $k\rightarrow \infty$ we have
   \begin{equation}
       \lambda_{k}(P)\sim  \left(\frac{k}{\beta}\right)^{\frac{1}{\mu}}, \qquad \beta=\Gamma(1+\mu)^{-1}A_{j_{0}}(P). 
   \end{equation}
It follows that $\lambda_{k}(P^{- \frac{d+2}{2m}})=\op{O}(k^{-\frac{1}{2}-\delta})$, with 
   $\delta=\frac{1}{2}(\frac{1}{\mu}\frac{d+2}{m}-1)>0$. In particular, we have $\sum_{k\geq 0}\lambda_{k}(P^{- \frac{d+2}{2m}})^{2}<\infty$, that is,
 $P^{- \frac{d+2}{2m}}$ is a Hilbert-Schmidt operator on $L^{2}(M,\cE)$. 
   
   In addition, observe that $P^{- \frac{d+2}{2m}}$ is a \psivdo\ of order $-\frac{(d+2)}{2}$ and that any operator $Q \in \pvdo^{-\frac{(d+2)}{2}}(M,\cE)$ can 
   be written as 
   \begin{equation}
       Q=\Pi_{0}(P)Q +(1-\Pi_{0}(P))P^{-\frac{d+2}{2m}}P^{\frac{d+2}{2m}}Q,
        \label{eq:Spectral.Hilbert-Schmidt}
   \end{equation}
   where $\Pi_{0}(P)$ denotes the orthogonal projection onto $\ker P$. Recall that the space of Hilbert-Schmidt operators is a two-sided ideal of 
   $\cL(L^{2}(M<\cE)$. Observe that in~(\ref{eq:Spectral.Hilbert-Schmidt}) the projection $\Pi_{0}(P)$ is a smoothing operator, hence is a Hilbert-Schmidt 
   operator. Moreover,  $P^{\frac{d+2}{2m}}Q$ is a bounded operator on $L^{2}(M,\cE)$, since this is 
   a zero'th order \psivdo. Therefore, we see that any $Q\in \pvdo^{-\frac{(d+2)}{2}}(M,\cE)$ is a Hilbert-Schmidt operator on $L^{2}(M,\cE)$.

    We now get a contradiction as follows. Let $\kappa: U \rightarrow V$ be a Heisenberg chart over which there is a trivialization 
    $\tau:\cE_{_{U}}\rightarrow U\times \C^{r}$ of $\cE$ 
    and such that the open $V\subset \Rd$ is bounded. 
    Let $\varphi \in C^{\infty}_{c}(\Rd)$ have non-empty support $L$, let $\chi \in C^{\infty}_{c}(V\times V)$ be such that $\chi(x,y)=1$ near $L\times L$, and
   let $Q:C_{c}^{\infty}(V)\rightarrow C^{\infty}(V)$ be given by the kernel,
    \begin{equation}
        k_{Q}(x,y)=|\varepsilon_{x}'|\varphi(x)\|\varepsilon_{x}(y)\|^{-\frac{d+2}{2}}\chi(x,y).
    \end{equation}
    
    The kernel $k_{Q}(x,y)$ has a compact support contained in $V\times V$ and, as $\varphi(x)(1-\chi(x,y))$ vanishes near the diagonal of $V\times V$, we have 
\begin{equation}
         k_{Q}(x,y)=|\varepsilon_{x}'|\varphi(x)\|-\varepsilon_{x}(y)\|^{-\frac{d+2}{2}} \bmod C^{\infty}(V\times V). 
\end{equation}
     Since $\varphi(x)\|y\|^{-\frac{d+2}{2}}$ belongs to $\cK^{-\frac{d+2}{2}}(V\times \Rd)$, 
     it follows from Proposition~\ref{prop:PsiVDO.characterisation-kernel2} 
    that $Q$ is a \psivdo\ of order $-\frac{d+2}{2}$. 
    
    Let $Q_{0}=\tau^{*}\kappa^{*}(Q\otimes 1)$. Then $Q_{0}$ belongs to $\pvdo^{-\frac{d+2}{2}}(M,\cE)$, hence is a Hilbert-Schmidt operator on 
    $L^{2}(M,\cE)$. This implies that  $Q$ is a Hilbert-Schmidt operator on $L^{2}(V)$, so by~\cite[p.~109]{GK:ITLNSO} its kernel lies in $L^{2}(V\times V)$. 
    This cannot be true, however, because we have 
   \begin{multline}
       \int_{V\times V}|k_{Q}(x,y)|^{2}dxdy\geq \int_{L\times L} |\varepsilon_{x}'|^{2}|\varphi(x)|^{2}\|\varepsilon_{x}(y)\|^{-(d+2)}dxdy \\  
       = \int_{L}|\varepsilon_{x}'||\varphi(x)|^{2}(\int_{\varepsilon_{x}(L)}\|y\|^{-(d+2)}dy)dx=\infty.
   \end{multline}
   This gives a contradiction, so we must have $A_{0}(P)\neq 0$. The proof is therefore complete. 
\end{proof}
 
\begin{example}
    Proposition~\ref{thm:Heat1.spectral-asymptotics} is valid for the following operators:\smallskip
  
   (a) Real selfadjoint sublaplacian $ \Delta=\nabla_{X_{1}}^{*}\nabla_{X_{1}}+\ldots+\nabla_{X_{m}}^{*}\nabla_{X_{m}}^{2}+\mu(x)$, 
   where $\nabla$ is connection on $\cE$, the vector fields $X_{1},\ldots,X_{m}$ span $H$ and $\mu(x)$ is a selfadjoint section of $\End \cE$, 
   provided that the Levi form of $(M,H)$ is non-vanishing; \smallskip 
  
   (b) The Kohn Laplacian on a compact CR manifold and acting on $(p,q)$-forms when the condition $Y(q)$ holds everywhere;\smallskip 
  
   (c) The horizontal sublaplacian on a compact Heisenberg manifold acting on horizontal forms of degree $k$ when the condition $X(k)$ holds 
   everywhere;\smallskip 
   
   (d) The horizontal sublaplacian on a compact CR manifold acting on $(p,q)$-forms when the condition $X(p,q)$ holds 
   everywhere;\smallskip 
   
 (e) The Gover-Graham operators $\boxdot_{\theta}^{(k)}$ on a strictly pseudoconvex CR manifold $M^{2n+1}$ with $k\neq n+1$;\smallskip 
 
   (f) The contact Laplacian on a compact orientable contact manifold.%\smallskip
% 
%    (f) The  of the horizontal sublaplacian acting on functions on a compact strictly pseudoconvex CR manifold.\smallskip 
\end{example}
 
Several authors have obtained Weyl asymptotics closely related to~(\ref{eq:Intro1.counting-function-asymptotics}) for bicharacteristic hypoelliptic 
 operators~(see~\cite{II:PDPEAABSFSP}, \cite{Me:HOCVC2WE}, \cite{MS:ECH2}), including sublaplacians on Heisenberg manifolds, and for more general hypoelliptic 
 operators~(see~\cite{Mo:ESOHCM1}, \cite{Mo:ESOHCM2}) using different approaches involving other pseudodifferential calculi. 
 
 These authors obtained their results in settings more general than the Heisenberg setting but, as far as the Heisenberg setting is concerned, our approach using the 
 Volterra-Heisenberg calculus presents the following advantages:\smallskip
 
 (i) The pseudodifferential analysis is somewhat simpler, since the Volterra-Heisenberg calculus yields for free 
 a heat kernel asymptotics, once the principal symbol of the heat operator is shown to be invertible, for which it is enough to use the Rockland condition 
 in many cases;\smallskip 
 
 (ii) As the Volterra-Heisenberg calculus takes fully into account the underlying Heisenberg geometry of the manifold and is invariant by change of 
 Heisenberg coordinates, we can explicitly compute the coefficient $\nu_{0}(P)$ in~(\ref{eq:Heat1.counting-function-asymptotics}) 
 for operators admitting a normal form. 
 This is illustrated  below by Proposition~\ref{prop:Spectral.normal-form-nu0P}, which will be used in the next two sections to give geometric 
 expressions for the Weyl asymptotics~(\ref{eq:Heat1.counting-function-asymptotics}) for geometric operators on CR and contact manifolds.\smallskip

Next, prior to dealing with operators admitting a normal form we need the following lemma.
 \begin{lemma}\label{lem:Spectral.a0(P)}
    Let $P:C^{\infty}(M,\cE)\rightarrow C^{\infty}(M,\cE)$ be a selfadjoint differential operator of Heisenberg order $m$ such that $P$ is bounded 
    from below and the principal symbol of $P+\partial_{t}$ is an invertible Volterra-Heisenberg symbol. Then, for any $a\in M$, the following 
    hold.\smallskip
    
    1) The model operator $P^{a}+\partial_{t}$ admits a unique fundamental solution $K^{a}(x,t)$ which, in addition, belongs  to the class 
    $\cK_{\op{v},-(d+2)}(G_{a}M\times\R_{(m)})$.\smallskip
    
    2) In Heisenberg coordinates centered at $a$ we have
    \begin{equation}
        a_{0}(P)(0)=K^{a}(0,1).
    \end{equation}
\end{lemma}
\begin{proof}
    First, it is enough to prove the result in a trivializing Heisenberg chart $U$ near $a$. Furthermore, we may as well assume that $P$ is scalar 
    because the proof for systems follows along similar lines.
      
    Let  $Q\in \pvhdo^{-m}(U\times\R_{(m)})$ be a parametrix for $P+\partial_{t}$ on $U\times \R$. Let $\sigma_{-m}(Q)\in S_{\op{v},-m}(\fg^{*}U\times 
    \R_{(m)})$ be its global principal symbol, so that $\sigma_{-m}(Q)$
    the inverse of $\sigma_{m}(P)+i\tau$, and let $Q^{a}$ be the left-invariant Volterra-\psivdo\ on $G_{a}U\times \R$ with symbol 
    $\sigma_{-m}(Q)(a,.,.)$. In particular, the operator $Q^{a}$ is the inverse of $P^{a}+\partial_{t}$ on $\cS(G_{a}U\times \R)$ and it agrees 
    with the left-convolution by $K^{a}=[\sigma_{-m}(a,.,t)]_{(\xi,\tau)\rightarrow (y,t)}^{\vee}$. 
   Therefore,  the left-convolution operator with $[(P^{a}+\partial_{t})K^{a}]$ agrees with $(P^{a}+\partial_{t})Q^{a}=1$. Thus, 
    \begin{equation}
        (P^{a}+\partial_{t})K^{a}(y,t)=\delta(y,t),
    \end{equation}
    that is, $K^{a}(y,t)$ is a fundamental solution for $P^{a}+\partial_{t}$. 
    
    Let $K\in \cS'(G_{a}U\times \R)$ be another fundamental solution for $P^{a}+\partial_{t}$.
    Then the left-convolution operator $Q$ by $K$ 
   is a right-inverse for $P^{a}+\partial_{t}$, so agrees with $Q^{a}$. Hence $K=K^{a}$, which shows that $K^{a}$ is the unique fundamental solution
  of  $P^{a}+\partial_{t}$.
  
  In addition, $K^{a}$ is in the class $\cK_{\op{v},-(d+2)}(G_{a}M\times\R_{(m)})$ because this the inverse 
  Fourier transform of a symbol in $S_{\op{v},-m}(\fg^{*}U\times \R_{(m)})$. 

   Finally, let $q_{-m}\in S_{\op{v},-m}(\URd\times\R_{(m)})$ be the local principal symbol of the parametrix $Q$.  
  As $Q$ has global principal symbol $\sigma_{-m}(Q)$ it follows from Proposition~\ref{prop:Powers1.principal-symbol} 
  that $q_{-m}(0,.,.)=\sigma_{-m}(Q)(a,.,.)$. Moreover, 
   since we are in the Heisenberg coordinates centered at $a$, the map $\varepsilon_{0}$ to the Heisenberg coordinates centered at $0$ is just the 
   identity. Combining this with~(\ref{eq:Heat1.counting-function-asymptotics}) then gives
   \begin{equation}
       a_{0}(P)(0)=|\varepsilon_{0}'|\check{q}_{m}(0,0,1)=[\sigma_{-m}(Q)]_{(\xi,\tau)\rightarrow (y,t)}(a,0,1)=K^{a}(0,1).
   \end{equation}
   The proof is thus achieved.
 \end{proof}

Now, assume that the Levi form of $(M,H)$ has constant rank $2n$. Therefore, by Proposition~\ref{prop:Bundle.intrinsic.fiber-structure} 
the tangent Lie group bundle $GM$ is a fiber bundle with typical fiber $G=\bH^{2n+1}\times \R^{d-2n}$. 

Let $P:C^{\infty}(M,\cE)\rightarrow C^{\infty}(M,\cE)$ be a selfadjoint differential operator of Heisenberg order $m$ such that $P$ is bounded 
    from below and satisfies the Rockland condition at every point. 
    
    We further assume that the density $d\rho(x)$ of $M$ and the model operator of $P$ 
    at a point admit normal forms. By this it is meant that there exist $\rho_{0}>0$ and a differential operator $P_{0}:C^{\infty}(G,\C^{r})\rightarrow 
    C^{\infty}(G,\C^{r})$  such that for any point $a \in M$ there exist trivializing Heisenberg coordinates, herewith called 
    normal trivializing Heisenberg coordinates centered at $a$, with respect to which the following hold:\smallskip 
    
    (i) We have $d\rho(x)_{|_{x=0}}=[\rho_{0}dx]_{|_{x=0}}$;\smallskip 
    
    (ii) At $x=0$ the tangent group is $G$ and the model operator of $P$ is $P_{0}$.\smallskip
    
In the sequel, we let $\op{vol}_{\rho}M$ denote the volume of $M$ with respect to $\rho$, that is, 
\begin{equation}
    \op{vol}_{\rho}M=\int_{M}d\rho(x).
\end{equation}
    As it turns out the assumptions (i) and (ii) allows us to relate the Weyl asymptotics~(\ref{eq:Heat1.counting-function-asymptotics}) 
    for $P$ to $\op{vol}_{\rho}M$ as follows.

\begin{proposition} \label{prop:Spectral.normal-form-nu0P}
Under the assumptions (i) and (ii) as $\lambda \rightarrow \infty$ we have
\begin{equation}
  N(P;\lambda)\sim \frac{\nu_{0}(P_{0})}{\rho_{0}}(\op{vol}_{\rho}M)\lambda^{\frac{d+2}{m}}, 
  \quad \nu_{0}(P_{0})=\Gamma(1+\frac{d+2}{m})^{-1}\tr_{\C^{r}}K_{0}(0,1),
     \label{eq:Spectral.normal-form-nu0P}
\end{equation}
where $K_{0}(x,t)$ denotes the fundamental solution of $P_{0}+\partial_{t}$.
\end{proposition}
\begin{proof}
   By Proposition~\ref{thm:Heat1.spectral-asymptotics} as $\lambda \rightarrow \infty$ we have 
    \begin{equation}
        N(P;\lambda) \sim \nu_{0}(P)\lambda^{\frac{d+2}{m}}, \qquad \nu_{0}(P)=\Gamma(1+\frac{d+2}{m})^{-1}\int_{M}\tr_{\cE}a_{0}(P)(x).
        \label{eq:Heat1.counting-function-asymptotics2}
    \end{equation}

    On the other hand, by Lemma~\ref{lem:Spectral.a0(P)} in normal trivializing Heisenberg coordinates centered at point $a\in M$ we have
    \begin{equation}
      [ \tr_{\C^{r}} a_{0}(P)(x)dx]_{|_{x=0}}=[\tr_{\C^{r}}K_{0}(0,1)dx]_{|_{x=0}}=\tr_{\C^{r}}K_{0}(0,1)\rho_{0}^{-1}[d\rho(x)]_{|_{x=0}}.
    \end{equation}
  Hence $\tr_{\cE}a_{0}(P)(x)=\rho_{0}^{-1}\tr_{\C^{r}}K_{0}(0,1)d\rho(x)$. Thus, 
  \begin{equation}
   \nu_{0}(P)=\Gamma(1+\frac{d+2}{m})^{-1}\int_{M}\rho_{0}^{-1}\tr_{\C^{r}}K_{0}(0,1)d\rho(x)=\rho_{0}^{-1}\nu_{0}(P_{0})\op{vol}_{\rho}M,
  \end{equation}
  with $\nu_{0}(P_{0})=\Gamma(1+\frac{d+2}{m})^{-1}\tr_{\C^{r}}K_{0}(0,1)$. Combining this~(\ref{eq:Heat1.counting-function-asymptotics2}) 
  then proves the lemma.
\end{proof}

\section{Weyl asymptotics and CR geometry}
\label{sec:Spectral-CR}
The aim of this section is to express in geometric terms the Weyl asymptotics~(\ref{eq:Heat1.counting-function-asymptotics}) 
for the Kohn Laplacian, the horizontal sublaplacian and the Gover-Graham on CR manifolds with Levi metrics.  

Let $M^{2n+1}$ be a  $\kappa$-strictly pseudoconvex manifold,  $0\leq \kappa \leq \frac{n}{2}$,  equipped with a pseudohermitian contact form $\theta$. 
Let $T_{1,0}\subset T_{\C}M$ be the CR tangent bundle of $M$, set $T_{0,1}=\overline{T_{0,1}}$ and $H=\Re (T_{1,0}\oplus T_{0,1})$ and let  $L_{\theta}$ 
denote the Levi form on $T_{1,0}$ so that
\begin{equation}
     L_{\theta}(Z,W)=-id\theta (Z,\bar W)=i\theta([Z,\bar W]) \qquad \forall Z,W \in C^{\infty}(M,T_{1,0}).
     \label{eq:Heat1.Levi-form}
\end{equation}

Let $X_{0}$ be the Reeb vector field of $\theta$, so that $\imath_{X_{0}} \theta=1$ and $\imath_{X_{0}}d\theta=0$. This gives rise to the splitting, 
\begin{equation}
    T_{\C}M=T_{1,0}\oplus T_{0,1}\oplus (\C X_{0}).
%     \label{eq:¥}
\end{equation}
For $p,q=0,\ldots,n$ we let $\Lambda^{p,q}$ denote the bundle of $(p,q)$-forms associated to this splitting. Moreover, we have
\begin{equation}
    [Z,\bar W]= -iL_{\theta}(Z,W)T \quad \bmod T_{1,0}\oplus T_{0,1} \qquad \forall Z, W \in C^{\infty}(M,T_{1,0}).
     \label{eq:Heat1.commutators-metric}
\end{equation}

We endow $M$ with a Levi metric as follows (see also~\cite{FS:EDdbarbCAHG}). 
Let $\tilde{h}$ be a (positive definite) Hermitian metric on $T_{1,0}$. Then there exists a Hermitian-valued section $A$ of $\End_{\C}T_{1,0}$ such that 
\begin{equation}
    L_{\theta}(Z,W) =\tilde{h}(AZ,W) \qquad \forall Z,W\in C^{\infty}(M,T_{1,0}).
\end{equation}

Let $A=U|A|$ be the polar decomposition of $A$. Since by assumption the Hermitian form $L_{\theta}$ has signature $(n-\kappa,\kappa,0)$, hence is 
nondegenerate, the section $A$ is invertible and the section $U$ is an orthogonal matrix with only the eigenvalues $ 1$ and $-1$ with multiplicities $n-\kappa$ 
and $\kappa$ respectively. Therefore, we have the splitting,
\begin{equation}
    T_{1,0}=T_{1,0}^{+}\oplus T_{1,0}^{-}, \qquad T_{1,0}^{\pm}=\ker (U\mp 1), 
    \label{eq:Heat1.positive-negative-splitting}
\end{equation}
where the restriction of $L_{\theta}$ to $T_{1,0}^{+}$ (resp.~$T_{1,0}^{-}$) is positive definite (resp.~negative definite). We then get a Levi 
metric on $T_{1,0}$ by letting
\begin{equation}
    h(Z,W)=\tilde{h}(|A|Z,W), \qquad Z,W\in C^{\infty}(M,T_{1,0}).
\end{equation}
In particular, on the direct summands $T_{1,0}^{\pm}$ of the splitting~(\ref{eq:Heat1.positive-negative-splitting}) we have
\begin{equation}
    L_{\theta}(Z,W)=h(UZ,W)=\pm h(Z,W), \qquad Z,W\in C^{\infty}(M,T_{1,0}^{\pm}).
     \label{eq:Heat1.Levi-form.Levi-metric}
\end{equation}

We now extend $h$ into a Hermitian metric $h$ on $T_{\C}M$ by making the following requirements:
\begin{gather}
  h(X_{0},X_{0})=1, \qquad   h(Z,W)=\overline{h(\bar Z,\bar W)} \quad \forall Z, W\in T_{0,1},\\ 
  \text{The splitting $T_{\C}M=\C T\oplus T_{1,0}\oplus T_{0,1}$ is orthogonal with respect to $h$}.
\end{gather}

This allows us to express the Levi form $\cL_{\C}:(H\otimes \C)\times (H\otimes \C)\rightarrow T_{\C}M/(H\otimes \C)$ as follows. 
Since $\theta(T)=1$  we have 
\begin{equation}
    \cL_{\C}(X,Y)=\theta([X,Y])T=-d\theta(X,Y)T, \qquad X,Y\in C^{\infty}(M,H\otimes \C). 
\end{equation}
Therefore, if follows from~(\ref{eq:Heat1.Levi-form.Levi-metric}) that we have 
\begin{equation}
    \cL(Z,\bar W)=-iL_{\theta}(Z,W)=h(Z,iUW), \qquad Z,W\in C^{\infty}(M,T_{1,0}).
\end{equation}

Since $\cL$ is antisymmetric and the integrability condition $[T_{1,0},T_{1,0}]\subset T_{1,0}$ implies that $\cL_{\C}$ vanishes on $T_{1,0} \times 
T_{1,0}$ and on $T_{0,1} \times T_{0,1}$, we get
\begin{equation}
       \cL_{\C}(X,Y)=h(X,LY), \qquad X,Y\in  C^{\infty}(M,H\otimes \C),
\end{equation}
where $L$ is the antilinear antisymmetric section of $\End_{\R} (H\otimes \C)$ such that
\begin{equation}
    L(Z+\bar{W})=iUW-i(\overline{UZ}) \qquad \forall Z,W\in C^{\infty}(M,T_{1,0}).
     \label{eq:Spectral.Levi-form-Levi-metric2}
\end{equation}
In particular, since $U^{*}=U$ and $U^{2}=1$ we have $|L|=1$. 

Let $Z_{1},\ldots,Z_{n}$ be a local orthonormal frame for $T_{1,0}$ (with respect to $h$) and such that $Z_{1},\ldots,Z_{n-\kappa}$ 
span $T_{1,0}^{+}$ and $Z_{n-\kappa+1},\ldots,Z_{n}$ span $T_{1,0}^{-}$. Then $\{T,Z_{j},\bar Z_{j}\}$ is an orthonormal frame. In the sequel we 
will call such a frame an \emph{admissible orthonormal frame of} $T_{\C}M$. Then from~(\ref{eq:Heat1.Levi-form.Levi-metric}) we get: 
\begin{equation}
 L_{\theta}(Z_{j},\bar Z_{k})=\epsilon_{j}h(Z_{j},Z_{k})=\epsilon_{j}\delta_{jk},
\label{eq:Heat1.dtheta.admissible-frame}
\end{equation}
where $\epsilon_{j}=1$ for $j=1,\ldots,n-\kappa$ and $\epsilon_{j}=-1$ for $j=n-\kappa+1,\ldots,n$.

Let $\{\theta,\theta^{j},\theta^{\bar j}\}$ be the dual coframe of $T^{*}_{\C}$ associated to $\{T,Z_{j},\bar Z_{j}\}$. Then the volume form of 
$h$ is locally given by 
\begin{equation}
    \sqrt{h(x)} dx=i^{n}\theta\wedge \theta^{1}\wedge \theta^{\bar 1}\wedge \cdots \wedge \theta^{n}\wedge \theta^{\bar n}.
\end{equation}
Furthermore, because of~(\ref{eq:Heat1.dtheta.admissible-frame}) we have 
$d\theta =i \sum_{j=1}^{n} \epsilon_{j}\theta^{j}\wedge \theta^{\bar j} \bmod \theta \wedge T^{*}M$, so as $ \epsilon_{1}\ldots 
\epsilon_{n}=(-1)^{\kappa}$ we have $\theta \wedge d\theta^{n}= n! i^{n}(-1)^{\kappa}\theta\wedge \theta^{1}\wedge \theta^{\bar 1}\wedge \cdots \wedge 
    \theta^{n}\wedge \theta^{\bar n}$. 
%     is equal to
% \begin{equation}
% %     n! i^{n} \epsilon_{1}\ldots \epsilon_{n}\theta\wedge \theta^{1}\wedge \theta^{\bar 1}\wedge \cdots \wedge 
% %     \theta^{n}\wedge \theta^{\bar n}\\ = 
%     n! i^{n}(-1)^{\kappa}\theta\wedge \theta^{1}\wedge \theta^{\bar 1}\wedge \cdots \wedge 
%     \theta^{n}\wedge \theta^{\bar n}. 
% \end{equation}
Therefore, we get the following global formula for the volume form,
\begin{equation}
     \sqrt{h(x)} dx=\frac{(-1)^{\kappa}}{n!}  \theta \wedge d\theta^{n}. 
\end{equation}
In particular, we see that the volume form depends only on $\theta$ and not on the choice of the Levi metric. 

\begin{definition}
    The pseudohermitian volume of $(M,\theta)$ is 
    \begin{equation}
        \op{vol}_{\theta}M=\frac{(-1)^{\kappa}}{n!} \int_{M} \theta \wedge d\theta^{n}.
    \end{equation}
\end{definition}

We shall now relate the asymptotics~(\ref{eq:Heat1.counting-function-asymptotics}) for the Kohn Laplacian, the horizontal sublaplacian and the 
Gover-Graham operators to the volume 
$\op{vol}_{\theta}M$. To this end consider the Heisenberg group $\bH^{2n+1}=\R \times \R^{2n}$ together with the left-invariant basis of $\fh^{2n+1}$ 
given by
\begin{gather}
    X_{0}^{0}=\frac{\partial}{\partial x_{0}}, \qquad X_{j}^{0}=\frac{\partial}{\partial x_{j}}+x_{n+j}\frac{\partial}{\partial x_{0}},    
    \label{eq:Heat1.Heisenberg-vector-fields1}\\ 
     X_{n+j}^{0}=\frac{\partial}{\partial x_{n+j}}-x_{j}\frac{\partial}{\partial 
    x_{0}}, \quad  j=1,\ldots, n.
     \label{eq:Heat1.Heisenberg-vector-fields2}
\end{gather}

For $\mu \in \C$ let $\cL_{\mu}$ denote  the Folland-Stein sublaplacian, 
\begin{equation}
    \cL_{\mu}=-\frac{1}{2}(X_{1}^{0})^{2}+\ldots+(X_{2n}^{0})^{2})-i\mu X_{0}^{0}.
    \label{eq:Weyl.Folland-Stein-sublaplacian}
\end{equation}
For this sublaplacian the condition~(\ref{eq:Sublaplacian.condition.scalar}) reduces to $\mu\neq \pm n, \pm(n+2),\ldots$, 
so in this case the operators $\cL_{\mu}$ and $\cL_{\mu}+\partial_{t}$ are hypoelliptic and admit unique fundamental solutions since their symbols are invertible.
% Volterra-Heisenberg symbol (see~\cite{BGS:HECRM}). 

In the sequel it will be convenient to use the variable $x'=(x_{1},\ldots,x_{2n})$. 

\begin{lemma}
 For $|\Re \mu|<n$ the fundamental solution $ k_{\mu}(x_{0},x',t)$ of $\cL_{\mu}+\partial_{t}$ is equal to
    \begin{equation}
%        k_{\mu}(x_{0},x',t) = 
       \chi(t) (2\pi t)^{-n} \int_{-\infty}^{\infty}e^{ix_{0}\xi_{0}-\mu t\xi_{0}}(\frac{t\xi_{0}}{\sinh t\xi_{0}})^{n}
        \exp [-\frac{1}{t} \frac{t\xi_{0}}{\tanh t\xi_{0}}|x'|^{2}]d\xi_{0}.
        \label{eq:Heat1.fundamental-solution-cLlambda}
    \end{equation}
    where $\chi(t)$ denotes the characteristic function of $(0,\infty)$. 
\end{lemma}
\begin{proof}
  The fundamental solution $k_{\mu}(x_{0},x',t)$ is solution to the equation, 
   \begin{equation}
       (\cL_{\mu}+\partial_{t})k(x_{0},x',t)=\delta(x_{0})\otimes \delta(x')\otimes \delta(t).
        \label{eq:Heat1.heat-cLmu}
   \end{equation}
   Let us make a Fourier transform with respect to $x_{0}$, so that letting $\hat{k}(\xi_{0},x',t)=\hat{k}_{x_{0}\rightarrow 
   \xi_{0}}(\xi_{0},x',t)$ the equation~(\ref{eq:Heat1.heat-cLmu}) becomes
    \begin{gather}
        (\hat{\cL}_{0}+\mu \xi_{0}+\partial_{t})\hat{k}_{\mu}=\delta(x')\otimes \delta(t), 
        \label{eq:Heat1.heat-cLmu-Fourier}\\ 
       \hat{\cL}_{0}= -\frac{1}{2} \sum_{j=1}^{n}(\frac{\partial}{\partial x_{j}}-ix_{n+j}\xi_{0})^{2}- 
       \frac{1}{2} \sum_{j=1}^{n}(\frac{\partial}{\partial x_{n+j}}+ix_{j}\xi_{0})^{2}.
    \end{gather}
    
    Notice that $\hat{\cL}_{0}=\frac{1}{2} H_{A(\xi_{0})}$, where 
    $H_{A(\xi_{0})}=-\sum_{j=1}^{2n}(\partial_{j}-\sum_{k=1}^{2n}iA(\xi_{0})_{jk}x_{k})^{2}$ is the harmonic oscillator associated to the 
    real antisymmetric matrix, 
    \begin{equation}
        A(\xi_{0})=\xi_{0}J, \qquad 
        J=\left(  
        \begin{array}{cc}
            0 & -I_{n}  \\
            I_{n} & 0
        \end{array} \right).
    \end{equation}

    Therefore, the fundamental solution of $\hat{\cL}_{0}+\mu \xi_{0}+\partial_{t}$ is given by a version of the Mehler formula 
    (see, e.g.,~\cite{GJ:QPFIPV}, \cite[p.~225]{Po:CMP1}), that is, $ \hat{k}_{0}(\xi_{0},x',t)$ is equal to
    \begin{multline}
        \chi(t) (2\pi t)^{-n} \det{}^{\frac{1}{2}}( \frac{it A(\xi_{0})}{\sinh (itA(\xi_{0}))}) 
        \exp [ -\frac{1}{2t} \acou{ \frac{itA(\xi_{0})}{\tanh (itA(\xi_{0})}x'}{x'}],\\
        =\chi(t) (2\pi t)^{-n} (\frac{t\xi_{0}}{\sinh (t\xi_{0}}))^{n} \exp [-\frac{1}{2t} \frac{t\xi_{0}}{\tanh (t\xi_{0})}|x'|^{2}]. 
    \end{multline}

A solution of~(\ref{eq:Heat1.heat-cLmu-Fourier}) is now given by $\hat{k}_{\mu}(\xi_{0},x',t)=e^{-\mu\xi_{0}t}\hat{k}_{0}(\xi_{0},x',t)$. 
Moreover, as we have 
\begin{equation}
    |h_{0}(\xi_{0},x',t)|\leq \pi^{-n}|\xi_{0}|^{n}e^{-tn|\xi_{0}|}, 
\end{equation}
we see that for  $|\Re \mu| <n$ the function $e^{-t\mu 
    \xi_{0}}\hat{k}_{0}$ is integrable with respect to $\xi_{0}$. Since $k_{\mu}(x_{0},x',t)$ is  the inverse Fourier transform with respect to $\xi_{0}$ of 
    $\hat{k}_{\mu}(\xi_{0},x',t)$ it follows that $ k_{\mu}(x_{0},x',t) $ is equal to
     \begin{equation}
       \chi(t) (2\pi t)^{-n} \int_{-\infty}^{\infty}e^{ix_{0}\xi_{0}-\mu t\xi_{0}}(\frac{t\xi_{0}}{\sinh t\xi_{0}})^{n}
        \exp [-\frac{1}{t} \frac{t\xi_{0}}{\tanh t\xi_{0}}|x'|^{2}]d\xi_{0}.
    \end{equation}
The lemma is thus proved. 
\end{proof}

Next, for $|\Re \mu|<n$ we let 
\begin{equation}
    \nu(\mu)=\frac{1}{(n+1)!}k_{\mu}(0,0,1)= \frac{(2\pi)^{-(n+1)}}{(n+1)!} \int_{-\infty}^{\infty}e^{-\mu\xi_{0}}(\frac{\xi_{0}}{\sinh \xi_{0}})^{n}d\xi_{0}. 
\end{equation}

\begin{lemma}\label{lem:Heat1.alpha0-volume}
  Let $\Delta:C^{\infty}(M,\cE)\rightarrow C^{\infty}(M,\cE)$ be a selfadjoint sublaplacian  which is bounded from below and assume there exists $\mu 
  \in (-n,n)$ such that near any point of $M$ there is an admissible orthonormal frame $Z_{1},\ldots,Z_{n}$ of $T_{1,0}$ with respect to which 
  $\Delta$ takes the form, 
    \begin{equation}
        \Delta= - \sum_{j=1}^{n}(\overline{Z}_{j}Z_{j}+Z_{j}\overline{Z}_{j}) - i\mu X_{0}+  \op{O}_{H}(1).
         \label{eq:Heat1.principal-term-Deltamu}
    \end{equation}
Then as $\lambda \rightarrow \infty$ we have
     \begin{equation}
        N(\Delta; \lambda) \sim 2^{n} \nu(\mu) \rk \cE(\op{vol}_{\theta}M)\lambda^{n+1}. 
           \label{eq:Heat1.asymptotics-sublaplacian-mu}
     \end{equation}
\end{lemma}
\begin{proof}
Let $Z_{1},\ldots,Z_{n}$ be a local admissible orthonormal frame of $T_{1,0}$. Then from~(\ref{eq:Heat1.Levi-form}) and~(\ref{eq:Heat1.Levi-form.Levi-metric})  
we obtain 
\begin{equation}
    [Z_{j},\bar Z_{k}]= -L_{\theta}(Z_{j},Z_{k})T=-i\epsilon_{j}\delta_{jk}T \quad \bmod T_{1,0}\oplus T_{0,1}. 
     \label{eq:Heat1.commutators-admissible-frame}
\end{equation}
In addition, we let $X_{1},\ldots,X_{2n}$ be the vector fields in $H$ such that
 \begin{equation}
     Z_{j}=\left\{ 
      \label{eq:Spectral.convention-Xj-orientation}
 \begin{array}{ll}
     \frac{1}{2}(X_{j}-iX_{n+j})& \text{for $j=1,\ldots,n-\kappa$},\\
     \frac{1}{2}(X_{n+j}-iX_{j}) & \text{for $j=n-\kappa+1,\ldots,n$}.
 \end{array}\right.
 \end{equation}
Then $X_{1},\ldots,X_{2n}$ is a local frame of $H=\Re (T_{1,0}\oplus T_{0,1})$ and from~(\ref{eq:Heat1.commutators-admissible-frame}) we get 
\begin{gather}
    [X_{j},X_{n+k}]=-2\delta_{jk}X_{0} \ \bmod H, 
    \label{eq:Heat1.almost-Heisenberg-relations1}
    \\ [X_{0},X_{j}]=[X_{j},X_{k}]=[X_{n+j},X_{n+k}]=0\ \bmod H.
     \label{eq:Heat1.almost-Heisenberg-relations2}
\end{gather}

Moreover, in terms of the vector fields $X_{1},\ldots,X_{2n}$ the formula~(\ref{eq:Heat1.principal-term-Deltamu}) becomes
    \begin{equation}
        \Delta= - \frac{1}{2}(X_{1}^{2}+\ldots+X_{2n}^{2}) + i\mu X_{0}+  \op{O}_{H}(1).
         \label{eq:Heat1.principal-term-Deltamu2}
    \end{equation}
Combining this with~(\ref{eq:Spectral.Levi-form-Levi-metric2}) shows that the condition~(\ref{eq:Sublaplacian.condition}) for $\Delta$, 
is given in terms of the eigenvalues of $|L|=1$ and becomes 
\begin{equation}
    \mu \not \in  \{\pm(n+2k);\ k=0,1,2,\ldots\}.
\end{equation}
Since by assumption we have $\mu\in (-n,n)$, the condition~(\ref{eq:Sublaplacian.condition}) is fulfilled at every point, 
so by Proposition~\ref{thm:Powers1.heat-sublaplacians} the principal symbol of $\Delta+\partial_{t}$ is an invertible 
Volterra-Heisenberg symbol. This allows us to apply Proposition~\ref{thm:Heat1.spectral-asymptotics} to deduce that as $\lambda \rightarrow \infty$ we have
\begin{equation}
    N(\Delta; \lambda) \sim \nu_{0}(\Delta) \lambda^{n+1},
     \label{eq:Spectral.Weyl-Delta-lemma}
\end{equation}
where $\nu_{0}(P)$ is given by~(\ref{eq:Heat1.counting-function-asymptotics}). 
    
Let us now work in Heisenberg coordinates centered at a point $a\in M$ related to a local $H$-frame $X_{0},X_{1},\ldots,X_{2n}$ as above. 
Because of~(\ref{eq:Heat1.almost-Heisenberg-relations1}) and (\ref{eq:Heat1.almost-Heisenberg-relations2}) the model vector fields 
$X_{0}^{a},X_{1}^{a}\ldots,X_{2n}^{a}$ coincide with the vector fields~(\ref{eq:Heat1.Heisenberg-vector-fields1})--(\ref{eq:Heat1.Heisenberg-vector-fields2}),  
% the model vector fields of $X_{0},\ldots,X_{2n}$ at $a$ are
% \begin{gather}
%     X_{0}^{a}=\frac{\partial}{\partial x_{0}}, \qquad X_{j}^{a}=\frac{\partial}{\partial x_{j}}+x_{n+j}\frac{\partial}{\partial 
%     x_{0}},\\ X_{n+j}^{a}=\frac{\partial}{\partial x_{n+j}}-x_{j}\frac{\partial}{\partial 
%     x_{0}}, \quad j=1,\ldots, n.
% \end{gather}
% These are the vector fields 
so $G_{a}M$ agrees with the Heisenberg group $\bH^{2n+1}$ in the Heisenberg coordinates. In addition, using~(\ref{eq:Heat1.principal-term-Deltamu2}) 
we see that the model operator of $\Delta$ is 
\begin{equation}
    \Delta^{a}= - \frac{1}{2}((X_{1}^{a})^{2}+\ldots+(X_{2n}^{a})^{2}) - i\mu X_{0}^{a}=\cL_{\mu}.
\end{equation}
As $\mu\in (-n,n)$ it follows that $\Delta$ satisfies the Rockland condition at every point. 

On the other hand, since we are in Heisenberg coordinates we have $X_{j}(0)=\frac{\partial}{\partial x_{j}}$. In particular, we have 
$\theta(0)=dx_{0}$. Moreover, as $d\theta(X,Y)=-\theta([X,Y])$ we deduce 
from~(\ref{eq:Heat1.almost-Heisenberg-relations1})--(\ref{eq:Heat1.almost-Heisenberg-relations2}) that, for $j=0,1,\ldots,n$ and $k=1,\ldots,n$, we 
have 
\begin{equation}
    d\theta(X_{j},X_{n+k})=2\delta_{jk}, \qquad d\theta(X_{j},X_{k})=0.
\end{equation}
 It then follows that $d\theta(0)=2\sum_{j=1}^{n} dx_{j}\wedge dx_{n+j}$. Therefore, the volume form is given 
by
\begin{equation}
  \frac{(-1)^{\kappa}}{n!}\theta \wedge d\theta^{n}(0)  
  =(-1)^{\kappa}2^{n}dx_{0}\wedge dx_{1}\wedge dx_{n+1}\wedge \ldots \wedge dx_{n}\wedge dx_{2n}.
     \label{eq:Spectral.volume-form-Heisenberg-coordinates}
\end{equation}
Note also that because of~(\ref{eq:Spectral.convention-Xj-orientation}) 
the orientation of $M$ is that given by $i^{n}T\wedge Z_{1}\wedge\bar{Z}_{1}\wedge \ldots \wedge 
Z_{n}\wedge\bar{Z}_{n}=(-1)^{\kappa}T\wedge X_{1}\wedge X_{n+1}\wedge \ldots \wedge X_{n}\wedge X_{2n}$. 
Therefore, at $x=0$ the volume form of $M$ 
is in the same orientation class as $(-1)^{\kappa}dx_{0}\wedge dx_{1}\wedge dx_{n+1}\wedge \ldots \wedge dx_{n}\wedge dx_{2n}$. Thus, 
% 
% In terms of density, 
% together with~(\ref{eq:Spectral.volume-form-Heisenberg-coordinates}), this means that we have 
\begin{equation}
    \frac{(-1)^{\kappa}}{n!}\theta \wedge d\theta^{n}(0)  =2^{n} dx_{|_{x=0}}. 
\end{equation}

All this shows that in the above Heisenberg coordinates the volume form and the model operator of $\Delta$ have normal forms in the sense of 
Proposition~\ref{prop:Spectral.normal-form-nu0P} with $\rho_{0}=2^{n}$ and $P_{0}=\cL_{\mu}\otimes 1_{\C^{r}}$. Therefore, 
we may apply~(\ref{eq:Spectral.normal-form-nu0P}) to get
\begin{equation}
    \nu_{0}(\Delta)=2^{n}\nu_{0}(\cL_{\mu}\otimes 1_{\C^{r}})\int_{M}  \frac{(-1)^{\kappa}}{n!}\theta \wedge d\theta^{n}=2^{n} \nu_{0}(\cL_{\mu})\rk\cE 
    \op{vol}_{\theta}M,
\end{equation}
where  $\nu_{0}(\cL_{\mu})=\frac{1}{(n+1)!}k_{\mu}(0,0,1)=\nu(\mu)$. The lemma is thus proved.
\end{proof}

We are now ready to relate the asymptotics~(\ref{eq:Heat1.counting-function-asymptotics}) for the Kohn Laplacian to the pseudohermitian volume of $M$. 
\begin{theorem}
  Assume that $M$ is endowed with a Levi metric and consider the Kohn Laplacian $\Boxbpq:C^{\infty}(M,\Lambda^{p,q})\rightarrow C^{\infty}(M,\Lambda^{p,q})$ 
  acting on $(p,q)$-forms with  $q\neq \kappa$ and $q\neq n-\kappa$. Then  
    as $\lambda \rightarrow \infty$ we have 
\begin{equation}
        N(\Boxbpq ;\lambda) \sim \alpha_{n\kappa pq}(\op{vol}_{\theta}M)\lambda^{n+1},
%     \label{eq:}
\end{equation}
where $ \alpha_{n\kappa pq}$ is equal to  
\begin{equation}
   \sum_{\max(0,q-\kappa)\leq  k\leq \min(q,n-\kappa)} \frac{1}{2} \binom{n}{p} \binom{n-\kappa}{k}\binom{\kappa}{q-k} 
        \nu(n-2(\kappa-q+2k)).
         \label{eq:Heat1.alphapq}  
%     \label{eq:}
\end{equation}
In particular $\alpha_{n\kappa pq}$ is a universal constant depending only on $n$, $\kappa$, $p$ and $q$.
\end{theorem}
\begin{proof}
As $M$ is $\kappa$-strictly pseudoconvex the $Y(q)$-condition reduces to $q\neq \kappa$ and  $q \neq n-\kappa$, so in this case $\Boxbpq$ has an 
invertible principal symbol. Since $\Boxbpq$ is a positive it follows from Proposition~\ref{thm:Heat1.spectral-asymptotics} 
that as $\lambda \rightarrow \infty$ we have 
\begin{equation}
    N(\Boxbpq;\lambda)  \sim \nu_{0}(\Boxbpq) \lambda^{n+1}.
     \label{eq:Heat1.Weyl-asymptotic-Boxb1}
\end{equation}
It then remains to show that $\nu_{0}(\Boxbpq)=\alpha_{n\kappa pq}\op{vol}_{\theta}M$ with $\alpha_{n \kappa pq}$ given by~(\ref{eq:Heat1.alphapq}). 

Let $\{T,Z_{j},\bar Z_{j}\}$ be a local admissible orthonormal frame of $T_{\C}M$ and let $\{\theta,\theta^{j}, \theta^{\bar j}\}$ be the dual 
coframe of $T^{*}_{\C}M$. For ordered subsets $J=\{j_{1},\ldots,j_{p}\}$ and $K=\{k_{1},\ldots,k_{q}\}$ of $\{1,\ldots,n\}$ with $j_{1}<\ldots<j_{p}$ 
and $k_{1}<\ldots<k_{q}$ we let 
\begin{equation}
    \theta^{J,\bar K}=\theta^{j_{1}}\wedge \ldots \wedge \theta^{j_{p}}\wedge \theta^{\bar k_{1}}\wedge \ldots \wedge \theta^{\bar k_{q}}.
\end{equation}
Then $\{\theta^{J,\bar K}\}$ form an orthonormal frame of $\Lambda^{*,*}$ and, as shown in~\cite[Sect.~20]{BG:CHM}, with respect to this frame 
$\Box_{b}$ takes the form, 
\begin{gather}
    \Box_{b}=\op{diag}\{ \Box_{J\bar K}\} +   \op{O}_{H}(1),
    \label{eq:Heat1.Boxb-local-form1}\\
    \Box_{J\bar K}=-\frac{1}{2}\sum_{1\leq j \leq n} (Z_{j}\bar Z_{j}+\bar Z_{j} Z_{j}) +\frac{1}{2}\sum_{j \in K}[Z_{j},\bar Z_{j}] - 
    \frac{1}{2}\sum_{j \not \in K}[Z_{j},\bar Z_{j}].
      \label{eq:Heat1.Boxb-local-form2}
\end{gather}

Moreover, using~(\ref{eq:Heat1.commutators-admissible-frame}) 
we see that the leading part of $\Box_{J \bar K}$ is equal to
\begin{equation}
    -\frac{1}{2}\sum_{1\leq j \leq n} (Z_{j}\bar Z_{j}+\bar Z_{j} Z_{j}) -\frac{i}{2}\mu_{K}X_{0}, \qquad \mu_{K}=\sum_{j \in K} \epsilon_{j} 
    -\sum_{j\not\in K}\epsilon_{j}.
    \label{eq:Heat1.principal-term-BoxJK}
\end{equation}
Notice that since $\epsilon_{j}=1$ for $j=1,\ldots,n-\kappa$ and $\epsilon_{j}=-1$ for $j=n-\kappa+1,\ldots,n$ we have $\mu_{K}\in [-n,n]$ and 
$\mu_{K}=\pm n$ if, and only if, we have $K=\{1,\ldots,n-\kappa\}$ or $K=\{n-\kappa+1,\ldots,n\}$. Thus if $|K|=q$ with $q \not\in\{\kappa,n-\kappa\}$ 
then $\Box_{J \bar K}$ is two times a sublaplacian of the form~(\ref{eq:Heat1.principal-term-Deltamu}) with $\mu=\mu_{K}$ in $(-n,n)$. 

On the other hand, complex conjugation is an isometry, so from the orthogonal splitting~(\ref{eq:Heat1.positive-negative-splitting}) 
we get the orthogonal splitting $T_{0,1}=T_{0,1}^{+}\oplus T_{0,1}^{-}$ with $T_{0,1}^{\pm}=\overline{T_{1,0}^{\pm}}$. 
By duality these splittings give rise to the orthogonal decompositions,
\begin{equation}
        \Lambda^{1,0}=\Lambda^{1,0}_{+}\oplus \Lambda^{1,0}_{-}, \qquad \Lambda^{0,1}=\Lambda^{0,1}_{+}\oplus \Lambda^{0,1}_{-}. 
%     \label{eq:}
\end{equation}
Therefore, letting $\Lambda^{p;q,k}=\Lambda^{p,0}\wedge (\Lambda^{0,1}_{+})^{k}\wedge (\Lambda^{0,1}_{-})^{q-k}$ we have
\begin{equation}
         \Lambda^{p,q}=  \bigoplus_{\max(0,q-\kappa)\leq  k\leq \min(q,n-\kappa)} \Lambda^{p;q,k}. 
%          \qquad \Lambda^{p;q,k}=\Lambda^{p,0}\wedge (\Lambda^{0,1}_{+})^{k}\wedge (\Lambda^{0,1}_{-})^{q-k}.
          \label{eq:Heat1.splitting-pqk}
%     \label{eq:}
\end{equation}

Let  $\Pi_{p;qk}$ be the orthogonal projection onto $\Lambda^{p;q,k}$ and set $\Box_{pqk}=\Pi_{pqk}\Box_{b} \Pi_{pqk}$.  
As~(\ref{eq:Heat1.Boxb-local-form1}) and~(\ref{eq:Heat1.Boxb-local-form2}) show that $\Boxbpq$ is scalar up to first order terms we have
\begin{equation}
    \Boxbpq= \sum_{\max(0,q-\kappa)\leq  k\leq \min(q,n-\kappa)} \Box_{p;qk} + \op{O}_{H}(1). 
\end{equation}
In particular,  
 as $\nu_{0}(\Boxbpq)$ depends only on the principal 
symbol of $\Boxbpq$ we get
\begin{equation}
    \nu_{0}(\Boxbpq)=\sum_{\max(0,q-\kappa)\leq  k\leq \min(q,n-\kappa)}  \nu_{0}(\Box_{p;qk}).
     \label{eq:Heat1-splitting-alpha0}
\end{equation}
We are thus reduced to express each coefficient $\nu_{0}(\Box_{p;qk})$ in terms of $\op{vol}_{\theta}M$. 

Next, if $\theta^{J\bar K}$ is a section of $\Lambda^{p,q,k}$ then we have
\begin{gather}
    \# K\cap \{1,\ldots,n-\kappa\}=k, \quad  \#K\cap \{n-\kappa+1,\ldots,n\}=q-k,\\
   \#K^{c}\cap \{1,\ldots,n-\kappa\}=n-\kappa-k,\\  \#K^{c}\cap \{n-\kappa+1,\ldots,n\}=\kappa-q+k,
\end{gather}
from which we get $ \mu_{K}= k-(q-k)-[(n-\kappa-k)-(\kappa-q+k)]=n+2q-2\kappa-4k.$
Combining this with~(\ref{eq:Heat1.Boxb-local-form1})--(\ref{eq:Heat1.principal-term-BoxJK}) then gives
\begin{equation}
    \Box_{p;qk}=-\frac{1}{2}\sum_{1\leq j \leq n} (Z_{j}\bar Z_{j}+\bar Z_{j} Z_{j}) -\frac{i}{2}(n+2q-2\kappa -4k)X_{0}+   \op{O}_{H}(1).
     \label{eq:Heat1.Boxpqk-local}
\end{equation}

Now, we can  apply Lemma~\ref{lem:Heat1.alpha0-volume} to $2 \Box_{p;qk}$ to get that $ \nu_{0}( \Box_{p;qk})$ is equal to
\begin{multline}
2^{-(n+1)}\nu_{0}( 2\Box_{p;qk})=
 2^{-(n+1)}2^{n}(\rk \Lambda^{p;q,k}) \nu(n+2q-2\kappa-4k).\op{vol}_{\theta}M \\
    = \frac{1}{2}\binom{n}{p} \binom{n-\kappa}{k}\binom{\kappa}{q-k}\nu(n+2q-2\kappa-4k).\op{vol}_{\theta}M.
\end{multline}
Combining this with~(\ref{eq:Heat1.Weyl-asymptotic-Boxb1}) and~(\ref{eq:Heat1-splitting-alpha0}) then shows that as $\lambda \rightarrow \infty$ we have 
\begin{equation}
    N(\Boxbpq;\lambda) \sim \alpha_{n\kappa pq}(\op{vol}_{\theta}M)\lambda^{n+1},  
\end{equation}    
with $\alpha_{n\kappa pq}$ given by~(\ref{eq:Heat1.alphapq}). The proof is thus achieved.
\end{proof}

Next, we turn to the horizontal sublaplacian on $(p,q)$-forms.

\begin{theorem}\label{thm:Heat1.Deltab}
  Assume that $M$ is endowed with a Levi metric and consider the  horizontal sublaplacian 
  $\Delta_{b}:C^{\infty}(M,\Lambda^{p,q})\rightarrow C^{\infty}(M,\Lambda^{p,q})$ acting on $(p,q)$-forms with  
  $(p,q)\neq (\kappa,n-\kappa)$ and $(p,q)\neq (n-\kappa,\kappa)$. Then as $\lambda  \rightarrow \infty$ we have 
\begin{equation}
             N(\Delta_{b;p,q};\lambda) \sim \beta_{n\kappa pq}(\op{vol}_{\theta}M) \lambda^{n+1},
%     \label{eq:}
\end{equation}
 where $\beta_{n\kappa pq}$ is equal to
\begin{equation}
        \!  \!  \!  \!   \sum_{\substack{\max(0,q-\kappa)\leq  k\leq \min(q,n-\kappa)\\ \max(0,p-\kappa)\leq l\leq \min(p,n-\kappa)}}  \!  \!  \!  \!
        2^{n}\binom{n-\kappa}{l}\binom{\kappa}{p-l} \binom{n-\kappa}{k}\binom{\kappa}{q-k} 
        \nu(2(q-p)+4(l-k)).
         \label{eq:Heat1.betapq}
%     \label{eq:}
\end{equation}
In particular $\beta_{n\kappa pq}$  is a universal constant depending only on $n$, $\kappa$, $p$ and $q$.
\end{theorem}
\begin{proof}
Thanks to~(\ref{eq:Operators.Tanaka-Kohn}) we know that we have 
 \begin{equation}
     \Delta_{b}=\Boxbpq+ \overline{\Box}_{b;p,q},
 \end{equation}
 where $\overline{\Box}_{b}$ denotes the conjugate operator of $\Boxb$, that is, $\overline{\Box}_{b}\eta =\overline{\Boxb\overline{\eta }}$, 
 (or, equivalently, the Laplacian of the $\partial_{b}$-complex). 
 
As in~(\ref{eq:Heat1.splitting-pqk}) we have orthogonal splittings, 
\begin{gather}
    \Lambda^{p,q}=    \bigoplus_{\max(0,p-\kappa)\leq  l\leq \min(p,n-\kappa)} \Lambda^{p,l;q}  
    = \bigoplus_{\substack{\max(0,q-\kappa)\leq  k\leq \min(q,n-\kappa)\\ \max(0,p-\kappa)\leq  l\leq \min(p,n-\kappa)}} \Lambda^{p,l;q,k},\\
    \Lambda^{p,l;q}= (\Lambda^{1,0}_{+})^{l}\wedge (\Lambda^{1,0}_{-})^{p-l} \wedge \Lambda^{0,q}, \qquad 
     \Lambda^{p,l;q,k}= \Lambda^{p,l;0}\wedge \Lambda^{p,l;0}.
\end{gather}
Since $\Boxbpq$ is a scalar operator modulo lower order terms, the same is true for $\overline{\Box}_{b}$ and $\Delta_{b}$. Therefore,
on $\Lambda^{p,q}$ we can write
\begin{equation}
        \overline{\Box}_{b}= \sum_{\max(0,p-\kappa)\leq  l\leq \min(p,n-\kappa)} \overline{\Box}_{p,l;q} +   \op{O}_{H}(1), \qquad  \overline{\Box}_{p,l;q}= 
        \Pi_{p,l;q} \overline{\Box}_{b}\Pi_{p,l;q}, 
%     \label{eq:}
\end{equation}
where $\Pi_{p,l;q}$ denote the orthogonal projection onto $ \Lambda^{p,l;q}$. 
Therefore, if we let  $ \Pi_{p,l;q,k}$ denote the orthogonal projection onto $\Lambda^{p,l;q,k}$ and set  $\Delta_{p,l;q,k}=  
\Pi_{p,l;q,k}\Delta_{b}\Pi_{p,l;q,k}$, then we have

\begin{equation}
        \Delta_{b}= \sum_{\substack{\max(0,q-\kappa)\leq  k\leq \min(q,n-\kappa)\\ \max(0,p-\kappa)\leq l\leq \min(p,n-\kappa)}} \Delta_{p,l;q,k} 
        +   \op{O}_{H}(1).
%     \label{eq:}
\end{equation}
In particular, as in~(\ref{eq:Heat1-splitting-alpha02}) we have 
\begin{equation}
    \nu_{0}( \Delta_{b;p,q})=  \sum_{\substack{\max(0,q-\kappa)\leq  k\leq \min(q,n-\kappa)\\ \max(0,p-\kappa)\leq l\leq \min(p,n-\kappa)}} 
    \nu_{0}(\Delta_{p,l;q,k}).
     \label{eq:Heat1-splitting-alpha02}
\end{equation}

Next, let $\{X_{0},Z_{j},\bar Z_{j}\}$ be a local admissible orthonormal frame for $T_{\C}M$. Since 
$\overline{\Box}_{p,l;q}=\overline{\Box_{q;p,l}}$, using~(\ref{eq:Heat1.Boxpqk-local}) we see that on $\Lambda^{p,l;q}$ we have 
\begin{equation}
    \overline{\Box}_{p,l;q}= - \frac{1}{2}\sum_{1\leq j \leq n} (Z_{j}\bar Z_{j}+\bar Z_{j} Z_{j}) +\frac{i}{2}(n+2p-2\kappa -4l)X_{0}+  
     \op{O}_{H}(1).
\end{equation}
Therefore, on $\Lambda^{p,l;q,k}$ we can write
\begin{multline}
      \Delta_{p,l;q,k}= \Box_{p;q,k}+\overline{\Box}_{p,l;q}\\ = 
      -\sum_{1\leq j \leq n} (Z_{j}\bar Z_{j}+\bar Z_{j} Z_{j}) -i(2(q-p)+4(l-k))X_{0}+   \op{O}_{H}(1).
     \label{eq:Heat1.Deltaplqk}
\end{multline}

Thanks to~(\ref{eq:Heat1.Deltaplqk}) we can apply Lemma~\ref{lem:Heat1.alpha0-volume} to get
\begin{multline}
    \nu_{0}(\Delta_{p,l;q,k})= 2^{n}\rk \Lambda^{p,l;q,k}\nu(2(q-p)+4(l-k)) \op{vol}_{\theta}M\\ 
    = 2^{n}\binom{n-\kappa}{l}\binom{\kappa}{p-l} \binom{n-\kappa}{k}\binom{\kappa}{q-k}  \nu(2(q-p)+4(l-k)) \op{vol}_{\theta}M.
\end{multline}
Combining this with~(\ref{eq:Heat1-splitting-alpha02}) then shows that as $\lambda\rightarrow \infty$ we have 
   \begin{equation}
         N(\Delta_{b};\lambda) \sim \beta_{n\kappa pq}(\op{vol}_{\theta}M) \lambda^{n+1}, 
    \end{equation}
with $\beta_{n\kappa pq}$ given by~(\ref{eq:Heat1.betapq}). The theorem is thus proved.
\end{proof}

 Finally, suppose that $M$ is strictly pseudoconvex, i.e., we have $\kappa=0$, and for $k=1,\ldots,n+1,n+2, n+4,\ldots$ let 
 $\boxdot_{\theta}^{(k)}:C^{\infty}(M)\rightarrow C^{\infty}(M)$ be the Gover-Graham operator of order $k$. Then we have:
 
  \begin{theorem}
 Assume $k\neq n+1$. Then there exists a  universal constant $\nu_{n}^{(k)}>0$ depending only on $n$ and $k$  such that as $\lambda 
     \rightarrow \infty$ we have 
     \begin{equation}
          N(\boxdot_{\theta}^{(k)};\lambda) \sim  \nu_{n}^{(k)}(\op{vol}_{\theta}M) \lambda^{\frac{n+1}{k}}. 
     \end{equation}
  \end{theorem}
\begin{proof}
Let $Z_{1},\ldots,Z_{n}$ be an orthonormal frame of $T_{1,0}$ over an open $U\subset M$. Since $\kappa=0$ this is an admissible orthonormal frame. Therefore, 
as shown in the proof of Lemma~\ref{lem:Heat1.alpha0-volume}, if for $j=1,\ldots,n$ we let 
$Z_{j}=\frac{1}{2}(X_{j}-iX_{n+j})$ then $X_{0},X_{1},\ldots,X_{2n}$ is an $H$-frame such 
that in associated Heisenberg coordinates centered at a point $a\in U$ we have $G_{0}M=\bH^{2n+1}$ and $\theta\wedge d\theta^{n}(0)=2^{n}dx$ and the 
model operator at $x=0$ of $\Delta_{\mu}=-\sum (\overline{Z}_{j}Z_{j}+Z_{j}\overline{Z}_{j})-i\mu T$, $\mu \in \C$, is the Folland-Stein 
sublaplacian $\cL_{\mu}$ 
in~(\ref{eq:Weyl.Folland-Stein-sublaplacian}).

By Proposition~\ref{prop:Examples.Gover-Graham} the operator $\boxdot_{\theta}^{(k)}$ is of the form, 
       \begin{equation}
      (\Delta_{b;0}+i(k-1)X_{0})(\Delta_{b;0}+i(k-3)X_{0}) \cdots (\Delta_{b;0}-i(k-1)X_{0}) +\op{O}_{H}(2k-1). 
%           \label{eq:Examples.Gover-Graham-principal-part}
     \end{equation}
Moreover, by~(\ref{eq:Examples.Deltabpq}) we have
     $\Delta_{b;0}=-\sum (\overline{Z}_{j}Z_{j}+Z_{j}\overline{Z}_{j}) +\op{O}_{H}(1)$, so in Heisenberg coordinates as above 
     the model operator  at $x=0$ of $\Delta_{b;0}-i\mu X_{0}$, $\mu \in \C$,  is just $\cL_{\mu}$. Incidentally, the model operator of 
     $\boxdot_{\theta}^{(k)}$ is  
      $\cL_{-k+1}\cL_{-k+3}\cdots \cL_{k-1}$. 
     
     All this shows that $\boxdot_{\theta}^{(k)}$ admits a normal form in the sense of Proposition~\ref{lem:Spectral.a0(P)} 
     and that we can make use of this proposition as in the proof of Lemma~\ref{lem:Heat1.alpha0-volume} to 
    deduce that as $\lambda \rightarrow \infty$ we have 
     \begin{equation}
          N(\boxdot_{\theta}^{(k)};\lambda) \sim  \nu_{n}^{(k)}(\op{vol}_{\theta}M) \lambda^{\frac{n+1}{k}}, 
          \qquad \nu_{n}^{(k)}=2^{n}\Gamma(1+\frac{2n+2}{k})^{-1}K^{(k)}(0,1),
         %          \label{eq:¥}
     \end{equation}
     where $K^{(k)}(x,t)$ denotes the fundamental solution of $\cL_{-k+1}\cL_{-k+3}\cdots \cL_{k-1}+\partial_{t}$. In particular, we see that 
     $\nu_{n}^{(k)}$ is universal constant depending only on $k$ and $n$. Furthermore, it follows from 
     Proposition~\ref{thm:Heat1.spectral-asymptotics} that $\nu_{n}^{(k)}$ is~$>0$. 
 \end{proof}
 
\section{Weyl asymptotics and contact geometry}
\label{sec:Spectral-contact}
In this section we express in more geometric terms the Weyl asymptotics for the horizontal sublaplacian and the contact Laplacian on a compact 
orientable contact manifold $(M^{2n+1},H)$ equipped with a contact form $\theta$. 

Let $X_{0}$ be the Reeb vector field associated to $\theta$, so that $\imath_{X_{0}}d\theta=0$ and $\imath_{X_{0}}\theta=1$.  
Since $M$ is orientable $H$ admits a calibrated almost complex structure $J\in C^{\infty}(M,\End H)$, $J^{2}=-1$, so that for any nonzero section $X$ of $H$ 
we have $d\theta(X,JX)=-d\theta(JX,X)>0$. 
We then endow $M$ with the orientation  defined by $\theta$ and the almost complex structure $J$, so that  $\theta\wedge d\theta^{n}>0$, and we equip 
it with  the Riemannian metric, 
\begin{equation}
    g_{\theta, J}=d\theta(.,J.)+\theta^{2}.
\end{equation}

A local orthonormal frame $X_{1},\ldots, X_{2n}$ of $H$ will be called admissible when we have $X_{n+j}=JX_{j}$ for $j=1,\ldots,n$. If 
$\theta^{1},\ldots,\theta^{2n}$ denotes the dual frame then we have $d\theta=\sum_{j=1}^{n}\theta^{j}\wedge \theta^{n+j}$, so the  
volume form of $g_{\theta}$ is equal to 
\begin{equation}
    \theta^{1}\wedge \theta^{n+1}\wedge \ldots\wedge \theta^{n}\wedge \theta^{2n}\wedge \theta=\frac{1}{n!}d\theta^{n}\wedge \theta.
\end{equation}
In particular, the volume form is independent of the choice of the almost complex structure and depends only on the contact form.

\begin{definition}
    The contact volume of $(M^{2n+1},\theta)$ is given by 
    \begin{equation}
        \op{vol}_{\theta}M=\frac{1}{n!}\int_{M}d\theta^{n}\wedge \theta.
    \end{equation}
\end{definition}

\begin{lemma}\label{lem:Heat1.alpha0-volume-contact}
  Let $\Delta:C^{\infty}(M,\cE)\rightarrow C^{\infty}(M,\cE)$ be a selfadjoint sublaplacian  such that $\Delta$ is bounded from below. We further assume 
  that there exists $\mu 
  \in (-n,n)$ so that near any point of $M$ there is an admissible orthonormal frame $X_{1},\ldots,X_{2n}$ of $H$ with respect to which 
  $\Delta$ takes the form, 
    \begin{equation}
        \Delta= - (X_{1}^{2}+\ldots+X_{2n}^{2}) - i\mu X_{0}+  \op{O}_{H}(1).
         \label{eq:Heat1.principal-term-Deltamu-contact}
    \end{equation}
Then as $\lambda \rightarrow \infty$ we have
     \begin{equation}
        N(\Delta; \lambda) \sim 2^{n}\nu(\mu) \rk \cE (\op{vol}_{\theta}M)\lambda^{n+1}. 
           \label{eq:Heat1.asymptotics-sublaplacian-mu-contact}
     \end{equation}
\end{lemma}
\begin{proof}
  Let $X_{1},\ldots,X_{2n}$ be an admissible local orthonormal frame of $TM$ and for $j=1,\ldots,2n$ let $\tilde{X}_{j}=\sqrt{2}X_{j}$. Then 
  $X_{0},\tilde{X}_{1},\ldots,\tilde{X}_{2n}$ is a local $H$-frame of $TM$ with respect to which $\Delta$ takes the form,
  \begin{equation}
      \Delta=-\frac{1}{2}(\tilde{X}_{1}^{2}+\ldots+\tilde{X}_{2n}^{2})-i\mu(x)T+ \op{O}_{H}(1).
       \label{eq:Spectral.local-form-admissible}
  \end{equation}
  
 Moreover, for $j,k=1,\ldots,2n$ we have 
 \begin{equation}
     \theta([\tilde{X}_{j},J\tilde{X}_{k}])=-2d\theta(X_{j},JX_{k})=-2g_{\theta}(X_{j},X_{k})=-2\delta_{jk}.
 \end{equation}
  Therefore, for $j,k=1,\ldots,n$ we get:
\begin{gather}
    [X_{j},X_{n+k}]=-2\delta_{jk}X_{0} \ \bmod H,\\ 
    [T,X_{j}]=[X_{j},X_{k}]=[X_{n+j},X_{n+k}]=0\ \bmod H.
     \label{eq:Heat1.almost-Heisenberg-relations-contact}
\end{gather}
The equalities~(\ref{eq:Spectral.local-form-admissible})--(\ref{eq:Heat1.almost-Heisenberg-relations-contact}) 
are the same as~(\ref{eq:Heat1.almost-Heisenberg-relations1})--(\ref{eq:Heat1.principal-term-Deltamu2}) in the case $\kappa=0$. Therefore, 
along the same lines as that of the proof  Lemma~\ref{lem:Heat1.alpha0-volume} we get
\begin{equation}
    \nu_{0}(\Delta)=2^{n}\frac{\nu(\mu)}{n!}\int_{M}d\theta^{n}\wedge \theta=2^{n}\nu(\mu)\op{vol}_{\theta}M.
\end{equation}
Combining this with Proposition~\ref{thm:Heat1.spectral-asymptotics} then proves the asymptotics~(\ref{eq:Heat1.asymptotics-sublaplacian-mu-contact}).
\end{proof}

\begin{theorem}
    Let $\Delta_{b;k}:C^{\infty}(M,\Lambda^{k}_{\C}H^{*})\rightarrow C^{\infty}(M,\Lambda^{k+1}_{\C}H^{*})$ be the horizontal sublaplacian on $M$ in 
    degree $k$ with $k \neq n$. Then as $\lambda  \rightarrow \infty$ we have
    \begin{gather}
         N(\Delta_{b;k};\lambda) \sim \gamma_{nk} (\op{vol}_{\theta}M) \lambda^{n+1},\qquad 
          \label{eq:Spectral.Weyl-Deltab-contact}
 \gamma_{nk}=\sum_{p+q=k}2^{n}\binom{n}{p} \binom{n}{q}\nu(p-q).         
    \end{gather}
In particular $\gamma_{nk}$ is universal constant depending on $n$ and $k$ only.
\end{theorem}
\begin{proof}
   As explained in Section~\ref{sec:Operators} the almost complex structure of $H$ gives rise to an orthogonal decomposition 
   $\Lambda^{k}_{\C}H^{*}=\bigoplus_{p+q=k}\Lambda^{p,q}$.   If $X_{1},\ldots,X_{2n}$ is a local admissible orthonormal frame of $H$ then, as shown by 
   Rumin~\cite[Prop.~2]{Ru:FDVC}, on $\Lambda^{p,q}$ the operator $\Delta_{b}$ takes the form,
   \begin{equation}
       \Delta_{b}=-(X_{1}^{2}+\ldots+X_{2n}^{2})+i(p-q)X_{0}+ \op{O}_{H}(1).
        \label{eq:Spectral.contact-local-form-Deltapq}
   \end{equation}
   where the lower order part is not scalar. Therefore, modulo lower order terms, $\Delta_{b}$ preserves the bidegree. We thus may write
   \begin{equation}
       \Delta_{b;k}=\sum_{p+q=k} \Delta_{p,q} + \op{O}_{H}(1), \qquad \Delta_{p,q}=\Pi_{p,q}\Delta_{b}\Pi_{p,q},
   \end{equation}
   where $\Pi_{p,q}$ denotes the orthogonal projection of $\Lambda^{*}_{\C}H^{*}$ onto $\Lambda^{p,q}$. In particular, 
   \begin{equation}
       \nu_{0}(\Delta_{b;k})=\sum_{p+q=k}\nu_{0}(\Delta_{p,q}).
        \label{eq:Spectral.decomposition-nuDeltab-contact}
   \end{equation}
   
   Moreover, since $\Delta_{p,q}$ takes the form~(\ref{eq:Spectral.contact-local-form-Deltapq}) 
   with respect to  any admissible orthonormal frame of $H$, we may apply 
   Lemma~\ref{lem:Heat1.alpha0-volume-contact} to get
   \begin{equation}
       \nu_{0}(\Delta_{p,q})=2^{n}\sum_{p+q=k}\binom{n}{p} \binom{n}{q}\nu(p-q). 
   \end{equation}
   Combining this with~(\ref{eq:Spectral.decomposition-nuDeltab-contact}) then gives the asymptotics~(\ref{eq:Spectral.Weyl-Deltab-contact}). 
\end{proof}

Finally, in the case of the contact Laplacian we can prove:

\begin{theorem}
%     Let $\Delta_{R}:C^{\infty}(M,\Lambda^{*}\oplus \Lambda^{n}_{*})\rightarrow C^{\infty}(M,\Lambda^{*}\oplus \Lambda^{n}_{*})$ 
%     be the contact Laplacian on $M$.\smallskip 
%     
    1)  Let  $\Delta_{R;k}:C^{\infty}(M,\Lambda^{k})\rightarrow C^{\infty}(M,\Lambda^{k})$ be the contact Laplacian in degree $k \neq n$. Then  
    there exists a  universal constant $\nu_{nk}>0$ depending only on $n$ and $k$ 
    such that as $\lambda \rightarrow \infty$ we have
     \begin{equation}
           N(\Delta_{R;k})\sim \nu_{nk} (\op{vol}_{\theta}M)\lambda^{n+1}. 
                \label{eq:Spectral.Weyl-contact-Laplacian1}
      \end{equation}
    
    2) For $j=1,2$ consider the contact Laplacian $\Delta_{R;nj}:C^{\infty}(M,\Lambda^{n}_{j})\rightarrow C^{\infty}(M, \Lambda^{n}_{j})$.  
    Then there exists a  universal constant $\nu_{n,j}>0$ depending only on $n$ and $j$  such that as $\lambda \rightarrow \infty$ we have 
    \begin{equation}
      N(\Delta_{R;nj})\sim \nu_{n,j} (\op{vol}_{\theta}M)\lambda^{\frac{n+1}{2}}. 
         \label{eq:Spectral.Weyl-contact-Laplacian2}
    \end{equation}
\end{theorem}
\begin{proof}
Let $a\in M$ and consider a chart around $a$ together with an admissible orthonormal frame $X_{1},\ldots,X_{2n}$ of $H$. Since 
$T,X_{1},\ldots,X_{2n}$ form a $H$-frame this chart is a Heisenberg chart. Moreover, as shown in the proofs of Lemma~\ref{lem:Heat1.alpha0-volume} 
and Lemma~\ref{lem:Heat1.alpha0-volume-contact} the following hold:\smallskip

(i) We have $X_{j}^{a}=X_{j}^{0}$, where $X_{0}^{0},\ldots,X_{2n}^{0}$ denote the left-invariant vector 
fields~(\ref{eq:Heat1.Heisenberg-vector-fields1}) and~(\ref{eq:Heat1.Heisenberg-vector-fields2}) on $\bH^{2n+1}$. In particular, we have 
$G_{a}M=\bH^{2n+1}$ and $H_{a}=H^{0}_{0}$, where $H^{0}_{0}$ is the left-invariant Heisenberg structure of $\bH^{2n+1}$.\smallskip 

(ii) We have $\theta(0)=dx_{0}=\theta^{0}(0)$ and $d\theta(0)=2\sum_{j=1}^{n}dx_{j}\wedge dx_{n+j}=d\theta^{0}(0)$, where 
$\theta^{0}=dx_{0}+\sum_{j=1}^{n}(x_{j}dx_{n+j}-x_{n+j}dx_{j})$ is the standard left-invariant contact form of $\bH^{2n+1}$.\smallskip 

(iii) The density on $M$ given by the contact volume form $\frac{1}{n!}\theta\wedge d\theta^{n}$ agrees at $x=0$ with the density $2^{n}dx$ on 
$\Rd$.\smallskip 

On the other hand, for $k=0,1,\ldots,2n$ the fiber at $a$ of the bundle $\Lambda^{k}_{*}$ depends only on $H_{a}$ and on the values of $\theta$ and 
$d\theta$ at $a$. Therefore, it follows from the statements (i) and (ii) that  in the Heisenberg coordinates centered at $a$  
the fibers at $x=0$ of the bundles $\Lambda^{*}\oplus \Lambda^{n}_{*}$ of $M$ and $\bH^{2n+1}$ agree. 

Next, let $\Delta_{R}^{0}$ denote the contact Laplacian on $\cH^{2n+1}$. Then we have:
\begin{lemma}
       In the Heisenberg coordinates centered at $a$ the model operators of $\Delta_{R;k}^{a}$ and $\Delta_{R;nj}$ agree with $\Delta_{R;k}^{0}$ and 
       $\Delta_{R;nj}^{0}$ respectively.
\end{lemma}
\begin{proof}[Proof of the lemma]
First, note that in view of the formulas~(\ref{eq:Operators.equalities-DeltaR-Deltab.1})--(\ref{eq:Operators.equalities-DeltaR-Deltab.5}) 
for $\Delta_{R}$ and of Proposition~\ref{prop:PsiHDO.composition2} and Proposition~\ref{prop:PsiHDO.transpose-global}, we only have to show that in the 
Heisenberg coordinates centered at $a$ the   
model operators $d_{R;k}^{a}$ and $D^{a}_{R;n}$ agree with the operators $d_{R;k}^{0}$ and $D^{0}_{R;n}$ on $\bH^{2n+1}$.

Let $\theta^{1},\ldots,\theta^{2n}$ be the coframe of $H^{*}$ dual to $X_{1},\ldots,X_{2n}$. This coframe gives rise to a trivialization of 
$\Lambda^{*}_{\C}H^{*}$ over the chart, in which we have $d_{b}=\sum_{j=1}^{2n}\varepsilon(\theta^{j})X_{j}$. Furthermore, 
since $X_{j}(0)=\frac{d}{dx_{j}}$ we have $\theta_{j}(0)=dx_{j}$, so the model operator of $d_{b}$ is 
$d_{b}^{a}=\sum_{j=1}^{2n}\varepsilon(dx_{j})X_{j}^{0}=d_{b}^{0}$, where $d_{b}^{0}$ is the $d_{b}$-operator on $\bH^{2n+1}$.  
In particular, as for $k=n+1,\ldots,2n$ we have $d_{R;k}=d_{b;k}$ on $\Lambda^{k}$, we get 
$d_{R;k}^{a}=d_{b;k}^{a}=d_{b;k}^{0}=d_{R;k}^{0}$. 

On the other hand, by definition for $k=0,\ldots,n-1$ we have $d_{R;k}^{*}=d_{b;k}^{*}$ on $\Lambda^{k}$. Thus 
$(d_{R;k}^{a})^{*}=(d_{R;k}^{*})^{a}=(d_{b;k}^{*})^{a}=(d_{b;k}^{a})^{*}=(d_{b;k}^{0})^{*}=(d_{R;k}^{0})^{*}$, 
which by taking adjoints gives $d_{R;k}^{a}=d_{R;k}^{0}$. 

Finally, as  $D_{R;n}=\cL_{X_{0}}+d_{b;n-1}\varepsilon(d\theta)^{-1}d_{b;n}$ the model operator 
$D^{a}_{R;n}$  is equal to $\cL_{T^{a}}+d_{b;n-1}^{a}(\varepsilon(d\theta)^{a})^{-1}d_{b;n}^{a}= 
\cL_{T^{0}}+d_{b;n-1}^{0}\varepsilon(d\theta^{0})^{*}d_{b;n}^{0}=D^{0}_{R;n}$. 
The proof  is therefore achieved. 
\end{proof}

Thanks to the statements (i) and (iii) and the claim above we may apply Proposition~\ref{prop:Spectral.normal-form-nu0P}. Letting $K_{0;k}(x,t)$ be 
the fundamental solution of $\Delta_{R;k}^{0}+\partial_{t}$ we then deduce that as $\lambda \rightarrow \infty$ we have
     \begin{equation}
           N(\Delta_{R;k})\sim \nu_{nk} (\op{vol}_{\theta}M)\lambda^{n+1}, \qquad  
           \nu_{nk}=\frac{2^{n}}{(n+1)!}\tr_{\Lambda^{k,0}}K_{0;k}(0,1).
      \end{equation}
In particular, the constant $\nu_{nk}$ depends only on $n$ and $k$, hence is a universal constant, and it follows from
Proposition~\ref{thm:Heat1.spectral-asymptotics}  that $\nu_{nk}$ is a positive number. 

Similarly, let  $K_{0;nj}(x,t)$ be 
the fundamental solution of $\Delta_{R;nj}^{0}+\partial_{t}$. Then as $\lambda \rightarrow \infty$, we have 
    \begin{equation}
      N(\Delta_{R;nj})\sim \nu_{n,j} (\op{vol}_{\theta}M)\lambda^{\frac{n+1}{2}}, 
    \end{equation}
where $\nu_{n,j}=2^{n}\Gamma(1+\frac{n+1}{2})^{-1}\tr_{\Lambda^{n,0}_{j}}K_{0|_{\Lambda^{n,0}_{j}}}(0,1)$ . In particular $\nu_{n,j}$ 
 is a universal constant depending only on $n$ and $j$. The proof is thus complete.
\end{proof}

\appendix 
% \chapter{Appendix}

\numberwithin{equation}{chapter}

\chapter{Proof of Proposition 3.1.18}
\label{chap.Appendix-Invariance}
First, we need the lemma below.

\begin{lemma-appendix}\label{lem:Appendix.Heisenberg.invariance}
   For $\Re m>0$ we have $\cK^{m}(\URd)\subset C^{\infty}(U)\hotimes C^{[\frac{\Re m}{2}]}(\Rd)$.  
\end{lemma-appendix}
\begin{proof}  
 Let $N=[\frac{\Re m}{2}]$ and let $\alpha$ be a multi-order such that $|\alpha|\leq N$. As $\brak\alpha \leq 2|\alpha|\leq 
-(k+d+2)$ the multiplication by $\xi^{\alpha}$ maps $S^{\hat{m}}(\URd)$ to $C^{\infty}(U)\hotimes 
L^{1}(\Rd)$. Composing it with the inverse Fourier transform with respect to $\xi$ then shows that the map 
$p\rightarrow \partial_{y}^{\alpha}\check{p}_{\xiy}$ maps $S^{\hat{m}}(\URd)$ to $C^{\infty}(U)\hotimes 
C^{0}(\Rd)$. It then follows that for any $p \in S^{\hat{m}}(\URd)$ the transform
$\check{p}_{\xiy}(x,y)$  belongs to $C^{\infty}(U)\hotimes C^{N}(\Rd)$.
     
Now, if $K \in \cK^{m}(\URd)$ then by Lemma~\ref{lem:PsiHDO.characterization.Km} 
     there exists a symbol $p\in S^{\hat{m}}(\URd)$, $\hat{m}=-(m+d+2)$, such that $K(x,y)$ is equal to $\check{p}_{\xiy}(x,y)$ modulo 
     a smooth function. Hence $K(x,y)$ belongs to  
     $C^{\infty}(U)\hotimes C^{N}(\Rd)$. The lemma is thus proved.
 \end{proof}
% \begin{proof}  
%      Let $\alpha$ be a multi-order such that $|\alpha|\leq N$. Then we have $\brak\alpha \leq 2|\alpha|\leq 
% -(k+d+2)$, so the multiplication by $\xi^{\alpha}$ maps continuously $S_{||}^{k}(\URd)$ to $C^{\infty}(U)\hotimes 
% L^{1}(\Rd)$. Composing it with the inverse Fourier transform with respect to $\xi$ then shows that the map 
% $p\rightarrow \partial_{y}^{\alpha}\check{p}_{\xiy}$ is continuous from $S_{||}^{k}(\URd)$ to $C^{\infty}(U)\hotimes 
% C^{0}(\Rd)$. Thus the map $p\rightarrow \check{p}_{\xiy}$  is continuous from $S_{||}^{k}(\URd)$ to 
% $C^{\infty}(U)\hotimes C^{N}(\Rd)$.\smallskip  
%      
%      2) Let $K \in \cK^{m}(\URd)$. Then by Lemma~\ref{lem:PsiHDO.characterization.Km} 
%      there exists $p\in S^{\hat{m}}(\URd)$, $\hat{m}=-(m+d+2)$, such that $K(x,y)$ is equal to $\check{p}_{\xiy}(x,y)$ modulo 
%      a smooth function. Since $S^{\hat m}(\URd) \subset S_{||}^{\Re \hat m}(\URd)$ it follows from the first part that $K$ is in 
%      $C^{\infty}(U)\hotimes C^{[\frac{\Re m}{2}]}(\Rd)$. Thus $\cK^{m}(\URd)$ is contained in 
%      $C^{\infty}(U)\hotimes C^{[\frac{\Re m}{2}]}(\Rd)$. 
%  \end{proof}

 We are now ready to prove Proposition~\ref{prop:PsiHDO.invariance}. Let $\tilde{U}$ be an open subset of $\Rd$ together with a hyperplane bundle 
 $\tilde{H}\subset T\tilde{U}$ and a $\tilde{H}$-frame of $T\tilde{U}$ and let $\phi:(U,H)\rightarrow (\tilde{U},\tilde{H})$ be a Heisenberg 
 diffeormorphism. Let $\tilde{P}\in \Psi_{\tilde{H}}^{m}(V)$ and set $P=\phi^{*}\tilde{P}$. We need to show that $P$ is a \psivdo\ of order $m$ on $U$. 

 First, by Proposition~\ref{prop:PsiVDO.characterisation-kernel2}
the distribution kernel of $\tilde{P}$ takes the form,   
 \begin{equation}
      k_{\tilde{P}}(\tilde{x},\tilde{y})= |{\tilde{\varepsilon}_{\tilde{x}}'}| 
      K_{\tilde{P}}(\tilde{x},-{\tilde{\varepsilon}}_{\tilde{x}}(\tilde{y})) + \tilde{R}(\tilde{x},\tilde{y}),  
\label{eq:Appendix.kP}
 \end{equation}
with $K_{\tilde{P}}(\tilde{x},\tilde{y})$ in $\cK^{\hat m}(\tilde{U}\times\Rd)$ and  $\tilde{R}(\tilde{x},\tilde{y})$ in 
$C^\infty(\tilde{U}\times \tilde{U})$. 
Therefore,  the distribution kernel  of $P=\phi^{*}P$ is given by 
\begin{equation}
    k_{P}(x,y) = |\phi'(y)|  k_{\tilde{P}}(\phi(x),\phi(y))= 
     |\varepsilon'_{x}| K(x, -\varepsilon_{x}(y)) +  \tilde{R}(\phi(x),\phi(y)),
\end{equation}
where $K$ is the distribution on $\cU=\{ (x,y)\in \URd; \  \varepsilon_{x}^{-1}(-y)\in U\} \subset \URd$ given by 
 \begin{equation}
     K(x,y)=|\partial_{y}\Phi(x,y)| K_{\tilde{P}}(\phi(x), \Phi(x,y)), 
     \quad \Phi(x,y)= -\tilde{\varepsilon}_{\phi(x)}\circ \phi\circ \varepsilon_{x}^{-1}(-y). 
 \end{equation}

 Next, it follows from \cite[Props.~3.16, 3.18]{Po:Pacific1} that we have
\begin{equation}
    \Phi(x,y)=-\phi_{H}'(x)(-y)+\Theta(x,y)=\phi_{H}'(x)(y)+\Theta(x,y),
\end{equation}
where $\Theta(x,y)$ is a smooth function on $\cU$ with a behavior near $y=0$ of the form
\begin{equation}
    \Theta(x,y)=(\op{O(\|y\|^{3})},  \op O(\|y\|^{2}),\ldots, \op O(\|y\|^{2})).
     \label{eq:Appendix.invariance.Theta}
\end{equation}
Then a Taylor expansion around $\tilde{y}=\phi_{H}'(x)y$ gives
\begin{equation}
   K(x,y) =  \sum_{\brak\alpha<N}|\partial_{y}\Phi(x,y)| \frac{\Theta(x,y)^\alpha}{\alpha!} 
    (\partial^\alpha_{\tilde{y}}K_{\tilde{P}})(\phi(x),\phi_{H}'(x)y) + R_{N}(x,y), 
   \label{eq:PsivDO.invariance.TaylorK1}
\end{equation}
where $R_{N}(x,y)$ is equal to 
\begin{equation}
   R_{N}(x,y) =  \sum_{\brak\alpha=N} 
   |\partial_{y}\Phi(x,y) | \frac{\Theta(x,y)^\alpha}{\alpha!} \int_{0}^1 (t-1)^{N-1}
    \partial^\alpha_{\tilde{y}}K_{\tilde{P}}(\phi(x),\Phi_{t}(x,y) ) dt.
    \label{eq:PsivDO.invariance.TaylorK2} 
\end{equation}
and we have let $\Phi_{t}(x,y) =\phi_{H}'(x)y +t\Theta(x,y)$.

Set $f_{\alpha}(x,y)=|\partial_{y}\Phi(x,y)| \Theta(x,y)^\alpha$. Then~(\ref{eq:Appendix.invariance.Theta})  implies 
that near $y=0$ we have 
\begin{equation}
    f_{\alpha}(x,y)=\op{O}(\|y\|^{3\alpha_{0}+2(\alpha_{1}+\ldots+\alpha_{d})})=\op{O}(\|y\|^{\frac32\brak\alpha}).
     \label{eq:Appendix.behavior-f-alpha}
\end{equation}
Thus all the homogeneous components of degree $< \frac32\brak\alpha$ in  
the Taylor expansion for $f_{\alpha}(x,y)$ at $y=0$ must be zero. Therefore, we can write
\begin{equation}
    f_{\alpha}(x,y)= \sum_{\frac32\brak\alpha \leq \brak\beta < \frac{3}{2}N} \!\! f_{\alpha\beta}(x) \frac{y^\beta }{\beta !}+ 
    \sum_{\brak\beta\dot{=}\frac{3}{2}N}r_{N\alpha\beta}(x,y)y^{\beta},
     \label{eq:Appendix.behavior-Theta-alpha}
\end{equation}
where we have let $f_{\alpha\beta}(x)=\partial_{y}^{\beta}f_{\alpha}(x,0)$, the functions $r_{M\alpha\beta}(x,y)$ are in $C^{\infty}(\cU)$ 
and the notation $\brak 
\beta\dot{=}\frac{3}{2}N$ means that $\brak 
\beta$ is equal to $\frac{3}{2}N$ if $\frac{3}{2}N$ is an integer and to $\frac{3}{2}N+\frac{1}{2}$ otherwise. Thus,
\begin{equation}
    K(x,y)= \sum_{\brak\alpha<N} \sum_{\frac32\brak\alpha \leq \brak\beta < \frac{3}{2}N}\!\! K_{\alpha\beta}(x,y) + 
    \sum_{\brak\alpha<N} R_{N\alpha}(x,y) + R_{N}(x,y),
     \label{eq:Appendix.expansion-K}
\end{equation}
where we have let
\begin{gather}
     K_{\alpha\beta}(x,y)= f_{\alpha\beta}(x) y^\beta 
    (\partial^\alpha_{\tilde{y}}K_{\tilde{P}})(\phi(x),\phi_{H}'(x)y),\\ R_{N\alpha }(x,y) = 
    \sum_{\brak\beta\dot{=}\frac{3}{2}N} r_{M\alpha\beta}(x,y)y^\beta (\partial^\alpha_{\tilde{y}}K_{\tilde{P}})(\phi(x),\phi_{H}'(x)y). 
\end{gather}

As in the proof of Proposition~\ref{prop:PsiVDO.characterisation-kernel2} 
the smoothness of $\phi_{H}'(x)y$ and the fact that $\phi_{H}'(x)(\lambda.y)=\lambda.\phi_{H'}(x)y$ 
for any $\lambda 
\in \R$  imply that $K_{\alpha\beta}(x,y)$ belongs to $\cK^{\hat m -\brak\alpha+\brak\beta}(U\times\Rd)$. 
Notice that if  $\frac32\brak\alpha \leq \brak\beta \dot{=}\frac{3}{2}N$ then we have 
$\Re \hat m-\brak\alpha+\brak\beta\geq \Re\hat m +\frac{1}{3}\brak\beta\geq \Re\hat m + \frac{N}{2}$.
It then follows from Lemma~\ref{lem:Appendix.Heisenberg.invariance} that, for any integer $J$, the remainder term $R_{N\alpha}$ is 
in $C^{J}(\URd)$ as soon as $N$ is large enough.  

Let $\pi_{x}:\URd \rightarrow U$ denote the projection on the first coordinate. In the sequel we will say that a distribution $K(x,y)\in \cD'(\URd)$ is properly 
supported with respect to $x$ when $\pi_{x|_{\supp K}}$ is a proper map, i.e., for any compact $L \subset U$ the set $\supp K \cap (L\times \Rd)$ is compact.

In order to deal with the regularity of $R_{N}(x,y)$ in~(\ref{eq:PsivDO.invariance.TaylorK2}) 
we need the lemma below. 
\begin{lemma-appendix}\label{lem:Appendix.Theta-alpha}
  There exists a  function $\chi_{n}\in  C^{\infty}_{c}(\cU)$ properly supported with respect to $x$ such that   
  $\chi(x,y)=1$ near $y=0$ and, for any multi-order $\alpha$, we can write   
 \begin{equation}
      \chi(x,y)\Theta(x,y)^{\alpha}=\sum_{\brak \beta\dot{=}\frac{3}{2}\brak \alpha} \theta_{\alpha\beta}(t,x,y) \Phi_{t}(x,y)^{\beta} 
      \label{eq:Appendix.claim}
\end{equation}
where the functions $\theta_{\alpha\beta}(t,x,y)$ are in $C^{\infty}([0,1]\times U\times \Rd)$.
\end{lemma-appendix}
\begin{proof}[Proof of the lemma]
   Let $U'$ be a relatively compact open subset of $U$ and let $(t_{0},x_{0})\in [0,1]\times U'$. 
   Since $\Phi_{t_{0}}(x_{0},0)=0$ and $\partial_{y}\Phi_{t_{0}}(x_{0},0)= \phi'_{H}(x_{0})$ is invertible
  the implicit function theorem implies that there 
  exist an open interval $I_{x_{0}}$ containing $t_{x_{0}}$, an open subset $U_{x_{0}}$ of $U$ containing $x_{0}$, open subsets 
   $V_{x_{0}}$ 
   and $\tilde{V}_{x_{0}}$ of $\Rd$ containing $0$ and a smooth map $\Psi_{x_{0}}(t,x,\tilde{y})$ from 
   $I_{x_{0}}\times U_{x_{0}}\times \tilde{V}_{x_{0}}$ to 
   $V_{x_{0}}$ such that 
    $U_{x_{0}}\times V_{x_{0}}$ is contained in $\cU$ and  for any $(t,x,y)$ in $I_{x_{0}}\times U_{x_{0}}\times V_{x_{0}}$ and any 
    $\tilde{y}$ in $\tilde{V}_{x_{0}}$ we have 
    \begin{equation}
        \tilde{y}=\Phi_{t}(x,y) \Longleftrightarrow y=\Psi_{x_{0}}(t,x,\tilde{y}). 
         \label{eq:Appendix.definition-Psi}
    \end{equation}
   
    Since $[0,1]\times \overline{U'}$ is compact we can cover it by finitely many products $I_{x_{k}}\times U_{x_{k}}$, $k=1,..,p$, with 
    $(t_{k},x_{k})\in [0,1]\times \overline{U'}$. In particular, the sets $I=\cup_{k}I_{k}$ and $U''=\cup U_{k}$ are open neighborhoods of $I$ and 
    $\overline{U'}$ respectively. Thanks to~(\ref{eq:Appendix.definition-Psi}) 
    we have $\Psi_{x_{k}}=\Psi_{x_{l}}$ on 
    $(I_{x_{k}}\times U_{x_{k}}\times V_{x_{k}})\cap (I_{x_{l}}\times U_{x_{l}}\times V_{x_{l}})$. Therefore, 
    setting $V= \cap_{k}V_{k}$ and $\tilde{V}=\cap_{k}\tilde{V}_{k}$ we have $U''\times V\subset \cU$ and there exists a smooth map $\Psi$ from 
    $I\times U''\times \tilde{V}$ such that for any $(t,x,y)$ in $I\times U''\times V$ and any 
    $\tilde{y}$ in $\tilde{V}$ we have 
    \begin{equation}
        \tilde{y}=\Phi_{t}(x,y) \Longleftrightarrow y=\Psi(t,x,\tilde{y}). 
    \end{equation}
    
   Furthermore, as $\partial_{\tilde{y}} \Psi(t,x,0)=[\partial_{y}\Phi_{t}(x,0)]^{-1}=\phi_{H}'(x)^{-1}$ and for any $\lambda \in \R$ we have
   $\phi_{H}'(x)^{-1}(\lambda.y)=\lambda \phi_{H}'(x)^{-1}(y)$, the function $\Theta(x,\Psi(t, x,\tilde{y}))$ 
   behaves near 
   $\tilde{y}=0$ as in~(\ref{eq:Appendix.invariance.Theta}), so as 
   in~(\ref{eq:Appendix.behavior-f-alpha}) and~(\ref{eq:Appendix.behavior-Theta-alpha})  for any multi-order 
   $\alpha$ we can write 
   \begin{equation}
       \Theta(x,\Psi(t, x,\tilde{y}))^{\alpha}= \sum_{\brak\beta \dot{=} \frac32 \brak\alpha} 
          \tilde{\theta}_{\alpha\beta}(t,x,\tilde{y}) \tilde{y}^{\beta},
   \end{equation}
 for some functions  $\tilde{\theta}_{\alpha\beta}(t,x,\tilde{y})\in C^{\infty}(I\times 
   U''\times \tilde{V})$.   Setting $\tilde{y}=\Phi_{t}(x,y)$ then gives 
  \begin{equation}
       \Theta(x,y)^{\alpha}=\sum_{\brak  \beta\dot{=}\frac{3}{2}\brak\alpha} \theta_{\alpha\beta}(t,x,y) 
       \Phi_{t}(x,y)^{\beta},
       \label{eq:Appendix.lemma.form-Theta-alpha}
  \end{equation}
    for some functions  $\theta_{\alpha\beta}(t,x,y)\in C^{\infty}(I\times U''\times V)$.
 
    All this allows us to construct locally finite coverings $(U'_{n})_{n\geq 0}$ and $(U''_{n})_{n\geq 0}$ of $U$ by relatively compact open 
   subsets in such way that, for each integer $n$, the open $U''_{n}$ contains $\overline{U'_{n}}$ and there exists an open $V_{n}\subset \Rd$ containing $0$
so that, for any multiorder $\alpha$, on $[0,1]\times U''_{n}\times V_{n}$ we have   
     \begin{equation}
       \Theta(x,y)^{\alpha}= \sum_{\brak \beta\dot{=}\frac{3}{2}\brak\alpha} \theta_{\alpha\beta}^{(n)}(t,x,y) \Phi_{t}(x,y)^{\beta}, %\qquad 
      \label{eq:Appendix.lemma.form-Theta-alpha-n}
   \end{equation}
for some functions $ \theta_{\alpha\beta}^{(n)}(t,x,y)\in C^{\infty}([0,1]\times U''_{n}\times V_{n})$.

For each $n$ let $\varphi_{n}\in C^{\infty}_{c}(U_{n}'')$ be such that $\varphi_{n}=1$ on $U_{n}'$ and  let 
  $\psi_{n}\in C^{\infty}_{c}(V_{n})$ be such that $\psi_{n}=1$ on a neighborhood $V'_{n}$ of $0$. Then we   
  construct a locally  finite family $(\chi_{n})_{n\geq 0} \subset C^{\infty}_{c}(\cU)$ as follows: for $n=0$ we set  
  $\chi_{0}(x,y)=\varphi_{0}(x)\psi_{0}(y)$ and for $n\geq 1$ we let
\begin{equation}
   \chi_{n}(x,y)= (1-\varphi_{0}(x)\psi_{0}(y))\ldots (1-\varphi_{n-1}(x)\psi_{n-1}(y))\varphi_{n}(x)\psi_{n}(y) . 
\end{equation}
  Then $\chi =  \sum_{n\geq 0} \chi_{n}$ is  a well defined smooth function on $\URd$ 
  supported on $\cup_{n\geq 0} (U_{n}''\times V_{n}) \subset \cU$, hence properly supported with respect to $x$. Also, as $\chi(x,y)$ is equal to $1$ on each product 
  $U_{n}'\times V_{n}'$ we see that $\chi(x,y)=1$ on a neighborhood of $U \times \{0\}$. In addition, 
  thanks to~(\ref{eq:Appendix.lemma.form-Theta-alpha-n}) 
  on $[0,1]\times\cU$ we have 
  \begin{equation}
      \chi(x,y)\Theta(x,y)^{\alpha}=\sum_{n\geq 0} \chi_{n}(x,y)\Theta(x,y)^{\alpha}=  \sum_{\brak \beta\dot{=}\frac{3}{2}\brak\alpha}  
      \theta_{\alpha\beta}(t,x,y) \Phi_{t}(x,y)^{\beta},
  \end{equation}
where $\theta_{\alpha\beta}(t,x,y) :=\sum_{n} \chi_{n}(x,y) \theta_{\alpha\beta}^{(n)}(t,x,y)$ belongs to 
$C^{\infty}([0,1]\times\URd)$. The lemma is thus proved.
 \end{proof}

Let us go back  to the proof of Proposition~\ref{prop:PsiHDO.invariance}. Thanks to~(\ref{eq:PsivDO.invariance.TaylorK2}) 
and~(\ref{eq:Appendix.claim})  we see that $ \chi(x,y) R_{N}(x,y)$ is equal to
\begin{equation}
   \sum_{\brak\alpha=N}\sum_{\brak\beta\dot{=}\frac{3}{2}N} \int_{0}^{1} 
    r_{\alpha\beta}(t,x,y)(\tilde{y}^{\beta}\partial^\alpha_{\tilde{y}}K_{\tilde{P}})(\phi(x),\Phi_{t}(x,y))dt, 
\end{equation}
with $r_{\alpha\beta}(t,x,y)$ in $C^{\infty}([0,1\times U\times \Rd)$. Since  
$\tilde{y}^{\beta}\partial^\alpha_{\tilde{y}}K_{\tilde{P}}$ is in $\cK^{\hat{m}+N/2}(\tilde{U}\times \Rd)$ it follows from 
Lemma~\ref{lem:Appendix.Heisenberg.invariance} that, for any integer $J\geq 0$, as soon as 
$N$ is taken large enough $\tilde{y}^{\beta}\partial^\alpha_{\tilde{y}}K_{\tilde{P}}$ is in $C^{J}(\tilde{U}\times \Rd)$ and so 
$\chi(x,y) R_{N}(x,y)$ is in $C^{J}(\URd)$.

On the other hand, set $K_{P}(x,y)=\chi(x,y)K(x,y)=\sum \chi_{n}(x,y)K(x,y)$. Since $\chi(x,y)$ is supported in $\cU$ and is properly supported with 
respect to $x$ this defines a distribution on $U\times \Rd$. Moreover, using~(\ref{eq:Appendix.expansion-K}) 
we get
\begin{equation}
    K_{P}(x,y)= \sum_{\brak \alpha<N} \sum_{\frac32\brak\alpha \leq \brak\beta 
    < \frac{3}{2}N} K_{\alpha\beta}(x,y)+\sum_{j=1}^{3}R_{N}^{(j)},
        \label{eq:Appendix.expansion-KP1}
 \end{equation}   
 where the remainder terms $R_{N,z}^{(j)}$, $j=1,2,3$ are given by 
 \begin{gather}
R_{N}^{(1)}=\chi(x,y)R_{N}(x,y),\quad  R_{N}^{(2)}= \sum_{\brak \alpha <N} \chi(x,y) R_{N\alpha}(x,y),\\
   R_{N}^{(2)}(x,y)=  \sum_{\brak \alpha<N} \sum_{\frac32\brak\alpha \leq \brak\beta 
    <\frac{3}{2}N}  (1-\chi(x,y))K_{\alpha\beta}(x,y).
\end{gather}
Each term $K_{\alpha\beta}(x,y)$ belongs to $\cK^{\hat 
m-\brak\alpha+\brak \beta}(\URd)$ and, as $\hat{m}+\brak\beta-\brak\alpha=\hat{m}+j$ 
and $\frac32\brak\alpha \leq \brak\beta$ imply $\brak\alpha\leq 2j$ and $\brak\beta \leq \frac{4}{3}j$, 
in the r.h.s.~(\ref{eq:Appendix.expansion-KP1}) there are only  finitely 
many such distributions in a given space $\cK^{\hat m+j}(\URd)$ as $\alpha$ and $\beta$ range over all multi-orders such that 
$\frac32\brak\alpha \leq \brak\beta$.    

Furthermore, the reminder term $R_{N}^{(3)}$ is smooth and the other remainder terms $R_{N}^{(j)}$, $j=1,2$, are in $C^{J}(\URd)$ as soon as $N$ is 
large enough.  Thus, 
\begin{equation}
    K_{P}(x,y)  \sim\sum_{\frac32\brak\alpha \leq \brak\beta} K_{\alpha\beta}(x,y), 
    \label{eq:Appendix.asymptotic-expansion-KP'}
\end{equation}
which implies that $K_{P}$ belongs to $\cK^{\hat m}(\URd)$ and satisfies~(\ref{eq:Appendix.asymptotic-expansion-KP'}).

Finally, from~(\ref{eq:Appendix.kP}) and the very definition of $\Phi(x,y)$ on $U\times U$, we deduce that the distribution kernel of $P$ differs from 
$ |\varepsilon_{x}'| K_{P}(x,-\varepsilon_{x}(y)$ by the smooth function
% of the form, 
% is equal to
% \begin{multline}
%   k_{P}(x,y)=   
\begin{equation}
      [1-\chi(x,\varepsilon_{x}(y))]  
           |\tilde{\varepsilon}'_{\phi(x)}| K_{\tilde{P}}(\phi(x),-\tilde{\varepsilon}_{\phi(x)}\circ\phi(y))  + \tilde{R}(\phi(x),\phi(y)).%\\
%     \label{eq:}
\end{equation}
%        =  |\varepsilon_{x}'| K_{P}(x,-\varepsilon_{x}(y)) \quad \bmod C^{\infty}(U\times U).
% \end{multline}
Combining this with Proposition~\ref{prop:PsiVDO.characterisation-kernel2} and the fact that $K_{P}(x,y)$ 
satisfies~(\ref{eq:Appendix.asymptotic-expansion-KP'}) proves Proposition~\ref{prop:PsiHDO.invariance}.

\chapter{Proof of Proposition~3.1.21}
\label{chap:Appendix-transpose}
Let $P:C^{\infty}_{c}(U)\rightarrow C^{\infty}(U)$ be a  \psivdo\ of order $m$ and let us show that its transpose operator 
$P^{t}:C^{\infty}_{c}(U)\rightarrow C^{\infty}(U)$ is a \psivdo\ of order $m$. 
By Proposition~\ref{prop:PsiVDO.characterisation-kernel2} the distribution kernel of $P$ is of the form, 
\begin{equation}
    k_{P}(x,y)=|\varepsilon_{x}'|K_{P}(x,-\varepsilon_{x}(y)) +R(x,y),
\end{equation}
with $K_{P}(x,y)$ in $\cK^{m}(\URd)$ and $R(x,y)$ in $C^{\infty}(U\times U)$. Thus the distribution kernel $k_{P^{t}}(x,y)=k_{P}(y,x)$ of 
$P^{t}$ is equal to  
\begin{equation}
   |\varepsilon_{y}'|K_{P}(y,-\varepsilon_{y}(x)) +R(y,x) = |\varepsilon_{x}'|K(x,-\varepsilon_{x}(y)) +R(y,x),
     \label{eq:Appendix.kPt-K}
\end{equation}
where $K$ is the distribution on the open $\cU=\{(x,y); \ \varepsilon_{x}^{-1}(-y)\in U\}$ given by
\begin{equation}
    K(x,y)=|\varepsilon_{x}'|^{-1}|\varepsilon_{y}'|K_{P}(\varepsilon_{x}^{-1}(-y), -\varepsilon_{\varepsilon_{x}^{-1}(-y)}(x)).
\end{equation}

\begin{lemma-appendix}
    On $\cU$ we have
    \begin{equation}
        \varepsilon_{\varepsilon_{x}^{-1}(-y)}(x)=y-\Theta(x,y),
    \end{equation}
 where $\Theta: \cU \rightarrow \Rd$ is a smooth map with a behavior near $y=0$ of the form~(\ref{eq:Appendix.invariance.Theta}). 
\end{lemma-appendix}
\begin{proof}
    Let $(x,y)\in \cU$ and $Y\in G_{x}U$ let $\lambda_{y}(Y)=y.Y$, that is,  $\lambda_{y}$ is the left multiplication by $y$ on $G_{x}U$. 
    Then by~\cite[Eq.~(3.32)]{Po:Pacific1} for $Y$ small enough we have
    \begin{equation}
        \lim_{t\rightarrow 0} \varepsilon_{x}\circ \varepsilon_{\varepsilon_{x}^{-1}(t.-y)}^{-1}(t.Y)=\lambda_{-y}(Y)=\lambda_{y}^{-1}(Y).
    \end{equation}
    Since $\varepsilon_{x}\circ \varepsilon_{\varepsilon_{x}^{-1}(t.-y)}^{-1}(t.Y)$ is a smooth function of $(t,Y)$ near $(0,0)$, it 
    follows from the implicit function theorem that for $Y$ small enough we have 
    \begin{equation}
        \lim_{t \rightarrow 0} t^{-1}.\varepsilon_{\varepsilon_{x}^{-1}(t.-y)} \circ \varepsilon_{x}^{-1}(Y)=\lambda_{y}(Y).
    \end{equation}
    In particular, for $Y=0$ we get 
    \begin{equation}
         \lim_{t \rightarrow 0} t^{-1}.\varepsilon_{\varepsilon_{x}^{-1}(t.-y)}(x)=y.
         \label{eq:Appendix.inverse-asymptotic}
    \end{equation}
    
    Now, the function $ \varepsilon_{\varepsilon_{x}^{-1}(-y)}(x)$ depends smoothly on $(x,y)\in \cU$, so~(\ref{eq:Appendix.inverse-asymptotic}) 
    allows us to put it into the form,
    \begin{equation}
        \varepsilon_{\varepsilon_{x}^{-1}(-y)}(x)=y-\Theta(x,y),
    \end{equation}
 where $\Theta=(\Theta_{0},\ldots,\Theta_{d})$ is smooth map from $\cU$ to $\Rd$ with a behavior near $y=0$ of the form
 \begin{equation}
     \Theta_{0}(x,y)=\op{O}(|y_{0}|^{2}+|y_{0}||y|+|y|^{3}), \qquad \Theta_{j}(x,y)=\op{O}(|y|^{2}), \quad j=1,\ldots,d. 
 \end{equation}
 In particular, near $y=0$ the map $\Theta$ has a behavior of the form~(\ref{eq:Appendix.invariance.Theta}). %The lemma is thus proved. 
\end{proof}

Next, a Taylor expansion around $(\varepsilon_{x}^{-1}(-y),-y)$ gives
\begin{equation}
        K(x,y)=\sum_{\brak\alpha<N} |\varepsilon_{x}'|^{-1}|\varepsilon_{y}'|\frac{\theta(x,y)^{\alpha}}{\alpha!} 
        (\partial_{y}^{\alpha}K_{P})(\varepsilon_{x}^{-1}(-y),-y) +R_{N}(x,y),
%     \label{eq:}
\end{equation}
  where $ R_{N}(x,y)$ is equal to
\begin{equation}
       \sum_{\brak\alpha=N}|\varepsilon_{x}'|^{-1}|\varepsilon_{y}'|\frac{\theta(x,y)^{\alpha}}{\alpha!} \int_{0}^{1}(1-t)^{N-1} 
        \partial_{y}^{\alpha}K_{P})(\varepsilon_{x}^{-1}(-y),\Phi_{t}(x,y)),
        \label{eq:AppendixB.RN}
\end{equation}
and we have let $\Phi_{t}(x,y)=-y+t\Theta(x,y)$. 

Let $a_{\alpha}(x,y)= 
|\varepsilon_{x}'|^{-1}|\varepsilon_{y}'|\frac{\theta(x,y)^{\alpha}}{\alpha!}$. Thanks to~(\ref{eq:Appendix.invariance.Theta}) 
the same arguments used to prove~(\ref{eq:Appendix.behavior-Theta-alpha}) show 
that there exist functions $r_{N\alpha}(x,y)\in C^{\infty}(\cU)$, $\brak \beta \dot{=}\frac{3}{2}N$, so that 
\begin{equation}
    a_{\alpha}(x,y)=\sum_{\frac{3}{2}\brak\alpha \leq \brak \beta < \frac{3}{2}N}a_{\alpha\beta}(x)y^{\beta} + 
    \sum_{\brak \beta \dot{=}\frac{3}{2}N}r_{N\alpha}(x,y)y^{\beta},
\end{equation}
where we have let $a_{\alpha \beta}(x)=\frac{1}{\beta!}\partial^{\beta}a_{\alpha}(x,0)$. Therefore, we get %$K(x,y)$ is equal to
\begin{multline}
K(x,y)=    \sum_{\brak\alpha<N} \sum_{\frac{3}{2}\brak\alpha \leq \brak \beta < \frac{3}{2}N} a_{\alpha\beta}(x)y^{\beta} 
    (\partial_{y}^{\alpha}K_{P})(\varepsilon_{x}^{-1}(-y),-y) \\ +  \sum_{\brak\alpha<N}R_{N\alpha}(x,y) +R_{N}(x,y),
\end{multline}
where we have let
\begin{equation}
        R_{N\alpha}(x,y)= \sum_{\brak \beta \dot{=}\frac{3}{2}N}r_{N\alpha}(x,y)y^{\beta} (\partial_{y}^{\alpha}K_{P})(\varepsilon_{x}^{-1}(-y),-y).
        \label{eq:AppendixB.RNalpha}
%     \label{eq:}
\end{equation}

Next, a further Taylor  expansion gives
\begin{multline}
(\partial_{y}^{\alpha}K_{P})(\varepsilon_{x}^{-1}(-y),-y)= \sum_{|\gamma|<N}\frac{1}{\gamma!}(\varepsilon_{x}^{-1}(-y)-x)^{\gamma} 
  (\partial^{\gamma}_{x}\partial_{y}^{\alpha}K_{P})(x,-y)\\ + \sum_{|\gamma|=N} \int_{0}^{1}(1-t)^{N-1} 
  (\partial^{\gamma}_{x}\partial_{y}^{\alpha}K_{P})(\varepsilon_{t}(x,y),-y),
\end{multline}
where we have let $\varepsilon_{t}(x,y)=x+t(\varepsilon_{x}^{-1}(-y)-x)$. Since $\varepsilon_{x}^{-1}(-y)-x$ is polynomial in $y$ of degree 2 and 
vanishes for $y=0$, we can write
\begin{equation}
    \frac{1}{\gamma!}(\varepsilon_{x}^{-1}(-y)-x)^{\gamma} =\sum_{|\gamma|\leq |\delta \leq 2|\gamma|} b_{\gamma\delta}(x)y^{\delta},
\end{equation}
where we have let $b_{\gamma\delta}(x)=\frac{1}{\gamma!\delta!}[\partial_{y}(\varepsilon_{x}^{-1}(-y)-x)^{\gamma}] _{y=0}$. 
Thus, 
% we can put $K(x,y)$ into the form, 
\begin{multline}
    K(x,y)= \sum_{\alpha,\beta,\gamma,\delta}^{(N)}
    K_{\alpha\beta\gamma\delta}(x,y) + \sum_{\brak\alpha<N} \sum_{\frac{3}{2}\brak\alpha \leq \brak \beta < \frac{3}{2}N} R_{N\alpha\beta}(x,y)\\ +
    \sum_{\brak\alpha<N} R_{N\alpha}(x,y)  +R_{N}(x,y),
\end{multline}
where the first summation goes over all the multi-orders $\alpha$, $\beta$, $\gamma$ and $\delta$ such that $\brak\alpha<N$, $\frac{3}{2}\brak\alpha \leq 
\brak \beta < \frac{3}{2}N$ and $|\gamma|\leq |\delta| \leq 2|\gamma|<2N$, and 
% \begin{gather}
\begin{equation}
         K_{\alpha\beta\gamma\delta}(x,y)=f_{\alpha\beta\gamma\delta}(x) y^{\beta+\delta}  (\partial^{\gamma}_{x}\partial_{y}^{\alpha}K_{P})(x,-y), 
%     \label{eq:}
\end{equation}
with $f_{\alpha\beta\gamma\delta}(x)=a_{\alpha\beta}(x) b_{\gamma\delta}(x)$ and $R_{N\alpha \beta}(x,y)$ is equal to 
\begin{equation}
    \sum_{|\gamma|=N}\sum_{N\leq |\delta|\leq 2N} a_{\alpha\beta\gamma\delta}(x)y^{\beta+\delta} \int_{0}^{1}(1-t)^{N-1} 
      (\partial^{\gamma}_{x}\partial_{y}^{\alpha}K_{P})(\varepsilon_{t}(x,y),-y).
      \label{eq:AppendixB.RNalphabeta}
%     \label{eq:}
\end{equation}
% \end{gather}

Now, the distribution $y^{\beta}K_{P}(x,-y)$ belongs to $\cK^{\hat m-\brak \alpha+\brak \beta}(\URd)$. In particular, 
if $\frac{3}{2}\brak \alpha \leq \brak \beta\dot{=}\frac{3}{2}N$ then $\Re \hat{m}-\brak\alpha+\brak \beta\geq \Re \hat{m}+\frac{1}{3}\brak \beta \geq  \Re 
\hat{m}+\frac{1}{2}N$. Therefore, for any given integer $J$  Lemma~\ref{lem:Appendix.Heisenberg.invariance} tells us that
$y^{\beta}K_{P}(x,-y)$ is in $C^{J}(\URd)$ as soon as $N$ is large enough. 
It follows that all the remainder terms 
$R_{N\alpha}(x,y)$, $\brak \alpha<N$, belong to $C^{J}(\cU)$ for $N$ large enough. 

Similarly, if $\frac{3}{2}\brak \alpha \leq \brak \beta $ and $|\gamma|=N\leq |\delta|\leq 2N$ then $\Re \hat{m}-\brak\alpha +\brak \beta+\brak 
\delta\geq \Re \hat{m}+\brak{\delta}\geq \Re \hat{m}+\frac{1}{2}N$, so using Lemma~\ref{lem:Appendix.Heisenberg.invariance} 
we see that 
$y^{\beta+\delta} (\partial^{\gamma}_{x}\partial_{y}^{\alpha}K_{P})(x,-y)$ is in $C^{J}(\URd)$ for $N$ large enough. It then follows that for $N$ 
large enough the remainder terms $R_{N\alpha\beta}(x,y)$  with $\brak \alpha<N$ and $\frac{3}{2}\brak \alpha \leq \brak \beta \dot{=}\frac{3}{2}N$ are all 
in $C^{J}(\cU)$ as soon as $N$ is chosen large enough. 

In order to deal with the last remainder term $R_{N}(x,y)$ notice that, along the same lines as that of the proof of 
Lemma~\ref{lem:Appendix.Theta-alpha},  
one can show that there exists a $\chi \in C^{\infty}(\cU)$ properly supported with respect to $x$ such that 
  $\chi(x,y)=1$ near $y=0$ and, for any multi-order $\alpha$, we can write   
 \begin{equation}
      \chi(x,y)\Theta(x,y)^{\alpha}=\sum_{\brak \beta\dot{=}\frac{3}{2}\brak \alpha} \theta_{\alpha\beta}(t,x,y) \Phi_{t}(x,y)^{\beta},
\end{equation}
where the functions $\theta_{\alpha\beta}(t,x,y)$ are in $C^{\infty}([0,1]\times U\times \Rd)$. Therefore $ \chi(x,y) R_{N}(x,y)$ is equal to
\begin{equation}
%     \chi(x,y) R_{N}(x,y)=
    \sum_{\brak\alpha=N}\sum_{\brak \beta\dot{=}\frac{3}{2}\brak \alpha} 
    |\varepsilon_{x}'|^{-1}|\varepsilon_{y}'| \int_{0}^{1}r_{N\alpha\beta}(t,x,y)
    (y^{\beta}\partial_{y}^{\alpha}K_{P})(\varepsilon_{x}^{-1}(-y),\Phi_{t}(x,y)),
\end{equation}
for some functions $r_{N\alpha\beta}(t,x,y)$ in $C^{\infty}([0,1]\times \URd)$. Since $(y^{\beta}\partial_{y}^{\alpha}K_{P})$ is in  
$\cK^{\hat{m}-\brak \alpha+\brak\beta}(\URd)$ and we have $\Re \hat{m}-\brak\alpha+\brak\beta \geq \Re \hat{m}+\frac{1}{2}\brak\alpha=\Re 
\hat{m}+\frac{1}{2}N$, using Lemma~\ref{lem:Appendix.Heisenberg.invariance} we see 
that for $N$ large enough $ \chi(x,y) R_{N}(x,y)$ is in $C^{J}(\cU)$. As $ \chi(x,y) R_{N}(x,y)$ is a properly supported  with respect to 
$x$ this shows that it belongs to $C^{J}(\URd)$.

Let $K_{P^{t}}(x,y)=\chi(x,y)K(x,y)$. This defines a distribution on $\URd$ since $\chi$ is properly supported with respect to $x$. Moreover, we have
\begin{equation}
  K_{P^{t}}(x,y)= \sum_{\alpha,\beta,\gamma,\delta}^{(N)}
    K_{\alpha\beta\gamma\delta}(x,y) + \sum_{j=1}^{4}R_{N}^{(j)}(x,y), 
\end{equation}
where the remainder terms $R_{N}^{(j)}$, $j=1,\ldots,4$, are given by the formulas,
\begin{gather}
    R_{N}^{(1)}=\chi(x,y)R_{N}(x,y), \qquad 
    R_{N}^{(2)}= \sum_{\brak\alpha<N}\chi(x,y)R_{N\alpha}(x,y),\\
    R_{N}^{(3)}= \sum_{\brak\alpha<N} \sum_{\frac{3}{2}\brak\alpha \leq \brak \beta < \frac{3}{2}N}\chi(x,y)R_{N\alpha\beta}(x,y),\\
    R_{N}^{(4)}= \sum_{\alpha,\beta,\gamma,\delta}^{(N)}(1-\chi(x,y))K_{\alpha\beta\gamma\delta}(x,y).
\end{gather}

Note that $K_{\alpha\beta\gamma\delta}(x,y)$ belongs to $\cK^{\hat{m}-\brak\alpha+\brak\gamma+\brak\delta}(\URd)$ and there are finitely many terms 
of a given order as $\alpha$, $\beta$, $\gamma$ and $\delta$ ranges over all the multi-orders such that $\frac{3}{2}\brak\alpha \leq 
\brak \beta$ and $|\gamma|\leq |\delta| \leq 2|\gamma|$. 

On the other hand, the remainder term $R_{N}^{(4)}$ is smooth and it follows from the above observations that the other remainder terms are in 
$C^{J}(\URd)$ as soon  as 
$N$ is large enough. Thus,
\begin{equation}
    K_{P^{t}}(x,y) \sim \sum_{\frac{3}{2}\brak\alpha \leq \brak \beta} \sum_{|\gamma|\leq |\delta| \leq 2|\gamma|}K_{\alpha\beta\gamma\delta}(x,y),
     \label{eq:Appendix.expansionKPt}
\end{equation}
which incidentally shows that $K_{P^{t}}(x,y)$ belongs to $\cK^{\hat{m}}(\URd)$. 

Finally, thanks to~(\ref{eq:Appendix.kPt-K}) we can put the kernel of $P^{t}$ into the form,
\begin{multline}
    k_{P^{t}}(x,y)=|\varepsilon_{x}'|K_{P}(x,-\varepsilon_{x}(y))+|\varepsilon_{x}'|[(1-\chi)K](x,-\varepsilon_{x}(y))+R(y,x)\\ 
    = |\varepsilon_{x}'|K_{P}(x,-\varepsilon_{x}(y)) \quad \bmod C^{\infty}(U\times U).
\end{multline}
It then follows from Proposition~\ref{prop:PsiVDO.characterisation-kernel2} that $P^{t}$ is a \psivdo\ of order $m$. Moreover, working out the expression for 
$K_{\alpha\beta\gamma\delta}$ shows that the asymptotics expansion~(\ref{eq:Appendix.expansionKPt}) 
reduces to~(\ref{eq:PsiHDO.transpose-expansion-kernel}). The proof of Proposition~\ref{prop:PsiHDO.transpose-chart} 
is thus achieved.

\backmatter

\end{document}